\documentclass[12pt,a4paper,twoside]{book}

\usepackage{amsmath, exscale, epsfig, amssymb}

\small
\normalsize
\usepackage{amsfonts}
\usepackage{amsbsy}
\usepackage{fancyhdr}
\usepackage{amscd}
\usepackage{theorem}
\usepackage{epsf}
\usepackage{psfrag}
\usepackage{multicol}
\allowdisplaybreaks

\def\build#1_#2^#3{\mathrel{\mathop{\kern 0pt#1}\limits_{#2}^{#3}}}

\def\rems{\noindent{\bf Remarks.} }

\def\supp{{\rm supp\,}}
\def\MM{\hbox{\bf M}}

\def\cq{\hfill\cqfd\par\smallskip}
\def\da{\downarrow}
\def\ua{\uparrow}
\def\ov{\overline}
\def\ba{\begin{eqnarray*}}
\def\ea{\end{eqnarray*}}
\def\cg{\Big[}
\def\cd{\Big]}
\def\pg{\Big(}
\def\pd{\Big)}
\def\proof{\noindent{\bf Proof.} }
\def\rem{\noindent{\bf Remark.} }
\def\d{\partial}

\def\g{{\cal G}}

\def\f{{\cal F}}

\def\b{{\cal B}}

\def\m{{\cal M}}

\def\h{{\cal H}}

\def\e{{\cal E}}
\def\t{{\cal T}}
\def\v{{\cal V}}

\def\z{{\cal Z}}
\def\n{{\cal N}}

\def\w{{\cal W}}

\def\W{{\cal W}}
\def\TT{{\bf T}}
\def\EE{\hbox{\bf E}}
\def\QQ{\hbox{\bf Q}}
\def\wt{\widetilde}

\newcommand{\Q}{\mathbb{Q}}

\newcommand{\K}{\mathbb{K}}
\newcommand{\D}{\mathbb{D}}

\newcommand{\M}{\mathbb{M}}
\newcommand{\T}{\mathbb{T}}

\def\la{\longrightarrow}
\newcommand{\beq}{\begin{eqnarray*}}

\newcommand{\be}{\begin{equation}}
\newcommand{\eeq}{\end{eqnarray*}}
\newcommand{\ben}{\begin{enumerate}}
\newcommand{\een}{\end{enumerate}}
\newcommand{\beqs}{\begin{eqnarray*}&\displaystyle}
\newcommand{\eeqs}{&\end{eqnarray*}}
\newcommand{\wh}{\widehat}

\newcommand{\card}{{\rm Card}}

\newtheorem{theorem}{Theorem}[section]
\newtheorem{lemma}[theorem]{Lemma}
\newtheorem{proposition}[theorem]{Proposition}
\newtheorem{corollary}[theorem]{Corollary}
{\theorembodyfont{\rmfamily}}
{\theorembodyfont{\rmfamily}}
{\theorembodyfont{\rmfamily}\newtheorem{definition}{Definition}[section]}

\newcommand{\E}{\mathbb{E}}

\newcommand{\N}{\mathbb{N}}
\newcommand{\bP}{\mathbb{P}}

\newcommand{\R}{\mathbb{R}}

\newcommand{\Z}{\mathbb{Z}}
\newcommand{\PP}{{\bf P}}

\newcommand{\cT}{\mathcal{T}}
\newcommand{\cU}{\mathcal{U}}

\def\cqfd{ \hfill $\blacksquare$ }

\begin{document}

\setcounter{chapter}{0}
\setcounter{section}{0}

\title{ \bf Random Trees, L\'evy Processes and Spatial Branching Processes.}

\author{
Thomas {\sc Duquesne,} \thanks{Universit\'e Paris 11, Math\'ematiques, 91405 Orsay Cedex, France } 
\and  Jean-Fran\c cois {\sc Le Gall} 
\thanks{D.M.A., Ecole normale sup\'erieure, 45 rue d'Ulm, 75005 
Paris, France} }
\vspace{4mm}
\date{2002}
\maketitle

\tableofcontents

\chapter*{Introduction}

The main goal of this work is to investigate the genealogical structure
of continuous-state branching processes in connection with limit
theorems for discrete Galton-Watson trees. Applications are also given
to the construction and various properties of spatial branching processes
including a general class of superprocesses.

\smallskip
Our starting point is the recent work of Le Gall and Le Jan \cite{LGLJ1}
who proposed a coding of the genealogy of general continuous-state branching processes
via a real-valued random process called the height process. Recall that 
continuous-state branching processes are the continuous analogues of
discrete Galton-Watson branching processes, and that the law of any such process
is characterized by a real function $\psi$ called the branching mechanism. 
Roughly speaking, the height process is a continuous analogue of the contour
process of a discrete branching tree, which is easy to visualize (see Section 0.1,
and note that the previous informal interpretation of the height process is
made mathematically precise by the results of Chapter 2). In the important special case of
the Feller branching diffusion ($\psi(u)=u^2$), the height process is
reflected linear Brownian motion: This unexpected connection between 
branching processes and Brownian motion, or random walk in a discrete setting
has been known for long and exploited by a number of authors (see e.g. 
\cite{Al2}, \cite{Dw}, \cite{Ge}, \cite{NP}, \cite{Ro}). The key contribution of 
\cite{LGLJ1} was to observe that for a general subcritical continuous-state
branching process, there is an explicit formula expressing the height process
as a functional of a spectrally positive L\'evy process
whose Laplace exponent $\psi$ is precisely the branching mechanism. This suggests that 
many problems concerning the genealogy of continuous-state branching processes
can be restated and solved in terms of spectrally positive L\'evy processes, for which a 
lot of information is available (see e.g. Bertoin's recent monograph \cite{Be}).
It is the principal aim of the present work to develop such applications.

In the first two sections below, we briefly describe the objects of interest
in a discrete setting. In the next sections, we outline the main contributions
of the present work.

\section{Discrete trees}

Let
$${\cal U}=\bigcup_{n=0}^\infty \N^n  $$
where $\N=\{1,2,\ldots\}$ and by convention $\N^0=\{\emptyset\}$. If
$u=(u_1,\ldots,u_n)\in\N^n$, we set $|u|=n$, so that 
$|u|$ represents the ``generation'' of $u$. If
$u=(u_1,\ldots u_m)$ and 
$v=(v_1,\ldots, v_n)$ belong to $\cal U$, we write $uv=(u_1,\ldots u_m,v_1,\ldots ,v_n)$
for the concatenation of $u$ and $v$. In particular $u\emptyset=\emptyset u=u$.

\smallskip
A (finite) rooted ordered tree $\cT$ is a finite subset of
$\cU$ such that:
\begin{description}
\item{(i)} $\emptyset\in \cT$.

\item{(ii)} If $v\in \cT$ and $v=uj$ for some $u\in {\cal U}$ and
$j\in\N$, then $u\in\cT$.

\item{(iii)} For every $u\in\cT$, there exists a number $k_u(\cT)\geq 0$
such that $uj\in\cT$ if and only if $1\leq j\leq k_u(\cT)$.
\end{description}
We denote by ${\bf T}$ the set of all rooted ordered trees. In what follows, we see each vertex of the
tree $\cT$ as an individual of a population  whose $\cT$ is the family tree. The cardinality $\#(\cT)$ of $\cT$
is the total progeny.

If $\cT$ is a tree and $u\in \cT$, we define the shift of $\cT$ at $u$
by $\theta_u\cT=\{v\in U:uv\in\cT\}$.
Note that $\theta_u\tau\in {\bf T}$. 

We now introduce the (discrete) {\bf height function} associated with a tree $\cT$. 
Let us denote by $u(0)=\emptyset, u(1), u(2),\ldots 
,u(\#(\cT) -1 ) $ the elements of $\cT$ listed in lexicographical order. The height function
$H(\cT)=(H_n ( \cT) ;0\leq n < \#(\cT)) $ is defined by 
$$ H_n (\cT)  =|u(n)|, \quad 0\leq n < \#(\cT) .$$
The height function is thus the sequence of the generations of the individuals of $\cT$, when these individuals
are
visited in the lexicographical order
(see Fig.1 for an example).  It is easy to check that $H(\cT) $ characterizes the tree $\cT$.

\begin{center}
\unitlength=1.1pt
\begin{picture}(300,140)

\thicklines \put(50,0){\line(-1,2){20}}
\thicklines \put(50,0){\line(1,2){20}}
\thicklines \put(30,40){\line(0,1){40}}
\thicklines \put(30,40){\line(-1,2){20}}
\thicklines \put(30,40){\line(1,2){20}}
\thicklines \put(30,80){\line(-1,2){20}}
\thicklines \put(30,80){\line(1,2){20}}
\put(40,-5){$\emptyset$}
\put(20,35){$1$}
\put(60,35){$2$}
\put(-3,75){$11$}
\put(17,75){$12$}
\put(37,75){$13$}
\put(-7,115){$121$}
\put(31,115){$122$}
\put(30,-25){tree $\cT$}
\thinlines \put(100,0){\line(1,0){120}}
\thinlines \put(100,0){\line(0,1){130}}
\thicklines \put(100,0){\line(1,5){8}}
\thicklines \put(108,40){\line(1,5){8}}
\thicklines \put(116,80){\line(1,-5){8}}
\thicklines \put(124,40){\line(1,5){8}}
\thicklines \put(132,80){\line(1,5){8}}
\thicklines \put(140,120){\line(1,-5){8}}
\thicklines \put(148,80){\line(1,5){8}}
\thicklines \put(156,120){\line(1,-5){8}}
\thicklines \put(164,80){\line(1,-5){8}}
\thicklines \put(172,40){\line(1,5){8}}
\thicklines \put(180,80){\line(1,-5){8}}
\thicklines \put(188,40){\line(1,-5){8}}
\thicklines \put(196,0){\line(1,5){8}}
\thicklines \put(204,40){\line(1,-5){8}}

\thinlines \put(108,0){\line(0,1){2}}
\thinlines \put(116,0){\line(0,1){2}}
\thinlines \put(124,0){\line(0,1){2}}
\thinlines \put(100,40){\line(1,0){2}}
\thinlines \put(100,80){\line(1,0){2}}
\thinlines \put(212,0){\line(0,1){2}}

\put(125,-25){contour function}
\put(240,-25){height function}
\thinlines\put(230,0){\line(1,0){80}}
\thinlines\put(230,0){\line(0,1){130}}
\thicklines \put(230,0){\line(1,4){10}}
\put(237.5,38){$\bullet$}
\thicklines \put(240,40){\line(1,4){10}}
\put(247.5,78){$\bullet$}
\thicklines\put(250,80){\line(1,0){10}}
\put(257.5,78){$\bullet$}
\thicklines\put(260,80){\line(1,4){10}}
\put(267.5,118){$\bullet$}
\thicklines\put(270,120){\line(1,0){10}}
\put(277.5,118){$\bullet$}
\thicklines\put(280,120){\line(1,-4){10}}
\put(287.5,78){$\bullet$}
\thicklines\put(290,80){\line(1,-4){10}}
\put(297.5,38){$\bullet$}

\thinlines \put(240,0){\line(0,1){2}}
\thinlines \put(250,0){\line(0,1){2}}
\thinlines \put(260,0){\line(0,1){2}}
\thinlines \put(300,0){\line(0,1){2}}
\thinlines \put(230,40){\line(1,0){2}}
\thinlines \put(230,80){\line(1,0){2}}

\put(107,-6){\scriptsize 1}
\put(115,-6){\scriptsize 2}
\put(123,-6){\scriptsize 3}
\put(208,-6){$\scriptstyle \zeta(\cT)$}
\put(239,-6){\scriptsize 1}
\put(249,-6){\scriptsize 2}
\put(259,-6){\scriptsize 3}
\put(290,-6){$\scriptstyle\#(\cT)-1$}

\put(95,39){\scriptsize 1}
\put(95,79){\scriptsize 2}
\put(225,39){\scriptsize 1}
\put(225,79){\scriptsize 2}

\end{picture}

\vskip 15mm

Figure 1
\end{center}

The {\bf contour function} gives another way of characterizing the tree, which is easier to visualize on a
picture (see Fig.1). Suppose that the tree is embedded in the half-plane in 
such a way that edges have length one. Informally, we imagine the motion of a
particle that starts at time 
$t=0$ from the root of the tree and then
explores the tree from the left to the right, moving continuously along the edges 
at unit speed, until it comes back to its starting point. Since it is clear
that each edge will be crossed twice in this evolution, the total time
needed to explore the tree is $\zeta(\cT):=2(\#(\cT)-1)$.
The value $C_t$ 
of the contour function at time $t$ is the distance (on the tree) between the position of the particle
at time $t$ and the
root. By convention $C_t=0$ if $t\geq \zeta(\cT)$.
Fig.1 explains the definition of the contour function better than a formal definition.

\section{Galton-Watson trees}

Let $\mu$ be a critical or subcritical offspring distribution. This means that
$\mu$ is a probability measure on $\Z_+$ such that
$$\sum_{k=0}^\infty k\mu(k)\leq 1.$$
We exclude the trivial case where $\mu(1)=1$.

There is a unique probability distribution ${\bf Q}_\mu$ on ${\bf T}$ such that

\begin{description}
\item{(i)} ${\bf Q}_\mu(k_\emptyset=j)=\mu(j)$,\quad
$j\in \Z_+$.

\item{(ii)} For every $j\geq 1$ with $\mu(j)>0$, the shifted trees
$\theta_1\cT,\ldots,\theta_j\cT$
are independent under the conditional probability
${\bf Q}_\mu(\cdot\mid k_\emptyset =j)$
and their conditional distribution is ${\bf Q}_\mu$.
\end{description}

A random tree with distribution ${\bf Q}_\mu$ is called a Galton-Watson tree
with offspring distribution $\mu$, or in short a $\mu$-Galton-Watson tree.

\smallskip
Let $\cT_1,\cT_2,\ldots$ be a sequence of independent $\mu$-Galton-Watson trees.
We can associate with this sequence a {\bf height process} obtained by concatenating the
height functions of each of the trees $\cT_1,\cT_2,\ldots$. More precisely,
for every $k\geq 1$, we set
$$H_n=H_{n-(\#(\cT_1)+\cdots+\#(\cT_{k-1}))}(\cT_k)\ \; {\rm if}\ \#(\cT_1)+\cdots+\#(\cT_{k-1})
\leq n<\#(\cT_1)+\cdots+\#(\cT_{k}).$$
The process $(H_n,n\geq 0)$ codes the sequence of trees. 

Similarly, we define 
a {\bf contour process} $(C_t,t\geq 0)$ coding the sequence of trees by concatenating
the contour functions $(C_t(\cT_1),t\in [0,\zeta(\cT_1)+2])$, $(C_t(\cT_2),t\in[0,\zeta(\cT_2)+2])$, etc.
Note that $C_t(\cT_n)=0$ for $t\in [\zeta(\cT_n),\zeta(\cT_n)+2]$, and that we are concatenating
the functions $(C_t(\cT_n),t\in [0,\zeta(\cT_n)+2])$ rather than the functions
$(C_t(\cT_n),t\in [0,\zeta(\cT_n)])$. This is a technical trick that will be useful in Chapter 2 below.
We may also observe that the process obtained by concatenating the functions $(C_t(\cT_n),t\in [0,\zeta(\cT_n)])$
would not determine the sequence of trees.

There is a simple relation between the height process and the contour process: See Section 2.4
in Chapter 2 for more details.

Although the height process is not a Markov process, except in very particular cases,
it turns out to be a simple functional of a Markov chain, which is even a random walk.
The next lemma is taken from \cite{LGLJ1},
but was obtained independently by other authors: See \cite{BV} and \cite{BK}.

\medskip
\noindent{\bf Lemma}
{\it
Let $\cT_1,\cT_2,\ldots$ be a sequence of independent $\mu$-Galton-Watson trees,
and let $(H_n,n\geq 0)$ be the associated height process. There exists a 
random walk $V$ on $\Z$ with initial value $V_0=0$ and jump distribution $\nu(k)=\mu(k+1)$, for
$k=-1,0,1,2,\ldots$, such that for every $n\geq 0$,
\begin{equation} 
\label{key}
H_n=\card\{k\in\{0,1,\ldots,n-1\}:V_k=\inf_{k\leq j\leq n}V_j\}.
\end{equation}
}

\medskip
A detailed proof of this lemma would be cumbersome, and we only
explain the idea. 
By definition, $H_n$ is the generation of the individual visited at time $n$,
for a particle that visits the different vertices of the
sequence of trees one tree after another and in lexicographical order
for each tree. Write $R_n$ for the quantity equal to the number of younger brothers
(younger means greater in the lexicographical order) of the individual visited at time $n$
plus the number of younger brothers of his father, plus the number of younger 
brothers of his grandfather etc. Then the random walk that appears
in the lemma may be defined by
$$V_n=R_n - (j-1)\qquad\hbox{if } \#(\cT_1)+\cdots+\#(\cT_{j-1})
\leq n<\#(\cT_1)+\cdots+\#(\cT_{j}).$$
To verify that $V$ is a random walk with jump distribution $\nu$, note that
because of the lexicographical order of visits, we have at time $n$ no
information on the fact that the individual visited at that time has children or not.
If he has say $k\geq 1$ children, which occurs with probability $\mu(k)$,
then the individual visited at time $n+1$ will be the first of these
children, and our definitions give
$R_{n+1}=R_n+(k-1)$ and
$V_{n+1}=V_n+(k-1)$. On the  other hand if he has no child, which occurs with probability $\mu(0)$,
then the individual visited at time $n+1$ is the first of the
brothers counted in the definition of $R_n$ (or the ancestor of the next
tree if $R_n=0$) and we easily see that
$V_{n+1}=V_n-1$. We thus get exactly the transition mechanism of the random walk with 
jump distribution $\nu$.

Let us finally explain formula (\ref{key}). From our definition of $R_n$
and $V_n$, it is easy to see that the condition $n<\inf\{j>k:V_j<V_k\}$
holds iff the individual visited at time $n$ is a descendant of the individual
visited at time $k$ (more precisely, $\inf\{j>k:V_j<V_k\}$ is the time of the first 
visit after $k$ of an individual that is not a descendant of individual $k$).
Put in a different way, the condition $V_k=\inf_{k\leq j\leq n}V_j$
holds iff the individual visited at time $k$ is an ascendant of the individual
visited at time $n$. It is now clear that the right-hand side of (\ref{key})
just counts the number of ascendants of the individual visited at time $n$,
that is the generation of this individual.

\section{The continuous height process}

To define the height process in a continuous setting, we use an analogue of
the discrete formula (\ref{key}). The role of the random walk $V$
in this formula is played by a L\'evy process $X=(X_t,t\geq 0)$ without negative jumps.
We assume that $X$ does not drift to $+\infty$
(this corresponds to the subcriticality of $\mu$ in the discrete setting), and that
the paths of $X$ are of infinite variation a.s.: The latter assumption
implies in particular that the process $X$
started at the origin will immediately hit both $(0,\infty)$ and $(-\infty,0)$.
The law of $X$ can be characterized by its Laplace functional $\psi$, which 
is the nonnegative function on $\R_+$ defined by
$$E[\exp(-\lambda X_t)]=\exp(t\psi(\lambda)).$$
By the L\'evy-Khintchine formula and our special assumptions on $X$, the function
$\psi$ has to be of the form
$$\psi(\lambda)=\alpha\lambda+\beta\lambda^2+\int \pi(dr)\,(e^{-\lambda r}-1+\lambda r),$$
where $\alpha,\beta\geq 0$ and $\pi$ is a $\sigma$-finite mesure on $(0,\infty)$
such that $\int \pi(dr)(r\wedge r^2)<\infty$. We write
$$S_t=\sup_{s\leq t} X_s\ ,\qquad I_t=\inf_{s\leq t}X_s\,.$$

By analogy with the discrete case, we would like to define $H_t$ as the ``measure''
of the set
\begin{equation} 
\label{minimaset}
\{s\leq t:X_s=\inf_{s\leq r\leq t} X_r\}.
\end{equation}
However, under our assumptions on $X$, the Lebesgue measure of this set
is always zero, and so we need to use some sort of local time. The key idea is
to introduce for every fixed $t>0$ the time-reversed process
$$\wh X^{(t)}_s=X_t-X_{(t-s)-}\ ,\quad 0\leq s\leq t\,,$$
and its associated supremum
$$\wh S^{(t)}_s=\sup_{r\leq s} \wh X^{(t)}_r\,.$$

We observe that via time-reversal $s\rightarrow  t-s$, the set (\ref{minimaset})
corresponds to $\{s\leq t:\wh S^{(t)}_s=\wh X^{(t)}_s\}$. This leads to the
rigorous definition of $H$: $H_t$ is defined as the local time at level $0$,
at time $t$ of the process $\wh S^{(t)}-\wh X^{(t)}$. This definition makes sense
because $\wh S^{(t)}-\wh X^{(t)}$ has the same law over $[0,t]$ as the 
so-called reflected process $S-X$ for which $0$ is a regular point under our assumptions.
Note that the normalization of local time has to be specified in some way:
See Section 1.1. The process $(H_t,t\geq 0)$ is called the 
$\psi$-{\bf height process}, or simply the height process.

Why is the $\psi$-height process $H$ an interesting object of study ? In the
same way as the discrete height process codes the genealogy of a sequence of 
independent Galton-Watson trees, we claim that the continuous height process represents the
genealogical structure of continuous-state branching processes, which
are the continuous analogues of Galton-Watson processes. This informal claim is at the heart
of the developments of the present work. Perhaps the best justification for it can be
found in the limit theorems of Chapter 2 that relate the discrete and continuous height
processes (see Section 0.4 below). Another justification is the Ray-Knight theorem for the
height process that will be discussed below.

The goal of Chapter 1 is to
present a self-contained construction and to derive several new properties
of the $\psi$-height process. Although there is some overlap with \cite{LGLJ1},
our approach is different and involves new approximations. It is important to realize that
$H_t$ is defined as the local time at time $t$ of a process which itself depends on $t$. 
For this reason, it is not clear whether the paths of $H$ have any regularity properties.
Also $H$ is not Markov, except in the very special case where $X$ has no jumps.
To 
circumvent these difficulties, we rely on the important tool of the 
exploration process: For every $t\geq 0$, we define a random measure
$\rho_t$ on $\R_+$ by setting
\begin{equation} 
\label{Irho}
\langle \rho_t,f\rangle=\int_{[0,t]} d_sI^s_t\,f(H_s)
\end{equation}
where 
$$I^s_t=\inf_{s\leq r\leq t} X_r$$
and the notation $d_sI^s_t$ refers to integration with respect to the nondecreasing function 
$s\rightarrow  I^s_t$. The {\bf exploration process} $(\rho_t,t\geq 0)$ is a Markov process with values in the
space $M_f(\R_+)$ of finite measures on $\R_+$. It was introduced and studied
in \cite{LGLJ1}, where its definition was motivated by a 
model of a LIFO queue (see \cite{Li} for some applications to queuing theory).

The exploration process has several interesting properties.
In particular it is c\` adl\` ag (right-continuous with left limits) and it has an explicit
invariant measure in terms of the subordinator with Laplace exponent $\psi(\lambda)/\lambda$
(see Proposition \ref{mesinv}). 
Despite its apparently complicated definition, the exploration process is the crucial
tool that makes it possible to answer most questions concerning the height process. A first illustration of
this is the choice of a ``good'' lower-semicontinuous modification of $H_t$, which is obtained by considering for
every 
$t\geq 0$ the supremum of the support of the measure $\rho_t$ (beforehand, to make
sense of the definition of $\rho_t$, one needs to use a first version of $H$ that can be defined by suitable
approximations of local times).

An important feature of both the height process and the 
exploration process is the fact that both $H_t$ and $\rho_t$ depend only on the values
of $X$, or of $X-I$, on the excursion interval of $X-I$ away from $0$ that straddles $t$. 
For this reason, it is possible to define and to study both the height process and the 
exploration process under the excursion measure of $X-I$ away from $0$. This excursion
measure, which is denoted by $N$, plays a major role throughout this work, and many
results are more conveniently stated under $N$. Informally, the height process under $N$
codes exactly one continuous tree, in the same way as each excursion away from $0$
of the discrete height process corresponds to one Galton-Watson tree in the sequence (cf Section 0.2).

 As a typical application of the exploration process, we introduce and study the
local times of the height process, which had not been considered in earlier
work. These local times play an important role in the sequel, in particular
in the applications to spatial branching processes.
The local time of $H$ at level $a\geq 0$ and at time $t$ is denoted by $L^a_t$ and
these local times can be defined through
the approximation
$$\lim_{\varepsilon\rightarrow  0}E\Big[\sup_{s\leq t}\Big|\varepsilon^{-1}\int_0^s 1_{\{a<H_r<a+\varepsilon\}}dr
-L^a_s\Big|\Big]=0$$
(Proposition \ref{LTapprox}). The proof of this approximation depends in a crucial
way on properties of the exploration process derived in Section 1.3: Since $H$ is in 
general not Markovian nor a semimartingale, one cannot use the standard methods of construction
of local time. 

The Ray-Knight theorem for the height process states that if $T_r=\inf\{t\geq 0:X_t=-r\}$,
for $r>0$, the process $(L^a_{T_r},a\geq 0)$ is a continuous-state branching process with
branching mechanism $\psi$ (in short a $\psi$-CSBP) started at $r$. 
Recall that 
the $\psi$-CSBP is the Markov process $(Y_a,a\geq 0)$ with values in $\R_+$
whose transition kernels are characterized by their
Laplace transform: For $\lambda>0$ and $b>a$,
$$E[\exp-\lambda Y_b\mid Y_a]=\exp(-Y_a\,u_{b-a}(\lambda)),$$
where $u_t(\lambda)$, $t\geq 0$ is the unique nonnegative solution
of the differential equation
$$
{\partial u_t(\lambda)\over \partial t}=-\psi(u_t(\lambda))\ ,\quad u_0(\lambda)=\lambda.
$$
By analogy with the discrete setting, we can think of 
$L^a_{T_r}$ as ``counting'' the number of individuals at generation $a$ in
a Poisson collection of continuous trees (those trees coded by the excursions of 
$X-I$ away from $0$ before time $T_r$). The Ray-Knight theorem corresponds to the intuitive fact that the population
at generation $a$ is a branching process.

The previous Ray-Knight theorem had already been derived in \cite{LGLJ1} although in a less precise
form (local times of the height process had not been constructed). An important consequence of the Ray-Knight
theorem, also derived in \cite{LGLJ1}, is a criterion for the path continuity of $H$: $H$
has continuous sample paths iff
\begin{equation} 
\label{contcond}
\int_1^\infty {d\lambda\over \psi(\lambda)}<\infty.
\end{equation}
This condition is in fact necessary and sufficient for the a.s. extinction of the $\psi$-CSBP. 
If it does not hold, the paths of $H$ have a very wild behavior: The values of $H$ over
any nontrivial interval $[s,t]$ contain a half-line $[a,\infty)$. On the other hand, 
(\ref{contcond}) holds if $\beta>0$, and in the stable case $\psi(\lambda)=c\lambda^\gamma$,
$1<\gamma\leq 2$ (the values $\gamma\in(0,1]$ are excluded by our assumptions).

In view of applications in Chapter 4, we derive precise information about the H\" older continuity of
$H$. We show that if
$$\gamma=\sup\{r\geq 0:\lim_{\lambda\rightarrow  \infty} \lambda^{-r}\psi(\lambda)=+\infty\},$$
then the height process $H$ is a.s. H\"older continuous with exponent $r$
for any $r\in(0,1-\gamma^{-1})$, and a.s. not H\" older continuous with exponent $r$
if $r>1-\gamma^{-1}$.

\section{From discrete to continuous trees}

Chapter 2 discusses limit theorems for rescaled Galton-Watson trees. These results demonstrate
that the $\psi$-height process is the correct continuous analogue of the discrete height process
coding Galton-Watson trees. 

It is well known \cite{La1} that continuous-state branching processes
are the only possible scaling limits of discrete-time Galton-Watson branching processes.
One may then ask for finer limit theorems involving the genealogy. 
Precisely, starting from a sequence of rescaled Galton-Watson processes
that converge in distribution towards a continuous-state branching process,
can one say that the corresponding discrete Galton-Watson trees also converge,
in some sense, towards a continuous genealogical structure ? The results of
Chapter 2 show that the answer is yes.

To be specific, consider a sequence $(\mu_p)$ of (sub)critical offspring distributions. For every
$p\geq 1$, let $Y^p$ be a (discrete-time) Galton-Watson process with offspring distribution $\mu_p$
started at $Y^p_0=p$. Suppose that the processes $Y^p$ converge after rescaling towards a
$\psi$-CSBP, where $\psi$ satisfies the conditions introduced in Section 0.3. Precisely, we assume that
there is a sequence $\gamma_p\uparrow \infty$ such that
\begin{equation} 
\label{convresGW}
(p^{-1}Y^p_{[\gamma_pt]},t\geq 0)
\build{\la}_{p\rightarrow \infty}^{\rm (d)} (Y_t,t\geq 0),
\end{equation}
where $Y$ is a $\psi$-CSBP, and the symbol (d) indicates convergence in distribution in the
Skorokhod space. Let $H^p$ be the discrete height process associated with $\mu_p$ in the sense
of Section 0.2. Then Theorem \ref{marginconv} shows that
\begin{equation} 
\label{convHmargin}
(\gamma_p^{-1}H^p_{[p\gamma_pt]},t\geq 0)\build{\la}_{p\rightarrow \infty}^{\rm (fd)} (H_t,t\geq 0),
\end{equation}
where $H$ is the $\psi$-height process and (fd) indicates convergence of finite-dimensional marginals.
A key ingredient of the proof is the observation due to Grimvall \cite{Gr} that the convergence
(\ref{convresGW}) implies the convergence in distribution (after suitable rescaling) of the random walks 
$V^p$ with jump distribution $\nu_p(k)=\mu_p(k+1)$, $k=-1,0,1,\ldots$, towards the L\'evy process with
Laplace exponent $\psi$. The idea is then to pass to the limit in the formula for $H^p$ in terms of $V^p$,
recalling that the $\psi$-height process is given by an analogous formula in terms of the L\'evy process $X$.
In the special case $\beta=0$ and under more restrictive assumptions, the convergence (\ref{convHmargin})
had already appeared in \cite{LGLJ1}.

In view of applications, the limit (\ref{convHmargin}) is not satisfactory because the convergence
of finite-dimensional marginals is too weak. In order to reinforce (\ref{convHmargin}) to a 
functional convergence, it is necessary to assume some regularity of the paths of $H$. We assume that condition
(\ref{contcond}) ensuring the path continuity of $H$ holds (recall that if this condition does not
hold, the paths of $H$ have a very wild behavior). Then, we can prove (Theorem \ref{functionalconv})
that the convergence (\ref{convHmargin}) holds in the sense of weak convergence on the
Skorokhod space, provided that the following condition is satisfied: For every $\delta>0$
\begin{equation} 
\label{Itech}
\liminf_{p\rightarrow \infty}P[Y^p_{[\delta\gamma_p]}=0]>0.
\end{equation}
Roughly speaking this means that the rescaled Galton-Watson process $(p^{-1}Y^p_{[\gamma_pt]},t\geq 0)$
may die out at a time of order $1$, as its weak limit $Y$ does (recall that we are
assuming (\ref{contcond})). The technical condition (\ref{Itech}) is both necessary and sufficient for the
reinforcement of (\ref{convHmargin}) to a functional convergence. Simple examples show that this condition cannot be
omitted in general.

However, in the important special case where $\mu_p=\mu$ for every $p$, we are able to show 
(Theorem \ref{stableconv}) that the technical condition (\ref{Itech}) is always satisfied . In that case,
$\psi$ must be of the form $\psi(u)=c\,u^\gamma$ with $1<\gamma\leq 2$, so that obviously (\ref{contcond})
also holds. Thus when $\mu_p=\mu$ for every $p$, no extra condition is needed to get a functional
convergence. 

In Section 2.4, we show that the functional convergence derived for rescaled discrete height processes 
can be stated as well in terms of the contour processes (cf Section 0.1). Let 
$C^p=(C^p_t,t\geq 0)$ be the contour process for a sequence of independent $\mu_p$-Galton-Watson trees.
Under the assumptions that warrant the functional convergence in (\ref{convHmargin}), Theorem
\ref{contour-conv} shows that we have also
$$
(p^{-1}C^p_{p\gamma_pt},t\geq 0)\build{\la}_{p\rightarrow \infty}^{\rm (d)} (H_{t/2},t\geq 0).
$$
Thus scaling limits are the same for the discrete height process and for the contour process.

In the remaining part of Chapter 2, we give applications of (\ref{convHmargin}) assuming that the
functional convergence holds. In particular, rather than considering a sequence of $\mu_p$-Galton-Watson
trees, we discuss the height process associated with a single tree conditioned to be large. 
Precisely, let $\wt H^p$ be the height process for one $\mu_p$-Galton-Watson tree conditioned to
non-extinction at generation $[\gamma_pT]$, for some fixed $T>0$. Then, Proposition \ref{tree-condit-height}
gives
$$(\gamma_p^{-1}\wt H^p_{[p\gamma_pt]},t\geq 0)\build{\la}_{p\rightarrow \infty}^{\rm (d)} (\wt H_t,t\geq 0),$$
where the limiting process is an excursion of the $\psi$-height process conditioned to hit level $T$.
This is of course reminiscent of a result of Aldous \cite{Al2} who proved that in the case of
a critical offspring distribution $\mu$ with finite variance, the contour process
of a $\mu$-Galton-Watson tree conditioned to have exactly $p$ vertices converges after a suitable rescaling
towards a normalized Brownian excursion (see also 
\cite{Git} and \cite{MaMo} for related results including the
convergence of the height process in Aldous' setting). Note that in Aldous' result, the conditioning becomes 
degenerate in the limit, since the ``probability'' that a Brownian excursion has length exactly one is zero.
This makes it more difficult to derive this result from our approach, although it seems very
related to our limit theorems. See however Duquesne \cite{Du2} for an extension of Aldous' theorem
to the stable case
using the tools of the present work (a related result in the stable case was obtained 
by Kersting \cite{Ker}).

The end of Chapter 2 is devoted to reduced trees. We consider again a single Galton-Watson tree conditioned 
to non-extinction at generation $[\gamma_pT]$. For every $k<[\gamma_pT]$, we denote by 
$Z^{(p),[\gamma_pT]}_k$ the number of vertices at generation $k$ that have descendants at generation
$[\gamma_pT]$. Under the assumptions and as a consequence of Proposition \ref{tree-condit-height},
we can prove that
$$(Z^{(p),[\gamma_pT]}_{[\gamma_pt]},0\leq t<T)\build{\la}_{p\rightarrow \infty}^{\rm (fd)} (Z^T_t,0\leq t<T)$$
where the limit $Z^T$ has a simple definition in terms of $\wt H$: $Z^T_t$ is the number of
excursions of $\wt H$ above level $t$ that hit level $T$. Thanks to the properties of the 
height process and the exploration process that have been derived in Chapter 1, it is possible to calculate the
distribution of the time-inhomogeneous branching process $(Z^T_t,t\geq 0)$. This distribution is derived in
Theorem \ref{lawreduced}. Of course in the stable case, corresponding to $\mu_p=\mu$ for every $p$, the
distribution of $Z^T$ had been computed previously. See in particular Zubkov \cite{Zu}
and Fleischmann and Siegmund-Schultze \cite{FlSG}.

\section{Duality properties of the exploration process}

In the applications developed in Chapters 3 and 4, a key role is played by the
duality properties of the exploration process $\rho$. We first observe that
formula (\ref{Irho}) defining the exploration process can be rewritten in the 
following equivalent way
$$\rho_t(dr)=\beta 1_{[0,H_t]}(r)\,dr+\sum_{s\leq t,X_{s-}<I^s_t}(I^s_t-X_{s-})\delta_{H_s}(dr)$$
where $\delta_{H_s}$ is the Dirac measure at $H_s$, and we recall that 
$I^s_t=\inf_{s\leq r\leq t}X_r$. We then define another measure $\eta_t$ by setting
$$\eta_t(dr)=\beta 1_{[0,H_t]}(r)\,dr+\sum_{s\leq t,X_{s-}<I^s_t}(X_s-I^s_t)\delta_{H_s}(dr).$$
To motivate this definition, we may come back to the discrete setting of Galton-Watson trees. 
In that setting, the discrete height process $H_n$ gives the generation of the 
$n$-th visited vertex by a ``particle'' that visits vertices in lexicographical order one tree
after another, and the analogue of $\rho_t$ gives for every $k\leq H_n$ the number of younger
(i.e. coming after in the lexicographical order) brothers of the ancestor at generation $k$
of the $n$-the visited vertex. Then the analogue of $\eta_t$ gives for every $k\leq H_n$
the number of older brothers of the ancestor at generation $k$
of the $n$-the visited vertex.

It does not seem easy to study directly the Markovian properties or the regularity of paths of the
process $(\eta_t,t\geq 0)$. The right point of view is to
consider the pair $(\rho_t,\eta_t)$, which is easily seen to be a Markov process in
$M_f(\R_+)^2$. The process $(\rho_t,\eta_t)$ has an invariant measure $\M$ determined 
in Proposition \ref{invariant-rho-eta}. The key result (Theorem \ref{duality-rho})
then states that the Markov processes $(\rho,\eta)$ and $(\eta,\rho)$ are in duality under $\M$. 
A consequence of this is the fact that $(\eta_t,t\geq 0)$ also has a 
c\` adl\` ag modification. More importantly, we obtain a crucial time-reversal property: Under the excursion
measure $N$ of $X-I$, the processes $(\rho_s,\eta_s;0\leq s\leq \sigma)$ and $(\eta_{(\sigma-s)-},\rho_{(\sigma-s)-};
0\leq s\leq \sigma)$ have the same distribution (here $\sigma$ stands for the duration of the excursion under
$N$). This time-reversal property plays a major role in many
subsequent calculations. It implies in particular that the law of $H$ under $N$ is invariant under
time-reversal. This property is natural in the discrete setting, if we think of the contour 
process of a Galton-Watson tree, but not obvious in the continuous case.

\section{Marginals of trees coded by the height process}

Let us explain more precisely how an excursion of the $\psi$-height process codes
a continuous branching structure. We consider first a deterministic continuous function
$e:\R_+\la \R_+$ such that $e(t)>0$ iff $0<t<\sigma$, for some $\sigma=\sigma(e)>0$. 
For any $s,s'\geq 0$, set
$$m_e(s,s')=\inf_{s\wedge s'\leq t\leq s\vee s'}e(t).$$
Then $e$
codes a continuous genealogical structure via the following simple prescriptions:

\begin{description}
\item{(i)} To each $s\in[0,\sigma]$ corresponds a vertex at generation $e(s)$.

\item{(ii)} Vertex $s$ is an ancestor of vertex $s'$ if $e(s)=m_e(s,s')$. In general, 
$m_e(s,s')$ is the
generation of the last common ancestor to $s$ and $s'$.

\item{(iii)} We put $d(s,s')=e(s)+e(s')-2m_e(s,s')$ and identify $s$ and $s'$ 
($s\sim s'$) if $d(s,s')=0$.

\end{description}

Formally, the tree coded by $e$ can be defined as the quotient set $[0,\sigma]/\sim$,
equipped with the distance $d$ and the genealogical relation specified in (ii). 

With these definitions, the line of ancestors of a vertex $s$ is isometric to the segment
$[0,e(s)]$. If we pick two vertices $s$ and $s'$, their lines of ancestors share a common
part
isometric to $[0,m_e(s,s')]$, and then become distinct. In general, if we consider
$p$ instants $t_1,\ldots,t_p$ with $0\leq t_1\leq\cdots\leq t_p\leq \sigma$, we can associate
with these $p$ instants a genealogical tree $\theta(e,t_1,\ldots,t_p)$, which consists
of a discrete rooted ordered tree with $p$ leaves, denoted by $\cT(e,t_1,\ldots,t_p)$
and marks $h_v(e,t_1,\ldots,t_p)\geq 0$ for $v\in\cT(e,t_1,\ldots,t_p)$, that correspond to the
lifetimes of vertices in $\cT (e,t_1,\ldots,t_p)$. See subsection 3.2.1 for a precise definition.

In the second part of Chapter 3, we use the duality results proved in the first part to
calculate the distribution of the tree $\theta(H,\tau_1,\ldots,\tau_p)$ under certain
excursion laws of $H$ and random choices of the instants $\tau_1,\ldots,\tau_p$. We assume that the
continuity condition (\ref{contcond}) holds. We first consider Poissonnian marks with intensity $\lambda$,
and the height process $H$ under the excursion measure $N$ of $X-I$. Let $\tau_1,\ldots,\tau_M$
be the marks that fall into the duration interval $[0,\sigma]$ of the excursion. Theorem
\ref{tree-Poisson} shows that under the probability measure $N(\cdot\mid M\geq 1)$, the tree
$\theta(H,\tau_1,\ldots,\tau_M)$ is distributed as the family tree of a continuous-time 
Galton-Watson process starting with one individual at time $0$ and where

\begin{description}
\item{$\bullet$} lifetimes have exponential distribution with parameter $\psi'(\psi^{-1}(\lambda))$;
\item{$\bullet$} the offspring distribution is the law of the variable $\xi$ with generating function
$$E[r^\xi]=r+{\psi((1-r)\psi^{-1}(\lambda))\over \psi^{-1}(\lambda)\psi'(\psi^{-1}(\lambda))}.$$
\end{description}

In the quadratic case, we get a critical binary branching $E[r^\xi]={1\over 2}(1+r^2)$. The result in that
case had been obtained by Hobson \cite{Ho}. 

We finally specialize to the stable case $\psi(\lambda)=\lambda^\gamma$, $\gamma\in(1,2]$. 
By scaling arguments, we can then
make sense of the law $N_{(1)}=N(\cdot \mid \sigma=1)$ of the normalized excursion of $H$.
Using the case of Poissonnian marks, we compute explicitly the law of the tree
$\theta(H,t_1,\ldots,t_p)$ under $N_{(1)}$, when $(t_1,\ldots,t_p)$ are
chosen  independently and uniformly over $[0,1]^p$. In the quadratic case $\psi(u)=u^2$, 
$H$ is under $N_{(1)}$ a normalized Brownian excursion, and the corresponding tree
is called the continuum random tree (see Aldous \cite{Al1},\cite{Al91},\cite{Al2}). By analogy,
in our more general case $\psi(u)=u^\gamma$, we may call the tree coded by $H$ under $N_{(1)}$
the stable continuum random tree. Our calculations give what Aldous calls the finite-dimensional
marginals of the tree. In the case $\gamma=2$, these marginals were computed by Aldous
(see also Le Gall \cite{LG99} for a different approach closer to the present work).
In that case, the discrete skeleton $\cT(H,t_1,\ldots,t_p)$ is uniformly distributed over all
binary rooted ordered trees with $k$ leaves. When $\gamma<2$, things become different as we can get
nonbinary trees (the reason why we get only binary trees in the Brownian case is the fact that
local minima of Brownian motion are distinct). Theorem \ref{stable-marginal} shows in particular
that if $\cT$ is a tree with $p$ leaves such that $k_u(\cT)\not =1$ for every $u\in\cT$
(this condition must be satisfied by our trees $\cT(e,t_1,\ldots,t_p)$) then the probability
that $\cT(H,t_1,\ldots,t_p)=\cT$ is
$${p!\over (\gamma-1)(2\gamma-1)\cdots((p-1)\gamma-1)}\
\prod_{v\in \n_\cT}{|(\gamma-1)(\gamma-2)\cdots(\gamma-k_v+1)|\over k_v!}$$
where $\n_\cT=\{v\in \cT:k_v>0\}$ is the set of nodes of $\cT$. It would be interesting to know whether
this distribution on discrete trees has occurred in other settings.

\section{The L\'evy snake}

Chapters 1 -- 3 explore the continuous genealogical structure coded by the $\psi$-height process $H$.
In Chapter 4, we examine the probabilistic objects obtained by combining this branching structure
with a spatial motion given by a c\` adl\` ag Markov process $\xi$ with state space $E$. Informally,
``individuals'' do not only reproduce themselves, but they also move in space independently
according to the law of $\xi$. 
The $(\xi,\psi)$-superprocess is then
a Markov process taking values in the space of finite measures on $E$,
whose value at time $t$ is a random measure putting mass on the
set of positions of ``individuals'' alive at time $t$. Note that the
previous description is very informal since in the continuous branching
setting there are no individual particles but rather a continuum of infinitesimal
particles. Recent accounts of the theory of superprocesses can be found in
Dynkin \cite{Dy01}, Etheridge \cite{Et} and Perkins \cite{Pe}.

Our coding of the genealogy by the height process leads to
introducing a Markov process whose values will give the 
historical paths followed by the ``individuals'' in the population. This a generalization
of the L\'evy snake introduced in \cite{LG93} and studied in particular in \cite{LG99}.
To give a precise definition, fix a starting point $x\in E$, consider the $\psi$-height process 
$(H_s,s\geq 0)$ and recall the notation $m_H(s,s')=\inf_{[s,s']}H_r$ for $s\leq s'$. We assume 
that the continuity condition (\ref{contcond}) holds. Then conditionally on
$(H_s,s\geq 0)$ we consider a time-inhomogeneous Markov process $(W_s,s\geq 0)$
whose distribution is described as follows:

\begin{description}
\item{$\bullet$} For every $s\geq 0$, $W_s=(W_s(t),0\leq t<H_s)$ is a path of $\xi$
started at $x$ and with finite lifetime $H_s$.

\item{$\bullet$} If we consider two instants $s$ and $s'$, the corresponding paths
$W_s$ and $W_{s'}$ are the same up to time $m_H(s,s)$ and then behave independently.

\end{description}

The latter property is consistent with the fact that in our coding of the
genealogy, vertices attached to $s$ and $s'$ have the same ancestors up to
generation $m_H(s,s')$. See Section 4.1 for a more precise definition.

The pair $(\rho_s,W_s)$ is then a Markov process with values in the product space
$M_f(\R_+)\times \W$, where $\W$ stands for the set of all finite c\` adl\` ag paths in $E$. 
This process is called the {\bf L\'evy snake} (with initial point $x$). It was introduced and studied in
\cite{LGLJ2}, where a  form of its connection with superprocesses was established. Chapter 4 gives much more
detailed information about its properties. In particular, we prove the strong Markov property 
of the L\'evy snake (Theorem \ref{strongMarkovsnake}), which plays a crucial role in
several applications.
 
We also use the local times of the
height process to give a nicer form of the connection with superprocesses. Write $\wh W_s$
for the left limit of $W_s$ at its lifetime $H_s$ (which exists a.s. for each fixed $s$), and recall the
notation
$T_r=\inf\{t\geq 0:X_t=-r\}$. For every $t\geq 0$, we can define a random measure $Z_t$ on $E$ by setting
$$\langle Z_t,\varphi\rangle=\int_0^{T_r} d_sL^t_s\,\varphi(\wh W_s)\,.$$
Then $(Z_t,t\geq 0)$ is a $(\xi,\psi)$-superprocess with initial value $r\delta_x$.
This statement is in fact a special case of Theorem \ref{super} which constructs  
a $(\xi,\psi)$-superprocess with an arbitrary initial value. For this more general statement,
it is necessary to use excursion measures of the L\'evy snake: Under the {\bf excursion measure} $\N_x$,
the process $(\rho_s,s\geq 0)$ is distributed according to its excursion measure $N$,
and $(W_s,s\geq 0)$ is constructed by the procedure explained above, taking 
$x$ for initial point.

As a second application, we use local time techniques to construct exit measures
from an open set and to establish the
integral equation satisfied by the Laplace functional of exit measures
(Theorem \ref{integrexit}). Recall that exit
measures of superprocesses play a fundamental role in the
connections with partial differential equations studied recently by Dynkin and Kuznetsov (a detailed account of
these connections can be found in the forthcoming book
\cite{Dy01}). 

We then study the continuity of the path-valued process $W_s$
with respect to the uniform topology on paths. This question is closely related 
to the compact support property for superprocesses. In the case when $\xi$
is Brownian motion in $\R^d$, Theorem \ref{compactness-Brownian} shows that the condition
$$\int_1^\infty \Big(\int_0^t \psi(u)\,du\Big)^{-1/2} dt<\infty$$
is necessary and sufficient
for $W_t$ to be continuous with respect to the uniform topology on paths. The proof relies on 
connections of the exit measure with partial differential equations and 
earlier work
of Sheu \cite{Sh1}, who was interested in the compact support
property for superprocesses. More generally, assuming only that $\xi$ has
H\"older continuous paths, we use the continuity properties of $H$ derived in Chapter 1 to
give simple sufficient conditions ensuring that the same conclusion holds.
 
Although we do not develop such applications in the present work, we expect that the L\'evy snake will be a powerful
tool to study connections with partial differential equations, in the spirit of \cite{LG95}, as well as 
path properties of superprocesses (see \cite{LGPer1995} for a typical application
of the Brownian snake to super-Brownian motion). 

In the last two sections of Chapter 4, we compute certain 
explicit distributions related to the L\'evy snake and the 
$(\xi,\psi)$-superprocess, under the excursion
measures $\N_x$. We assume that the path-valued process $W_s$
is continuous with respect to the uniform topology on paths, and then
the value $W_s(H_s)$ can be defined as a left limit at the lifetime,
simultaneously for all $s\geq 0$. If $D$ is an open set in $E$
such that $x\in D$, we consider the first exit time
$$T_D=\inf\{s\geq 0:\tau(W_s)<\infty\}$$
where $\tau(W_s)=\inf\{t\in[0,H_s]:W_s(t)\notin D\}$. Write 
$u(y)=\N_y(T_D<\infty)<\infty$ for every $y\in D$. Then the distribution
of $W_{T_D}$ under $\N_x(\cdot\cap\{T_D<\infty\})$ is characterized
by the function $u$ and the distribution $\Pi_x$
of $\xi$ started at $x$ via the formula: For every $a\geq 0$
\ba
&&\N_x\Big(1_{\{T_D<\infty\}}1_{\{a<H_{T_D}\}}F(W_{T_D}(t),0\leq t\leq a)\Big)\\
&&\qquad=\Pi_x\Big[1_{\{a< \tau\}}u(\xi_a)F(\xi_r,0\leq r\leq a)\exp\Big(-\int_0^a
\wt\psi(u(\xi_r))dr\Big)\Big],
\ea
where $\tau$ stands for the first exit time of $\xi$ from $D$, and $\wt \psi(r)=\psi(r)/r$.
Theorem \ref{law-exit} gives more generally the law of the pair $(W_{T_D},\rho_{T_D})$
under $\N_x(\cdot\cap\{T_D<\infty\})$. In the special case when $\xi$
is Brownian motion in $\R^d$, the function $u$ can be identified as 
the maximal 
nonnegative solution of ${1\over 2}\Delta u=\psi(u)$ in $D$, and the law of 
$W_{T_D}$ is that of a Brownian motion with drift ${\nabla u/u}$ up to its
exit time from $D$. This considerably extends a result of \cite{LG94}
proved in the quadratic branching case by a very different method.

The last section of Chapter 4 investigates reduced spatial trees,
again under the assumption that the path-valued process $W_s$
is continuous with respect to the uniform topology on paths. We consider 
a spatial open set $D$ with $x\in D$, and the L\'evy snake under its excursion
measure $\N_x$ (in the superprocess setting this means that we are looking at all
historical paths corresponding to one ancestor at time $0$). 
We condition on the event that $\{T_D<\infty\}$, that is one at least of the paths
$W_s$ exits $D$, and we want to describe the spatial structure of all the paths that
exit $D$, up to their respective exit times. This is an analogue (and in fact a generalization)
of the reduced tree problem studied in Chapter 2. In the spatial situation, all paths $W_s$
that exit $D$ will be the same up to a certain time $m_D$ at which there is a branching
point with finitely many branches, each corresponding to an excursion of the
height process $H$ above level $m_D$, in which the L\'evy snake exits $D$.  
In each such excursion the paths $W_s$ that exit $D$ will be the same up to a level 
strictly greater than $m_D$, at which there is another branching point, and so on.

To get a full description of the reduced spatial tree, one only needs to compute the
joint distribution of the path $W^D_0=W_{T_D}(.\wedge m_D)$, that is the common
part to all paths that do exit $D$, and the number $N_D$ of branches at the first branching point.
Indeed, conditionally on the pair $(W^D_0,N_D)$, the ``subtrees of paths'' that originate
from the first branching point will be independent and distributed according to the 
full reduced tree with initial point $\wh W^D_0=W^D_0(m_D)$ (see Theorem \ref{reduced-domain}
for more precise statements). Theorem \ref{reduced-domain} gives explicit formulas for the 
joint distribution of $(W^D_0,m_D)$, again in terms of the function $u(y)=\N_y(T_D<\infty)<\infty$.
Precisely, the law of the ``first branch'' $W^D_0$ is given by
\ba
&&\N_x(1_{\{T_D<\infty\}}F(W^D_0))\\
&&\ =\int_0^\infty db\,\Pi_x\Big[1_{\{b<\tau\}}u(\xi_b)\,\theta(u(\xi_b))\,
\exp\Big(-\int_0^b \psi'(u(\xi_r))dr\Big) F(\xi_r,0\leq r\leq b)\Big],
\ea
where $\theta(r)=\psi'(r)-\wt \psi(r)$. Furthermore the conditional distribution
of $N_D$ given $W^D_0$ depends only on the branching point $\wh W^D_0$ and is given by
$$\N_x[r^{N_D}\mid T_D<\infty, W^D_0]=r\,{\psi'(U)-\gamma_\psi(U,(1-r)U)\over
\psi'(U)-\gamma_\psi(U,0)}\ ,\qquad 0\leq r\leq 1,$$
where $U=u(\wh W^D_0)$ and $\gamma_\psi(a,b)={\psi(a)-\psi(b)\over a-b}$. In the
stable case
$\psi(u)=u^\gamma$, the variable $N_D$ is independent of $W^D_0$ and its generating function
is $(\gamma-1)^{-1}((1-r)^\gamma-1+\gamma r)$.

\medskip
\noindent{\it Acknowledgment}. We would like to thank Yves Le Jan for allowing
us to use several ideas that originated in his work in collaboration
with one of us.

\chapter{The height process}

\section{Preliminaries on L\'evy processes}

\subsection{Basic assumptions}

In this section, we introduce the class of L\'evy 
processes that will be relevant to our study and
we record some of their basic properties. 
For almost all facts about L\'evy processes that we need, we
refer to the recent book of Bertoin \cite{Be} (especially
Chapter VII). 

We consider a L\'evy process $X$ on the real 
line. It will be convenient to assume that
$X$ is the canonical process on the Skorokhod
space $\D(\R_+,\R)$ of c\` adl\` ag (right-continuous
with left limits) real-valued paths. The canonical filtration
will be denoted by $(\g_t,t\in[0,\infty])$. Unless otherwise noted,
the underlying probability measure $P$ is the law of the
process started at $0$.

We assume that the following three properties hold a.s.:

\smallskip

(H1) $X$ has no negative jumps.

\smallskip
(H2) $X$ does not drift to $+\infty$.

\smallskip

(H3) The paths of $X$ are of infinite variation.

\smallskip
Thanks to (H1), the ``Laplace transform'' 
$E[\exp-\lambda X_t]$ is well defined 
for every $\lambda\geq 0$ and $t\geq 0$, and can be written
as
$$E[\exp-\lambda X_t]=\exp(t\psi(\lambda)),$$
with a function $\psi$ of the form 
$$\psi(\lambda)=\alpha_0 \lambda
+\beta \lambda^2+\int_{(0,\infty)}
\pi(dr)\,(e^{-\lambda r}-1+1_{\{r<1\}}\lambda r),$$
where $\alpha_0\in\R$, $\beta\geq 0$ and the
L\'evy measure $\pi$ is a Radon measure 
on $(0,\infty)$ such that
$\int_{(0,\infty)} (1\wedge
r^2)\,\pi(dr)<\infty$.

Assumption (H2) then holds iff $X$ has first moments
and $E[X_1]\leq 0$. The first moment assumption is equivalent to saying that
$\pi$ satisfies the stronger integrability
condition
$$\int_{(0,\infty)} (r\wedge
r^2)\,\pi(dr)<\infty.$$
Then $\psi$ can be written in the form
\begin{equation} 
\label{form-Laplace}
\psi(\lambda)=\alpha\lambda
+\beta \lambda^2+\int_{(0,\infty)}
\pi(dr)\,(e^{-\lambda r}-1+\lambda r),
\end{equation}
 Note that $\psi$ is then convex and that we have
$E[X_t]=-t\,\psi'(0)=-t\alpha$. The condition $E[X_1]\leq 0$ thus holds iff $\alpha\geq 0$. 
The process $X$ is recurrent or drifts to 
$-\infty$ according as $\alpha=0$ or $\alpha>0$.

Finally, according to \cite{Be}
(Corollary VII.5), assumption (H3)
holds iff at least one of the following two
conditions is satisfied: $\beta>0$, or
$$\int_{(0,1)} r\,\pi(dr)=\infty.$$

Summarizing, we assume that $X$ is a L\'evy process with no negative jumps,
whose Laplace exponent $\psi$ has the form (\ref{form-Laplace}),
where $\alpha\geq 0$, $\beta\geq 0$ and $\pi$ is a $\sigma$-finite measure
on $(0,\infty)$ such that $\int (r\wedge r^2)\pi(dr)<\infty$,
and we exclude the case where both $\beta=0$ and 
$\int_{(0,1)} r\,\pi(dr)<\infty$.

\medskip
\rem Only assumption (H1) is crucial to the connections
with branching processes that are presented in this
work. Assumption (H2) means that we restrict our
attention to the critical or subcritical case. 
We impose assumption (H3) in order to concentrate
on the most interesting cases: A 
simpler parallel theory can be developed in the finite
variation case, see Section 3 of \cite{LGLJ1}.

\medskip

We will use the notation $T_y=\inf\{t\geq 0:X_t=-y\}$ for $y\in \R$.
By convention $\inf\emptyset=+\infty$.

Under our assumptions, the point $0$ is regular
for $(0,\infty)$ and for $(-\infty,0)$, meaning that
$\inf\{t>0 : X_t>0\}=0$ and $\inf\{t>0 : X_t<0\}=0$
a.s. (see \cite{Be}, Theorem VII.1 and Corollary VII.5).
We sometimes use this property in connection
with the so-called duality property: For every $t>0$, define a 
process $\wh X^{(t)}=(\wh X^{(t)}_s,0\leq s\leq t)$
by setting
$$\wh X^{(t)}_s=X_t-X_{(t-s)-}\ ,\quad {\rm if} 
\ 0\leq s<t,$$
and $\wh X^{(t)}_t=X_t$. Then  $(\wh X^{(t)}_s,0\leq s\leq
t)$ has the same law as
$(X_s,0\leq s\leq t)$. 

If we combine the duality property with the
regularity of $0$ for both $(0,\infty)$ and $(-\infty,0)$, we easily get that the set
$$\{s>0:X_{s-}=I_s\ {\rm or}\ X_{s-}=S_{s-}\}$$
almost surely does not intersect $\{s\geq 0:\Delta X_s\not =0\}$. This property
will be used implicitly in what follows.

\subsection{Local times at the maximum and the minimum}

For every $t\geq 0$, set
$$S_t=\sup_{s\leq t} X_s\ ,\quad I_t=\inf_{s\leq t}X_s.$$

Then both processes $X-S$ and $X-I$ are strong Markov
processes, and the results recalled at the end of the previous
subsection imply that the point $0$ is regular 
for itself with respect to each of these two 
Markov processes.
We can thus define the corresponding 
Markovian local times and excursion measures, which both
play a fundamental role in this work.

\medskip
Consider first $X-S$.
We denote by $L=(L_t,t\geq 0)$
a local time at $0$ for $X-S$. Observe
that $L$ is only defined up to a positive multiplicative
constant, that will be specified later.
Let $N^*$  be
the associated excursion measure, which is a 
$\sigma$-finite measure on $\D(\R_+,\R)$. It will 
be important for our purposes
to keep track of the final jump under
$N^*$. This can be achieved by the 
following construction. 
Let $(a_j,b_j),\ j\in J$ be the excursion
intervals of $X-S$ away from $0$. In the 
transient case
($\alpha>0$), there is exactly one value
$j\in J$ such that $b_j=+\infty$. For every
$j\in J$ let $\omega^j\in\D(\R_+,\R)$ 
be defined by
$$\omega^j(s)=X_{(a_j+s)\wedge b_j}-X_{a_j}\ ,
\qquad s\geq 0.$$
Then the point measure
$$\sum_{j\in J} \delta_{(L_{a_j},\omega^j)}$$
is distributed as $1_{\{l\leq \eta\}}
\n(dld\omega)$, where $\n$
denotes a Poisson point measure 
with intensity $dl\,N^*(d\omega)$,
and $\eta=\inf\{l:\n([0,l]\times
\{\sigma=+\infty\})\geq 1\}$, if
$$\sigma(\omega)=\inf\{t>0:\omega(r)
=\omega(t)\ {\rm for}\ {\rm every}\ r\geq t\}$$
stands for the duration of the excursion $\omega$. This
statement characterizes the
excursion measure
$N^*$, up to the multiplicative 
constant already mentioned. Note that $X_0=0$ and
$X_t=X_\sigma\geq 0$
for $t\geq \sigma$, $N^*$ a.e.

\smallskip
Consider then $X-I$.  
It is easy to verify that the continuous increasing process
$-I$ is a local time at $0$ for the Markov
process $X-I$. We will denote by $N$ the associated
excursion measure, which can be characterized
in a way similar to $N^*$ (with the difference that
we have always  $-I_\infty=+\infty$ a.s., in contrast
to the property $L_\infty<\infty$ a.s. in the transient
case).  We already noticed that
excursions
of $X-I$ cannot start with a jump. Hence, $X_0=0$,
$N$ a.e. It is also clear from our assumptions on $X$
that $\sigma<\infty$, $X_t>0$ for every $t\in(0,\sigma)$
and $X_{\sigma-}=0$, $N$ a.e.

\medskip
We will now specify the normalization of 
$N^*$, or equivalently of $L$. Let 
$m$ denote Lebesgue measure on $\R$.

\medskip
\begin{lemma}
\label{invarmeas}
We can fix the normalization of $L$, or
equivalently of $N^*$, so that, for every
Borel subset $B$ of $(-\infty,0)$,
\begin{equation}
\label{inv}
N^*\left(\int_0^\sigma ds\,1_B(X_s)\right)=m(B).
\end{equation}
\end{lemma}

\proof
For every $x\in\R$, write $P_x$ for the law of the
L\'evy process started at $x$. Also set
$\tau=\inf\{s\geq 0:X_s\geq 0\}$ and recall that
$(X_t,t>0)$ is Markovian under $N^*$
with the transition kernels of the
underlying L\'evy process stopped when hitting
$[0,\infty)$.
Thanks to this observation, it is enough to prove that,
for every $\varepsilon>0$,
there exists a constant $c(\varepsilon)$
such that for every Borel subset $B$ of
$(-\infty,-\varepsilon)$,
$$E_{-\varepsilon}\cg\int_0^\tau ds\,1_B(X_s)\cd
=c(\varepsilon)m(B).$$

Consider first the transient case. By applying 
the strong Markov property at hitting times of negative values, it is easy
to verify that the measure on $(-\infty,-\varepsilon)$
defined by
$$B\la E_{-\varepsilon}\cg\int_0^\infty ds\,1_B(X_s)\cd$$
must be a multiple of Lebesgue measure. However,
writing $T^\varepsilon_0=0,T^\varepsilon_1,\ldots,
T^\varepsilon_n$, etc. for the successive visits
of $-\varepsilon$ via $[0,\infty)$, we have
$$E_{-\varepsilon}\cg\int_0^\infty ds\,1_B(X_s)\cd
=\sum_{i=0}^\infty E_{-\varepsilon}
\cg 1_{\{T_i^\varepsilon<\infty\}}
\int_{T_i^\varepsilon}^{T^{\varepsilon}_{i+1}}
ds\,1_B(X_s)\cd={
E_{-\varepsilon}\cg\int_0^\tau ds\,1_B(X_s)\cd\over
P_{-\varepsilon}[\tau=\infty]} .$$
The desired result follows.

In the recurrent case, the ergodic theorem gives
$${1\over n}\int_0^{T^\varepsilon_n}ds\,1_B(X_s)
\build\la_{n\rightarrow \infty}^{\rm a.s.}
E_{-\varepsilon}\cg\int_0^\tau ds\,1_B(X_s)\cd,$$
whereas the Chacon-Ornstein ergodic theorem implies
$${\int_0^{T^\varepsilon_n}ds\,1_B(X_s)
\over
\int_0^{T^\varepsilon_n}ds\,
1_{(-2\varepsilon,-\varepsilon)}(X_s)}
\build\la_{n\rightarrow \infty}^{\rm a.s.} {m(B)\over \varepsilon}.$$
The conclusion easily follows. \cq

\smallskip

In what follows we always assume that the normalization
of $L$ or of $N^*$ is fixed as in Lemma \ref{invarmeas}.

Let $L^{-1}(t)=\inf\{s,L_s>t\}$. By convention,
$X_{L^{-1}(t)}=+\infty$ if $t\geq L_\infty$. The 
process $(X_{L^{-1}(t)},t\geq 0)$ is a 
subordinator (the so-called
ladder height process) killed at an 
independent exponential time in the transient
case. 

\begin{lemma}
\label{Laplaceladder}
For every $\lambda> 0$,
$$E[\exp-\lambda
X_{L^{-1}(t)}]=\exp(-t\wt\psi(\lambda)),$$
where 
$$\wt\psi(\lambda)={\psi(\lambda)\over\lambda}
=\alpha+\beta\lambda+
\int_0^\infty (1-e^{-\lambda r})\pi([r,\infty))
\,dr.$$
\end{lemma}

\proof
By a classical result of fluctuation theory (see
e.g. \cite{Bi1} Corollary p.724), we have
$$E[\exp-\lambda
X_{L^{-1}(t)}]=\exp(-ct\wt\psi(\lambda)),$$
where $c$ is a positive constant. We have to verify that
$c=1$ under our normalization.

Suppose first that 
$\pi\not =0$.  Then notice that the L\'evy
measure $c\pi([r,\infty))dr$ of $X_{L^{-1}(t)}$
is the ``law'' of $X_\sigma$ under
$N^*(\cdot\cap \{X_\sigma>0\})$. However,
for any nonnegative measurable function
$f$ on $[0,\infty)^2$, we get by 
a predictable projection
\begin{eqnarray*}
N^*\cg f(\Delta X_\sigma,X_\sigma)\,1_{\{
X_\sigma>0\}}\cd
&=&N^*\cg \sum_{0<s\leq \sigma}f(\Delta
X_s,X_s)\,1_{\{
X_s>0\}}\cd\\
&=&N^*\cg\int_0^\sigma ds \int \pi(dx)\,
f(x,X_{s-}+x)\,1_{\{X_{s-}+x>0\}}\cd\\
&=&\int_{-\infty}^0 dy\int\pi(dx)\,f(x,y+x)\,
1_{\{y+x>0\}},
\end{eqnarray*}
using Lemma \ref{invarmeas} in the last 
equality. It follows that
\begin{equation}\label{joint}
N^*\cg f(\Delta X_\sigma,X_\sigma)\,1_{\{
X_\sigma>0\}}\cd
=\int \pi(dx)\int_0^x dz\,f(x,z),
\end{equation} 
and we get $c=1$ by comparing with the 
L\'evy measure of $X_{L^{-1}(t)}$. 

In the case $\pi=0$,
$X$ is a scaled linear Brownian motion 
with drift, and the same conclusion follows
from direct computations. \cq

\medskip
Note that we have in particular
$P[L^{-1}(t)<\infty]=e^{-\alpha t}$, which shows that
$L_\infty$ has an exponential distribution with parameter
$\alpha$ in the transient case.

When $\beta>0$, we can get a simple expression for
$L_t$. From well-known results on subordinators,
we have a.s. for every $u\geq 0$,
$$m(\{X_{L^{-1}(t)};t\leq u,L^{-1}(t)<\infty\})
=\beta(u\wedge L_\infty).$$
Since the sets $\{X_{L^{-1}(t)};t\leq u,
L^{-1}(t)<\infty\}$ and $\{S_r;r\leq L^{-1}(u)\}$
coincide except possibly for a countable set,
we have also
\begin{equation}\label{localLeb}
m(\{S_r;r\leq t\})=\beta\,L_t
\end{equation}
for every $t\geq 0$ a.s.

The next lemma provides a useful approximation of
the local time $L_t$.

\begin{lemma}
\label{localstrip}
For every $x>0$,
$$\lim_{\varepsilon\da 0}
{1\over \varepsilon}
\int_0^{L^{-1}(x)}
1_{\{S_s-X_s<\varepsilon\}}ds=x\wedge L_\infty$$
in the $L^2$-norm. Consequently, for every $t\geq 0$,
$$\lim_{\varepsilon\da 0}
{1\over \varepsilon}
\int_0^{t}
1_{\{S_s-X_s<\varepsilon\}}ds=L_t$$
in probability.
\end{lemma}

\proof It is enough to prove the first assertion. Let $\n$
be as previously a Poisson point measure
on $\R_+\times \D(\R_+,\R)$ 
with intensity $dl\,N^*(d\omega)$, and $\eta=\inf\{l:\n([0,l]\times
\{\sigma=+\infty\})\geq 1\}$. For 
every $x>0$ set
$$J_\varepsilon(x)=
{1\over \varepsilon}
\int \n(dld\omega)\,1_{\{l\leq x\}}
\int_0^{\sigma(\omega)}
1_{(-\varepsilon,0]}(\omega(s))\,ds.$$
Then,
$$E[J_\varepsilon(x)]
={x\over \varepsilon}
N^*\cg \int_0^\sigma 1_{(-\varepsilon,0]}(X_s)
ds\cd=x$$
by (\ref{inv}). Furthermore,
$$E[J_\varepsilon(x)^2]
=(E[J_\varepsilon(x)])^2
+x\varepsilon^{-2}
N^*\cg \pg\int_0^\sigma 1_{(-\varepsilon,0]}(X_s)
ds\pd^2\cd,$$
and
\begin{eqnarray*}
N^*\cg \pg\int_0^\sigma 1_{(-\varepsilon,0]}(X_s)
ds\pd^2\cd
&=&2\,N^*\cg \int_{0\leq s\leq t\leq \sigma}
1_{(-\varepsilon,0]}(X_s) 
1_{(-\varepsilon,0]}(X_t)dsdt\cd\\
&=&2\,N^*\cg \int_0^\sigma 
ds\,1_{(-\varepsilon,0]}(X_s)
E_{X_s}\cg\int_0^{\tau} dt\,
1_{(-\varepsilon,0]}(X_t)\cd\cd\\
&\leq &2\varepsilon
\sup_{0\geq y>-\varepsilon}
E_y\cg\int_0^{\tau} dt\,
1_{(-\varepsilon,0]}(X_t)\cd,
\end{eqnarray*}
using the same notation 
${\tau}=\inf\{t\geq 0:X_t\geq 0\}$ as previously. We
then claim that
\begin{equation}\label{techlocaltime}
\sup_{0\geq y>-\varepsilon}
E_y\cg\int_0^{\tau} dt\,
1_{(-\varepsilon,0]}(X_t)\cd=o(\varepsilon)
\end{equation}
as $\varepsilon\rightarrow  0$. Indeed, by applying the
strong Markov property at $T_y$, we have for $y>0$,
$$N^*[T_y<\infty]\,E_{-y}\cg\int_0^
{\tau}
dt\, 1_{(-\varepsilon,0]}(X_t)\cd
\leq N^*\cg\int_0^\sigma dt\,
1_{(-\varepsilon,0]}(X_t)\cd=\varepsilon,$$
and the claim follows since $N^*[T_y<\infty]\ua
+\infty$ as $y\da 0$. From (\ref{techlocaltime})
and the preceding calculations,
we get
$$\lim_{\varepsilon\rightarrow  0}
E[(J_\varepsilon(x)-x)^2]=0.$$
By Doob's inequality (or a monotonicity argument), 
we have also
$$\lim_{\varepsilon\rightarrow  0}
E[\sup_{0\leq z\leq
x}(J_\varepsilon(z)-z)^2]=0.$$
The lemma now follows, since the pair
$$\pg{1\over \varepsilon}
\int_0^{L^{-1}(x)}
1_{\{S_s-X_s<\varepsilon\}}ds,L_\infty\pd$$
has the same distribution as
$(J_\varepsilon(x\wedge \eta),\eta)$.
\cq

\smallskip
As a consequence of Lemma \ref{localstrip}, we 
may choose a sequence $(\varepsilon_k,k=1,2,\ldots)$
of positive real numbers decreasing to $0$, such that
\begin{equation}\label{approxLT}
 L_t=\lim_{k\rightarrow \infty}
{1\over \varepsilon_k}
\int_0^{t}
1_{\{S_s-X_s<\varepsilon_k\}}ds\;,\qquad P\hbox{ a.s.}
\end{equation}
Using monotonicity arguments and a diagonal
subsequence, we may and will assume that the
previous convergence holds simultaneously for
every $t\geq 0$ outside a single set of
zero probability. In particular, if
we set for $\omega\in\D([0,t],\R)$,
$$\Phi_t(\omega)
=\liminf_{k\rightarrow \infty}
{1\over \varepsilon_k}
\int_0^{t}
1_{\{\sup_{[0,s]}\omega(r)-\omega(s)<\varepsilon_k\}}ds,$$
we have $L_t=\Phi_t(X_s,0\leq s \leq t)$, for every
$t\geq 0$, $P$ a.s.

Recall the notation $\wh X^{(t)}$ for the process $X$
time-reversed at time $t$.

\begin{proposition}
\label{keyinv}
For any nonnegative measurable functional
$F$ on the Skorokhod space $\D(\R_+,\R)$,
$$N\cg\int_0^\sigma dt\,F(\wh
X^{(t)}_{s\wedge t},s\geq 0)\cd =E\cg
\int_0^{L_\infty} dx\,F(X_{s\wedge
L^{-1}(x)},s\geq 0)\cd.$$
\end{proposition}

\proof We may assume that $F$ is bounded and
continuous. Fix $t>0$ and if $\omega\in\D([0,t],\R)$,
set $T_{\rm max}(\omega)=\inf\{s\in[0,t]:
\sup_{[0,s]}\omega(r)=\sup_{[0,t]}\omega(r)\}$
and let $\theta\omega\in \D(\R_+,\R)$
be defined by $\theta\omega(t)=\omega(t\wedge
T_{\rm max}(\omega))$. 
Let $z>0$. Excursion theory for
$X-I$ shows that, for every $\varepsilon>0$,
$$N\cg\int_0^\sigma dt\,1_{\{\Phi_t(\hat
X^{(t)})\leq z\}}\,F(\wh X^{(t)}_{s\wedge
t},s\geq 0)\cd ={1\over \varepsilon}
E\cg\int_0^{T_{\varepsilon}}
dt\,1_{\{\Phi_t(\hat
X^{(t)})\leq z\}}\,F\circ\theta(\wh
X^{(t)})\cd.$$
In deriving this equality, we also apply to the
time-reversed process $\wh X^{(t)}$ the fact that
the local time $L_s$ does not increase after
the (first) time of the maximum over $[0,t]$. Then,
\begin{eqnarray*}
&&{1\over \varepsilon}E\cg\int_0^{T_{\varepsilon}}
dt\,1_{\{\Phi_t(\hat
X^{(t)})\leq z\}}\,F\circ\theta(\wh
X^{(t)})\cd\\
&&\qquad={1\over \varepsilon}
E\cg\int_0^\infty
dt\,1_{\{I_t>-\varepsilon\}}\,1_{\{\Phi_t(\hat
X^{(t)})\leq z\}}\,F\circ\theta(\wh
X^{(t)})\cd\\
&&\qquad={1\over \varepsilon}
\int_0^\infty dt\,
E[1_{\{S_t-X_t<\varepsilon\}}\,1_{\{L_t\leq z\}}
\,F\circ\theta(X_s,s\leq t)]\\
&&\qquad=E\cg{1\over \varepsilon}
\int_0^{L^{-1}(z)} dt\,
1_{\{S_t-X_t<\varepsilon\}}
F\circ\theta(X_s,s\leq t)\cd.
\end{eqnarray*}
We then take $\varepsilon=\varepsilon_k$ and pass to the
limit
$k\rightarrow \infty$, using the $L^2$ bounds provided by Lemma \ref{localstrip}. Note that
the measures 
$${1\over \varepsilon_k}\,1_{[0,L^{-1}(z)]}(t)\,1_{\{S_t-X_t
<\varepsilon_k\}}dt$$
converge weakly to the finite
measure $1_{[0,L^{-1}(z)]}(t)dL_t$. Furthermore, 
$\theta(X_s,s\leq t)=(X_{s\wedge t},s\geq 0)$,
$dL_t$ a.e., a.s., and it is easy to verify
that the mapping $t\la F\circ \theta(X_s,s\leq
t)$
is continuous on a set of full $dL_t$-measure.
We conclude that
\begin{eqnarray*}N\cg\int_0^\sigma dt\,1_{\{\Phi_t(\hat
X^{(t)})\leq z\}}\,F(\wh X^{(t)}_{s\wedge
t},s\geq 0)\cd 
&=&E\cg\int_0^{L^{-1}(z)}
F(X_{s\wedge t},s\geq 0)\,dL_t\cd\\
&=&
E\cg\int_0^{z\wedge L_\infty} 
F(X_{s\wedge L^{-1}(x)},s\geq 0)\,dx\cd,
\end{eqnarray*}and the desired result follows by letting 
$z\rightarrow \infty$.\cq

\section{The height process and the exploration process}

We write $\wh S^{(t)}_s=\sup_{[0,s]}\wh X^{(t)}_r$
($0\leq s\leq t$) for the supremum process of
$\wh X^{(t)}$. 

\begin{definition}
The height process is the real-valued process
$(H_t,t\geq 0)$ defined as follows. First $H_0=0$
and for every $t>0$, $H_t$ is the local time at 
level $0$ at time $t$ of the process $\wh
X^{(t)}-\wh S^{(t)}$.
\end{definition}

The normalization of local time is of course that
prescribed by Lemma \ref{invarmeas}. 

Note that
the existence of a modification of the process
$(H_t,t\geq 0)$ with good continuity
properties is not clear from the previous definition.
When $\beta>0$ however, we can use (\ref{localLeb})
to see that 
\begin{equation}\label{heightLeb}
H_t={1\over \beta}\,m(\{I^s_t;s\leq t\})\,,
\end{equation}
where for $0\leq s\leq t$,
$$I^s_t=\inf_{s\leq r\leq t} X_r.$$
Clearly the right-hand side of (\ref{heightLeb}) 
gives a continuous modification of $(H_t,t\geq 0)$. 
When $\beta=0$, this argument does not apply and we 
will see later that there may exist no continuous
(or even c\` adl\` ag) modification of $(H_t,t\geq 0)$. 

At the present stage, we will use
the measurable modification of $(H_t,t\geq 0)$
with values in $[0,\infty]$ obtained 
by taking
\begin{equation}\label{approxH}
H^o_t=\Phi_t(\wh X^{(t)}_s,0\leq s\leq t)
=\liminf_{k\rightarrow  \infty}
{1\over \varepsilon_k}\int_0^t
1_{\{X_s<I^s_t+\varepsilon_k\}}\,ds,
\end{equation}

The liminf in (\ref{approxH}) is a limit
(and is finite) a.s. for 
every fixed $t\geq 0$. The following lemma shows that
more is true.

\begin{lemma}
\label{techheight}
Almost surely for every $t\geq 0$, we
have
$$H^o_s=\lim_{k\rightarrow  \infty}
{1\over \varepsilon_k}\int_0^s
1_{\{X_r<I^r_s+\varepsilon_k\}}\,dr 
<\infty,$$
for every $s<t$ such that $X_{s-}\leq I^s_t$,
and for $s=t$ if $\Delta X_t>0$.
\end{lemma}

\proof Let $s$ and $t$ be as in the statement. Then there
must exist a rational $u\in(s,\infty)$ such that $X_{s-}\leq
I^s_u$. We can then apply to the time-reversed process
$\wh X^{(u)}$ the approximation result (\ref{approxLT})
at times $u$ and $u-s$ respectively. The desired result
follows. \cq

\medskip
We denote by $M_f(\R_+)$ the space of all
finite measures on $\R_+$, which is equipped 
with the topology of weak convergence.

\begin{definition}
\label{explordef}
The exploration process is the process
$(\rho_t,t\geq 0)$ with values in $M_f(\R_+)$
defined as follows. For every nonnegative
measurable function $f$,
\begin{equation}\label{rhodef}
\langle \rho_t,f\rangle
=\int_{[0,t]} d_sI^s_t\,f(H^o_s)
\end{equation}
where the notation $d_sI^s_t$ refers to integration
with respect to the nondecreasing function $s\rightarrow  I^s_t$.
\end{definition}

Since we did not exclude the value $+\infty$
for $H^o_t$ (as defined by (\ref{approxH})), it may not be
obvious that the measure $\rho_t$
is supported on $[0,\infty)$. However,
this readily follows from the previous lemma
since the measure $d_sI^s_t$ is supported 
on the set $\{s<t:X_{s-}\leq I^s_t\}$ (to
which we need to add the point $t$ if $\Delta X_t>0$).

Notice that if $u$ and $v$ belong to the
set $\{s\leq t:X_{s-}\leq I^s_t\}$, and if $u\leq v$, then 
for every $r\in[0,u)$ the condition $X_{r}<I^r_u+\varepsilon_k$ implies 
$X_r< I^r_v+\varepsilon_k$, and by construction it follows that
$H^o_u\leq H^o_v$. Using the previous remark on the support 
of the measure $d_sI^s_t$, we see that the measure $\rho_t$
is supported on $[0,H^o_t]$, for every $t\geq 0$, a.s.

\medskip
The total mass of $\rho_t$ is
$$\langle \rho_t,1\rangle=I^t_t-I^0_t=X_t-I_t.$$
In particular $\rho_t=0$ iff $X_t=I_t$. 

\medskip
It will be useful to rewrite the definition of
$\rho_t$ in terms of the time-reversed process $\wh X^{(t)}$.
Denote by $\wh L^{(t)}=(\wh L^{(t)}_s,0\leq s\leq t)$ the local time
at $0$ of $\wh X^{(t)}-\wh S^{(t)}$ (in particular 
$H^o_t=\wh L^{(t)}_t$). Note that for $t\geq 0$ fixed,
we have $H^o_s=\wh L^{(t)}_t-\wh L^{(t)}_s$ for every 
$s\in[0,t]$ such that $X_{s-}\leq I^s_t$, a.s. (compare (\ref{approxLT})
and (\ref{approxH})).
Hence,
\begin{equation}\label{rhoreversed}
\langle \rho_t,f\rangle
=\int_{[0,t]} d\wh S^{(t)}_s\,f(\wh L^{(t)}_t-\wh L^{(t)}_s).
\end{equation}

If $\mu$ is a nonzero measure in $M_f(\R_+)$,
we write $\supp \mu$ for the topological support
of $\mu$ and set $H(\mu)=\sup(\supp\mu)$. By
convention $H(0)=0$. By a preceding remark, $H(\rho_t)\leq H^o_t$
for every $t\geq 0$, a.s.

\begin{lemma}
\label{techrho}
For every $t\geq 0$, $P[H(\rho_t)=H^o_t]=1$. Furthermore, almost
surely for every $t> 0$, we have
\begin{description}
\item{\rm (i)} $\rho_t(\{0\})=0$ ;
\item{\rm (ii)} $\supp \rho_t=[0,H(\rho_t)]$ if $\rho_t\not =0$ ;
\item{\rm (iii)} $H(\rho_s)=H^o_s$ for every $s\in[0,t)$ such that
$X_{s-}\leq I^s_t$ and for $s=t$ if $\Delta X_t>0$.
\end{description}
\end{lemma}

\proof It is well known, and easy to prove from the strong Markov property,
that the two random measures $dS_s$ and $dL_s$ have the same support
a.s. Then (\ref{rhoreversed}) implies that $\supp \rho_t=[0,H^o_t]$
a.s. for every fixed $t>0$. In particular, $P[H^o_t=H(\rho_t)]=1$.
Similarly (\ref{rhoreversed}) implies that $P[\rho_t(\{0\})>0]=0$
for every fixed $t>0$. However, if we have $\rho_t(\{0\})>0$
for some $t\geq 0$, our definitions and the right-continuity of paths
show that the same
property must hold for some rational
$r>t$. Property (i) follows.

Let us now prove (ii), which is a little more delicate. We already noticed
that (ii) holds for every fixed $t$, a.s., hence for every rational
outside a set of zero probability. Let $t>0$ with $X_t>I_t$, and set
$$\gamma_t=\sup\{s<t:I^s_t<X_t\}.$$
We consider two different cases.

(a) Suppose first that $X_{\gamma_t-}<X_t$, which holds in particular if $\Delta X_t>0$.
Then note that
$$\langle\rho_t,f\rangle=\int_{[0,\gamma_t)} d_sI^s_t\,f(H^o_s)+(X_t-X_{\gamma_t-})f(H^o_t).$$
Thus we can
find a rational $r>t$ sufficienty close to $t$, so that $\rho_r$ and 
$\rho_t$ have the same restriction to $[0,H^o_t)$. The fact that 
property (ii) holds for $r$ implies that it holds for $t$, and
we see also that $H^o_t=H(\rho_t)$ in that case.

(b) Suppose that $X_{\gamma_t-}=X_t$. Then we set for every $\varepsilon>0$,
$$\langle \rho^\varepsilon_t,f\rangle
=\int_{[0,t]} d_sI^s_t\,1_{\{I^s_t<X_t-\varepsilon\}}\,f(H^o_s).$$
From the remarks following the definition of $\rho_t$, it is 
clear that there exists some $a\geq 0$ such that $\rho^\varepsilon_t$ 
is bounded below by the restriction of $\rho_t$ to $[0,a)$, and
bounded above by the restriction of $\rho_t$ to $[0,a]$.
Also note that $\rho_t=\lim\uparrow \rho^\varepsilon_t$ as
$\varepsilon\downarrow 0$. Now, for every $\varepsilon>0$, we can pick 
a rational $r>t$ so that $I^t_r>X_t-\varepsilon$, and we have by
construction
$$\rho^\varepsilon_t=\rho^{\varepsilon+X_r-X_t}_r.$$
From the rational case, the support of $\rho^{\varepsilon+X_r-X_t}_r$ must be an 
interval $[0,a]$, and thus the same is true for $\rho^\varepsilon_t$.
By letting $\varepsilon\da 0$, we
get (ii) for $t$.

We already obtained (iii) for $s=t$ when $\Delta X_s>0$
(see (a) above). If $s\in (0,t)$ is such that $X_{s-}\leq I^s_t$,
we will have also $X_{s-}\leq I^s_r$ for any rational $r\in(s,t)$.
Then $H^o_s\leq H^o_r=H(\rho_r)$, and on the other hand, it is
clear that 
the measures $\rho_s$ and $\rho_r$ have the same restriction
to $[0,H^o_s)$. Thus the desired result follows from (ii). \cq

\begin{proposition}
\label{explor}
The process $(\rho_t,t\geq 0)$ is a c\` adl\` ag strong Markov
process in $M_f(\R_+)$. 
\end{proposition}
\rem The proof will show that $(\rho_t,t\geq 0)$ is even c\` adl\` ag with
respect to the variation distance on finite measures.

\medskip
\proof 
We first explain how to define the process
$\rho$ started at an arbitrary initial value
$\mu\in M_f(\R_+)$. 
To this end, we introduce some
notation. Let $\mu\in M_f(\R_+)$
and $a\geq 0$. If $a\leq \langle\mu,1\rangle$, we let
$k_a\mu$ be the unique finite measure on $\R_+$ such that,
for every
$r\geq 0$,
$$k_a\mu([0,r])=\mu([0,r])\wedge
(\langle\mu,1\rangle-a).$$ In particular,
$\langle k_a\mu,1\rangle=\langle\mu,1\rangle-a$. If $a\geq
\langle\mu,1\rangle$, we take
$k_a\mu=0$.

If $\mu\in M_f(\R_+)$ has compact support 
and $\nu\in M_f(\R_+)$, we define the
concatenation $[\mu,\nu]\in M_f(\R_+)$ by the formula
$$\int [\mu,\nu](dr)\,f(r)=\int \mu(dr)\,f(r)+\int
\nu(dr)\,f(H(\mu)+r).$$ 

With this notation, the law of the process $\rho$ started
at $\mu\in M_f(\R_+)$ is the distribution 
of the process $\rho^\mu$ defined by
\begin{equation}\label{rhoinitial}
\rho^\mu_t=[k_{-I_t}\mu,\rho_t]\ ,\qquad
t>0.
\end{equation} Note that this definition makes sense because
$k_{-I_t}\mu$ has compact support, for every $t>0$ a.s.

We then verify that the process $\rho$ has the stated
properties. For simplicity, we consider 
only the case when the initial value 
is $0$, that is when $\rho$ is defined as 
in Definition \ref{explordef}.
The
right-continuity of paths is straightforward from the
definition since the measures
$1_{[0,t']}(s)d_sI^s_{t'}$ converge to $1_{[0,t]}(s)d_sI^s_{t}$ 
in the variation
norm as $t'\downarrow t$. Similarly, we get
the existence of left-limits from the fact that
the measures $1_{[0,t']}(s)d_sI^s_{t'}$
converge to $1_{[0,t)}(s)d_sI^s_{t}$ in the variation
norm as $t'\uparrow t$, $t'<t$. 
We see in particular that $\rho$ and $X$ have the same
discontinuity times and that
\begin{equation} 
\label{discontH}
\rho_t=\rho_{t-}+\Delta X_t\,\delta_{H^o_t}.
\end{equation}

We now turn to the strong Markov 
property. Let $T$ be a stopping time of the
canonical 
filtration. We will express $\rho_{T+t}$ in terms
of $\rho_T$ and the shifted process
$X^{(T)}_t=X_{T+t}-X_T$. We claim that a.s. for every $t>0$
\begin{equation}\label{strongMarkov}
\rho_{T+t}=[k_{-I^{(T)}_t}\rho_T,\rho^{(T)}_t]
\end{equation}
where $\rho^{(T)}_t$ and $I^{(T)}_t$ obviously denote the
analogues of
$\rho_t$ and $I_t$ when $X$ is replaced by $X^{(T)}$. When
we have  proved (\ref{strongMarkov}), the strong Markov
property of the process $\rho$ follows 
by standard arguments, using also our definition
of the process started at a general initial value.

\smallskip
For the proof of (\ref{strongMarkov}), 
write 
$$\langle \rho_{T+t},f\rangle
=\int_{[0,T]}d_sI^s_{T+t} f(H^o_s)
+\int_{(T,T+t]}d_sI^s_{T+t} f(H^o_s).
$$
We consider separately each term in the right-hand side.
Introduce $u=\sup\{r\in(0,T]:X_{r-}<I^T_{T+t}\}$,
with the usual convention $\sup\emptyset=0$.
We have $I^s_{T+t}=I^s_T$ for $s\in [0,u)$
and $I^s_{T+t}=I^T_{T+t}$ for $s\in [u,T]$. Since
$X_{T}-I^T_{T+t}=-I^{(T)}_t$, it
then follows from our definitions that
\begin{equation}\label{sMarkov1}
\int_{[0,T]}d_sI^s_{T+t}\, f(H^o_s)
=\int_{[0,u)}d_sI^s_T\,f(H^o_s)+1_{\{u>0\}}(I^T_{T+t}-X_{u-})f(H^o_u)
=\langle k_{-I^{(T)}_t}\rho_T,f\rangle.
\end{equation}
Also notice that the measures $\rho_u$ and $k_{-I^{(T)}_t}\rho_T$
coincide, except possibly at the point $H^o_u$. In any case,
$H(\rho_u)=H(k_{-I^{(T)}_t}\rho_T)$, and we have also
$H^o_u=H(\rho_u)$ by Lemma \ref{techrho} (iii).

Now observe that for $d_sI^s_{T+t}$ almost every 
$s\in(T,T+t]$, we have $H^o_{s}=H^o_u+H^{o,(T)}_{s-T}$,
with an obvious notation. To see this, pick a rational
$r>T+t$ such that $I^{T+t}_r>X_{s-}$ and argue on the time-reversed
process $\widehat X^{(r)}$ as in the proof of Lemma 
\ref{techheight}. Hence,
\begin{equation}\label{sMarkov2}
\int_{(T,T+t]}d_sI^s_{T+t}\,f(H^o_s)=
\int_{(T,T+t]}d_sI^s_{T+t}\,f(H^o_u+H^{o,(T)}_{s-T})=
\int \rho^{(T)}_t(dx)\,f(H^o_u+x).
\end{equation}
Formula (\ref{strongMarkov}) follows by combining 
(\ref{sMarkov1}) and (\ref{sMarkov2}).\cq

\medskip
We now come back to the problem of finding a modification 
of the height process with good continuity properties.
By the first assertion of Lemma \ref{techrho},
$(H(\rho_t),t\geq 0)$ is a modification of $(H^o_t,t\geq 0)$.
{\bf From now on, we will
systematically use this modification} and write 
$H_t=H(\rho_t)$. 
From Lemma \ref{techrho} (iii), we see that
formula (\ref{rhodef}) defining $(\rho_t,t\geq 0)$ remains true 
if $H^o_s$ is replaced by $H_s$. The same applies to
formula (\ref{discontH}) giving the jumps of $\rho$.
Furthermore, the continuity properties of the process
$\rho_t$ (and especially the form of its jumps) imply that
the mapping $t\la H(\rho_t)=H_t$ is lower semicontinuous
a.s.

\medskip
Let us make an important remark at this point. Write 
$$g_t=\sup\{s\leq t:X_{s}=I_s\}$$
for the beginning of the excursion of $X-I$ that straddles $t$.
Then a simple time-reversal argument shows that a.s. for
every $t$ such that $X_t-I_t>0$, we have
$$\lim_{k\rightarrow \infty}{1\over \varepsilon_k}\int_0^{g_t}
1_{\{X_s<I^s_t+\varepsilon_k\}}\,ds=0$$
and thus we can replace (\ref{approxH})
by
$$H^o_t=\liminf_{k\rightarrow  \infty}
{1\over \varepsilon_k}\int_{g_t}^t
1_{\{X_s<I^s_t+\varepsilon_k\}}\,ds.$$
Recalling (\ref{rhodef}), we see that, a.s. for every 
$t\geq 0$ such that $X_t-I_t>0$, we can write $\rho_t$ and $H_t$ as
measurable functions of the excursion of $X-I$ that straddles $t$
(and of course $\rho_t=0$ and $H_t=0$ if $X_t=I_t$).
We can thus define both the height process and the
exploration process under the excursion measure $N$. 
Furthermore, if $(\alpha_j,\beta_j)$, $j\in J$,
denote the excursion intervals of $X-I$, and if $\omega_j$, $j\in J$,
denote the corresponding excursions, we have 
$H_t=H_{t-\alpha_j}(\omega_j)$ and 
$\rho_t=\rho_{t-\alpha_j}(\omega_j)$ for every $t\in(\alpha_j,\beta_j)$
and $j\in J$, a.s.

\medskip
Since $0$ is a regular point for $X-I$, we also see that
the measure $0$ is a regular point for the exploration
process $\rho$. It is immediate from the previous
remark that the excursion measure of $\rho$
away from $0$ is the ``law'' of $\rho$ under $N$.
Similarly, the process $-I$, which is the local time
at $0$ for $X-I$, is also the local time at $0$ for $\rho$.

\smallskip
We now state and prove a useful technical lemma about the
process $H$.

\begin{lemma}
\label{IVP}
{\rm (Intermediate values property)} Almost surely for every $t<t'$, the 
process $H$ takes all values between $H_t$ and $H_{t'}$ on the time interval $[t,t']$.
\end{lemma}

\proof First consider the case when $H_t>H_{t'}$. By using the
lower semi-continuity of $H$, we may assume that $t$ is rational. From
(\ref{strongMarkov}), we have $\rho_{t+s}=[k_{-I^{(t)}_s}\rho_t,\rho^{(t)}_s]$
for every $s>0$, a.s. Hence, if
$$\gamma_r=\inf\{s>0:I^{(t)}_s=-r\},$$
we have $\rho_{t+\gamma_r}=k_r\rho_t$,
and so $H_{t+\gamma_r}=H(k_r\rho_t)$ for every $r\geq 0$, a.s. However,
Lemma \ref{techrho} (ii) implies that the mapping $r\rightarrow  H(k_r\rho_t)$
is continuous. Now note that for $r=0$, $H_{t+\gamma_r}=H_t$,
whereas for $r=X_t-I^t_{t'}=-I^{(t)}_{t'-t}$ we have $t+\gamma_r\leq t'$ and
our definitions easily imply $\rho_{t+\gamma_r}\leq \rho_{t'}$ and
$H_{t+\gamma_r}\leq H_{t'}$.

Consider then the case when $H_t<H_{t'}$. By lower semi-continuity again, we
may assume that $t'$ is rational. In terms of the process
time-reversed at time $t'$, we have $H_{t'}=\wh L^{(t')}_{t'}$. Set
$$\sigma_r=\inf\{s\geq 0:\wh S^{(t')}_s\geq r\},$$
which is well defined for $0\leq r\leq X_{t'}-I_{t'}$. Since the
subordinator $S_{L^{-1}(t)}$ is strictly increasing, we see
that the mapping $r\rightarrow  \wh L^{(t')}_{\sigma_r}$ is continuous for
$r\in [0,X_{t'}-I_{t'}]$, a.s. Now note that
$$H_{t'-\sigma_r}=\wh L^{(t')}_{t'}-\wh L^{(t')}_{\sigma_r}$$
for every $r\in [0,X_{t'}-I_{t'}]$, a.s. For $r=X_{t'}-I^t_{t'}=\wh S^{(t')}_{t'-t}$,
we have 
$t'-\sigma_r\geq t$ and $H_{t'-\sigma_r}\leq H_t$ by construction.
The desired result follows. \cq

\smallskip
The next proposition is a corollary of 
Proposition \ref{keyinv}.
We denote by 
$U$ a subordinator defined on an auxiliary probability space, with Laplace exponent 
$$E[\exp-\lambda U_t]=\exp\pg-t\pg\beta\lambda+
\int_0^\infty (1-e^{-\lambda r})\,\pi([r,\infty))\,dr\pd\pd
=\exp(-t(\wt\psi(\lambda)-\alpha)).$$
For every $a\geq 0$, we let $J_a$ be the random 
element of $M_f(\R_+)$ defined by $J_a(dr)=1_{[0,a]}(r)\,dU_r$.

\begin{proposition}
\label{mesinv}
For every nonnegative measurable function $\Phi$ on $M_f(\R_+)$,
$$N\cg\int_0^\sigma dt\,\Phi(\rho_t)\cd
=\int_0^\infty da\,e^{-\alpha a}\,E[\Phi(J_a)].$$
Let $b\geq 0$. Then $\rho_s(\{b\})=0$
for every $s\geq 0$, $N$ a.e. or $P$ a.s.
\end{proposition}

\proof We have $\rho_t=\Sigma(\wh X^{(t)}_{s\wedge t},s\geq 0)$,
with a functional $\Sigma$ that is made explicit in (\ref{rhoreversed}).
We then apply Proposition \ref{keyinv} to obtain
$$
N\cg\int_0^\sigma dt\,\Phi(\rho_t)\cd=
E\cg\int_0^{L_\infty} da\,\Phi\circ\Sigma(X_{s\wedge L^{-1}(a)},s\geq 0)\cd.$$
However, for $a<L_\infty$,
$$\langle\Sigma(X_{s\wedge L^{-1}(a)},s\geq 0),\varphi\rangle
=\int_0^{L^{-1}(a)} dS_s\,\varphi(a-L_s).$$
The first assertion is now a consequence of 
Lemma \ref{Laplaceladder}, which shows that
$P[a<L_\infty]=e^{-\alpha a}$ and that, conditionally
on $\{L^{-1}(a)<\infty\}$, $S_{L^{-1}(r)}=X_{L^{-1}(r)}$, $0\leq r\leq a$
has the same distribution as $U$.

\smallskip
Consider now the second assertion. Note that the case $b=0$ is a consequence of
Lemma \ref{techrho} (i). So we may assume that $b>0$ and it is enough
to prove the result under the excursion measure $N$. However, 
since $b$ is a.s. not a jump time
of $U$, the right side of the formula of the proposition vanishes 
for $\Phi(\mu)=\mu(\{b\})$. The desired result follows, using
also the fact that $\rho$ is c\` adl\` ag in the variation norm. \cq

\medskip
We denote by $\MM$ the measure on ${M}_f({\R}_+)$
defined by:
$$\langle \MM,\Phi\rangle=\int_0^\infty da\,e^{-\alpha a}\,E[\Phi(J_a)].$$
Proposition \ref{mesinv} implies that the measure 
$\MM$ is invariant for $\rho$. 

The last proposition
of this section describes the potential kernel 
of the exploration process killed when it hits $0$.
We fix $\mu\in M_f(\R_+)$ and let $\rho^\mu$
be as in (\ref{rhoinitial}) the exploration process
started at $\mu$. We use the notation introduced in the 
proof of Proposition \ref{explor}.

\begin{proposition}
\label{potker}
Let 
$\tau_0=\inf\{s\geq 0,\,\rho^\mu_s=0\}$. Then,
$$E\Big[ \int_0^{\tau_0} ds\,\Phi(\rho^\mu_s)\Big]
=\int_0^{<\mu,1>} dr\int \MM(d\theta)\,\Phi([k_r\mu,\theta]).$$
\end{proposition}

\proof First note that $\tau_0=T_{<\mu,1>}$
by an immediate application of the definition of $\rho^\mu$. Then,
denote by $(\alpha_j,\beta_j)$, $j\in J$ the excursion intervals of $X-I$
away from $0$ before time $T_{<\mu,1>}$, and by $\omega_j$, $j\in J$
the corresponding excursions. 
As we observed before Proposition \ref{mesinv}, we have
$\rho_t=\rho_{t-\alpha_j}(\omega_j)$ for every
$t\in (\alpha_j,\beta_j)$, $j\in J$, a.s. 
Since $\{s\geq 0:X_s=I_s\}$ has zero Lebesgue measure a.s., it follows that
$$E\Big[\int_0^{\tau_0}ds\,\Phi(\rho^\mu_s)\Big] =E\Big[\sum_{j\in
J}\int_0^{\beta_j-\alpha_j} dr\,
\Phi([k_{-I_{\alpha_j}}\mu,\rho_r(\omega_j)])\Big].$$
By excursion theory, the point measure
$$\sum_{j\in J} \delta_{I_{\alpha_j},\omega_j}(dude)$$
is a Poisson point measure with intensity
$1_{[-<\mu,1>,0]}(u)du\,N(d\omega)$. Hence,
$$E\Big[\int_0^{\tau_0}ds\,\Phi(\rho^\mu_s)\Big]
=\int_0^{<\mu,1>}du\,N\cg\int_0^\sigma dr\,
\Phi([k_u\mu,\rho_r])\cd,$$
and the desired result follows from Proposition \ref{mesinv}. \cq

\section{Local times of the height process}

\subsection{The construction of local times}

Our goal is to construct a
local time process for $H$ at each level $a\geq 0$.
Since $H$ is in general not Markovian (and not a
semimartingale) we cannot apply a general theory, but
still we will use certain ideas which are familiar
in the Brownian motion setting. In the case $a=0$,
we can already observe that $H_t=0$ iff $\rho_t=0$
or equivalently $X_t-I_t=0$. Therefore the 
process $-I $ is the natural candidate for the 
local time of $H$ at level $0$.

\medskip
Let us fix $a\geq 0$. Since $t\rightarrow  \rho_t$ is c\` adl\` ag
in the variation norm, it follows that the mapping $t\rightarrow 
\rho_t((a,\infty))$ is c\` adl\` ag. Furthermore, it follows from 
(\ref{discontH})
that the
discontinuity times of this mapping are exactly
those times $t$ such that $\Delta X_t>0$ and
$H_t>a$, and the corresponding jump is $\Delta
X_t$.

Set
$$\tau^a_t=\inf\{s\geq 0:\int_0^s
1_{\{H_r>a\}}\,dr>t\}
=\inf\{s\geq 0:\int_0^s
1_{\{\rho_r((a,\infty))>0\}}\,dr>t\}.$$
From Proposition \ref{mesinv}, we get that $\int_0^\infty
1_{\{H_r>a\}}\,dr=\infty$ a.s., so that the
random times $\tau^a_t$ are a.s. finite.

When $a>0$, we also set  
$$\wt\tau^a_t=\inf\{s\geq 0:\int_0^s
1_{\{H_r\leq a\}}\,dr>t\}$$
and we let $\h_a$ be the $\sigma$-field
generated by the c\` adl\` ag process
$(X_{\tilde\tau^a_t},\rho_{\tilde\tau^a_t};t\geq 0)$
and the class of $P$-negligible sets of $\g_\infty$.
We also define $\h_0$ as
the $\sigma$-field
generated by the class of $P$-negligible sets of
$\g_\infty$.

\begin{proposition}
\label{reflec}
For every $t\geq 0$, let $\rho^a_t$ be the 
random measure on $\R_+$ defined by
$$\langle\rho^a_t,f\rangle 
=\int_{(a,\infty)}
\rho_{\tau^a_t}(dr)\,f(r-a).$$
The process $(\rho^a_t,t\geq 0)$ has
the same distribution as $(\rho_t,t\geq 0)$
and is independent of $\h_a$.
\end{proposition}

\proof {\it First step.} We first verify
that the process $(\langle\rho^a_t,1\rangle ,t\geq 0)$
has the same distribution
as $(\langle\rho_t,1\rangle ,t\geq 0)$.

Let $\varepsilon>0$. We introduce two
sequences of stopping times $S^k_\varepsilon$,
$T^k_\varepsilon$, $k\geq 1$, defined inductively
as follows:
\begin{eqnarray*}&&S^1_\varepsilon=\inf\{s\geq 0:\rho_s((a,\infty))\geq \varepsilon\},\\
&&T^k_\varepsilon=\inf\{s\geq
S^k_\varepsilon:\rho_s((a,\infty))=0\},\\
&&S^{k+1}_\varepsilon=\inf\{s\geq
T^k_\varepsilon:\rho_s((a,\infty))\geq
\varepsilon\}.
\end{eqnarray*}It is easy to see that these stopping times
are a.s. finite, and $S^k_\varepsilon\ua \infty$,
$T^k_\varepsilon\ua\infty$ as $k\ua \infty$.

From (\ref{strongMarkov})
applied with $T=S_\varepsilon^k$, we 
obtain that,
for every $k\geq 1$,
\begin{equation}\label{reftech1}
T^k_\varepsilon=\inf\{s>S^k_\varepsilon:
X_s=X_{S^k_\varepsilon}-\rho_{S^k_\varepsilon}
((a,\infty))\}.
\end{equation}
Formula (\ref{strongMarkov})
also implies that, for every
$s\in[0,T^k_\varepsilon -S^k_\varepsilon]$,
\begin{eqnarray}\label{reftech2}
\rho_{S^k_\varepsilon+s}((a,\infty))
&=&(\rho_{S^k_\varepsilon}((a,\infty))
+I^{(S^k_\varepsilon)}_s)
+(X^{(S^k_\varepsilon)}_s-
I^{(S^k_\varepsilon)}_s)\\
&=&
X_{S^k_\varepsilon+s}
-(X_{S^k_\varepsilon}-\rho_{S^k_\varepsilon}((a,\infty))).\nonumber
\end{eqnarray}
We set
$$Y^{k,\varepsilon}_s=\rho_{(S^k_\varepsilon+s)
\wedge T^k_\varepsilon}((a,\infty)).$$
As a straightforward consequence of
(\ref{reftech1}) and (\ref{reftech2}),
conditionally on
$\g_{S^k_\varepsilon}$, the process
$Y^{k,\varepsilon}$ is distributed as the 
underlying L\'evy process started at
$\rho_{S^k_\varepsilon}((a,\infty))$
and stopped at its first hitting time of $0$.

We then claim that, for every $t\geq 0$,
\begin{equation}\label{reflec1}
\lim_{\varepsilon\rightarrow  0}
\sup_{\{k\geq 1,S^k_\varepsilon\leq t\}}
\rho_{S^k_\varepsilon}((a,\infty))=0,
\qquad{\rm a.s.}
\end{equation}
Indeed, by previous observations about
the continuity properties of the mapping
$s\la \rho_s((a,\infty))$, we have
$$\sup_{\{k\geq 1,S^k_\varepsilon\leq t\}}
\rho_{S^k_\varepsilon}((a,\infty))
\leq \varepsilon+\sup\{\Delta X_s;s\leq t,
H_s>a,\rho_s((a,\infty))\leq \varepsilon\}.$$
However, the sets $\{s\leq t;\Delta X_s>0,
H_s>a,\rho_s((a,\infty))\leq \varepsilon\}$
decrease to $\emptyset$ as $\varepsilon\da 0$,
and so
$$\lim_{\varepsilon\rightarrow  0}\pg\sup\{\Delta X_s;s\leq t,
H_s>a,\rho_s((a,\infty))\leq \varepsilon\}\pd=0,$$
a.s., which yields the desired claim.

Set
$$Z^{k,\varepsilon}_s=
Y^{k,\varepsilon}_s-\inf_{0\leq r\leq s}
Y^{k,\varepsilon}_r\leq Y^{k,\varepsilon}_s.$$
Then, conditionally on $\g_{S^k_\varepsilon}$,
$Z^{k,\varepsilon}$ is distributed as an
independent copy 
of the
reflected process $X-I$, stopped when its local
time at $0$ hits
$\rho_{S^k_\varepsilon}((a,\infty))$.

Denote by $U^\varepsilon=(U^\varepsilon_s,s\geq
0)$ the process obtained by patching together
the paths $(Z^{k,\varepsilon}_s,0\leq s\leq
T^k_\varepsilon-S^k_\varepsilon)$. By the 
previous remarks, $U^\varepsilon$ is distributed
as the reflected L\'evy process $X-I$.

We then set
$$\tau^{a,\varepsilon}_s
=\inf\{t\geq 0,
\int_0^t \sum_{k=1}^\infty 1_{[S^k_\varepsilon,
T^k_\varepsilon]}(r)\,dr>s\}.$$
Observe that the time-changed process
$(\rho_{\tau^{a,\varepsilon}_s}((a,\infty)),
s\geq 0)$ is obtained by patching together
the paths $(Y^{k,\varepsilon}_s,0\leq s
\leq T^k_\varepsilon-S^k_\varepsilon)$. 
Moreover, we have for every $k\geq 1$,
$$\sup_{0\leq s\leq
T^k_\varepsilon-S^k_\varepsilon}
(Y^{k,\varepsilon}_s-Z^{k,\varepsilon}_s)
=
\rho_{S^k_\varepsilon}((a,\infty))
=Y^{k,\varepsilon}_0.$$
From (\ref{reflec1}), we conclude that for every
$t\geq 0$,
\begin{equation}\label{reflec2}
\lim_{\varepsilon\rightarrow  0}
\pg\sup_{s\leq t}\, |U^\varepsilon_s
-\rho_{\tau^{a,\varepsilon}_s}((a,\infty))|\pd
=0.
\end{equation}
Notice that $\tau^{a,\varepsilon}_s\da \tau^a_s$
as $\varepsilon\da 0$ and recall that
for every $\varepsilon>0$, $U^\varepsilon$ is distributed
as the reflected L\'evy process $X-I$.
We then get from
(\ref{reflec2}) that the process 
$\langle\rho^a_s,1\rangle =
\rho_{\tau^{a}_s}((a,\infty))$
is distributed as the reflected proces
$X-I$, which completes the first step.

{\it Second step.} We will now verify that
$\rho^a$ can be obtained as a functional 
of the total mass process
$\langle\rho^a,1\rangle $ in the same way as
$\rho$ is obtained from $\langle\rho,1\rangle $.
It will be enough to argue on one
excursion of $\langle\rho^a,1\rangle $ away from $0$.
Thus, let $(u,v)$ be the interval 
corresponding to one such excursion.
We can associate with $(u,v)$ a unique 
connected component $(p,q)$ of the
open set $\{s,H_s>a\}$, such that
$\tau^a_{u+r}=p+r$ for every $r\in[0,v-u)$,
and $q=\tau^a_{v-}$. By the intermediate values
property, we must have $H_p=H_q=a$. 

We also claim
that $X_r>X_p$ for every $r\in (p,q)$. If this were not the case,
we could find $r\in(p,q)$ such that $X_{r}=\inf\{X_s,p\leq s\leq r\}$,
which forces $H_r\leq H_p=a$ and gives a contradiction. 

The previous
observations and the definition of the process $\rho$
imply that, for every $r\in(p,q)$,
the restriction of $\rho_r$ to $[0,a]$ is exactly $\rho_p=\rho_q$.
Define
$$\omega(r)=X_{(p+r)\wedge q}-X_p
=\langle\rho_{(p+r)\wedge q},1\rangle -\langle\rho_p,1\rangle 
=\langle\rho^a_{(u+r)\wedge v},1\rangle ,$$
so that $\omega$ is the excursion of $\langle\rho^a,1\rangle $
corresponding to $(u,v)$. The construction of the
process $\rho$ implies that, for $0<r<q-p=v-u$,
$$\rho_{p+r}=[\rho_p,\rho_r(\omega)],$$
and so, for the same values of $r$,
$$\rho^a_{u+r}=\rho_{r}(\omega).$$
This completes the second step of the proof.

\smallskip
{\it Third step.} It remains to prove that
$\rho^a$ is independent of the $\sigma$-field
$\h_a$. For $\varepsilon>0$,
denote by $\h^\varepsilon_a$
the $\sigma$-field generated by the 
processes
$$(X_{(T^k_\varepsilon+s)\wedge
S^{k+1}_\varepsilon},s\geq 0)$$
for $k=0,1,\ldots$ (by convention
$T^0_\varepsilon=0$),
and the negligible sets of $\g_\infty$. From our
construction
(in particular the fact that $X_s>X_{T^k_\varepsilon}$
for $s\in [S^k_\varepsilon,T^k_\varepsilon)$), it is easy to verify that
the processes $(\rho_{(T^k_\varepsilon+s)
\wedge S^{k+1}_\varepsilon},s\geq 0)$
are measurable with respect to
$\h^\varepsilon_a$, and since $H_t>a$
for $t\in (S^k_\varepsilon,T^k_\varepsilon)$,
it follows that $\h_a\subset\h^\varepsilon_a$.

By the arguments of the first step, the
processes $Z^{k,\varepsilon}$ are 
independent conditionally on $\h^\varepsilon_a$,
and the conditional law of $Z^{k,\varepsilon}$ 
is the law of an
independent copy 
of the
reflected process $X-I$, stopped 
when its local time at $0$
hits
$\rho_{S^k_\varepsilon}((a,\infty))$.
It follows that the process $U^\varepsilon$ of
the first step is independent of
$\h^\varepsilon_a$, hence also of $\h_a$.
By passing to the limit $\varepsilon\rightarrow  0$,
we obtain that the total mass process
$\langle\rho^a,1\rangle $ is independent of $\h_a$.
As we know that
$\rho^a$ can be reconstructed as a 
measurable functional of its total mass process,
this completes the proof.
\cq

\smallskip
We let $l^a=(l^a(s),s\geq 0)$ be the local time
at $0$ of $\langle\rho^a,1\rangle $, or equivalently
of $\rho^a$.

\begin{definition}
\label{localtimedef}
The 
local time at level $a$ and at time $s$ of the
height process $H$ is defined by
$$L^a_s=l^a\pg\int_0^s 1_{\{H_r>a\}}\,dr\pd.$$
\end{definition}

This definition will
be justified below: see in particular Proposition \ref{LTapprox}.

\subsection{Some properties of local times}

The next lemma can be seen as dual to
Lemma \ref{localstrip}.

\begin{lemma}
\label{heightinfimum}
For every $t\geq 0$,
$$\lim_{\varepsilon\rightarrow  0}
{1\over \varepsilon}\int_0^t
1_{\{H_s\leq \varepsilon\}}\,ds=-I_t,$$
in the $L^1$-norm.
\end{lemma}

\proof We use arguments similar to the proof
of Lemma \ref{localstrip}. We first establish
that for every $x>0$,
\begin{equation}\label{heighttech1}
\lim_{\varepsilon\rightarrow  0}
{1\over \varepsilon}\int_0^{T_{x}}
1_{\{H_s\leq \varepsilon\}}\,ds=x,
\end{equation}
in probability. Note that Proposition 
\ref{mesinv} gives for any nonnegative
measurable function $g$
$$N\cg\int_0^\sigma ds\,g(H_s)\cd
=\int_0^\infty da\,e^{-\alpha a}\,g(a).$$
Let $\omega^j$, $j\in J$ denote the 
excursions of $X-I$ away from $0$ and let
$(\alpha_j,\beta_j)$ be the corresponding time 
intervals. We already noticed that
$H_s=H_{s-\alpha_j}(\omega^j)$
for $s\in (\alpha_j,\beta_j)$. Hence, using
also the previous displayed formula, we have
\begin{equation}\label{heighttech11}
E\cg{1\over \varepsilon}\int_0^{T_{x}}
1_{\{H_s\leq\varepsilon\}}\,ds\cd
={x\over \varepsilon}N\cg
\int_0^\sigma 1_{\{H_s\leq \varepsilon\}}\,ds\cd
={x\over \varepsilon}\pg
{1-e^{-\alpha\varepsilon}\over \alpha}\pd\leq x,
\end{equation}
and in particular,
\begin{equation}\label{heighttech2}
\lim_{\varepsilon\rightarrow  0}
E\cg{1\over \varepsilon}\int_0^{T_{x}}
1_{\{H_s\leq\varepsilon\}}\,ds\cd=x.
\end{equation}
We then want to get a second moment estimate.
To this end, it is necessary to introduce 
a suitable truncation. Fix $K>0$. A slight
modification of the proof of (\ref{heighttech2})
gives
\begin{equation}\label{heighttech3}
\lim_{\varepsilon\rightarrow  0}
E\cg{1\over \varepsilon}\int_0^{T_{x}}
1_{\{H_s\leq\varepsilon\}}\,
1_{\{X_s-I_s\leq K\}}\,ds\cd=x.
\end{equation}

If $H^{(s)}$ denotes the
height process for the shifted process
$X^{(s)}_t=X_{s+t}-X_s$, the
bound $H^{(s)}_{t-s}\leq H_t$ 
(for $0\leq s\leq t$) is obvious
from our construction. We can use this simple
observation to bound
\begin{eqnarray*}&&N\cg\pg\int_0^{\sigma}
1_{\{H_s\leq\varepsilon\}}\,
1_{\{X_s\leq K\}}\,ds\pd^2\cd\\
&&\quad \leq 2
N\cg\int_0^{\sigma}ds\,
1_{\{H_s\leq\varepsilon\}}\,
1_{\{X_s\leq K\}}
\int_s^\sigma dt\,1_{\{H_t\leq\varepsilon\}}\cd\\
&&\quad \leq 2
N\cg\int_0^{\sigma}ds\,
1_{\{H_s\leq\varepsilon\}}\,
1_{\{X_s\leq K\}}
\int_s^\sigma
dt\,1_{\{H^{(s)}_{t-s}\leq\varepsilon\}}\cd\\ 
&&\quad 
= 2  N\cg\int_0^{\sigma}ds\,
1_{\{H_s\leq\varepsilon\}}\,
1_{\{X_s\leq K\}}
\,E_{X_s}\cg\int_0^{T_{0}}dt\,
1_{\{H_t\leq\varepsilon\}}\cd\cd\\
&&\quad \leq 2\varepsilon
N\cg\int_0^{\sigma}ds\,
1_{\{H_s\leq\varepsilon\}}\,
1_{\{X_s\leq K\}}
\,X_s\cd\qquad\qquad\qquad({\rm by}\
(\ref{heighttech11}))\\ 
&&\quad = 2\varepsilon
\int_0^\varepsilon dy\,
E[X_{L^{-1}(y)}\,1_{\{L^{-1}(y)<\infty,
X_{L^{-1}(y)}\leq K\}}]
\qquad({\rm Proposition\ \ref{keyinv}})\\
&&\quad \leq 2\varepsilon^2
E[X_{L^{-1}(\varepsilon)}\wedge K]\\
&&\quad =o(\varepsilon^2),
\end{eqnarray*}by dominated convergence. 
As in the proof of Lemma \ref{localstrip},
we can conclude from (\ref{heighttech2})
and the previous estimate that
$$\lim_{\varepsilon\rightarrow  0}
{1\over \varepsilon}
\int_0^{T_{x}}
1_{\{H_s\leq\varepsilon\}}\,
1_{\{X_s-I_s\leq K\}}\,ds=x$$
in the $L^2$-norm. Since this holds for
every $K>0$, (\ref{heighttech1}) follows.

From (\ref{heighttech1}) and a monotonicity
argument we deduce that the convergence 
of Lemma \ref{heightinfimum} holds in
probability. To get $L^1$-convergence, 
we need a few other estimates. First 
observe that
\begin{equation}\label{heighttech4}
E\cg\int_0^t 1_{\{H_s\leq\varepsilon\}}ds\cd
=\int_0^tds P[H_s\leq\varepsilon]
=\int_0^tds P[L_s\leq\varepsilon]
=E[L^{-1}(\varepsilon)\wedge t]
\leq C\,\varepsilon,
\end{equation}
with a constant $C$ depending only on $t$
(in the last bound we applied to $L^{-1}$
an estimate valid for any subordinator).
Then,
\begin{eqnarray*}E\cg\pg\int_0^t 1_{\{H_s\leq\varepsilon\}}ds
\pd^2\cd
&=&2\,E\cg\int\!\!\int_{\{0<r<s<t\}} drds
1_{\{H_r\leq\varepsilon\}}
1_{\{H_s\leq\varepsilon\}}\cd\\
&\leq& 2\,E\cg\int\!\!\int_{\{0<r<s<t\}} drds
1_{\{H_r\leq\varepsilon\}}
1_{\{H^{(r)}_{s-r}\leq\varepsilon\}}\cd\\
&=&2\,E\cg\int_0^t dr\,
1_{\{H_r\leq\varepsilon\}}
E\cg\int_0^{t-r}
ds\,1_{\{H_s\leq\varepsilon\}}\cd\cd\\
&\leq&2\pg E\cg\int_0^t dr\,
1_{\{H_r\leq\varepsilon\}}\cd\pd^2.
\end{eqnarray*}As a consequence of the last estimate and
(\ref{heighttech4}), the 
variables $\varepsilon^{-1}\int_0^t
1_{\{H_s\leq\varepsilon\}}ds$, $\varepsilon>0$
are bounded in $L^2$. This completes the
proof of Lemma \ref{heightinfimum}.\cq

\smallskip
We can now give a useful approximation result
for local times of the height process.

\begin{proposition}
\label{LTapprox}
For every $t\geq 0$, 
$$\lim_{\varepsilon\rightarrow  0}
\sup_{a\geq 0}E\cg
\sup_{s\leq t}\Big|\varepsilon^{-1}
\int_0^s 1_{\{a<H_r\leq a+\varepsilon\}}\,dr-L^a_s
\Big|\cd=0.$$
Similarly, for every $t\geq 0$,
$$\lim_{\varepsilon\rightarrow  0}
\sup_{a\geq \varepsilon}E\cg
\sup_{s\leq t}\Big|
\varepsilon^{-1}\int_0^s 1_{\{a-\varepsilon<H_r
\leq
a\}}\,dr-L^a_s
\Big|\cd=0.$$
There exists a jointly measurable
modification of the collection
$(L^a_s,a\geq 0,s\geq 0)$, which is 
continuous and nondecreasing in the
variable $s$, and such that, a.s. for any
nonnegative measurable function $g$
on $\R_+$ and any $s\geq 0$,
\begin{equation}\label{densityoccup}
\int_0^s
g(H_r)\,dr=\int_{\R_+}g(a)\,L^a_s\,da.
\end{equation}
\end{proposition}

\proof First consider the case $a=0$. Then,
$\rho^0=\rho$ and $L^0_t=l^0(t)=-I_t$.
Lemma \ref{heightinfimum} and a simple
monotonicity argument, using the continuity of $L^0_s$, give
\begin{equation}\label{approx1}
\lim_{\varepsilon\rightarrow  0}
E\cg
\sup_{s\leq t}\Big|\varepsilon^{-1}
\int_0^s 1_{\{0<H_r\leq\varepsilon\}}\,dr-L^0_s
\Big|\cd=0.
\end{equation}
For $a>0$, set $A^a_t=\int_0^t
1_{\{H_s>a\}}\,ds$. Note that
$\{a<H_s\leq a+\varepsilon\}
=\{\rho_s((a,\infty))>0\}\cap
\{\rho_s((a+\varepsilon,\infty))=0\}$,
and so
\begin{eqnarray*}\int_0^s 1_{\{a<H_r\leq a+\varepsilon\}}\,dr
&=&\int_0^t 1_{\{\rho_s((a,\infty))>0\}\cap
\{\rho_s((a+\varepsilon,\infty))=0\}}\,ds\\
&=&\int_0^{A^a_t} 1_{\{\rho^a_r((\varepsilon,
\infty))=0\}}dr\\
&=&\int_0^{A^a_t} 1_{\{0<H^a_r\leq \varepsilon\}}
dr,
\end{eqnarray*}where $H^a_t=H(\rho^a_t)$. The first convergence
of the proposition then follows from
(\ref{approx1}), the trivial bound $A^a_t\leq t$
and the fact that $\rho^a$ has the same
distribution as $\rho$.

The second convergence is easily derived 
from the first one by elementary arguments.
Let us only sketch the method. 
For any fixed $\delta>0$, we can choose 
$\varepsilon_0>0$
sufficiently small so that for every
$a\geq 0$, $\varepsilon\in(0,\varepsilon_0]$,
\begin{equation}\label{approx2}
E\cg
\Big|\varepsilon^{-1}
\int_0^t 1_{\{a<H_r\leq a+\varepsilon\}}\,dr-L^a_s
\Big|\cd\leq \delta.
\end{equation}
Then, if $0<\varepsilon<\varepsilon_0\wedge a$,
$$E\cg\Big|\varepsilon^{-1}
\int_0^t 1_{\{a-\varepsilon<H_r\leq
a\}}\,dr-
\varepsilon_0^{-1}
\int_0^t 1_{\{a-\varepsilon<H_r\leq
a-\varepsilon+\varepsilon_0\}}\,dr\Big|\cd
<2\delta.$$
However, if $\varepsilon$ is very small
in comparison with $\varepsilon_0$, one 
also gets the bound
$$E\cg\Big|
\varepsilon_0^{-1}
\int_0^t 1_{\{a-\varepsilon<H_r\leq
a-\varepsilon+\varepsilon_0\}}\,dr
-(\varepsilon_0-\varepsilon)^{-1}
\int_0^t 1_{\{a<H_r\leq a+\varepsilon_0
-\varepsilon\}}\,dr\Big|\cd\leq \delta.$$
We get the desired result by combining the
last two bounds and (\ref{approx2}).

The existence of a jointly measurable 
modification of the process $(L^a_s,a\geq 0,s\geq 0)$
that satisfies the density of occupation
time formula (\ref{densityoccup}) follows 
from the first assertion of the proposition
by standard arguments.
\cq

From now on, we will only deal with the
jointly measurable modification of the
local times $(L^a_s,a\geq 0,s\geq 0)$
given by Proposition \ref{LTapprox}.
We observe that it is easy to extend
the definition of these local times 
under the excursion measure $N$. 
First notice that, as an obvious
consequence of the first assertion 
of Proposition \ref{LTapprox}, we have also
for every $a\geq 0$, $t\geq 0$
\begin{equation}\label{LTapproxbis}
\lim_{\varepsilon\rightarrow  0}
E\cg
\sup_{0\leq r\leq s\leq t}\Big|\varepsilon^{-1}
\int_r^s 1_{\{a<H_u\leq
a+\varepsilon\}}\,du-(L^a_s-L^a_r)
\Big|\cd=0.
\end{equation}
Then, let $V$ be a
measurable subset of $\D(\R_+,\R)$
such that $N[V]<\infty$. For instance we
may take $V=\{\sup_{s\geq 0}X_s>\delta\}$
for $\delta>0$. By considering the
first excursion of $X-I$ that belongs to $V$,
and then using (\ref{LTapproxbis}), we
immediately obtain the existence under $N$
of a continuous increasing process,
still denoted by $(L^a_s,s\geq 0)$,
such that
$$\sup_{s\leq t}\Big|\varepsilon^{-1}
\int_0^s 1_{\{a<H_r\leq a+\varepsilon\}}\,dr-L^a_s\Big|
\build{\la}_{\varepsilon\rightarrow  0}^{} 0$$
in $N$-measure. More precisely, for any $V$ such that $N[V]<\infty$,
\begin{equation}\label{LTapproxexc}
\lim_{\varepsilon\rightarrow  0}
N\cg 1_V\,
\sup_{s\leq t}\Big|\varepsilon^{-1}
\int_0^s 1_{\{a<H_r\leq a+\varepsilon\}}\,dr-L^a_s
\Big|\cd=0.
\end{equation} 
The next corollary will now be 
a consequence of Proposition \ref{keyinv}.
We use the notation introduced before 
Proposition \ref{mesinv}.

\begin{corollary}
\label{localinvar}
For any nonnegative measurable function
$F$ on $\D(\R_+,\R)$, and every $a\geq 0$,
$$N\cg \int_0^\sigma
dL^a_s\,F(\wh X^{(s)}_{r\wedge s},r\geq 0)\cd
=E[1_{\{L^{-1}(a)<\infty\}}\,F(X_{r\wedge
L^{-1}(a)},r\geq 0)].$$
In particular, for any nonnegative measurable function $F$
on $M_f(\R_+)$,
$$N\cg \int_0^\sigma
dL^a_s\,F(\rho_s)\cd
=e^{-\alpha a}\,E[F(J_a)].$$
\end{corollary}

\proof We may assume that $F$ is bounded and
continuous. Then let $h$ be a 
nonnegative continuous
function on $\R_+$, which vanishes 
outside $[\delta,A]$, for some
$0<\delta<A<\infty$. For the first identity, it is enough to 
prove that
$$N\cg \int_0^\sigma
dL^a_s\,h(s)\,F(\wh X^{(s)}_{r\wedge s},r\geq 0)\cd
=E[1_{\{L^{-1}(a)<\infty\}}\,
h(L^{-1}(a))\,F(X_{r\wedge
L^{-1}(a)},r\geq 0)].$$
Notice that the mapping
$s\la (\wh X^{(s)}_t,t\geq 0)$
is continuous except possibly on a
countable set that is not charged by
the measure $dL^a_s$. From
(\ref{LTapproxexc}), applied with
$V=\{\omega,\sigma(\omega)>\delta\}$,
and then Proposition \ref{keyinv},
we get
\begin{eqnarray*}&&N\cg \int_0^\sigma
dL^a_s\,h(s)\,F(\wh X^{(s)}_{r\wedge s},r\geq 0)\cd\\
&&\qquad=\lim_{\varepsilon\rightarrow  0}
N\cg{1\over \varepsilon}\int_0^\sigma
ds\,1_{\{a<H_s\leq a+\varepsilon\}}\,
h(s)\,F(\wh X^{(s)}_{r\wedge s},r\geq 0)\cd\\
&&\qquad=\lim_{\varepsilon\rightarrow  0}
{1\over \varepsilon}E\cg
\int_{a\wedge L_\infty}^{(a+\varepsilon)
\wedge L_\infty}dx\,h(L^{-1}(x))\,
F(X_{t\wedge L^{-1}(x)},t\geq 0)\cd\\
&&\qquad=E[1_{\{L^{-1}(a)<\infty\}}\,
h(L^{-1}(a))\,F(X_{r\wedge
L^{-1}(a)},r\geq 0)],
\end{eqnarray*}which completes the proof of the first assertion. The second
assertion follows from the first one in the same way
as Proposition \ref{mesinv} was derived from Proposition
\ref{keyinv}.
\cq

\medskip
We conclude this section with some remarks that will be useful in the 
applications developed below. Let $x>0$ and let $(\alpha_j,\beta_j)$, resp.
$\omega_j$, $j\in J$, denote the excursion intervals, resp.
the excursions of $X-I$ before time $T_x$. For every $a>0$, we have $P$ a.s.
\begin{equation} 
\label{decompLT}
L^a_{T_x}=\sum_{j\in J} L^a_{\sigma(\omega_j)}(\omega_j).
\end{equation}
A first inequality is easily derived by writing
$$L^a_{T_x}\geq \int_0^{T_x}d_sL^a_s\,1_{\{X_s>I_s\}}
=\sum_{j\in J}(L^a_{\beta_j}-L^a_{\alpha_j})=\sum_{j\in J} L^a_{\sigma(\omega_j)}(\omega_j)$$
where the last equality follows from the approximations of local time. The converse inequality seems
to require a different argument in our general setting. Observe that, by excursion theory
and then Proposition \ref{mesinv},
\ba E[L^a_{T_x}]&\leq&\liminf_{k\rightarrow \infty} E\Big[{1\over \varepsilon_k}\int_0^{T_x} ds
\,1_{\{a<H_s<a+\varepsilon_k\}}\Big]\\
&=&\liminf_{k\rightarrow \infty} E\Big[\sum_{j\in J}{1\over \varepsilon_k}\int_0^{\sigma(\omega_j)} ds
\,1_{\{a<H_s(\omega_j)<a+\varepsilon_k\}}\Big]\\
&=&\liminf_{k\rightarrow \infty}x\,N\Big({1\over \varepsilon_k}\int_0^{\sigma} ds
\,1_{\{a<H_s<a+\varepsilon_k\}}\Big)\\
&=&\liminf_{k\rightarrow \infty}{x\over \varepsilon_k}\int_a^{a+\varepsilon_k} db\,e^{-\alpha b}\\
&=&x\,e^{-\alpha a}
\ea
whereas Corollary \ref{localinvar} (with $F=1$) gives $E[\sum_{j\in J} L^a_{\sigma(\omega_j)}]
=x\,N(L^a_\sigma)=x\,e^{-\alpha a}$. This readily yields (\ref{decompLT}).

Let us finally observe that we can extend the definition of
the local times $L^a_s$ to the process $\rho$ started at
a general initial value $\mu\in M_f(\R_+)$. In view of
forthcoming applications consider the case when 
$\mu$ is supported on $[0,a)$, for $a>0$. Then, the previous
method can be used to construct a continuous increasing 
process $(L^a_s(\rho^\mu),s\geq 0)$ such that
$$L^a_s(\rho^\mu)=\lim_{\varepsilon\rightarrow  0}
{1\over \varepsilon}\int_0^s dr\,1_{\{a<H(\rho^\mu_r)<a+\varepsilon\}}$$
in probability (or even in the $L^1$-norm). Indeed the arguments of the
proof of Proposition \ref{reflec} remain valid when $\rho$ is 
replaced by $\rho^\mu$, and the construction and approximation
of $L^a_s(\rho^\mu)$ follow. Recall the notation 
$\tau_0=\inf\{s\geq 0:\rho^\mu_s=0\}$ and observe that $\tau_0=T_x$
if $x=\langle \mu,1\rangle$. Let $(\alpha_j,\beta_j)$, 
$\omega_j$, $j\in J$ be as above and set $r_j=H(k_{-I_{\alpha_j}}\mu)$. Then we have
\begin{equation}\label{localshift}
L^a_{\tau_0}(\rho^\mu)=\sum_{j\in J} L^{a-r_j}_{\beta_j-\alpha_j}(\omega_j).
\end{equation}
The proof is much similar to that of (\ref{decompLT}):
The fact that the left side 
of (\ref{localshift}) is greater than the right side is easy from our
approximations of local time. The equality is then obtained
from a first-moment argument, using Proposition \ref{potker} and
Fatou's lemma to handle the left side.  

\section{Three applications}

\subsection{The Ray-Knight theorem}

Recall that the $\psi$-continuous-state branching process (in short
the $\psi$-CSBP) is the Markov process $(Y_a,a\geq 0)$ with values in $\R_+$
whose transition kernels are characterized by their
Laplace transform: For $\lambda>0$ and $b>a$,
$$E[\exp-\lambda Y_b\mid Y_a]=\exp(-Y_a\,u_{b-a}(\lambda)),$$
where $u_t(\lambda)$, $t\geq 0$ is the unique nonnegative solution
of the integral equation
\begin{equation}\label{integralE}
u_t(\lambda)+\int_0^t \psi(u_s(\lambda))\,ds=\lambda.
\end{equation}

\begin{theorem}
\label{RK} 
Let $x>0$. The process $(L^a_{T_x},a\geq 0)$ is 
a $\psi$-CSBP started at $x$. 
\end{theorem}

\proof First observe that $L^a_{T_x}$ is $\h_a$-measurable. This is trivial
for $a=0$ since $L^0_{T_x}=x$. For $a>0$, note that,
if
$$T^a_{x}=\inf\{s\geq 0:X_{\wt\tau^a_s}=-x\},$$
we have
$$\int_0^{T_x}ds\,1_{\{a-\varepsilon
<H_s\leq a\}}=
\int_0^{T^a_x}ds\,1_{\{a-\varepsilon
<H_{\tilde\tau^a_s}\leq a\}}.$$
and the right-hand side is measurable with
respect to the $\sigma$-field $\h_a$. The 
measurability of $L^a_{T_{x}}$ with 
respect to $\h_a$ then follows from the 
second convergence of Proposition \ref{LTapprox}.

\medskip
We then verify that the function 
$$u_a(\lambda)=N[1-e^{-\lambda L^a_\sigma}]\quad (a>0),\qquad u_0(\lambda)=\lambda$$
solves equation (\ref{integralE}). From the strong Markov property
of $\rho$ under the excursion measure $N$, we get for $a>0$
$$
u_a(\lambda)=\lambda N\cg \int_0^\sigma dL^a_s\,e^{-\lambda(L^a_\sigma-L^a_s)}\cd
=\lambda N\cg \int_0^\sigma dL^a_s\,F(\rho_s)\cd,$$
where, for $\mu\in M_f(\R_+)$,
$F(\mu)=E[\exp(-\lambda L^a_{\tau_0}(\rho^\mu))]$. By 
Corollary \ref{localinvar}, we can concentrate on the case when 
$\mu$ is supported on $[0,a)$, and then 
(\ref{localshift}) gives
\begin{eqnarray*}F(\mu)&=&\exp\pg -\int_0^{\langle\mu,1\rangle} du\,
N[1-\exp(-\lambda L^{a-H(k_{-u}\mu)}_\sigma)]\pd\\
&=&\exp\pg -\int \mu(dr)\,N[1-\exp(-\lambda L^{a-r}_\sigma)]\pd.
\end{eqnarray*}Hence, using again Corollary \ref{localinvar},
\begin{eqnarray*}u_a(\lambda)&=&\lambda\,N\cg \int_0^\sigma dL^a_s 
\exp( -\int \rho_s(dr)u_{a-r}(\lambda))\cd\\
&=&\lambda\,e^{-\alpha a}
E\cg\exp\pg -\int J_a(dr)\,u_{a-r}(\lambda)\pd\cd\\
&=&\lambda\,\exp\pg-\int_0^a dr\,\wt\psi(u_{a-r}(\lambda))\pd.
\end{eqnarray*}It is a simple matter to verify that (\ref{integralE})
follows from this last equality.

\medskip
By (\ref{decompLT}) and excursion theory, we have
\begin{equation}\label{OnedimLT}
E[\exp(-\lambda\,L^a_{T_x})]=\exp(-x\,N[1-\exp(-\lambda L^a_\sigma)])
=\exp(-x\,u_a(\lambda)).
\end{equation}
To complete the proof, it is enough to
show that for $0<a<b$,
\begin{equation}\label{MarkovLT}
E[\exp(-\lambda L^b_{T_x})\mid \h_a]=\exp(-u_{b-a}(\lambda)\,L^a_{T_x}).
\end{equation}
Recall the notation $\rho^a$ from Proposition \ref{reflec}, and
denote by $\wt L^c_s$ the local times of $H^a_s=H(\rho^a_s)$. 
From our approximations of local times, it is 
straightforward to verify that
$$L^b_{T_x}=\wt L^{b-a}_{A^a_{T_x}},$$
where $A^a_s=\int_0^s dr\,1_{\{H_r>a\}}$ as previously.  
Write $U=L^a_{T_x}$ to simplify notation. If 
$T^a_r=\inf\{t\geq 0: l^a(t)>r\}$, we have $A^a_{T_x}=T^a_U$
(note that $l^a({A^a_{T_x}})=U$ by construction,
and that the strong Markov property of $X$ at time $T_x$
implies $l^a(t)>l^a(A^a_{T_x})$ for every $t>A^a_{T_x}$). Hence,
$$E[\exp(-\lambda L^b_{T_x})\mid \h_a]
=E[\exp(-\lambda \wt L^{b-a}_{T^a_U})\mid \h_a]
=E[\exp(-\lambda \wt L^{b-a}_{T^a_u})]_{u=U},$$
where in the second equality, we use the fact that the process $(\wt
L^{b-a}_{T^a_u},u\geq 0)$ is
a functional of $\rho^a$, and is thus independent of $\h_a$
(Proposition \ref{reflec}), whereas $U=L^a_{T_x}$
is $\h_a$-measurable. Since $\rho^a$ has the same distribution as $\rho$,
$\wt
L^{b-a}_{T^a_u}$ and $L^{b-a}_{T_u}$ also have the same law, and
the desired result (\ref{MarkovLT})
follows from (\ref{OnedimLT}). \cq 

\begin{corollary}
\label{extinction}
For every $a\geq 0$, set 
$$v(a)=N\cg \sup_{0\leq s\leq \sigma} H_s > a\cd.$$
Then,
\begin{description}
\item{\rm (i)} If $\int_1^\infty {du\over\psi(u)}=\infty$, we have $v(a)=\infty$
for every $a>0$.

\item{\rm(ii)} If $\int_1^\infty {du\over\psi(u)}<\infty$, the function
$(v(a),a>0)$ is determined by
$$\int_{v(a)}^\infty {du\over\psi(u)}=a.$$
\end{description}
\end{corollary}

\proof By the lower semi-continuity of $H$, the condition $\sup_{0\leq s\leq \sigma}
H_s > a$ holds iff $A^a_\sigma>0$, and our construction shows that this is the
case iff $L^a_\sigma>0$. Thus,
$$v(a)=N[L^a_\sigma>0]=\lim_{\lambda\rightarrow \infty}N[1-e^{-\lambda L^a_\sigma}]
=\lim_{\lambda\rightarrow \infty} u_a(\lambda),$$
with the notation of the proof of Theorem \ref{RK}. From (\ref{integralE}), we have
$$\int_{u_a(\lambda)}^\lambda {du\over\psi(u)}=a,$$
and the desired result follows. \cq

\subsection{The continuity of the height process}

We now use Corollary \ref{extinction} to give a necessary and sufficient condition for
the path continuity of the height process $H$. 

\begin{theorem}
\label{continuityheight}
The process $H$ has continuous sample paths $P$ a.s. iff $\int_1^\infty
{du\over\psi(u)}<\infty$.
\end{theorem}

\proof By excursion theory, we have
$$P\cg \sup_{0\leq s\leq T_x} H_s>a\cd =1-\exp(-xv(a)).$$
By Corollary \ref{extinction} (i), we see that 
$H$ cannot have continuous paths if
$\int_1^\infty {du\over\psi(u)}=\infty$.

Assume that $\int_1^\infty {du\over\psi(u)}<\infty$.
The previous formula and the property $v(a)<\infty$ imply that
\begin{equation}\label{limitzero}
\lim_{t\downarrow 0} H_t=0\qquad P\ {\rm a.s.}
\end{equation}

The path continuity of $H$ will follow from Lemma \ref{IVP} if
we can show that for every fixed interval $[a,a+h]$, $h>0$, the 
number of upcrossings of $H$ along $[a,a+h]$ is a.s. finite
on every finite time interval.
Set $\gamma_0=0$ and define by induction for every $n\geq 1$,
\begin{eqnarray*}
&&\delta_n=\inf\{t\geq \gamma_{n-1}:H_t\geq a+h\},\\
&&\gamma_n=\inf\{t\geq \delta_n:H_t\leq a\}.
\end{eqnarray*}
Note that $H_{\gamma_n}\leq a$ by the lower semi-continuity of $H$. 
On the other hand, as a straightforward consequence of
(\ref{strongMarkov}), we have a.s. for every $t\geq 0$,
$$H_{\gamma_n+t}\leq H_{\gamma_n}+H^{(\gamma_n)}_t.$$
Therefore $\delta_{n+1}-\gamma_n\geq \kappa_n$, if
$\kappa_n=\inf\{t\geq 0:H^{(\gamma_n)}_t\geq h\}$. The strong Markov
property implies that the variables $\kappa_n$
are i.i.d. . Furthermore, $\kappa_n>0$ a.s. by (\ref{limitzero}). 
It follows that $\delta_n\ua\infty$ as $n\ua\infty$, which completes the
proof. \cq

\medskip
It is easy to see that the condition $\int_1^\infty
{du\over\psi(u)}<\infty$ is also necessary and sufficient for $H$ to
have continuous sample paths $N$ a.e. On the other hand, we may
consider the process $\rho$ started at an arbitrary initial value
$\mu\in M_f(\R^d)$, as defined by formula 
(\ref{rhoinitial}), and ask about the sample path continuity of
$H(\rho_s)$. Clearly, the answer will be no if the support 
of $\mu$ is not connected. For this reason, we introduce the set
$M^0_f$ which consists of all measures $\mu\in M_f(\R_+)$
such that $H(\mu)<\infty$ and $\supp \mu=[0,H(\mu)]$. By convention the
zero measure also belongs to $M^0_f$. 

From (\ref{rhoinitial}) and Lemma \ref{techrho}, it is easy to verify that
the process $\rho$ started at an initial value $\mu\in M^0_f$ will
remain forever in $M^0_f$, and furthermore $H(\rho_s)$
will have continuous sample paths a.s. Therefore, we may restrict the
state space of $\rho$ to $M^0_f$. This restriction will be
needed in Chapter 4. 

\subsection{H\"older continuity of the height process}

In view of applications in Chapter 4, we now discuss the H\"older
continuity properties of the height process. We assume that the
condition $\int_1^\infty du/\psi(u)<\infty$ holds
so that $H$ has continuous sample paths by Theorem \ref{continuityheight}.
We set
$$ \gamma  =\sup \{ r \geq 0  : 
\lim_{\lambda\rightarrow \infty}\lambda^{-r} \psi (\lambda )  =+\infty \}.$$
The convexity of $\psi$ implies that $\gamma\geq 1$.

\begin{theorem}
\label{Holder-height}
The height process $H$ is $P$-a.s. locally H\"older continuous with exponent $\alpha$
for any $\alpha\in(0,1-1/\gamma)$, and is $P$-a.s. not locally H\" older continuous with
exponent $\alpha$ if $\alpha>1-1/\gamma$. 
\end{theorem}

\proof We rely on the following key lemma. Recall the notation
$\wh L^{(t)}$ for the local time at $0$ of $\wh X^{(t)}-\wh S^{(t)}$
(cf Section 1.2).

\begin{lemma}
\label{Hold-heitec}
Let $t\geq 0$ and $s>0$. Then $P$ a.s.,
\ba
&&H_{t+s}-\inf_{r\in[t,t+s]}H_r=H(\rho^{(t)}_s)\,,\\
&&H_{t}-\inf_{r\in[t,t+s]}H_r=\wh L^{(t)}_{t\wedge R}\,,
\ea
where $R=\inf\{r\geq 0:\wh X^{(t)}_r>-I^{(t)}_s\}$ $(\inf\emptyset =\infty)$.
\end{lemma}

\proof From (\ref{strongMarkov}), we get, a.s. for every $r>0$,
\begin{equation} 
\label{strongMarkov-bis}
H_{t+r}=H(k_{-I^{(t)}_r}\rho_t)+H(\rho^{(t)}_r).
\end{equation}
From this it follows that
$$\inf_{r\in[t,t+s]}H_r=H(k_{-I^{(t)}_s}\rho_t)$$
and the minimum is indeed attained at the (a.s.\,unique) time $v\in[t,t+s]$
such that $X_v=I^t_{t+s}$. The first assertion of the lemma now follows by combining the
last equality with (\ref{strongMarkov-bis}) written with $r=s$.

Let us turn to the second assertion. If $I_t\geq I^t_{t+s}$, 
then on one hand $X_v=I_v$ and
$\inf_{r\in[t,t+s]}H_r=H_v=0$, on the other hand, $R=\infty$, and the second
assertion reduces to $H_t=\wh L^{(t)}_t$ which is the definition of $H_t$. 
Therefore we can assume that $I_t<I^t_{t+s}$. Let
$$u=\sup\{r\in[0,t]:X_{r-}<I^t_{t+s}\}.$$
As in the proof of Proposition \ref{explor}, we have
$$H_u=H(k_{-I^{(t)}_s}\rho_t)=\inf_{r\in[t,t+s]}H_r.$$
On the other hand, the construction of the height process
shows that the equality $H_r=\wh L^{(t)}_t-\wh L^{(t)}_{t-r}$
holds simultaneously for all $r\in[0,t]$ such that $X_{r-}\leq I^r_t$
(cf Lemma \ref{techheight}). In particular for $r=u$ we get
$$H_t-\inf_{r\in[t,t+s]}H_r=H_t-H_u=\wh L^{(t)}_t-(\wh L^{(t)}_t-\wh 
L^{(t)}_{t-u})=\wh L^{(t)}_{t-u}.$$
To complete the proof, simply note that we have $t-u=R$ on the event 
$\{I_t<I^t_{t+s}\}$.\cq

\medskip
To simplify notation we set $\varphi(\lambda)=\lambda/\psi^{-1}(\lambda)$.
The right-continuous inverse $L^{-1}$
of $L$ is a subordinator with Laplace exponent $\varphi$: See Theorem VII.4 (ii)
in \cite{Be}, and note that the constant $c$ in this statement is equal to $1$
under our normalization of local time (compare with Lemma \ref{Laplaceladder}).

\begin{lemma}
\label{moments-height}
For every $t\geq 0$, $s>0$ and $q>0$, 
$$ E [ |H_{t+s} - \inf_{r\in[t, t+s]} H_r |^q  ] \leq C_q\, \varphi (1/s )^{-q }\; ,$$
and 
$$ E[ |H_t - \inf_{r\in[t, t+s]} H_r |^q ]  \leq C_q\, \varphi (1/s )^{-q} \; ,$$
where $C_q= e \Gamma (q +1 ) $ is a finite constant depending only on $q $.
\end{lemma}

\proof Recall that $H(\rho_s)=H_s\build=_{}^{\rm (d)}L_s$. 
From Lemma \ref{Hold-heitec} we have
$$ E[| H_{t+s } - \inf_{r\in[t, t+s ]} H_r |^q ]= 
E[L^q_s] = q \int_0^{+\infty } x^{q-1} P [
L_s > x ]dx  \; . $$
However,
$$ P[L_s > x ]= P[s > L^{-1 }(x) ] \leq e\,E[ \exp ( - L^{-1}(x) /s ) ] =
e \exp ( -x \varphi (1/s ) ) . $$
\noindent Thus 
$$ E[ | H_{t+s } - \inf_{[t, t+s ]} H |^q ] \leq  eq 
\int_0^{+\infty } x^{q-1} \exp (-x \varphi (1/s ) )\,dx =
C_q \varphi ( 1/s )^{-q } \; . $$
This completes the proof 
of the first assertion. 
 
In order to prove the second one, first note that $I^{(t)}_s$ is independent
of ${\cal G}_t$ and therefore also of the time-reversed process $\wh X^{(t)}$.
Writing $\tau_a=\inf\{r\geq 0:S_r>a\}$, we get from the second
assertion of Lemma \ref{Hold-heitec}
$$ E[ | H_t - \inf_{r\in[t, t+s ]} H_r |^q ]\leq \int_{ [0, +\infty ) } P[-I_s \in da]\,
E[L^q_{\tau_a } ]. $$
Note that
$$ E[L^q_{\tau_a  }] = q\int_0^{+\infty } x^{q-1 } P[L_{\tau_a } >x ]\,dx \,,$$
and that $P[L_{\tau_a}>x]=P[S_{L^{-1}(x)}<a]$. It follows that
$$  E[| H_t - \inf_{r\in[t, t+s ]} H_r |^q ]\leq q\int_{ [0, +\infty ) }P[-I_s \in da
]
\int_0^{+\infty } dx \,x^{q-1 } P[S_{ L^{-1 }(x)} < a ]\; . $$
An integration by parts leads to
$$  E[| H_t - \inf_{r\in[t, t+s ]} H_r |^q ]\leq q \int_0^{+\infty } dx\,x^{q-1 } 
\int_{ [0, +\infty ) } P [S_{ L^{-1 }(x)}\in db ]
P[ -I_s >b]\; . $$
However 
$$ P[-I_s >b ]= P[T_{b }<s ] \leq 
e\,E[\exp (-T_{b }/s ) ]=e\, \exp (-b \psi^{-1} 
(1/s ) ) $$
since we know (\cite{Be} Theorem VII.1) that $(T_b,b\geq 0)$ 
is a subordinator with exponent $\psi^{-1}$.
Recalling Lemma \ref{Laplaceladder}, we get
\ba
E[| H_t - \inf_{r\in[t, t+s ]} H_r |^q ]&\leq  &
eq\int_0^{+\infty } dx\,x^{q-1 } E[\exp ( -\psi^{-1 } ( 1/s ) S_{L^{-1 }(x) } ) ]\\
&=&eq\int_0^{+\infty } dx\,x^{q-1 }\,\exp(-x/(s\psi^{-1}(1/s)))  \\
&=&  C_q
\varphi (1/s )^{-q}.
\ea
This completes the proof of Lemma \ref{moments-height}. \cq

\medskip
\noindent{\bf Proof of Theorem \ref{Holder-height}.} From Lemma \ref{moments-height} and
an elementary inequality, we get for every $t\geq 0$, $s>0$ and $q>0$
$$E[|H_{t+s}-H_t|^q]\leq 2^{q+1}C_q\,\varphi(1/s)^q.$$
Let $\alpha\in(0,1-1/\gamma)$. Then $(1-\alpha)^{-1}<\gamma$ 
and thus $\lambda^{-(1-\alpha)^{-1}}\psi(\lambda)$
tends to $\infty$ as $\lambda\rightarrow \infty$. It easily follows 
that $\lambda^{\alpha-1}\psi^{-1}(\lambda)$
tends to $0$ and so $\lambda^{-\alpha}\varphi(\lambda)$ 
tends to $\infty$ as $\lambda\rightarrow \infty$.
The previous bound then yields the existence 
of a constant $C$ depending on $q$
and $\alpha$ such that for every $t\geq 0$ and $s\in(0,1]$,
$$E[|H_{t+s}-H_t|^q]\leq C\,s^{q\alpha}.$$
The classical Kolmogorov lemma gives the first assertion of the theorem.

To prove the second assertion, observe that for every $\alpha>0$ and $A>0$,
$$ P[H_s < A s^{\alpha } ] = P[ L_s <  A s^{\alpha } ]
=P[ s < L^{-1 }(As^{\alpha }) ] \; . $$
Then use the elementary inequality
$$P[s < L^{-1 }(As^{\alpha } )] \leq {e\over e-1}\,
E[ 1- \exp (-L^{-1 }(As^{\alpha } ) /s ) ],$$
which leads to
$$P[H_s \leq  A s^{\alpha } ]\leq {e\over e-1}\, (1- \exp (-As^{\alpha } \varphi (1/s ) ))
.$$
If  $ \alpha >1-1/\gamma $, we 
can find a sequence $(s_n )$ decreasing to zero such $s_n^{\alpha } 
\varphi (1/s_n )$ tends to $0$. Thus, for 
any $A>0$ 
$$ \lim_{n \rightarrow  \infty } P[ H_{s_n } \leq A s_n^{\alpha } ] =0  \; ,$$
and it easily follows that $\displaystyle{\limsup_{s\rightarrow  0} s^{-\alpha}H_s=\infty}$, $P$
a.s.,  which completes the proof.\cq

\chapter{Convergence of Galton-Watson trees}

\section{Preliminaries}

Our goal in this chapter is to study the convergence
in distribution of Galton-Watson trees, under the
assumption that the associated Galton-Watson processes,
suitably rescaled, converge in distribution to
a continuous-state branching process. To give a precise
meaning to the convergence of trees, we will code
Galton-Watson trees by a discrete height process,
and we will establish the convergence 
of these (rescaled) discrete processes to the 
continuous height process of the previous chapter.
We will also prove that similar convergences
hold when the discrete height processes are
replaced by the contour processes of the trees.

\smallskip
Let us introduce the basic objects considered
in this chapter. For every $p\geq 1$, let $\mu_p$
be a subcritical or critical offspring distribution. That is,
$\mu_p$ is a probability distribution on $\Z_+=\{0,1,\ldots\}$
such that
$$\sum_{k=0}^\infty k\,\mu_p(k)\leq 1.$$
We systematically exclude the trivial cases where $\mu_p(1)=1$ or $\mu_p(0)=1$.
We also define another probability measure $\nu_p$ 
on $\{-1,0,1,\ldots\}$ by setting $\nu_p(k)=\mu_p(k+1)$
for every $k\geq -1$. 

We denote by $V^p=(V^p_k,k=0,1,2,\ldots)$ a discrete-time random walk
on $\Z$ with jump distribution $\nu_p$ and started at $0$. We also
denote by $Y^p=(Y^p_k,k=0,1,2,\ldots)$
a Galton-Watson branching process with offspring distribution
$\mu_p$ started at $Y^p_0=p$. 

Finally, we consider a L\'evy process $X=(X_t,t\geq 0)$ 
started at the origin and satisfying assumptions
(H1) and (H2) of Chapter 1. As in Chapter 1, we write $\psi$
for the Laplace exponent of $X$. We denote by
$Y=(Y_t,t\geq 0)$ a $\psi$-continuous-state branching
process started at $Y_0=1$.

The following variant of a result due to Grimvall \cite{Gr} plays an
important role in our approach. Unless otherwise specified the convergence
in distribution of processes is  in the functional sense, that is in the
sense of the weak  convergence of the laws of the processes on the
Skorokhod space
$\D(\R_+,\R)$. We will use the notation $\build{\la}_{}^{\rm(fd)}$
to indicate weak convergence of finite-dimensional marginals.

For $a\in \R$, $[a]$ denotes the integer part of $a$.

\begin{theorem}
\label{Grim}
Let $(\gamma_p,,p=1,2,\ldots)$ be a nondecreasing sequence of
positive integers converging to $\infty$.
The convergence in distribution
\begin{equation}
\label{Grim1}
\left(p^{-1}Y^p_{[\gamma_pt]}\,,\,t\geq 0\right)
\build{\longrightarrow}_{p\rightarrow \infty}^{{\rm (d)}} (Y_t,t\geq 0)
\end{equation}
holds if and only if
\begin{equation}
\label{Grim2}
\left(p^{-1}V^p_{[p\gamma_pt]}\,,\,t\geq 0\right)
\build{\longrightarrow}_{p\rightarrow \infty}^{{\rm (d)}} (X_t,t\geq 0).
\end{equation}
\end{theorem}

\proof By standard results on the convergence of 
triangular arrays (see e.g. Theorem 2.7 in Skorokhod \cite{Sko}),
the functional convergence (\ref{Grim2})
holds iff
\begin{equation}\label{Grim2bis}
p^{-1}V^p_{p\gamma_p}
\build{\longrightarrow}_{p\rightarrow \infty}^{{\rm (d)}} X_1.
\end{equation}Fix any sequence $p_1<p_2<\cdots<p_k<\cdots$ such that
$\gamma_{p_1}<\gamma_{p_2}<\cdots$. If $j=\gamma_{p_k}$
for some $k\geq 1$, set $c_j=p_k$, $V^{(j)}=V^{p_k}$ and let $\theta_j$ be
the  probability measure on $\R$ defined by $\theta_j({n\over
c_j})=\nu_{p_k}(n)$ for every integer $n\geq -1$. 
Then (\ref{Grim2bis}) is equivalent to saying that
$${1\over c_j}V^{(j)}_{jc_j}
\build{\longrightarrow}_{j\rightarrow \infty}^{{\rm (d)}} X_1$$
for any choice of the sequence $p_1<p_2<\cdots$. 
Equivalently the convolutions $(\theta_j)^{*jc_j}$ converge 
weakly to the law of $X_1$. By Theorems 3.4 and 3.1 of
Grimvall \cite{Gr}, this property holds iff the
convergence (\ref{Grim1}) holds along the sequence $(p_k)$. 
(Note that condition (b) in Theorem 3.4 of \cite{Gr}
is automatically satisfied here since we restrict our attention
to the (sub)critical case.) This completes the proof. \cq

\section{The convergence of finite-dimensional marginals}

From now on, we suppose that assumption (H3) holds
in addition to (H1) and (H2). Thus we can consider the 
height process $H=(H_t,t\geq 0)$ of Chapter 1.

For every $p\geq 1$, let $H^p=(H^p_k,k\geq 0)$
be the discrete height process associated with
a sequence of independent Galton-Watson trees
with offspring distribution $\mu_p$ (cf Section 0.2).
As was observed in Section 0.2, we may and will assume
that the processes $H^p$ and $V^p$ are related by the
formula
\begin{equation}\label{heightwalk}
H^p_k={\rm
Card}\{j\in\{0,1,\ldots,k-1\}:V^p_j=\inf_{j\leq l\leq k}V^p_l\}.
\end{equation}The following theorem sharpens a result of \cite{LGLJ1}.

\begin{theorem}
\label{marginconv}
Under either of the convergences {\rm (\ref{Grim1})}
or {\rm(\ref{Grim2})}, we have also
\begin{equation}
\label{heightconv}
\left({1\over \gamma_p}H^p_{[p\gamma_pt]}\,,\,t\geq 0\right)
\build{\longrightarrow}_{p\rightarrow \infty}^{{\rm (fd)}} (H_t,t\geq 0).
\end{equation}
\end{theorem}

\proof
Let $f_0$ be a truncation function, that is a 
bounded continuous function from $\R$ into $\R$
such that $f_0(x)=x$ for every $x$ belonging to
a neighborhood of $0$. By standard results on
the convergence of rescaled random walks
(see e.g. Theorem II.3.2 in \cite{Jacod}), the
convergence (\ref{Grim2}) holds iff the following 
three conditions are satisfied:
\begin{eqnarray*}
({\rm C1})\qquad&&\lim_{p\rightarrow \infty}p\gamma_p \sum_{k
=-1}^\infty f_0({k\over p})\,\nu_p(k)=-\alpha+
\int_0^\infty (f_0(r)-r)\,\pi(dr)\\
({\rm C2})
\qquad&& \lim_{p\rightarrow \infty}p\gamma_p \sum_{k
=-1}^\infty f_0({k\over p})^2\,\nu_p(k)=2\beta+
\int_0^\infty f_0(r)^2\,\pi(dr)\\
({\rm C3})
\qquad&& \lim_{p\rightarrow \infty} p\gamma_p \sum_{k=
-1}^\infty \varphi({k\over
p})\,\nu_p(k)=
\int_0^\infty \varphi(r)\,\pi(dr),
\end{eqnarray*}
for any bounded continuous function $\varphi$ on $\R$
that vanishes on a neighborhood of $0$.

\smallskip

By (\ref{heightwalk}) and time-reversal, $H^p_k$ has the
same distribution as 
$$\Lambda^{(p)}_k={\rm Card}\{j\in\{1,\ldots,k\}:V^p_j=\sup_{0\leq
l\leq j}V^p_l\}.$$
Without loss of generality, the Skorokhod representation
theorem allows us to assume that the convergence
\begin{equation} 
\label{DonsLevy}
\left(p^{-1}V^p_{[p\gamma_pt]}\,,\,t\geq 0\right)
\build{\longrightarrow}_{p\rightarrow \infty}^{} (X_t,t\geq 0)
\end{equation}
holds a.s. in the sense of Skorokhod's topology. Suppose we can prove
that for every $t>0$, 
\begin{equation}
\label{locconv}
\lim_{p\rightarrow \infty} {\gamma_p^{-1}}\,\Lambda^{(p)}_{[p\gamma_pt]}=L_t
\end{equation}
in probability. (Here $L=(L_t,t\geq 0)$ is the local time process
of $X-S$ at level $0$ as in Chapter 1.) Then a simple 
time-reversal argument implies that $\gamma_p^{-1}H^p_{[p\gamma_pt]}$
also converges in probability to $\widehat L^{(t)}_t=H_t$,
with the notation of Chapter 1. Therefore the proof of Theorem
\ref{marginconv} reduces to showing that (\ref{locconv}) holds.

We first consider the case where $\int_{(0,1)} r\pi(dr)=\infty$. We introduce the stopping times
$(\tau^p_k)_{k\geq 0}$  defined recursively as follows:
\begin{eqnarray*}
\tau^p_0
&=&0,\\
\tau^p_{m+1}&=&\inf\{n>\tau^p_m:V^p_n
\geq V^p_{\tau^p_m}\}.
\end{eqnarray*}
Conditionally on the event $\{\tau^p_m<\infty\} $, the random variable
$1_{\{\tau^p_{m+1}<\infty\}}(V^p_{\tau^p_{m+1}}-V^p_{\tau^p_m})$ 
is independent of the past of $V^p$ up to time $\tau^p_m$ and
has the same law as
$1_{\{\tau^p_1<\infty\}}V^p_{\tau^p_1}$. Also recall the 
classical equality (cf (5.4) in \cite{LGLJ1}):
\begin{equation}
\label{WH}
P[\tau^p_1<\infty,\,V^p_{\tau^p_1}=j]=\nu_p([j,\infty))\ ,\qquad j\geq 0.
\end{equation}

\smallskip
  For every $u>\delta>0$, set:
\begin{eqnarray*}
\kappa(\delta ,u)&=&\int_0^\infty
\pi(dr)\int_0^r dx\,1_{(\delta,u]}(x)=\int_0^\infty 
\pi (dr)\,\Big((r-\delta)^{+}\wedge 
(u-\delta)\Big)\,,\\
\kappa_p(\delta, u)&=&{{\displaystyle
\sum_{p\delta< j\leq pu} 
\nu_p([j,\infty))}\over
{\displaystyle \sum_{j\geq 0}
\nu_p([j,\infty))}}=
P[p\delta< V^p_{\tau^p_1}\leq pu\mid
\tau^p_1<\infty]\,,\\
L^{\delta ,u}_t &=&{\rm Card}\{s\leq t,
\Delta S_{s}\in (\delta 
,u] \}\,,\\
l^{p,\delta ,u}_k &=&{\rm
Card}\{j<k,\,\overline{V}^p_j+p\delta
< V^p_{j+1}\leq\overline{V}^p_j+
pu\}\,,\\
\end{eqnarray*}
where  $\overline{V}^p_j=\sup\{V^p_i,0
\leq i\leq j\}$. Note that $\kappa(\delta,u)\uparrow\infty$
as $\delta\downarrow 0$, by our assumption $\int_{(0,1)} r\pi(dr)=\infty$.  From the a.s. convergence of the 
processes $p^{-1}V_{[p\gamma_pt]}$, we have
\begin{equation}\label{margintech1}
\lim_{p\rightarrow \infty} l^{p,\delta ,u
}_{[p\gamma_pt]}=L^{\delta ,u}_t \ ,\quad {\rm a.s.}
\end{equation}
(Note that $P[\Delta S_s=a\hbox{ for some }s>0]=0$ for every fixed $a>0$, by (\ref{joint}).)
By applying excursion
theory to the process
$X-S$ and using formula (\ref{joint}), one easily gets
for every $u>0$
\begin{equation}
\label{margintech2}
\lim_{\delta\rightarrow  0}\kappa(
\delta ,u)^{-1}
L^{\delta, u}_t=L_t \ ,\quad {\rm a.s.}
\end{equation}\smallskip
We claim that we have also
\begin{equation} \label{margintech3}
\lim_{p\rightarrow \infty}{\gamma_p}\,\kappa_p(\delta ,u)
=\int_0^\infty \Big((r-\delta)^{+}\wedge
(u-\delta )\Big)
\,\pi(dr)=\kappa(\delta ,u).\end{equation}
To get this convergence, first apply (C3)
to the function
$\varphi(x)=(x-\delta)^{+}\wedge
(u-\delta )$. It follows that 
$$\lim_{p\rightarrow \infty}p\,\gamma_p
\sum_{k=-1}^{\infty}\nu_p(k)\,\Big(\big(\frac{k}{p}
-\delta)^+\wedge (u-\delta )\Big)
=\kappa(\delta ,u).$$
On the other hand, it is elementary to verify that
$$\left|p\,\gamma_p
\sum_{k=-1}^{\infty}\nu_p(k)\,\Big(\big(\frac{k}{p}
-\delta)^+\wedge (u-\delta )\Big)-
\gamma_p\sum_{p\delta< j\leq pu} 
\nu_p([j,\infty))\right|
\leq \gamma_p\sum_{k\geq \delta p}\nu_p(k)$$
and the right-hand side tends to $0$ by (C3). 
Thus we get
$$\lim_{p\rightarrow \infty}
\gamma_p\sum_{p\delta< j\leq pu}
\nu_p([j,\infty))=\kappa(\delta ,u).$$
Furthermore, as a simple consequence of (C1) and the
(sub)criticality of $\mu_p$, we have also
$$\sum_{j=0}^\infty \nu_p([j,\infty))
=1+\sum_{k=-1}^\infty k\nu_p(k)\build{\la}_{p\rightarrow \infty}^{}
1.$$
(This can also be obtained from (\ref{WH}) and the
weak convergence (\ref{Grim2}).)
Our claim (\ref{margintech3}) now follows.

Finally, we can also obtain a relation between
$l^{p,\delta,u}_k$ and $\Lambda^{(p)}_{[p\gamma_pt]}$. Simply 
observe that
conditional on $\{\tau^p_k<\infty\}$, $l^{p,\delta,u}_{\tau^p_k}$ is the sum of 
$k$ independent Bernoulli variables with parameter
$\kappa_p(\delta,u)$. Fix an integer $A>0$ and set $A_p=\gamma_pA+1$.
From Doob's inequality, we easily get (see \cite{LGLJ1},
p.249 for a similar estimate)
$$E\left[\sup_{0\leq j\leq
\tau^p_{A_p}}\Big|{1\over \gamma_p}
\big(\Lambda^{(p)}_j-\kappa_p(\delta ,u)^{-1}l^{p,\delta
,u}_j\big)\Big|^2\right]\leq
{8(A+1)\over \gamma_p}\,\kappa_p(\delta ,u)^{-1}.
$$
Hence, using (\ref{margintech3}), we have
\begin{equation}\label{margintech4}
\limsup_{p\rightarrow \infty}
E\Big[\sup_{j\leq \tau^p_{A_p}}\Big| {1\over \gamma_p}
\big(\Lambda^{(p)}_j-\kappa_p(\delta ,u)^{-1}l^{p,\delta ,u}_j\big)
\Big|^2\Big]
\leq {8(A+1)\over \kappa(\delta ,u)}.
\end{equation}To complete the proof, let $\varepsilon>0$ and first choose
$A$ large enough so that $P
[L_t\geq A-3\varepsilon]<\varepsilon$. If $u>0$ is fixed, we can 
use (\ref{margintech2})
and (\ref{margintech4}) to pick $\delta>0$ small enough
and then $p_0=p_0(\delta)$ so that
\begin{equation} 
\label{margintech17}
P \Big[\Big|\kappa(\delta ,u)^{-1}L^{
\delta ,u}_t-L_t\Big|>\varepsilon\Big]<
\varepsilon
\end{equation}
and 
\begin{equation} \label{margintech5}
P\left[\sup_{j\leq \tau^p_{A_p}}\Big| {1\over \gamma_p}
\big(\Lambda^{(p)}_j-\kappa_p(\delta ,u)^{
-1}l^{p,\delta ,u}_j\big)\Big|>
\varepsilon\right]<\varepsilon,
\qquad {\rm if}\ p\geq p_0.
\end{equation}From (\ref{margintech1}) and (\ref{margintech3}), we can also
find $p_1(\delta)$ so that for every $p\geq p_1$,
\begin{equation} \label{margintech6}
P\left[\ \Big|{1\over
\gamma_p\kappa_p(\delta ,u)}\,l^{p,
\delta ,u}_{[p\gamma_pt]}
-\kappa(\delta ,u)^{-1}L^{\delta ,u}_t
\Big|>\varepsilon\right]<\varepsilon.
\end{equation}By combining the previous estimates (\ref{margintech17}),
(\ref{margintech5}) and (\ref{margintech6}), we get
for $p\geq p_0\vee p_1$  
\begin{equation} \label{margintech7}
P \Big[\Big|{1\over \gamma_p}\,
\Lambda^{(p)}_{[p\gamma_pt]}-L_t\Big|>3
\varepsilon\Big]
\leq 3\varepsilon+P[[p\gamma_pt]>
\tau^p_{A_p}].
\end{equation}Furthermore, by using (\ref{margintech5}) and
then (\ref{margintech17}) and (\ref{margintech6}), we have 
for $p$ sufficiently large,
\begin{eqnarray*}
P[ \tau^p_{A_p} < [p\gamma_pt ] ] & \leq & \varepsilon+
P \Big[{
1\over \gamma_p\kappa_p(\delta ,u)}\,l^{p,\delta ,u}_{[p\gamma_p t]}
\geq A-\varepsilon\Big] \\
&\leq &3\varepsilon+P[L_t\geq A-3
\varepsilon]\\
&\leq &4\varepsilon
\end{eqnarray*}
from our choice of $A$. Combining this estimate with 
(\ref{margintech7}) completes the proof of 
(\ref{locconv})
in the case $\int_0^1 r\pi(dr)=\infty$.

It remains to treat the case where $\int_0^1 r\pi(dr)<\infty$. In that case,
(H3) implies that $\beta>0$, and we know from (\ref{localLeb})
that 
$$ L_t = {1\over \beta} m( \{ S_s ; s\leq t \} ) \; .$$
Furthermore, (\ref{joint}) and the assumption $\int_{(0, 1)} r \pi (dr ) <\infty $ imply
that for any 
$t >0 $ , 
$$\card \{ s\in [0, t] ;\; \Delta S_s >0 \}
< \infty \; , \quad {\rm a.s.} $$
For every $\delta>0$ and $t\geq 0$, we set
$$ \widetilde{S}^{\delta}_t = S_t - \sum_{ s\in [0, t ] } 
{\bf 1}_{(\delta, \infty )} (\Delta S_s ) \Delta S_s \; .$$
By the previous remarks, we have a.s. for $\delta$ small enough,
\begin{equation}
\label{ecart}
\widetilde{S}^{\delta}_t=S_t - \sum_{ s\in [0, t ] } 
\Delta S_s=m( \{ S_s ; s\leq t \} )= \beta L_t \; .
\end{equation}
Let us use the same notation $\tau^p_m$, $\ov{V}^p_j$ as in the 
case $\int_0^1 r\pi(dr)=\infty$, and also set
for any $m \geq 1 $, 
$$ d^p_m ={\bf 1}_{ \{ \tau^p_{m}<\infty \} } 
( V^p_{ \tau^p_{m} }- V^p_{ \tau^p_{m-1} } ) $$
and 
$$ \widetilde{S}^{\delta, p }_m = \sum_{n \geq 1 } d^p_n 
{\bf 1}_{ \{ d^p_n \leq p\delta ,\; \tau^p_{n} \leq  m\} } \; . $$
The convergence (\ref{DonsLevy}) implies that 
\begin{equation} 
\label{betatech8}
\left( {1\over p} 
\overline{V}^p_{[p\gamma_p s]} ,{s \geq 0 }\right)\ 
\build{\la}_{p\rightarrow  \infty}^{}\ (S_s ,s\geq 0)\; ,  \quad {\rm a.s.},
\end{equation}
and, for every $t\geq 0$,
$${1\over p}  \sum_{n \geq 1 } d^p_n
{\bf 1}_{ \{ d^p_n > p\delta,\; \tau^p_{n} \leq  
[p\gamma_p t] \} 
}\ \build{\la}_{p\rightarrow  \infty}^{}\ 
\sum_{ s\in [0, t ] } 
{\bf 1}_{(\delta, \infty )} (\Delta S_s ) \Delta S_s  \; , \quad {\rm a.s.}$$
Thus we have 
\begin{equation}
\label{convps}
\lim_{ p\rightarrow  \infty } {1\over p} \widetilde{S}^{\delta, p }_{[p\gamma_p t]}
= \widetilde{S}^{\delta}_t 
\quad{\rm a.s.} 
\end{equation}
The desired convergence (\ref{locconv}) is then a consequence
of (\ref{ecart}), (\ref{convps}) and the following result: For every $\varepsilon>0$,
\begin{equation}
\label{identif}
\lim_{\delta \rightarrow  0} \; 
\limsup_{p\rightarrow  \infty } \;
P [|  {1\over p} 
\widetilde{S}^{\delta, p }_{[p\gamma_p t]} \; -\; 
{\beta\over\gamma_p} \Lambda^p_{[p\gamma_pt]}| > \varepsilon
]= 0 \;.
\end{equation}

To prove (\ref{identif}), set
\begin{eqnarray*}
\alpha_1 (p, \delta ) &=& 
E \left[ d^p_1 \,{\bf 1}_{\{d^p_1 \leq \delta p , \tau^p_{1} <\infty\}}
\right] \; ; \\
\alpha_2 (p ,\delta ) &=& 
E \left[ (d^p_1)^2\, {\bf 1}_{\{d^p_1 \leq \delta p ,\, \tau^p_{1} <\infty \}}
\right] \; .
\end{eqnarray*}
Observe that 
$$ E \left[ (d^p_1 -\alpha_1 (p, \delta) )^2 \,{\bf 1}_{\{d^p_1 \leq \delta p ,\, \tau^p_{1} <\infty \}}
\right] \leq \alpha_2 (p ,\delta) \; .$$
Let $A>0$ be an integer and let $A_p=\gamma_pA+1$ as above. By Doob's inequality,
$$
E [ \sup_{1\leq m \leq A_p } | 
\widetilde{S}^{\delta, p }_{\tau^p_{m} } - m \alpha_1 (p,\delta)|^2 ] \leq 4 A_p \alpha_2 (p, \delta) \; . 
$$
Since
$$
\sup_{1\leq m \leq A_p } | 
\widetilde{S}^{\delta, p }_{\tau^p_{m} } - m \alpha_1 (p, \delta) |
= \sup_{1\leq j\leq \tau^p_{A_p }} | 
\widetilde{S}^{\delta, p }_{j } - \alpha_1 (p, \delta) \Lambda^p_{j} | \; .
$$
we have 
\begin{equation}
\label{doob}
E\Big[ \sup_{0\leq j\leq \tau^p_{A_p }}| 
{1\over p}
\widetilde{S}^{\delta, p }_{j} - 
{ \alpha_1 (p, \delta)\over p} \Lambda^p_{
j} |^2 
\Big] \leq {4A_p \over p^2} \alpha_2 (p, \delta) \; .
\end{equation}

We now claim that

\begin{equation}
\label{betatech1}
\lim_{ p\rightarrow  \infty } {\gamma_p \over p}
\alpha_1 (p, \delta) = \beta + {1\over 2}
\int_{(0, \infty)} (r\wedge \delta)^2\pi (dr)  \ 
\build{\la}_{\delta\rightarrow  0}^{}\ \beta  \; , 
\end{equation}
and
\begin{equation}
\label{betatech2}
 \lim_{\delta\rightarrow  0 } \; 
\limsup_{ p \rightarrow  \infty } \; {\gamma_p \over p^2} \; 
\alpha_2 (p ,\delta) =  0 \; . 
\end{equation}
To verify (\ref{betatech1}), note that, by (\ref{WH}),
\begin{equation} 
\label{betatech3}
{\gamma_p\over p} \alpha_1(p,\delta)={\gamma_p\over p}\sum_{j=0}^{[\delta p]}j\,\nu_p([j,\infty))
={p\gamma_p\over 2}\sum_{k=0}^\infty \nu_p(k)\Big({k\over p}\wedge {[\delta p]\over p}\Big)
\Big({k\over p}\wedge {[\delta p]\over p}+{1\over p}\Big)
\end{equation}
We now apply (C1) and (C2) with the truncation function $f_0(x)=(x\wedge \delta)\vee (-\delta)$. 
Multiplying by $p^{-1}$ the convergence in (C1) and adding the one in (C2), we get
$$\lim_{p\rightarrow \infty} p\gamma_p\sum_{k=0}^\infty \nu_p(k)
\Big({k\over p}\wedge {\delta}\Big)
\Big({k\over p}\wedge {\delta}+{1\over p}\Big)=2\beta+\int_{(0,\infty)} (r\wedge \delta)^2\pi(dr).$$
Comparing with (\ref{betatech3}) we immediately get (\ref{betatech1}). The proof of 
(\ref{betatech2}) is analogous.

By (\ref{doob}) and an elementary inequality, we have
$$
E \Big[ \sup_{0\leq j\leq \tau^p_{A_p }} | 
{1\over p}
\widetilde{S}^{\delta, p }_{j } - 
{\beta \over\gamma_p} 
\Lambda^p_{j}|^2 
\Big] 
\leq {8A_p\over p^2} \alpha_2 (p, \delta) + 2 (
{A_p\over \gamma_p })^2 \,( \beta - {\gamma_p\over p}
\alpha_1 (p, \delta) )^2 \; .
$$
Thus, (\ref{betatech1}) and (\ref{betatech2}) imply that for any $A>0 $ 
\begin{equation}
\label{premierpt}
\lim_{\delta \rightarrow  0 } 
\limsup_{p\rightarrow  \infty}
E \Big[ \sup_{0\leq j\leq \tau^p_{A_p }} | 
{1\over p}
\widetilde{S}^{\delta, p }_{j } - 
{\beta\over \gamma_p} 
\Lambda^p_{j} |^2 
\Big] = 0 \; .
\end{equation}
It follows that
$$\lim_{\delta\rightarrow  0} \limsup_{p\rightarrow  \infty }
P \Big[|  {1\over p} \widetilde{S}^{\delta, p }_{[p\gamma_p t]} -
{\beta\over \gamma_p} \Lambda^p_{[p\gamma_pt]} | > \varepsilon
\Big] \leq  \limsup_{p\rightarrow  \infty } P [
\tau^p_{A_p } < [p \gamma_p t]  ]
\; . $$
However,
$$P [\tau^p_{A_p } < [p \gamma_p t]  ]\leq P[{1\over p}\ov{V}^p_{[p\gamma_pt]}\geq {1\over p}
\wt S^{\delta,p}_{\tau^p_{A_p}}]$$
and by (\ref{premierpt}) the right side is bounded above for $p$ large by $P[{1\over p}\ov{V}^p_{[p\gamma_pt]}\geq
{\beta\over \gamma_p}A_p-1]+\varepsilon_\delta$, where $\varepsilon_\delta\rightarrow  0$ as $\delta\rightarrow  0$. 
In view of (\ref{betatech8}), this is enough to conclude that
$$\lim_{A \rightarrow  \infty }
\limsup_{p\rightarrow  \infty } P [\tau^p_{A_p } < [p \gamma_p t]  ]= 0 \; , $$
and the desired result (\ref{identif}) follows. 
This completes the proof of 
(\ref{locconv}) and of Theorem \ref{marginconv}. \cq

\section{The functional convergence}

Our goal is now to discuss conditions that ensure that
the convergence of Theorem \ref{marginconv} holds in a 
functional sense. We assume that the function
$\psi$ satisfies the condition
\begin{equation}\label{continuitycond}
\int_1^\infty {du\over \psi(u)}<\infty.
\end{equation}By Theorem \ref{continuityheight}, this implies that the height process 
$(H_t,t\geq 0)$ has continuous sample paths. On the
other hand, if this condition does not hold, the paths
of the height process do not belong to any of the usual
functional spaces.

For every $p\geq 1$, we denote by $g^{(p)}$ the generating
function of $\mu_p$, and by $g^{(p)}_n=g^{(p)}\circ\cdots
\circ g^{(p)}$ the $n$-th
iterate  of $g^{(p)}$.

\begin{theorem}
\label{functionalconv}
Suppose that the convergences {\rm(\ref{Grim1})} and {\rm(\ref{Grim2})}
hold and that the continuity condition {\rm(\ref{continuitycond})}
is satisfied. Suppose in addition that for every $\delta>0$,
\begin{equation} \label{extinctech}
\liminf_{p\rightarrow \infty} g^{(p)}_{[\delta \gamma_p]}(0)^{p}>0.\end{equation}Then,
\begin{equation}
\label{functionconv}
\left(\gamma_p^{-1}H^p_{[p\gamma_pt]}\,,\,t\geq 0\right)
\build{\longrightarrow}_{p\rightarrow \infty}^{{\rm (d)}} (H_t,t\geq 0)
\end{equation}
in the sense of weak convergence on $\D(\R_+,\R_+)$.
\end{theorem}

Let us make some important remarks. Condition (\ref{extinctech})
can be restated in probabilistic terms as follows: For every $\delta>0$,
$$\liminf_{p\rightarrow \infty} P[Y^p_{[\delta \gamma_p]}=0]>0.$$
(As will follow from our results, this implies that the
extinction time of $Y^p$, scaled by $\gamma_p^{-1}$, converges in
distribution to the extinction time of $Y$, which is finite a.s. under
(\ref{continuitycond}).) It is easy to see that the condition
(\ref{extinctech}) is necessary for the conclusion (\ref{functionconv})
to hold. Indeed, suppose that (\ref{extinctech}) fails, so that 
there exists $\delta>0$ such that
$P[Y^p_{[\delta\gamma_p]}=0]$ converges to $0$ as $p\rightarrow \infty$, at least
along a suitable subsequence. Clearly, this convergence also holds (along
the same subsequence) if $Y^p$ starts at $[ap]$ instead 
of $p$, for any fixed $a>0$. From the definition
of the discrete height process, we get that
$$P\Big[\sup_{k\leq T^p_{[ap]}}H^p_k\geq [\delta\gamma_p]\Big]
\build{\la}_{p\rightarrow \infty}^{} 1\ ,$$
where $T^p_j=\inf\{k\geq 0:V^p_k=-j\}$. From (\ref{Grim2}), we know 
that $(p\gamma_p)^{-1}T^p_{[ap]}$ converges in distribution
to $T_a$. Since $T_a\da 0$ as $a\da 0$, a.s., we easily
conclude that, for every $\varepsilon>0$, 
$$P\Big[\sup_{t\leq \varepsilon}\gamma_p^{-1} H^p_{[p\gamma_pt]}\geq
{[\delta\gamma_p]\over \gamma_p}
\Big]\build{\la}_{p\rightarrow \infty}^{} 1\ ,$$
and thus (\ref{functionconv}) cannot hold.

\smallskip
On the other hand, one might think that the condition
(\ref{extinctech}) is automatically satisfied under 
(\ref{Grim1}) and (\ref{continuitycond}). Let us explain 
why this is not the case. Suppose for simplicity 
that $\psi$ is of the type
$$\psi(\lambda)=\alpha \lambda+\int_{(0,\infty)}\pi(dr)\,(e^{-\lambda
r}-1+\lambda r),$$
and for every $\varepsilon>0$ set
$$\psi_\varepsilon(\lambda)=\alpha
\lambda+\int_{(\varepsilon,\infty)}\pi(dr)\,(e^{-\lambda r}-1+\lambda r).$$
Note that $\psi_\varepsilon(\lambda)\leq C_\varepsilon\lambda$ and so 
$\int_1^\infty \psi_\varepsilon(\lambda)^{-1}d\lambda=\infty$. Thus, if
$Y^\varepsilon$ is a 
$\psi_\varepsilon$-CSBP started at $1$, we have $Y^\varepsilon_t>0$ for every
$t>0$ a.s. (Grey \cite{Grey}, Theorem 1). It is easy to verify that
$$(Y^\varepsilon_t,t\geq 0)\build{\la}_{\varepsilon\rightarrow 0}^{} (Y_t,t\geq 0)$$
at least in the sense of the weak convergence of finite-dimensional marginals.
Let us fix a sequence $(\varepsilon_k)$ decreasing to $0$. Then for every 
$k$, we can find a subcritical or critical offspring distribution $\nu_k$,
and two positive integers $p_k\geq k$ and $\gamma_k\geq k$, in such a way
that if $Z^k=(Z^k_j,j\geq 0)$ is a Galton-Watson process with
offspring distribution $\nu_k$ started at $Z^k_0=p_k$, the law of the
rescaled process
$$Z^{(k)}_t=(p_k)^{-1}Z^k_{[\gamma_kt]}$$
is arbitrarily close to the law of $Y^{\varepsilon_k}$. In particular,
we may assume that
$P[Z^{(k)}_k>0]>1-2^{-k}$, and that the rescaled processes 
$Z^{(k)}$ converge in distribution to $Y$. However, the extinction 
time of $Z^{(k)}$ converges in probability to $+\infty$, and
so the condition (\ref{extinctech}) cannot hold.

\smallskip
There is however a very important special case where 
(\ref{extinctech}) holds.

\begin{theorem}
\label{stableconv}
Suppose 
that $\mu_p=\mu$ for every $p$ and that
the convergence {\rm (\ref{Grim1})} holds. Then the condition {\rm(\ref{extinctech})}
is automatically satisfied and the conclusion
of Theorem \ref{functionalconv} holds.
\end{theorem}

As we will see in the proof, under the assumption of Theorem
\ref{stableconv}, the process $X$ must be stable with index
$\alpha\in(1,2]$.  Clearly the condition (\ref{continuitycond}) holds in
that case.

\medskip
\noindent{\bf Proof of Theorem \ref{functionalconv}.} To simplify notation,
we set $H^{(p)}_t=\gamma_p^{-1}H^p_{[p\gamma_pt]}$ and 
$V^{(p)}_t=p^{-1}V^p_{[p\gamma_pt]}$. In view of
Theorem \ref{marginconv}, the proof of Theorem \ref{functionalconv}
reduces to checking that the laws of the processes
$(H^{(p)}_t,t\geq 0)$ are tight in the set of 
probability measures on $\D(\R_+,\R)$. By standard results
(see e.g. Corollary 3.7.4 in \cite{EK}), it is enough to
verify the following two properties:

\smallskip
(i) For every $t\geq 0$ and $\eta>0$, there exists a
constant $K\geq 0$ such that
$$\liminf_{p\rightarrow  \infty }\ P [ H_{t}^{(p)}\leq K] \geq 1-\eta .$$

(ii) For every $T>0$ and $\delta>0$,
$$\lim_{n\rightarrow  \infty }\ \limsup_{p\rightarrow  \infty }\ P \Big[ \sup_{1\leq i\leq
2^{n}} \ \sup_{t\in [
(i-1)2^{-n}T ,{i}{2^{-n}T} ]} 
\ |H_{t}^{(p)}-H_{(i-1){2^{-n}T}}^{(p)} |>\delta \Big] =0.$$

Property (i) is immediate from the convergence of
finite-dimensional marginals. Thus the real problem
is to prove (ii). We fix $\delta>0$ and $T>0$ and first observe that
\begin{eqnarray}
\label{formu1}
&&P \Big[ \sup_{1\leq i\leq
2^{n}} \ \sup_{t\in [
(i-1)2^{-n}T ,{i}{2^{-n}T} ]} 
\ |H_{t}^{(p)}-H_{(i-1){2^{-n}T}}^{(p)} |>\delta \Big] \\
&&\quad \leq A_1(n,p)+A_2(n,p)+A_3(n,p)\nonumber
\end{eqnarray}
where
\begin{eqnarray*}
A_1(n,p)&=&P \Big[ \sup_{1\leq i\leq 2^{n}}
|H_{i2^{-n}T}^{(p)}-H_{(i-1){2^{-n}T}}^{(p)} | >\frac{
\delta}{5}
 \Big] \\
A_2(n,p)&=&P \Big[
\sup_{t\in [
(i-1)2^{-n}T ,i2^{n}T ]} 
\ H_{t}^{(p)}>H_{(i-1)2^{-n}T}
^{(p)}
+\frac{4\delta }{5}
\ {\rm for \ some\ } 1\leq i\leq 2^n\Big] \\
A_3(n,p)&=&P \Big[
\inf_{t\in [
(i-1)2^{-n}T ,i2^{n}T ]} 
\ H_{t}^{(p)}< H_{i2^{n}T}^{(p)} 
-\frac{4\delta }{5}\ {\rm for \ some\ } 1\leq i\leq 2^n\Big]
\end{eqnarray*}
The term $A_1$ is easy to bound. By the 
convergence of finite-dimensional marginals, we have
$$\limsup_{p\rightarrow  \infty } 
A_1(n,p)\leq 
P \Big[ \sup_{1\leq i\leq 2^{n}}
|H_{i2^{-n}T} -H_{(i-1){2^{-n}T}}|\geq \frac{
\delta}{5}
\Big]$$
and the path continuity of the process $H$ ensures that the
right-hand side tends to $0$ as $n\rightarrow  \infty$.

To bound the terms $A_2$ and $A_3$, we introduce the
stopping times $\tau^{(p)}_k$, $k\geq 0$ defined by
induction as follows:
\ba
&&\tau_{0}^{(p)}=0\\
&&\tau_{k+1}^{(p)}=\inf \{ t\geq \tau^{(p)}_k : H^{(p)}_{t} >
\inf_{\tau^{(p)}_k\leq r\leq t}H^{(p)}_r+\frac{\delta}{5}
\}.
\ea
Let $i\in\{1,\ldots,2^n\}$ be such that
\begin{equation}\label{functiontech1}
\sup_{t\in [
(i-1)2^{-n}T ,i2^{n}T ]} 
\ H_{t}^{(p)}>H_{(i-1)2^{-n}T}
^{(p)}
+\frac{4\delta }{5}.
\end{equation}Then it is clear that the interval $[
(i-1)2^{-n}T ,i2^{n}T ]$ must contain at least one of the
random times $\tau^{(p)}_k$, $k\geq 0$. Let
$\tau^{(p)}_j$ be the first such time. By construction we have
$$\sup_{t\in[(i-1)2^{-n}T,\tau^{(p)}_j)}H^{(p)}_t\leq 
H^{(p)}_{(i-1)2^{-n}T}+{\delta\over 5},$$
and since the positive jumps of $H^{(p)}$
are of size $\gamma_p^{-1}$, we get also
$$H^{(p)}_{\tau^{(p)}_j}\leq 
H^{(p)}_{(i-1)2^{-n}T}+{\delta\over 5}+\gamma_p^{-1}< 
H^{(p)}_{(i-1)2^{-n}T}+{2\delta\over 5}$$
provided that $\gamma_p>5/\delta$. From (\ref{functiontech1}), we have then
$$\sup_{t\in [
\tau^{(p)}_j ,i2^{n}T ]} 
\ H_{t}^{(p)}>H_{\tau^{(p)}_j}
^{(p)}
+\frac{\delta }{5},$$
which implies that $\tau^{(p)}_{j+1}\leq i2^{-n}T$.
Summarizing, we get for $p$ large enough so that $\gamma_p>5/\delta$
\begin{equation}\label{functiontech2}
A_2(n,p)\leq P\Big[\tau^{(p)}_k<T\ {\rm and}\ \tau^{(p)}_{k+1}
-\tau^{(p)}_k<2^{-n}T\ {\rm for\ some\ } k\geq 0\Big].
\end{equation}A similar argument gives exactly the same bound for the quantity
$A_3(n,p)$.

The following lemma is directly inspired from \cite{EK} p.134-135.

\begin{lemma}
\label{EKtech}
For every $x>0$ and $p\geq 1$, set 
$$G_p(x)=P\Big[\tau^{(p)}_k<T\hbox{ and }\tau^{(p)}_{k+1}
-\tau^{(p)}_k<x\ \hbox{for some } k\geq 0\Big]$$
and
$$F_p(x)=\sup_{k\geq 0}P\Big[\tau^{(p)}_k<T \hbox{ and}\ \tau^{(p)}_{k+1}
-\tau^{(p)}_k<x\Big].$$
Then, for every integer $L\geq 1$,
$$G_p(x)\leq L\,F_p(x)+L\,e^T\int_0^\infty dy\,e^{-Ly}\,F_p(y).$$
\end{lemma}

\proof For every integer $L\geq 1$, we have
\ba
G_p(x)&\leq&\sum_{k=0}^{L-1}P[\tau_{k}^{(p)}< T\ {\rm and}\ 
\tau_{k+1}^{(p)}-\tau_{k}^{(p)}<x]+P 
[\tau_{L}^{(p)}< T]\\
&\leq &LF_p(x)+e^{T} E\Big[{1}_{\{ \tau_{L
}^{(p)}< T\} }\exp \Big( 
-\sum_{k=0}^{L-1}(\tau_{k+1}^{(p)}-\tau_{k}^{(p)}) \Big)\Big]\\
&\leq &LF_p(x)+e^{T} \prod_{k=0}^{L-1} E\Big[{1}_{\{ \tau_{L}^{(p)}<T\}
} 
\exp (-L(\tau_{k+1}^{(p)}-\tau_{k}^{(p)}))\Big]^{1/L}\,.
\ea
Then observe that for every $k\in\{0,1,\ldots,L-1\}$,
\begin{eqnarray*}
E\Big[{1}_{\{ \tau_{L}^{(p)}<T\}} 
\exp (-L(\tau_{k+1}^{(p)}-\tau_{k}^{(p)}))\Big]
&\leq& E\Big[1_{\{ \tau_{k}^{(p)}<T\}}
\int_{\tau_{k+1}^{(p)}-\tau_{k}^{(p)}}^\infty dy\,Le^{-Ly}\Big]\\
&\leq& \int_0^\infty dy\,Le^{-Ly}\,F_p(y).
\end{eqnarray*}
The desired result follows.\cq

\medskip
Thanks to Lemma \ref{EKtech}, the limiting behavior of the right-hand
side of (\ref{functiontech2}) will be reduced to that of the
function $F_p(x)$. To handle $F_p(x)$, we use the next lemma.

\begin{lemma}
\label{functionkey}
The random variables $\tau^{(p)}_{k+1}-\tau^{(p)}_k$
are independent and identically distributed.
Under the assumptions of Theorem \ref{functionalconv}, we
have
$$\lim_{x\da 0}\Big(\limsup_{p\rightarrow \infty} P[\tau^{(p)}_1\leq x]\Big)=0.$$
\end{lemma}

We need a simple lemma. 

\begin{lemma}
\label{Markov-height}
Let $V$ be a random walk on $\Z$. For every $n\geq 0$, set
\begin{equation}
\label{keybis}
H^o_n=\card\{k\in\{0,1,\ldots,n-1\}:V_k=\inf_{k\leq j\leq n}V_j\}.
\end{equation}
Let $\tau$ be a stopping time of the filtration $({\cal F}^o_n)$ generated by $V$.
Then the process
$$\Big(H^o_{\tau+n}-\inf_{\tau\leq k\leq \tau+n}H^o_k,n\geq 0\Big)$$
is independent of ${\cal F}^o_\tau$ and has the same distribution as $(H^o_n,n\geq 0)$.
\end{lemma}

\proof By considering the first time after $\tau$
where the random walk $V$ attains its minimum over $[\tau,\tau+n]$, one easily gets
$$\inf_{\tau\leq k\leq \tau+n}H^o_k=\card\{k\in\{0,1,\ldots,\tau-1\}:V_k=\inf_{k\leq j\leq \tau+n}
V_j\}.$$ Hence,
\ba 
H^o_{\tau+n}-\inf_{\tau\leq k\leq \tau+n}H^o_k
&=&\card\{k\in\{\tau,\ldots,\tau+n-1\}:V_k=\inf_{k\leq j\leq \tau+n} V_j\}\cr
&=&\card\{k\in\{0,\ldots,n-1\}:V^\tau_k=\inf_{k\leq j\leq n} V^\tau_j\},
\ea
where $V^\tau$ denotes the shifted random walk $V^\tau_n=V_{\tau+n}-V_\tau$.
Since $V^\tau$ is independent of ${\cal F}_\tau$ and has the same distribution as
$V$, the desired result follows from the previous formula and (\ref{keybis}). \cq

\noindent{\bf Proof of Lemma \ref{functionkey}.} Fix $k\geq 1$ and set
for every $t\geq 0$,
$$
\wt H^{(p)}_t=H^{(p)}_{\tau^{(p)}_k+t}-\inf_{\tau^{(p)}_k\leq r\leq
\tau^{(p)}_k+t}H^{(p)}_r.
$$
As a consequence of Lemma \ref{Markov-height}, the process $(\wt H^{(p)}_t,t\geq 0)$ is
independent of the past of $V^{(p)}$ up to the stopping time $\tau^{(p)}_k$ and has the same
distribution as $(H^{(p)}_t,t\geq 0)$. Since by definition
$$\tau^{(p)}_{k+1}-\tau^{(p)}_k=\inf\{t\geq 0:\wt H^{(p)}_t>{\delta\over
5}\}$$
the first assertion of the lemma follows.

Let us turn to the second assertion. To simplify notation, we write 
$\delta'=\delta/5$. For every $\eta>0$, set
$$T^{(p)}_\eta=\inf\{t\geq 0:V^{(p)}_t=-{[p\eta]\over
p}\}.$$
Then,
$$P[\tau^{(p)}_1\leq x]
=P\Big[\sup_{s\leq x}H^{(p)}_s>\delta'\Big]
\leq P\Big[\sup_{s\leq T^{(p)}_\eta}H^{(p)}_s>\delta'\Big]
+P[T^{(p)}_\eta<x].$$
On one hand, 
$$\limsup_{p\rightarrow \infty} P[T^{(p)}_\eta<x]\leq P[T_\eta\leq x],$$
and for any choice of $\eta>0$, the right-hand side goes to zero as
$x\da 0$. On the other hand, the construction of the discrete
height process shows that the quantity
$$\sup_{s\leq T^{(p)}_\eta}H^{(p)}_s$$
is distributed as $\gamma_p^{-1}(M_p-1)$, where $M_p$ is the extinction
time of a Galton-Watson process with offspring distribution $\mu_p$,
started at $[p\eta]$. Hence,
$$P\Big[\sup_{s\leq T^{(p)}_\eta}H^{(p)}_s>\delta'\Big]
=1-g^{(p)}_{[\delta'\gamma_p]+1}(0)^{[p\eta]},$$
and our assumption (\ref{extinctech}) implies that
$$\lim_{\eta\rightarrow  0}\Big(\limsup_{p\rightarrow \infty}
P\Big[\sup_{s\leq T^{(p)}_\eta}H^{(p)}_s>\delta'\Big]\Big)=0.$$
The second assertion of the lemma now follows. \cq

\smallskip
We can now complete the proof of Theorem \ref{functionalconv}. 
Set:
$$F(x)=\limsup_{p\rightarrow \infty} F_p(x)\ ,\quad G(x)=\limsup_{p\rightarrow \infty} G_p(x).
$$
Lemma \ref{functionkey} immediately shows that $F(x)\da 0$ as 
$x\da 0$. On the other hand, we get from Lemma \ref{EKtech}
that for every integer $L\geq 1$,
$$G(x)\leq L\,F(x)+L\,e^T\int_0^\infty dy\,e^{-Ly}\,F(y).$$
It follows that we have also $G(x)\da 0$ as $x\da 0$. By 
(\ref{functiontech2}), this gives
$$\lim_{n\rightarrow  \infty}\Big(\limsup_{p\rightarrow \infty} A_2(n,p)\Big)=0,$$
and the same property holds for $A_3(n,p)$. This completes the 
proof of (ii) and of Theorem \ref{functionalconv}.\cq

\bigskip
\noindent{\bf Proof of Theorem \ref{stableconv}.}
We now assume that $\nu_p=\nu$ for every $p$ and so $g^{(p)}_n=g_n$.
We first observe that the process $X$ must be stable. This is not
immediate, since the convergence (\ref{Grim2}) 
a priori implies only that $\nu$ belongs to the domain of
partial attraction of the law of $X_1$, which is not enough to 
conclude. However, the conditions (C1) -- (C3), which are equivalent
to (\ref{Grim2}), immediately show that
the sequence 
$\gamma_{p}/\gamma_{p+1}$ converges to $1$ as $p\rightarrow \infty$. 
Then Theorem 2.3 in \cite{Mej} implies that $\nu$ belongs to the
domain of attraction of the law of $X_1$, and by classical results
the law of $X_1$ must be stable with index $\alpha\in (0,2]$. We can exclude
$\alpha\in(0,1]$ thanks to our assumptions (H2) and (H3) (the latter
is only needed to exclude the trivial case $\psi(\lambda)=c\lambda$). Thus 
$\alpha\in(1,2]$ and $\psi(\lambda)=c\,\lambda^\alpha$
for some $c>0$. As a consequence of (\ref{integralE}), we have
$E[e^{-\lambda Y_\delta}]=
\exp-(\lambda^{-\bar\alpha}+c\bar\alpha\delta)^{-1/\bar\alpha}$,
where $\bar \alpha=\alpha-1$. In particular,
$P[Y_\delta=0]=\exp-(c\bar\alpha\delta)^{-1/\bar\alpha}>0$.

Let $g=g_1$ be the generating function of $\mu$. We have
$g'(1)=\sum k\,\mu(k)=1$, because otherwise this would contradict
(\ref{Grim2}). From Theorem 2 in \cite{Fe}, p.577, the function 
$$\sum_{k\geq x} \mu(k)$$
must be regularly
varying as $x\rightarrow \infty$, with exponent $-\alpha$. Then note that
$$g(e^{-\lambda})-1+\lambda =
\sum_{k=0}^\infty \mu(k)\,(e^{-\lambda k}-1+\lambda k)
=\lambda\int_0^\infty dx(1-e^{-\lambda x})\sum_{k\geq x} \mu(k).$$
An elementary argument shows that $g(e^{-\lambda})-1+\lambda $
is also regularly varying as $\lambda\rightarrow  0$ with exponent $\alpha$. 
Put differently,
$$g(r)=r+(1-r)^\alpha L(1-r)\ ,\quad 0\leq r<1\ ,$$
where the function $L(x)$ is slowly varying as $x\rightarrow  0$. 
This is exactly what we need to apply a result of Slack \cite{Slack}.

Let $Z^{(p)}_1$ be a random variable distributed as 
$(1-g_{[\delta \gamma_p]}(0))$ times the value at time $[\delta\gamma_p]$
of a Galton-Watson process with offspring distribution $\mu$
started with one individual at time $0$ and
conditioned to be non-extinct at time $[\delta\gamma_p]$. Theorem 1
of \cite{Slack} implies that
$$Z^{(p)}_1\build{\la}_{p\rightarrow \infty}^{\rm (d)} U$$
where $U>0$ a.s. In particular, we can choose positive constants 
$c_0$ and $c_1$ so that $P[Z^{(p)}_1>c_0]>c_1$ for all 
$p$ sufficiently large. On the other hand, we have
$${1\over p}Y^p_{[\delta\gamma_p]}
\build{=}_{}^{\rm(d)}{1\over p(1-g_{[\delta \gamma_p]}(0))}
\Big(Z^{(p)}_1+\cdots+Z^{(p)}_{M_p}\Big)$$
where $Z^{(p)}_1,Z^{(p)}_2,\ldots$ are i.i.d., and $M_p$
is independent of the sequence $(Z^{(p)}_j)$ and has
a binomial ${B}(p,1-g_{[\delta \gamma_p]}(0))$ distribution.

It is now easy to obtain the condition (\ref{extinctech}). Fix $\delta>0$.
Clearly
it suffices to verify that the sequence $p(1-g_{[\delta \gamma_p]}(0))$
is bounded. If not the case, we can choose a sequence $(p_k)$
such that ${p_k}(1-g_{[\delta \gamma_{p_k}]}(0))$ converges to $\infty$. 
From the previous representation for the law
of ${1\over p}Y^p_{[\delta \gamma_p]}$, it then follows that
$$P\Big[{1\over {p_k}}Y^{p_k}_{[\delta \gamma_{p_k}]}> c_0c_1\Big]
\build{\la}_{k\rightarrow  \infty}^{} 1.$$
From  (\ref{Grim1}), we get that $P[Y_\delta\geq c_0c_1]=1$,
which gives a contradiction
since $P[Y_{\delta}=0]>0$. This completes
the proof of (\ref{extinctech}). 

Finally, since
(\ref{continuitycond}) holds, we can apply Theorem
\ref{functionalconv}.
\cq

\section{Convergence of contour processes}

In this section, we show that the limit theorems
obtained in the previous section for rescaled discrete height
processes can be formulated as well in terms of the contour processes
of the Galton-Watson trees. The proof relies on simple 
connections between the height process and the contour process
of a sequence of Galton-Watson trees.

To begin with, we consider a (subcritical or critical)
offspring distribution $\mu$, and a sequence of independent
$\mu$-Galton-Watson trees. Let $(H_n,n\geq 0)$ 
and $(C_t,t\geq 0)$ be respectively the 
height process and the contour process associated with 
this sequence of trees (see Section 0.2). We also set
$$K_n=2n-H_n.$$
Note that the sequence $K_n$ is strictly increasing and $K_n\geq n$.

Recall that the value at time $n$ of the height process corresponds to
the generation of the individual visited at time $n$, assuming that 
individuals are visited in lexicographical order one tree after another.
It is easily checked by induction on $n$ that $[K_n,K_{n+1}]$ is
exactly the time interval during which the contour process goes from
the individual $n$ to the individual $n+1$. From this observation, we get
$$\sup_{t\in[K_n,K_{n+1}]}|C_t-H_n|\leq |H_{n+1}-H_n|+1.$$
A more precise argument for this bound follows from the explicit
formula for $C_t$ in terms of the height process: For $t\in[K_n,K_{n+1}]$,
$$\begin{array}{lll}
&C_t=(H_n-(t-K_n))^+\qquad&\hbox{if }t\in[K_n,K_{n+1}-1],\\
\noalign{\smallskip}
&C_t=(H_{n+1}-(K_{n+1}-t))^+\qquad&\hbox{if }t\in[K_{n+1}-1,K_{n+1}].
\end{array}$$
These formulas are easily checked by induction on $n$.

Define a random function $f:\R_+\longrightarrow \Z_+$ by setting
$f(t)=n$ iff $t\in[K_n,K_{n+1})$. From the previous bound, we 
get for every integer $m\geq 1$,
\begin{equation}
\label{contour1}
\sup_{t\in[0,m]}|C_t-H_{f(t)}|\leq\sup_{t\in[0,K_m]}|C_t-H_{f(t)}|\leq 1+\sup_{n\leq m}|H_{n+1}-H_n|.
\end{equation}
Similarly, it follows from the definition of $K_n$ that
\begin{equation}
\label{contour2}
\sup_{t\in[0,m]}|f(t)-{t\over 2}|\leq\sup_{t\in[0,K_m]}|f(t)-{t\over 2}|
\leq {1\over 2}\sup_{n\leq m}H_n+1.
\end{equation}

We now come back to the setting of the previous sections, considering
for every $p\geq 1$
a sequence of independent Galton-Watson trees with offspring
distribution $\mu_p$. For every $p\geq 1$, we denote by 
$(C^p_t,t\geq 0)$ the corresponding contour process.

\begin{theorem}
\label{contour-conv}
Suppose that the convergences {\rm (\ref{Grim2})} and {\rm (\ref{functionconv})}
hold. Then,
\begin{equation}
\label{contour-convergence}
\left(\gamma_p^{-1}C^p_{p\gamma_pt}\,,\,t\geq 0\right)
\build{\longrightarrow}_{p\rightarrow \infty}^{{\rm (d)}} (H_{t/2},t\geq 0).
\end{equation}
In particular, {\rm(\ref{contour-convergence})} holds under the assumptions
of Theorem \ref{functionalconv} or those of Theorem \ref{stableconv}.
\end{theorem}

\noindent{\bf Proof.} For every $p\geq 1$, write $f_p$ for the analogue
of the function $f$ introduced above. Also set $\varphi_p(t)=(p\gamma_p)^{-1}
f_p(p\gamma_p t)$. By (\ref{contour1}), we have for every $m\geq 1$,
\begin{equation}
\label{cont-tech1}
\sup_{t\leq m}\Big|{1\over \gamma_p}C^p_{p\gamma_p t}-{1\over \gamma_p}H^p_{p\gamma_p \varphi_p(t)}\Big|
\leq{1\over \gamma_p}+{1\over \gamma_p}\sup_{n\leq mp\gamma_p}
|H^p_{n+1}-H^p_{n}|
\build{\longrightarrow}_{p\rightarrow \infty}^{} 0 
\end{equation}
in probability, by (\ref{functionconv}).

On the other hand, we get from (\ref{contour2})
\begin{equation}
\label{cont-tech3}
\sup_{t\leq m}|\varphi_p(t)-{t\over 2}|
\leq{1\over 2p\gamma_p}\sup_{k\leq mp\gamma_p}H^p_k
+{1\over p\gamma_p}
\build{\longrightarrow}_{p\rightarrow \infty}^{} 0 
\end{equation}
in probability, by (\ref{functionconv}).

The statement of the theorem now follows from (\ref{functionconv}), (\ref{cont-tech1})
and (\ref{cont-tech3}). \cq

\section{A joint convergence\\ and an application to conditioned trees}

The convergences in distribution (\ref{functionconv})  and (\ref{contour-convergence})
hold jointly with (\ref{Grim1}) and (\ref{Grim2}). This fact is useful in applications
and we state it here as a corollary.

As previously, we consider for every $p$ a sequence of independent $\mu_p$-Galton-Watson trees
and we denote by $(H^p_n,n\geq 0)$ the associated height process and by
$(C^p_t,t\geq 0)$ the associated contour process. The random walk $V^p$
with jump distribution $\nu_p(k)=\mu_p(k+1)$ is related to $H^p$ via
formula (\ref{heightwalk}). Finally, for every integer $k\geq 0$, we denote by
$Y^p_k$ the number of individuals at generation $k$ in the first $p$
trees of the sequence, so that, in agreement with the previous notation,
$(Y^p_n,n\geq 0)$ is a Galton-Watson process with offspring distribution $\mu_p$
started at $Y^p_0=p$.

Recall that $(L^a_t,a\geq 0,t\geq 0)$ denote the local times of the
(continuous-time) height process associated with the L\'evy process $X$.
From Theorem \ref{RK}, we know that $(L^a_{T_1},a\geq 0)$ is a $\psi$-CSBP
and thus has a c\` adl\` ag modification.

\begin{corollary}
\label{jointconv}
Suppose that the assumptions of Theorem \ref{functionalconv} are satisfied. Then,
$$\left(p^{-1}V^p_{[p\gamma_pt]},\gamma_p^{-1}H^p_{[p\gamma_pt]},
\gamma_p^{-1}C^p_{2p\gamma_pt}\,;\,t\geq 0\right)
\build{\longrightarrow}_{p\rightarrow \infty}^{{\rm (d)}} (X_t,H_t,H_t;t\geq 0)$$
in distribution in $\D(\R_+,\R^3)$. We have also
$$\left(p^{-1}Y^p_{[\gamma_pa]}\,,\,a\geq 0\right)
\build{\longrightarrow}_{p\rightarrow \infty}^{{\rm (d)}} (L^a_{T_1},a\geq 0)$$
in distribution in $\D(\R_+,\R)$. Furthermore, these two convergences hold jointly,
in the sense that, for any bounded continuous function $F$ on $\D(\R_+,\R^3)\times \D(\R_+,\R)$,
\ba
&&\lim_{p\rightarrow  \infty} E\Big[F\Big((p^{-1}V^p_{[p\gamma_pt]},\gamma_p^{-1}H^p_{[p\gamma_pt]},
\gamma_p^{-1}C^p_{2p\gamma_pt})_{t\geq 0},(p^{-1}Y^p_{[\gamma_pa]})_{a\geq 0}\Big)\Big]\\
&&\quad =E[F((X_t,H_t,H_t)_{t\geq 0},(L^a_{T_1})_{a\geq 0})].
\ea
\end{corollary}

\proof To simplify notation, write $V^{(p)}_t=p^{-1}V^p_{[p\gamma_pt]}$, $H^{(p)}_t=
\gamma_p^{-1}H^p_{[p\gamma_pt]}$, $C^{(p)}_t=\gamma_p^{-1}C^p_{2p\gamma_pt}$
and $Y^{(p)}_a=p^{-1}Y^p_{[\gamma_pa]}$. By (\ref{Grim2}), 
(\ref{functionconv})  and (\ref{contour-convergence}), we know that 
each of the three sequences of the laws of the processes $V^{(p)}$, $H^{(p)}$,
$C^{(p)}$ is tight, and furthermore $H^{(p)}$ and $C^{(p)}$ converge in
distribution towards a continuous process. By a standard result
(see e.g. Corollary II.3.33 in \cite{Ja}), we get that the laws of the
triples $(V^{(p)},H^{(p)},C^{(p)})$ are tight in $\D(\R_+,\R^3)$. Let
$(X,H^*,H^{**})$ be a weak limit point of this sequence of triples
(with a slight abuse of notation, we may assume that the first component
of the limiting triple is the underlying L\'evy process $X$). By the
Skorokhod representation theorem, we may assume that along a subsequence,
$$(V^{(p)},H^{(p)},C^{(p)})\la (X,H^*,H^{**})$$
a.s. in $\D(\R_+,\R^3)$. However, the convergence (\ref{DonsLevy}) and 
a time-reversal argument imply that 
$$\lim_{p\rightarrow \infty}H^{(p)}_t=\wh L^{(t)}_t=H_t$$
in probability. This is enough to conclude that $H^*_t=H_t$. Similarly, the proof
of Theorem \ref{contour-conv} shows that
$$\lim_{p\rightarrow \infty}(C^{(p)}_t-H^{(p)}_t)=0$$
in probability. This yields $H^{**}_t=H^*_t=H_t$ and  we see that the limiting 
triple is equal to $(X,H,H)$ and does not depend on the choice of the 
subsequence. The first convergence of the corollary now follows.

By (\ref{Grim1}),
we know that 
$$(Y^{(p)}_a,a\geq 0)\build{\longrightarrow}_{p\rightarrow \infty}^{{\rm (d)}}
(Y_a,a\geq 0)$$
where $Y$ is a $\psi$-CSBP started at $1$. Since we also know that
$(L^a_{T_1},a\geq 0)$ is a $\psi$-CSBP started at $1$, the second convergence
in distribution is immediate, and the point is to verify that this convergence
holds jointly with the first one. To this end, note that the laws of the
pairs $((V^{(p)},H^{(p)},C^{(p)}),Y^{(p)})$ are tight in the space of
probability measures on $\D(\R_+,\R^3)\times \D(\R_+,\R)$. By extracting
a subsequence and using the  
Skorokhod representation theorem, we may assume that
$$\Big((V^{(p)},H^{(p)},C^{(p)}),Y^{(p)}\Big)\build{\longrightarrow}_{p\rightarrow \infty}^{}
\Big((X,H,H),Z\Big),$$
a.s. in $\D(\R_+,\R^3)\times \D(\R_+,\R)$. The proof will be finished if we can verify that
$Z_a=L^a_{T_1}$, the local time of $H$ at level $a$ and time $T_1$. To this end,
let $g$ be a Lipschitz continuous function from $\R_+$ into $\R_+$ with
compact support. The preceding convergence implies
\begin{equation} 
\label{jointtech1}
\lim_{p\rightarrow \infty} \int_0^\infty g(a)Y^{(p)}_a\,da=\int_0^\infty g(a)Z_a\,da\ ,\qquad\hbox{a.s.}
\end{equation}
On the other hand, let $T^p_p$ be the hitting time of $-p$ by $V^p$. The convergence
of $V^{(p)}$ towards $X$ easily implies
\begin{equation} 
\label{jointtech2}
\lim_{p\rightarrow \infty}{1\over p\gamma_p}T^p_p=\inf\{t\geq 0:X_t=-1\}=T_1\ ,\quad\hbox{a.s.}
\end{equation}
Then, from the definition of the height process of a sequence of trees, we have
\ba
\int_0^\infty g(a)Y^{(p)}_a\,da&=&\int_0^\infty g(a)\,{1\over p}Y^p_{[\gamma_pa]}\,da\\
&=&{1\over p} \sum_{k=0}^\infty\int_{\gamma_p^{-1}k}^{\gamma_p^{-1}(k+1)}
g(a)\Big(\sum_{j=0}^{T^p_p-1}1_{\{H^p_j=k\}}\Big)da\\
&=&{1\over p}\sum_{j=0}^{T^p_p-1}\int_{\gamma_p^{-1}H^p_j}^{\gamma_p^{-1}(H^p_j+1)}
g(a)\,da\\
&=&{1\over p\gamma_p}\sum_{j=0}^{T^p_p-1}g(\gamma_p^{-1}H^p_j)+O({1\over p\gamma_p^2}T^p_p)\\
&=&\int_0^{(p\gamma_p)^{-1}T^p_p} g(\gamma_p^{-1}H^p_{[p\gamma_ps]})\,ds
+O({1\over p\gamma_p^2}T^p_p)
\ea
and in view of (\ref{jointtech2}) this converges to
$$\int_0^{T_1} g(H_s)\,ds=\int_0^\infty g(a)\,L^a_{T_1}\,da.$$
Comparing with (\ref{jointtech1}), we conclude that
$$\int_0^\infty g(a)Z_a\,da=\int_0^\infty g(a)\,L^a_{T_1}\,da.$$
This implies that $Z_a=L^a_{T_1}$ and completes the proof. \cq

As an application, we now discuss conditioned trees. Fix $T>0$ and on some probability space, consider 
a $\mu_p$-Galton-Watson tree conditioned on non-extinction at generation $[\gamma_pT]$,
which is denoted by $\wt T^p$. Let $\wt H^p=(\wt H^p_n,n\geq 0)$ be the associated
height process, with the convention that $\wt H^p_n=0$ for $n\geq \card(\wt T^p)$.

\begin{proposition}
\label{tree-condit-height}
Under the assumptions of Theorem \ref{functionalconv}, we have
$$\Big(\gamma_p^{-1}\wt H^p_{[p\gamma_pt]},t\geq 0\Big)
\build{\longrightarrow}_{p\rightarrow \infty}^{{\rm (d)}}
(\wt H_t,t\geq 0),$$
where the limiting process $\wt H$ is distributed as $H$ under $N(\cdot\mid\sup H_s\geq T)$.
\end{proposition}

\rem We could have stated a similar result for the contour process instead
of the discrete height process.

\smallskip
\proof Write $\wt H^{(p)}_s=\gamma_p^{-1}\wt H^p_{[p\gamma_ps]}$ to simplify
notation. Also let $H^{(p)}_s=\gamma_p^{-1} H^p_{[p\gamma_ps]}$ be as above the
rescaled height process for a sequence of independent $\mu_p$-Galton-Watson trees. Set
\ba
&&R^{(p)}_T=\inf\{s\geq 0 : H^{(p)}_s
={[\gamma_pT]\over \gamma_p}\},\\
&&G^{(p)}_T=\sup\{s\leq
R^{(p)}_T:H^{(p)}_s=0\},\\
&&D^{(p)}_T=\inf\{s\geq
R^{(p)}_T:H^{(p)}_s=0\}.
\ea
Then without loss of generality we may assume that
$$\wt H^{(p)}_s=H^{(p)}_{(G^{(p)}_T+s)\wedge D^{(p)}_T}\ , \qquad s\geq 0.$$
This is simply saying that the first tree with height at least $[\gamma_pT]$
in a sequence of independent $\mu_p$-Galton-Watson trees is a $\mu_p$-Galton-Watson tree
conditioned on non-extinction at generation $[\gamma_pT]$.

Set
\ba
&&R_T=\inf\{s\geq 0 : H_s=T\},\\
&&G_T=\sup\{s\leq R_T:H_s=0\},\\
&&D_T=\inf\{s\geq R_T:H_s=0\},
\ea
and note that we may take $\wt H_s=H_{(G_T+s)\wedge D_T}$, by excursion theory for $X-I$.

We now claim that the convergence in distribution of $\wt H^{(p)}$ towards $\wt H$ follows from
the previous corollary, and more precisely from the joint convergence
$$(V^{(p)},H^{(p)})\build{\longrightarrow}_{p\rightarrow \infty}^{{\rm (d)}}(X,H).$$
It is again convenient to use the Skorokhod representation theorem
and to assume that the latter convergence holds a.s. We can then prove that
$\wt H^{(p)}$ converges a.s. towards $\wt H$.

To this end we need a technical lemma about the height process. We state it in greater generality than
needed here in view of other applications.

\begin{lemma}
\label{extrema}
Let $b>0$. Then $P$ a.s. or $N$ a.e. $b$ is not a local maximum nor
a local minimum of the height process.
\end{lemma}

\proof
Let 
$$D=\{b>0: P[\sup_{p\leq s\leq q}H_s=b]>0\ {\rm for\ some\ rationals}
\ q>p\geq 0\}.$$
Clearly $D$ is at most countable. However, from Proposition \ref{reflec}
and the relation between the height process and the 
exploration process, it immediately follows that if $b\in D$
then $b-a\in D$ for every $a\in[0,b)$. This is only possible
if $D=\emptyset$. The case of local minima is treated
in the same way. \cq

It follows from the lemma that we have also $R_T=\inf\{s\geq 0:H_s>T\}$.
Then the a.s. convergence of $H^{(p)}$ towards $H$ easily implies that
$R^{(p)}_T$ converges to $R_T$ a.s., and that
$$\limsup_{p\rightarrow  \infty}G^{(p)}_T\leq G_T\ ,\quad 
\liminf_{p\rightarrow  \infty}D^{(p)}_T\geq D_T\ .$$
To get reverse inequalities, we may argue 
as follows. Recall that the support of the random
measure $dI_s$ is exactly the set $\{s:H_s=0\}$,
so that for every fixed $s\geq 0$, we have $I_{s}>I_{R_T}$
a.s. on the set $\{s<G_T\}$. If $I^{(p)}_s=\inf\{V^{(p)}_r\,,\,r\leq s\}$,
it readily follows that a.s. on the set $\{s<G_T\}$ we have
$I^{(p)}_s>I^{(p)}_{R^{(p)}_T}$ for all $p$
sufficiently large. Hence a.s. for $p$ large, we have
$s<G^{(p)}_T$ on the set $\{s<G_T\}$. We conclude that
$G^{(p)}_T\rightarrow  G_T$ a.s., and a similar argument
gives $D^{(p)}_T\rightarrow  D_T$. From the preceding formulas for $\wt H^{(p)}$ and $\wt H$, 
it follows that $\wt H^{(p)}\rightarrow \wt H$ a.s. This completes the
proof of the proposition.

\cq

\rem The methodology of proof of Proposition \ref{tree-condit-height}
could be applied to other conditioned limit theorems. For instance, we
could consider the rescaled height (or contour) process of the 
$\mu_p$-Galton-Watson tree conditioned to have at least $p\gamma_p$
vertices and derive a convergence towards the excursion of the height
process $H$ conditioned to have length greater than $1$. We will leave such
extensions to the reader. We point out here that it is much harder to 
handle degenerate conditionings. To give an important example, consider the case 
where $\mu_p=\mu$ for every $p$. It is natural to ask for a limit
theorem for the (rescaled) height or contour process of a $\mu$-Galton-Watson
tree conditioned to have a large fixed number of vertices. The previous
results strongly suggest that the limiting process should be
a normalized (i.e. conditioned to have a fixed length) excursion of
the height process $H$. This is indeed true under
suitable assumptions: When $\mu$ is critical with finite variance, this was
proved by Aldous \cite{Al2} in the case of the contour process and the limit
is a normalized Brownian excursion as expected.
Aldous' result has been extended by Duquesne 
\cite{Du2} to the case when $\mu$ is in the
domain of attraction of a stable law of index $\gamma\in(1,2]$.

\section{The convergence of reduced trees}

Consider a $\mu$-Galton-Watson tree, which describes the genealogy of
a Galton-Watson process with offspring distribution $\mu$ starting
with one individual at time $0$. For every integer $n\geq 1$,
denote by $P^{(n)}$ the conditional probability knowing that
the process is not extinct at time $n$, or equivalently the height of
the tree is at least $n$. Under $P^{(n)}$, we can consider the
reduced tree that consists only of those individuals
in the generations up to time $n$ that have descendants at 
generation $n$. The results of the previous sections can be used to
investigate the limiting behavior of these reduced trees when 
$n$ tends to $\infty$, even in the more general setting where 
the offspring distribution depends on $n$. 

Here, we will concentrate on the population 
of the reduced tree at every generation. For 
every $k\in\{0,1,\ldots,n\}$, we denote by $Z^n_k$
the number of individuals in the tree at generation $k$
which have descendants at generation $n$. Obviously, $k\rightarrow  Z^n_k$
is nondecreasing, $Z^n_0=1$ and $Z^n_n$ is equal to the 
number of individuals in the original tree at generation $n$. 
If $g$ denotes the generating function of $\mu$ and 
$g_n$, $n\geq 0$ are the iterates of $g$, it is easy to
verify that $(Z^n_k,0\leq k\leq n)$ is a time-inhomogeneous
Markov chain whose transition kernels are characterized by:
$$E^{(n)}[r^{Z^n_{k+1}}\mid Z^n_k]=\Big(
{g(r(1-g_{n-k-1}(0))+g_{n-k-1}(0))-g_{n-k}(0)\over
1-g_{n-k}(0)}\Big)^{Z^n_k}
\ ,\quad 0\leq k<n\,.$$
The process $(Z^n_k,0\leq k\leq n)$ (under the 
probability measure $P^{(n)}$)
will be called the 
reduced process of the $\mu$-Galton-Watson tree 
at generation $n$. It is easy to see that for every
$k\in\{0,1,\ldots,n-1\}$, $Z^n_k$ can be written as
a simple functional of the height process of the tree: $Z^n_k$
counts the number of excursions of the height process 
above level $k$ that hit level $n$.

Consider as in the previous sections a sequence
$(\mu_p,p=1,2,\ldots)$ of (sub)critical offspring 
distributions, and for every integer $n\geq 1$
let $Z^{(p),n}=(Z^{(p),n}_k,0\leq k\leq n)$ be
the reduced process of the $\mu_p$-Galton-Watson tree 
at generation $n$. For every
$T>0$, we denote by $N_{(T)}$ 
the conditional probability $N(\cdot\mid \sup\{H_s,s\geq 0\}\geq T)$
(this makes sense provided that the condition (\ref{continuitycond})
holds, cf Corollary \ref{extinction}).

\begin{theorem}
\label{reducedconv}
Suppose that the assumptions of Theorem \ref{functionalconv}
hold and let $T>0$. Then,
$$\Big(Z^{(p),[\gamma_pT]}_{[\gamma_pt]},0\leq t<T\Big)
\build{\la}_{p\rightarrow \infty}^{(fd)} (Z^T_t,0\leq t< T)\ ,$$
where the limiting process $(Z^T_t,0\leq t<T)$ is defined
under $N_{(T)}$
as follows: For every $t\in[0,T)$, $Z^T_t$ is the number of excursions
of $H$ above level $t$ that hit level $T$. 
\end{theorem}

A more explicit description of the limiting process and
of the associated tree will be given in the next section.

\smallskip
\proof We use the notation of the proof of Proposition \ref{tree-condit-height}.
In particular, the height process of the $\mu_p$-Galton-Watson tree conditioned
on non-extinction at generation $[\gamma_pT]$ is $(\wt H^p_k,k\geq 0)$ and the associated
rescaled process is $\wt H^{(p)}_s=\gamma_p^{-1}\wt H^p_{[p\gamma_ps]}$.
We may and will assume that $\wt H^{(p)}_s$ is given by the formula
$$\wt H^{(p)}_s=H^{(p)}_{(G^{(p)}_T+s)\wedge D^{(p)}_T}$$
and that $(\wt H^{(p)}_s,s\geq 0)$ converges a.s. in the sense of the
Skorokhod topology, towards the process $\wt H_s=H_{(G_T+s)\wedge D_T}$
whose law is the distribution of $H$ under $N_{(T)}$.

Now we observe that the reduced process $Z^{(p),[\gamma_pT]}_{[\gamma_pt]}$
can be expressed in terms of $\wt H^{(p)}$. More precisely,
it is clear by construction that
for every $k\in\{0,1,\ldots,[\gamma_pT]-1\}$, 
$Z^{(p),[\gamma_pT]}_k$ is the number of excursions 
of $\wt H^p$ above level $k$ that hit level $[\gamma_pT]$.
Equivalently, for every $t$ such that $[\gamma_pt]<[\gamma_pT]$,
$$\wt Z^{(p)}_t:=Z^{(p),[\gamma_pT]}_{[\gamma_pt]}$$
is the number of excursions of $\wt H^{(p)}$ above level
$[\gamma_pt]/\gamma_p$ that hit level $[\gamma_pT]/\gamma_p$.

Let $t>0$. Using the fact that $t$, resp. $T$, is a.s. not a
local minimum, resp. maximum, of $H$ (Lemma \ref{extrema}), 
it is easy to deduce from the
convergence
$\wt H^{(p)}\rightarrow \wt H$ that the number of excursions of
$\wt H^{(p)}$ above level $[\gamma_pt]/\gamma_p$ that hit level
$[\gamma_pT]/\gamma_p$ converges a.s. to the number of
excursions of $\wt H$ above level $t$ that hit level $T$.
In other words, $\wt Z^{(p)}_t$ converges a.s. to $Z^T_t$.
This completes the proof. \cq

\section{The law of the limiting reduced tree}

In this section, we will describe the law of the
process $(Z^T_t,0\leq t<T)$ of the previous section,
and more precisely the law of the underlying branching tree.
We suppose that the L\'evy process $X$ satisfies
(\ref{continuitycond}) in addition to (H1) -- (H3).
The random variable $Z^T_t$ (considered under the probability
measure $N_{(T)}$) counts the number of excursions of
$H$ above level $t$ that hit level $T$.

Before stating our result, we recall the notation of Section 1.4.
For every $\lambda>0$ and $t>0$, 
$$u_t(\lambda)=N(1-\exp(-\lambda L^t_\sigma))$$
solves the integral equation
$$u_t(\lambda)+\int_0^t \psi(u_s(\lambda))\,ds=\lambda$$
and 
$$v(t)=u_t(\infty)=N(L^t_\sigma>0)=N\Big(\sup_{s\geq 0}H_s>t\Big)$$
is determined by
$$\int_{v(t)}^\infty {dx\over \psi(x)}=t.$$
Note the composition property $u_t\circ u_s =u_{t+s}$,
and in particular $u_t(v(r))=v(t+r)$.

\begin{theorem}
\label{lawreduced}
Under $N_{(T)}$, the process $(Z^T_t,0\leq t<T)$ is a time-inhomogeneous
Markov process whose law is characterized by the
following identities: For every $\lambda>0$,
\begin{equation}
\label{marginalreduced}
N_{(T)}[\exp-\lambda Z^T_t]=1-{u_t((1-e^{-\lambda})v(T-t))\over v(T)}.
\end{equation}and if $0\leq t<t'<T$,
\begin{equation}\label{reducedMarkov}
N_{(T)}[\exp-\lambda Z^T_{t'}\mid Z^T_t]
=(N_{(T-t)}[\exp-\lambda Z^{T-t}_{t'-t}])^{Z^T_t}
\end{equation}Alternatively, we can describe the law of
the process $(Z^T_t,0\leq t<T)$ under $N_{(T)}$
by the following properties.
\begin{itemize}
\item $Z^T_r=1$ if and only if
$r\in[0,{\gamma_T})$, where the law of ${\gamma_T}$ is given by
\begin{equation}\label{lifetimereduced}
N_{(T)}[{\gamma_T}>t]={\wt\psi(v(T))\over \wt\psi(v(T-t))}
\ ,\qquad 0\leq t<T,
\end{equation}where $\wt\psi(x)=x^{-1}\psi(x)$.
\item The conditional distribution of $Z^T_{\gamma_T}$ 
knowing ${\gamma_T}$ is characterized by
\begin{equation}\label{offspringreduced}
N_{(T)}[r^{Z^T_{\gamma_T}}\mid {\gamma_T}=t]
=r\,{\psi'(U)-\gamma_\psi(U,(1-r)U)\over
\psi'(U)-\gamma_\psi(U,0)},\qquad 0\leq r\leq 1
\end{equation}where $U=v(T-t)$ and for every $a,b\geq 0$,
$$\gamma_{\psi }(a,b)= \left\{
\begin{array}{ll}
\left(\psi (a) -\psi ( b) \right) / (a-b)
\quad& {\it if } \quad a\neq b, \\
\psi'(a) \quad& {\it if } \quad a=b\; .\hfill\\ 
\end{array}
\right. $$
\item
Conditionally on ${\gamma_T}=t$ and $Z^T_{\gamma_T}=k$, the process
$(Z^T_{t+r},0\leq r<T-t)$ is distributed as the sum 
of $k$ independent copies of the process $(Z^{T-t}_{r},0\leq r<T-t)$
under $N_{(T-t)}$.
\end{itemize}
\end{theorem}

\proof One can give several approaches to Theorem
\ref{lawreduced}. In particular, the time-inhomoge\-neous Markov
property could be deduced from the analogous result
for discrete reduced trees by using Theorem \ref{reducedconv}.
We will prefer to give a direct approach relying on the
properties of the height process.

Before stating a key lemma, we introduce some notation.
We fix $t\in(0,T)$. Note that the definition
of $Z^T_t$ also makes sense under the conditional
probability $N_{(t)}$.
We denote by $(e^t_i,i=1,\ldots,Z^T_t)$
the successive excursions of $H$ above level $t$ that
hit level $T-t$, shifted in space and time so that each
starts from $0$ at time $0$. Recall the notation 
$L^a_s$ for the local times of the height process.
We also write $L^t_{(i)}$ for the local time of $H$
at level $t$ at the beginning of excursion $e^t_i$.

\begin{lemma}
\label{reduced-tech} Under $N_{(t)}$, conditionally
on the local time $L^t_\sigma$, the point measure
$$\sum_{i=1}^{Z^T_t} \delta_{(L^t_{(i)},e^t_i)}$$
is Poisson with intensity $1_{[0,L^t_\sigma]}(\ell)d\ell\,N(de\cap \{\sup H_s>T-t\})$.
In particular, under $N_{(t)}$ or under $N_{(T)}$,
conditionally on $Z^T_t$, the excursions $(e^t_i,i=1,\ldots,Z^T_t)$
are independent with distribution $N_{(T-t)}$.
\end{lemma}

\proof We rely on
Proposition \ref{reflec} and use the notation of Chapter 1. Under the probability
measure $P$, denote by $f^t_i$, $i=1,2,\ldots$ the successive excursions
of $H$ above level $t$ that hit $T$, and let $\ell^t_i$
be the local time of $H$ at level $t$ at the beginning (or the end)
of excursion $f^t_i$. 
Then the $f^t_i$'s are also the successive 
excursions of the process $H^t_s=H(\rho^t_s)$ that hit level $T-t$,
and the numbers $\ell^t_i$ are the corresponding local times (of $H^t$)
at level $0$. By Proposition \ref{reflec} and excursion theory, the point
measure
$$\sum_{i=1}^\infty \delta_{(\ell^t_i,f^t_i)}$$
is Poisson with intensity $d\ell\,N(df\cap \{\sup H_s>T-t\})$
and is independent of the $\sigma$-field ${\mathcal H}_t$. 

On the other hand, let $\lambda_1$ be the local time 
of $H$ at level $t$ at the end of the first excursion of $H$
away from $0$ that hits level $t$. From 
the approximation of local time provided by
Proposition \ref{LTapprox}, it is easy to see that $\lambda_1$ is ${\mathcal H}_t$-measurable.
By the same argument
as in the proof of Theorem \ref{reducedconv}, the law under $N_{(t)}$
of the pair
$$\Big(L^t_\sigma,\sum_{i=1}^{Z^T_t} \delta_{(L^t_{(i)},e^t_i)}\Big)$$
is the same as the law under $P$ of
$$\Big(\lambda_1,\sum_{\{i:\ell^t_i\leq \lambda_1\}}
\delta_{(\ell^t_i,f^t_i)}\Big).$$
The first assertion of the lemma now follows from the
preceding considerations.

The second assertion stated under $N_{(t)}$
is an immediate consequence of the first one. The statement
under $N_{(T)}$ follows since $N_{(T)}=N_{(t)}(\cdot\mid Z^T_t\geq 1)$. \cq

\medskip

We return to the proof of Theorem \ref{lawreduced}. Note that
(\ref{reducedMarkov}) is an immediate consequence of the second 
assertion of the lemma. Let us prove (\ref{marginalreduced}).
By the first assertion of the lemma, $Z^T_t$ is Poisson with
intensity $v(T-t)L^t_\sigma$, conditionally on $L^t_\sigma$, under $N_{(t)}$. Hence,
\begin{eqnarray*}
N_{(t)}[e^{-\lambda Z^T_t}]&=&N_{(t)}\Big[e^{-L^t_\sigma v(T-t)(1-e^{-\lambda})}\Big]\cr
&=&1-{1\over v(t)}N\Big(1-e^{-L^t_\sigma v(T-t)(1-e^{-\lambda})}\Big)\cr
&=&1-{1\over v(t)}u_t((1-e^{-\lambda})v(T-t)).
\end{eqnarray*}
Then observe that
$$N_{(t)}[1-e^{-\lambda Z^T_t}]={1\over v(t)}N(1-e^{-\lambda Z^T_t})={v(T)\over
v(t)}N_{(T)}[1-e^{-\lambda Z^T_t}].$$
Formula (\ref{marginalreduced}) follows immediately.

It is clear that there exists a random variable ${\gamma_T}$ such
that $Z^T_t=1$ iff $0\leq t<{\gamma_T}$, $N_{(T)}$ a.s. (${\gamma_T}$
is the minimum of the height process between the first and the
last hitting time of $T$). Let us prove (\ref{lifetimereduced}).
By (\ref{marginalreduced}), we have,
$$N_{(T)}[{\gamma_T}>t]=\lim_{\lambda\rightarrow \infty} e^\lambda N_{(T)}[e^{-\lambda Z^T_t}]
=\lim_{\lambda\rightarrow \infty} e^\lambda\Big(1-{u_t((1-e^{-\lambda})v(T-t))\over v(T)}\Big).$$
Recalling that $u_t(v(T-t))=v(T)$, we have as $\varepsilon\rightarrow  0$,
$$u_t((1-\varepsilon)v(T-t))=v(T)-\varepsilon v(T-t){\partial u_t\over
\partial\lambda}(v(T-t))+o(\varepsilon),$$
and it follows that
$$N_{(T)}[{\gamma_T}>t]={v(T-t)\over v(T)}\,{\partial u_t\over
\partial\lambda}(v(T-t)).$$
Formula (\ref{lifetimereduced}) follows from that identity and the fact that, for $\lambda>0$,
\begin{equation}\label{deriv-u}
{\partial u_t\over
\partial\lambda}(\lambda)={\psi(u_t(\lambda))\over \psi(\lambda)}.
\end{equation}To verify (\ref{deriv-u}), differentiate the integral equation for $u_t(\lambda)$:
$${\partial u_t\over
\partial\lambda}(\lambda)=1-\int_0^t {\partial u_s\over
\partial\lambda}(\lambda)\,\psi'(u_s(\lambda))\,ds$$
which implies
$${\partial u_t\over
\partial\lambda}(\lambda)=\exp\Big(-\int_0^t \psi'(u_s(\lambda))ds\Big).$$
Then note that 
${\partial\over \partial t}\log\psi(u_t(\lambda))=-\psi'(u_t(\lambda))$
and thus
$$\int_0^t \psi'(u_s(\lambda))ds=\log\psi(u_t(\lambda))-\log\psi(\lambda).$$
This completes the proof of (\ref{deriv-u}) and (\ref{lifetimereduced}).

\medskip
We now prove the last assertion of the theorem. 
Recall the notation introduced before Lemma \ref{reduced-tech}. Clearly it suffices to
prove that the following property holds:

\smallskip
(P) Under $N_{(T)}$, conditionally on ${\gamma_T}=t$ and $Z^T_{\gamma_T}=n$, the excursions
$e^{\gamma_T}_1,\ldots,e^{\gamma_T}_n$ are i.i.d. according to the distribution $N_{(T-t)}$.

\smallskip
We can deduce property (P) from Lemma \ref{reduced-tech} via an approximation
procedure. Let us sketch the argument. For any $p\geq 2$
and any bounded continuous functional $F$ on $\R_+\times C(\R_+,\R_+)^p$,
\begin{eqnarray}
\label{approx-decomp}
&&N_{(T)}[1_{\{Z^T_{\gamma_T}=p\}}F({\gamma_T},e^{\gamma_T}_1,\ldots,e^{\gamma_T}_p)]\nonumber\\
&&\ =\lim_{n\rightarrow \infty}
\sum_{j=1}^{n-1}N_{(T)}\Big[1_{\{Z^T_{jT/n}=p;(j-1)T/n<{\gamma_T}\leq
jT/n\}}F({jT\over n},e^{jT/n}_1,\ldots,e^{jT/n}_p)\Big].
\end{eqnarray}
Note that the event $\{{\gamma_T}\leq jT/n\}$ contains $\{Z^T_{jT/n}=p\}$. As a consequence 
of the second part of Lemma \ref{reduced-tech} (applied with $t=jT/n$) we have
\begin{eqnarray*}
&&N_{(T)}\Big[1_{\{Z^T_{jT/n}=p;{\gamma_T}\leq
jT/n\}}F({jT\over n},e^{jT/n}_1,\ldots,e^{jT/n}_p)\Big]\\
&&=N_{(T)}\Big[1_{\{Z^T_{jT/n}=p;{\gamma_T}\leq
jT/n\}}\\
&&\qquad \qquad\int N_{(T-jT/n)}(df_1)\ldots N_{(T-jT/n)}(df_p)\,F({jT\over n},f_1,\ldots,f_p)\Big].
\end{eqnarray*}

We want to get a similar identity where the event $\{{\gamma_T}\leq
jT/n\}$ is replaced by $\{{\gamma_T}\leq
(j-1)T/n\}=\{Z^T_{(j-1)T/n}\geq 2\}$.
A slightly more complicated 
argument (relying on two applications of Lemma \ref{reduced-tech}, the first one with
$t=(j-1)T/n$ and then with $t=T/n$)
shows similarly that
\begin{eqnarray*}
&&N_{(T)}\Big[1_{\{Z^T_{jT/n}=p;{\gamma_T}\leq
(j-1)T/n\}}F({jT\over n},e^{jT/n}_1,\ldots,e^{jT/n}_p)\Big]\\
&&\ =N_{(T)}\Big[1_{\{Z^T_{jT/n}=p;{\gamma_T}\leq
(j-1)T/n\}}\\
&&\qquad\times\int N_{(T-jT/n)}(df_1)\ldots N_{(T-jT/n)}(df_p)\,F({jT\over n},f_1,\ldots,f_p)\Big].
\end{eqnarray*}

By making the difference between the last two displays, we see that the sum
in the right side of (\ref{approx-decomp}) exactly equals
\begin{eqnarray*}
&&\sum_{j=1}^{n-1}N_{(T)}\Big[1_{\{Z^T_{jT/n}=p;(j-1)T/n<{\gamma_T}\leq
jT/n\}}\\
&&\qquad\times\int N_{(T-jT/n)}(df_1)\ldots N_{(T-jT/n)}(df_p)\,F({jT\over n},f_1,\ldots,f_p)\Big].
\end{eqnarray*}
Using an easy continuity property of the mapping $r\rightarrow  N_{(r)}$, we get from
this and (\ref{approx-decomp}) that
\begin{eqnarray*}
&&N_{(T)}[1_{\{Z^T_{\gamma_T}=p\}}F({\gamma_T},e^{\gamma_T}_1,\ldots,e^{\gamma_T}_p)]\\
&&\ =N_{(T)}\Big[1_{\{Z^T_{\gamma_T}=p\}}\int N_{(T-{\gamma_T})}(df_1)\ldots
N_{(T-{\gamma_T})}(df_p) F({\gamma_T},f_1,\ldots,f_p)\Big],
\end{eqnarray*}
which completes the proof of property (P) and of the last assertion
of the theorem.

\smallskip
We finally verify (\ref{offspringreduced}). First observe from (\ref{lifetimereduced})
that the density of the law of $\gamma_T$ under $N_{(T)}$ is given by
$$h_T(t)=\wt\psi(v(T))\,h(T-t)$$
where
$$h(t)={v(t)\psi'(v(t))\over \psi(v(t))}-1.$$
On the other hand, fix $\delta\in(0,T)$, and note that
$\{\gamma_T>\delta\}=\{Z^T_\delta=1\}$. By the last assertion
of Lemma \ref{reduced-tech} we have for any nonnegative function $f$,
$$N_{(T)}[f(\gamma_T,Z^T_{\gamma_T})\,1_{\{\gamma_T>\delta\}}\mid \gamma_T>\delta]
=N_{(T-\delta)}[f(\gamma_{T-\delta}+\delta,Z^{T-\delta}_{\gamma_{T-\delta}})].$$
Hence, if $(\theta^T_t(k),k=2,3,\ldots)$, $0<t<T$ denotes a regular version
of the conditional law of $Z^T_{\gamma_T}$ knowing that $\gamma_T=t$, we have
\begin{eqnarray*}
\int_\delta^T dt\,h(T-t)\sum_{k=2}^\infty \theta^T_t(k)f(t,k)
&=&\int_0^{T-\delta}dt\,h(T-\delta-t)\sum_{k=2}^\infty \theta^{T-\delta}_t(k)f(t+\delta,k)\\
&=&\int_\delta^{T}dt\,h(T-t)\sum_{k=2}^\infty \theta^{T-\delta}_{t-\delta}(k)f(t,k).
\end{eqnarray*}
This shows that we must have $\theta^T_t=\theta^{T-\delta}_{t-\delta}$ for a.a. $t\in(\delta,T)$.
By simple arguments, we can choose the regular versions $\theta^T_t(k)$ in such a
way that $\theta^T_t(k)=\theta_{T-t}(k)$ for every $k\geq 2$, $T>0$ and $t\in(0,T)$.

We can then compute $N_{(T)}[e^{-\lambda L^T_\sigma}]$ in two different ways. First,
$$N_{(T)}[e^{-\lambda L^T_\sigma}]
=1-{N(1-e^{-\lambda L^T_\sigma})\over v(T)}=1-{u_T(\lambda)\over v(T)}.$$
Then, using property (P) once again,
\begin{eqnarray*}
N_{(T)}[e^{-\lambda L^T_\sigma}]&=&N_{(T)}\Big[\Big(N_{(T-t)}[e^{-\lambda
L^{T-t}_\sigma}]_{t=\gamma_T}\Big)^{Z^T_{\gamma_T}}\Big]\\
&=&\int_0^T dt\,h_T(t)\sum_{k=2}^\infty \theta_{T-t}(k)\,\Big(1-{u_{T-t}(\lambda)\over v(T-t)}\Big)^k.
\end{eqnarray*}
By comparing with the previous display and using the formula for 
$h_T(t)$, we get
$$\int_0^T dt\Big({\psi'(v(t))v(t)\over \psi(v(t))}-1\Big)
\sum_{k=2}^\infty \theta_{t}(k)\,\Big(1-{u_{t}(\lambda)\over v(t)}\Big)^k
={v(T)-u_T(\lambda)\over \psi(v(T))}.$$
We can now differentiate with respect to $T$ (for a proper justification we should
argue that the mapping $t\rightarrow  \theta_t$ is continuous, but we omit details). It follows
that
\begin{eqnarray*}
&&\Big({\psi'(v(T))v(T)\over \psi(v(T))}-1\Big)
\sum_{k=2}^\infty \theta_{T}(k)\,\Big(1-{u_{T}(\lambda)\over v(T)}\Big)^k\\
&&\qquad=-1+{\psi'(v(T))v(T)\over \psi(v(T))}+{\psi(u_T(\lambda))-u_T(\lambda)\psi'(v(T))\over
\psi(v(T))}.
\end{eqnarray*}
Hence,
$$\sum_{k=2}^\infty \theta_{T}(k)\,\Big(1-{u_{T}(\lambda)\over v(T)}\Big)^k
=1-{\psi(u_T(\lambda))-u_T(\lambda)\psi'(v(T))\over \psi(v(T))-v(T)\psi'(v(T))}.$$
If we substitute $r=1-{u_T(\lambda)\over v(T)}$
in this last identity we get
$$\sum_{k=2}^\infty \theta_{T}(k)\,r^k=1-{\psi((1-r)v(T))-(1-r)v(T)\psi'(v(T))\over
\psi(v(T))-v(T)\psi(v(T)}.$$
Formula (\ref{offspringreduced}) follows after straightforward transformations
of the last expression.

The proof of Theorem \ref{lawreduced} is now complete. Observe that the (time-inhomoge\-neous)
Markov property of the process $(Z^T_t,0\leq t<T)$ is a consequence of the description
provided in the second part of the theorem, and in particular of
the special form of the law of $\gamma_T$ and the fact that the law of
$Z^T_{\gamma_T}$ under $N_{(T)}[\cdot \mid \gamma_T>\delta]$
coincides with the law of $Z^{T-\delta}_{\gamma_{T-\delta}}$ under $N_{(T-\delta)}$. \cq

\bigskip
Let us discuss special cases of the theorem. When $\psi(u)=cu^\alpha$, with 
$c>0$ and $1<\alpha\leq 2$,
we have $v(t)=(c(\alpha-1)t)^{-1/(\alpha-1)}$, and formula (\ref{lifetimereduced})
shows that the law of $\gamma_T$ is uniform over $[0,T]$. This is the only case where
this property holds: If we assume that $\gamma_T$ is uniform over $[0,T]$,
(\ref{lifetimereduced}) implies that $\wt\psi(v(t))=C/t$ for some $C>0$.
By differentiating $\log\,v(t)$, we then get that $v(t)=C't^{-C}$
and it follows that $\psi$ is of the desired form.

\smallskip
Also in the stable case $\psi(u)=cu^\alpha$, formula (\ref{offspringreduced})
implies that $Z^T_{\gamma_T}$ is independent of $\gamma_T$,
and that its distribution is characterized by 
$$N_{(T)}[r^{Z^T_{\gamma_T}}]={(1-r)^\alpha-1+\alpha r\over \alpha-1}.$$
Of course when $\alpha=2$, we recover the well known fact that
$Z^T_{\gamma_T}=2$. When $\alpha\in(1,2)$, we get
$$N_{(T)}[Z^T_{\gamma_T}=k]={\alpha(2-\alpha)(3-\alpha)\cdots(k-1-\alpha)\over k!}\ , \quad k\geq
2.$$

To conclude let us mention that limiting reduced trees have been studied extensively
in the literature. In the finite variance case, the uniform distribution for $\gamma_T$ appears in Zubkov \cite{Zu},
and the full structure of the reduced tree is derived by
Fleischmann and Siegmund-Schultze \cite{FlSG}. Analogous results in the stable case
(and in the more general setting of multitype branching processes) can be found in
Vatutin \cite{Va} and Yakymiv \cite{Ya}.

\chapter{Marginals of continuous trees}

\section{Duality properties of the exploration process}

In this section, we study certain duality properties of
the process $\rho$. In view of forthcoming applications,
the main result is the time-reversal property stated 
in Corollary \ref{reversal} below. However the intermediate
results needed to derive this property are of independent
interest.

We work in the general setting of Chapter 1. In particular, 
the L\'evy process $X$ satisfies assumptions (H1) -- (H3),
and starts at $0$ under the probability measure $P$. 
Since the subordinator $S_{L^{-1}(t)}$ has drift $\beta$
(Lemma \ref{Laplaceladder}), it readily follows from formula (\ref{rhoreversed})
that the continuous part of $\rho_t$ is $\beta 1_{[0,H_t]}(r)dr$.
We can thus rewrite Definition \ref{explordef} in an equivalent way as
follows:
\begin{equation}
\label{new-rho}
\rho_t(dr)=\beta 1_{[0,H_t]}(r)\,dr+\build{\sum_{0<s\leq t}}
_{X_{s-}<I^s_t}^{} (I^s_t-X_{s-})\,\delta_{H_s}(dr).
\end{equation}
We then introduce another measure-valued process $(\eta_t,t\geq 0)$
by setting
\begin{equation}
\label{new-eta}
\eta_t(dr)=\beta 1_{[0,H_t]}(r)\,dr+\build{\sum_{0<s\leq t}}
_{X_{s-}<I^s_t}^{} (X_s-I^s_t)\,\delta_{H_s}(dr).
\end{equation}
In the same way as $\rho_t$, the measure $\eta_t$ is supported
on $[0,H_t]$. We will see below that $\eta_t$ is a.s. a {\it finite}
measure, a fact that is not obvious from the previous formula.
In the queuing system interpretation of \cite{LGLJ1}, the process
$(\rho_t,t\geq 0)$ accounts for the remaining service times
for all customers present in the queue at time $t$. In this interpretation,
$\eta_t$ describes the services already accomplished for these
customers.

We will see that in some sense, the process $(\eta_t,t\geq 0)$ 
is the dual of $(\rho_t,t\geq 0)$. It turns out that the study of
$(\eta_t,t\geq 0)$ is significantly more difficult than that
of $(\rho_t,t\geq 0)$. We start with a basic lemma.

\begin{lemma}
\label{basic-eta}
For each fixed value of $t>0$, we have $\langle\eta_t,1\rangle<\infty$,
$P$ a.s. or $N$ a.e. The process $(\eta_t,t\geq 0)$, which takes
values in $M_f(\R_+)$, is right-continuous in probability under $P$
Similarly, $(\eta_t,t>0)$ is right-continuous in measure 
under $N$.
\end{lemma}

\proof Let us first prove that $\langle\eta_t,1\rangle<\infty$,
$P$ a.s. It is enough to verify that
$$\sum_{0<s\leq t} \Delta X_s\,{1}_{\{X_{s-}<I^s_t\}}<\infty$$
$P$ a.s. By time-reversal, this is equivalent to
\begin{equation}
\label{eta-tech1}
\sum_{0<s\leq t} \Delta X_s\,{1}_{\{X_{s}>S_{s-}\}}<\infty
\end{equation}
$P$ a.s. However, for every $a>0$,
\begin{eqnarray*}
E\Big[\sum_{0<s\leq L^{-1}(a)} (\Delta X_s\wedge 1)\,{1}_{\{X_{s}>S_{s-}\}}\Big]
&=&aN^*((\Delta X_\sigma)\wedge 1)1_{\{X_\sigma>0\}})\\
&=&a\int \pi(dx)\int_0^x dz(z\wedge 1)\\
&\leq&a\int\pi(dx)\,(x\wedge x^2)\\
&<&\infty
\end{eqnarray*}
using (\ref{joint}) in the second equality. This gives our claim (\ref{eta-tech1})
and the first assertion of the lemma under $P$. 
The property $\langle\eta_t,1\rangle<\infty$, $N$ a.e., then follows 
from arguments of excursion theory, using in particular the Markov
property of $X$ under $N$.

The preceding considerations also imply that
$$\lim_{t\da 0} \sum_{0<s\leq t} \Delta X_s\,{1}_{\{X_{s}>S_{s-}\}}=0$$
in $P$-probability. Via time-reversal, it follows that the 
process $\eta_t$ is right-continuous at $t=0$ in probability
under $P$. Then let $t_0>0$. We first observe that
$\eta_{t_0}(\{H_{t_0}\})=0$ $P$ a.s. This follows from the fact 
that there is a.s. no value of $s\in(0,t_0]$ with $X_s>S_{s-}$
and $L_s=0$. Then, for $t>t_0$,
write $u=u(t)$ for the (first) time of the minimum of $X$
over $[t_0,t]$. Formula (\ref{new-eta}) implies that 
$\eta_{t}$ is bounded below by the restriction of $\eta_{t_0}$
to $[0,H_u)$, and bounded above by $\eta_{t_0}+\widetilde \eta^{(t)}_{t-t_0}$,
where  $\langle\widetilde \eta^{(t)}_{t-t_0},1\rangle$ has the
same distribution as $\langle\eta_{t-t_0},1\rangle$ (more precisely,
$\widetilde \eta^{(t)}_{t-t_0}$ is distributed as $\eta_{t-t_0}$,
up to a translation by $H_u$). The right-continuity 
in $P$-probability of the mapping $t\rightarrow \eta_t$ at $t=t_0$
follows from this observation, the property $\eta_{t_0}(\{H_{t_0}\})=0$,
the a.s. lower semi-continuity of $H_t$, and the case $t_0=0$.

The right-continuity in measure under $N$ follows from the same arguments.
\cq

\smallskip
Rather than investigating the Markovian properties
of $(\eta_t,t\geq 0)$ we will consider
the pair $(\rho_t,\eta_t)$. 
We first introduce some notation.
Let $(\mu,\nu)\in M_f(\R_+)^2$, and let $a>0$. Recall the notation of Proposition
\ref{explor}.
In a way analogous
to Chapter 1, we define $k_a(\mu,\nu)\in M_f(\R_+)^2$ by setting
$$k_a(\mu,\nu)=(\ov \mu,\ov\nu)$$
where $\ov \mu=k_a\mu$ and the measure $\ov\nu$ is the unique
element of $M_f(\R_+)$ such that
$$(\mu+\nu)_{|[0,H(k_a\mu)]}=k_a\mu+\ov\nu.$$
Note that the difference $\mu_{|[0,H(k_a\mu)]}-k_a\mu$ is a nonnegative
multiple of the Dirac measure at $H(k_a\mu)$, so that
$\ov \nu$ and $\nu_{|[0,H(k_a\mu)]}$ may only differ at the point
$H(k_a\mu)$.

Then, if $\theta_1=(\mu_1,\nu_1)\in M_f(\R_+)^2$
and $\theta_2=(\mu_2,\nu_2)\in M_f(\R_+)^2$, and if 
$H(\mu_1)<\infty$, we define the concatenation
$[\theta_1,\theta_2]$ by
$$[\theta_1,\theta_2]=([\mu_1,\mu_2],\nu)$$
where $\langle \nu,f\rangle=\int \nu_1(ds)1_{[0,H(\mu_1)]}(s)f(s)
+\int \nu_2(ds)f(H(\mu_1)+s)$.

\begin{proposition}
\label{Markov-eta}
{\rm (i)} Let $s\geq 0$ and $t>0$. Then, for every
nonnegative measurable function $f$ on $M_f(\R_+)^2$,
$$E[f(\rho_{s+t},\eta_{s+t})\mid \f_s]=\Pi^0_tf(\rho_s,\eta_s)$$
where $\Pi^0_t((\mu,\nu),d\mu'd\nu')$ is the distribution of the pair
$$[k_{-I_t}(\mu,\nu),(\rho_t,\eta_t)]$$
under $P$. The collection $(\Pi^0_t,t>0)$
is a Markovian semigroup on $M_f(\R_+)^2$.

{\rm (ii)} Let $s> 0$ and $t>0$. Then, for every
nonnegative measurable function $f$ on $M_f(\R_+)^2$,
$$N(f(\rho_{s+t},\eta_{s+t})\,1_{\{s+t<\sigma\}}\mid \f_s)=
1_{\{s<\sigma\}}\,\Pi_tf(\rho_s,\eta_s)$$
where $\Pi_t((\mu,\nu),d\mu'd\nu')$ is the distribution of the pair
$$[k_{-I_t}(\mu,\nu),(\rho_t,\eta_t)]$$
under $P(\cdot\cap \{T_{<\mu,1>}>t\})$. The collection $(\Pi_t,t>0)$
is a submarkovian semigroup on $M_f(\R_+)^2$.
\end{proposition}

\proof (i) Recall the notation of the proof of
Proposition \ref{explor}, and in particular formula (\ref{strongMarkov}). According to
this formula, we have
\begin{equation}
\label{Marko-rho}
\rho_{s+t}=[k_{-I^{(s)}_t}\rho_s,\rho^{(s)}_t]
\end{equation}
where the pair $(I^{(s)}_t,\rho^{(s)}_t)$ is defined in terms of the
shifted process $X^{(s)}$, which
is independent of $\f_s$. We then
want to get an analogous expression for $\eta_t$. Precisely,
we claim that
\begin{equation}
\label{Marko-rho-eta}
(\rho_{s+t},\eta_{s+t})=[k_{-I^{(s)}_t}(\rho_s,\eta_s),(\rho^{(s)}_t,\eta^{(s)}_t)]
\end{equation}
with an obvious notation. Note that (\ref{Marko-rho}) is the equality of the
first components in (\ref{Marko-rho-eta}).

To deal with the second components, recall the definition of $\eta_{s+t}$
$$\eta_{s+t}(du)=\beta 1_{[0,H_{s+t}]}(u)\,du+\build{\sum_{0<r\leq s+t}}
_{X_{r-}<I^r_{s+t}}^{} (X_r-I^r_{s+t})\,\delta_{H_r}(du).$$
First consider the absolutely continuous part. By (\ref{Marko-rho}),
we have
$$H_{s+t}=H(k_{-I^{(s)}_t}\rho_s)+H(\rho^{(s)}_t)=
H(k_{-I^{(s)}_t}\rho_s)+H^{(s)}_t$$
and thus
\ba
&&\int du\,1_{[0,H_{s+t}]}(u)\,f(u)\\
&&\qquad=\int du\,1_{[0,H(k_{-I^{(s)}_t}\rho_s)]}(u)\,f(u)
+\int du\,1_{[0,H^{(s)}_t]}(u)\,f(H(k_{-I^{(s)}_t}\rho_s)+u).
\ea
This shows that the absolutely continuous part of $\eta_{s+t}$
is the same as that of the second component of the right side of
(\ref{Marko-rho-eta}). 

Then the singular part of $\eta_{s+t}$ is equal to
\begin{equation} 
\label{sing-eta}
\build{\sum_{0<r\leq s}}
_{X_{r-}<I^r_{s+t}}^{} (X_r-I^r_{s+t})\,\delta_{H_r}
+\build{\sum_{s<r\leq s+t}}
_{X_{r-}<I^r_{s+t}}^{} (X_r-I^r_{s+t})\,\delta_{H_r}.
\end{equation}
Note that, if $r\in(s,s+t]$ is such that $X_{r-}<I^r_{s+t}$,
we have $H_r=H(k_{-I^{(s)}_t}\rho_s)+H^{(s)}_{r-s}$ (see the proof
of Proposition \ref{explor}). Thanks to this remark, we see that the
second term of the sum in (\ref{sing-eta}) is the image of the
singular part of $\eta^{(s)}_t$ under the mapping 
$u\rightarrow  H(k_{-I^{(s)}_t}\rho_s)+u$.

To handle the first term of (\ref{sing-eta}), we consider two cases.
Suppose first that $I_s<I^s_{s+t}$. Then set
$$v=\sup\{r\in(0,s]:X_{r-}<I^s_{s+t}\}.$$
In the first term of (\ref{sing-eta}), we need only consider
values $r\in(0,v]$. Note that $H_v=H(k_{-I^{(s)}_t}\rho_s)$
and that the measures $\rho_v$ and $k_{-I^{(s)}_t}\rho_s$
are equal except possibly at the point $H_v$
(see again the proof
of Proposition \ref{explor}). Then,
$$\build{\sum_{0<r<v}}
_{X_{r-}<I^r_{s+t}}^{} (X_r-I^r_{s+t})\,\delta_{H_r}
=\build{\sum_{0<r<v}}
_{X_{r-}<I^r_{s}}^{} (X_r-I^r_{s})\,\delta_{H_r}$$
coincides with the restriction of the singular part of
$\eta_s$ to $[0,H_v)=[0,H(k_{-I^{(s)}_t}\rho_s))$.
On the other hand, $\eta_{s+t}(\{H_v\})$
is equal to
$$X_v-I^v_{s+t}
=\eta_s(\{H_v\})+\rho_s([0,H(k_{-I^{(s)}_t}\rho_s)])-\langle
k_{-I^{(s)}_t}\rho_s,1\rangle$$
since by construction
\begin{eqnarray*}
&&\eta_s(\{H_v\})=X_v-I^v_s\,,\\
&&\rho_s([0,H(k_{-I^{(s)}_t}\rho_s)])=\rho_s([0,H_v])=I^v_s-I_s\,,\\
&&\langle k_{-I^{(s)}_t}\rho_s,1\rangle=X_s-I_s+I^{(s)}_t=I^s_{s+t}-I_s
=I^v_{s+t}-I_s\,.
\end{eqnarray*}
By comparing with the definition of $k_a(\mu,\nu)$,
we see that the proof of (\ref{Marko-rho-eta})
is complete in the case $I_s<I^s_{s+t}$.

The case $I_s\geq I^s_{s+t}$ is easier. In that case
$k_{-I^{(s)}_t}\rho_s=0$, and even $k_{-I^{(s)}_t}(\rho_s,\eta_s)=(0,0)$
(note that $\eta_s$ gives no mass to $0$, a.s.). Furthermore,
the first sum in (\ref{sing-eta}) vanishes, and it immediately
follows that (\ref{Marko-rho-eta}) holds.

The first assertion in (i) is a consequence of (\ref{Marko-rho-eta})
and the fact that $X^{(s)}$ is independent of $\f_s$.

As for the second assertion, it is enough to verify that,
for every $s,t>0$ we have
\begin{equation} 
\label{Marko-rho-eta-bis}
[k_{-I^{(s)}_t}[k_{-I_s}(\mu,\nu),(\rho_s,\eta_s)],(\rho^{(s)}_t,\eta^{(s)}_t)]
=[k_{-I_{s+t}}(\mu,\nu),(\rho_{s+t},\eta_{s+t})].
\end{equation}
Note that the case $\mu=\nu=0$ is just (\ref{Marko-rho-eta}).
To prove (\ref{Marko-rho-eta-bis}), we consider the same two cases as previously.

If $I^s_{s+t}> I_s$, or equivalently $-I^{(s)}_t< \langle \rho_s,1\rangle$,
then $I_s=I_{s+t}$ and so $k_{-I_s}(\mu,\nu)=k_{-I_{s+t}}(\mu,\nu)$. 
Furthermore, it is easy to verify that a.s.
$$k_{-I^{(s)}_t}[k_{-I_s}(\mu,\nu),(\rho_s,\eta_s)]
=[k_{-I_s}(\mu,\nu),k_{-I^{(s)}_t}(\rho_s,\eta_s)].$$
Hence
\ba
[k_{-I^{(s)}_t}[k_{-I_s}(\mu,\nu),(\rho_s,\eta_s)],(\rho^{(s)}_t,\eta^{(s)}_t)]
&=&[[k_{-I_s}(\mu,\nu),k_{-I^{(s)}_t}(\rho_s,\eta_s)],
(\rho^{(s)}_t,\eta^{(s)}_t)]\\
&=&[k_{-I_s}(\mu,\nu),[k_{-I^{(s)}_t}(\rho_s,\eta_s),
(\rho^{(s)}_t,\eta^{(s)}_t)]],
\ea
and (\ref{Marko-rho-eta-bis}) follows from (\ref{Marko-rho-eta}).

Finally, if $I^s_{s+t}\leq I_s$, or equivalently $-I^{(s)}_t>\langle
\rho_s,1\rangle$, it easily follows from our definitions
(and from the fact that $\eta_s(\{0\})=0$ a.s.) that
$$k_{-I^{(s)}_t}[k_{-I_s}(\mu,\nu),(\rho_s,\eta_s)]
=k_{-I_{s+t}}(\mu,\nu)\,,\qquad\hbox{a.s.}$$
Furthermore, the property $I^s_{s+t}\leq I_s$ also implies that
$(\rho^{(s)}_t,\eta^{(s)}_t) =(\rho_{s+t},\eta_{s+t})$,
and this completes the proof of (\ref{Marko-rho-eta-bis}).

(ii) First note that, for $s,t>0$, the identity
(\ref{Marko-rho-eta}) also holds $N$ a.e. on $\{s+t<\sigma\}$
with the same proof (the argument is even simpler as we
do not need to consider the case $I^s_{s+t}\leq I_s$). Also
observe that $N$ a.e. on $\{s<\sigma\}$, the condition
$s+t<\sigma$ holds iff $-I^{(s)}_t<X_s=\langle \rho_s,1\rangle$,
or equivalently $t<T^{(s)}_{<\rho_s,1>}=\inf\{r\geq 0:
X^{(s)}_r=-\langle \rho_s,1\rangle\}$. The first assertion 
in (ii) follows from these observations and the Markov property
under $N$.

The second assertion in (ii) follows from (\ref{Marko-rho-eta-bis})
and the fact that
\ba
\{T_{<\mu,1>}>s+t\}&=&\{I_{s+t}>-\langle \mu,1\rangle\}\\
&=&\{I_s>-\langle \mu,1\rangle\}\cap \{I^{(s)}_t>-\langle \mu,1\rangle-X_s\}\\
&=&\{I_s>-\langle \mu,1\rangle\}\cap \{T^{(s)}_{<k_{-I_s}\mu+\rho_s,1>}>t\}\,.
\ea
\ \cq

\smallskip
The previous proposition shows that the process
$(\rho_s,\eta_s)$ is Markovian under $P$. We now proceed
to investigate its invariant measure.

Let $\n(dsd\ell dx)$ be a Poisson point measure on $(\R_+)^3$
with intensity
$$ds\,\pi(d\ell)\,1_{[0,\ell]}(x)dx.$$
For every $a>0$, we denote by $\M_a$ the law on $M_f(\R_+)^2$
of the pair $(\mu_a,\nu_a)$ defined by
\ba
&&\langle \mu_a,f\rangle=\int \n(dsd\ell dx)\,1_{[0,a]}(s)\,xf(s)
+\beta\int_0^a ds\,f(s)\\
&&\langle \nu_a,f\rangle=\int \n(dsd\ell dx)\,1_{[0,a]}(s)\,(\ell-x)f(s)
+\beta\int_0^a ds\,f(s).
\ea
Note that $\M_a$ is invariant under the symmetry $(\mu,\nu)\rightarrow (\nu,\mu)$.
We also set
$$\M=\int_0^\infty da\,e^{-\alpha a}\,\M_a.$$
The marginals of $\M$ coincide with the measure $\MM$ of Chapter 1.

\begin{proposition}
\label{invariant-rho-eta}
Let $\Phi$ be a nonnegative measurable function on $M_f(\R_+)^2$. Then,
$$N\Big(\int_0^\sigma dt\,\Phi(\rho_t,\eta_t)\Big)
=\int \M(d\mu d\nu)\,\Phi(\mu,\nu).$$
\end{proposition}

\proof This is an extension of Proposition \ref{mesinv} and the
proof is much analogous. Consider (under $P$) the countable 
collection of instants $s_i$, $i\in I$ such that
$X_{s_i}>S_{s_i-}$. It follows from (\ref{joint}) that
\begin{equation} 
\label{Poisson-max}
\Big(L_\infty,\sum_{i\in I} \delta_{(L_{s_i},\Delta X_{s_i},
X_{s_i}-S_{s_i-})}(dsd\ell dx)\Big)
\build{=}_{}^{\rm(d)} \Big(\zeta,1_{[0,\zeta]}(s)
\n(dsd\ell dx)\Big)
\end{equation}
where $\zeta$ is an exponential variable with parameter $\alpha$
independent of $\n$ ($\zeta=\infty$ if $\alpha=0$). Recall from Chapter 1 the definition of the
time-reversed process $\wh X^{(t)}$. As in (\ref{rhoreversed}), we can 
rewrite the definition of $\rho_t$ and $\eta_t$ in terms of the reversed process $\wh X^{(t)}$:
\ba
&&\langle\rho_t,f\rangle=\beta \int_0^{\wh L^{(t)}_t}dr\,f(r)
+\build{\sum_{0<s\leq t}}
_{\wh X^{(t)}_s>\wh S^{(t)}_{s-}}^{} (\wh X^{(t)}_s-\wh
S^{(t)}_{s-})\,f(\wh L^{(t)}_t-\wh L^{(t)}_s),\\
&&\langle\eta_t,f\rangle=\beta \int_0^{\wh L^{(t)}_t}dr\,f(r)
+\build{\sum_{0<s\leq t}}
_{\wh X^{(t)}_s>\wh S^{(t)}_{s-}}^{} (\wh S^{(t)}_{s-}-\wh X^{(t)}_{s-}
)\,f(\wh L^{(t)}_t-\wh L^{(t)}_s).
\ea
Hence we can write $(\rho_t,\eta_t)=\Gamma(\wh X^{(t)}_{s\wedge t},s\geq 0)$
with a measurable functional $\Gamma$ that is made explicit in the
previous formulas. Proposition \ref{keyinv} now gives
$$N\Big(\int_0^\sigma dt\,\Phi(\rho_t,\eta_t)\Big)
=E\Big[\int_0^{L_\infty} da\,\Phi\circ\Gamma(X_{s\wedge L^{-1}(a)},s\geq 0)\Big].$$
However, $\Gamma(X_{s\wedge L^{-1}(a)},s\geq 0)=(\ov \mu_a,\ov\nu_a)$,
with
\ba
&&\langle\ov\mu_a,f\rangle=
\beta\int_0^a dr\,f(r)+\sum_{i\in I}
1_{[0,a]}(L_{s_i})\,(X_{s_i}-S_{s_i-})\,f(a-L_{s_i})\\
&&\langle\ov\nu_a,f\rangle=
\beta\int_0^a dr\,f(r)+\sum_{i\in I}
1_{[0,a]}(L_{s_i})\,(\Delta X_{s_i}-(X_{s_i}-S_{s_i-}))\,f(a-L_{s_i}).
\ea
Now use (\ref{Poisson-max}) to complete the proof. \cq

\smallskip
For every $t>0$, we denote by $\wh \Pi_t$ the image of the kernel
$\Pi_t$ under the symmetry $(\mu,\nu)\rightarrow (\nu,\mu)$, that is
$$\wh \Pi_t\Phi(\mu,\nu)=\int \Pi_t((\nu,\mu),d\nu'd\mu')\,\Phi(\mu',\nu').$$

\begin{theorem}
\label{duality-rho}
The kernels $\Pi_t$ and $\wh \Pi_t$ are in duality under $\M$.
\end{theorem}

This means that for any nonnegative measurable functions
$\Phi$ and $\Psi$ on $M_f(\R_+)^2$,
$$\M(\Phi \Pi_t\Psi)=\M(\Psi \wh\Pi_t\Phi).$$

\proof We first consider the potential kernels
$$U=\int_0^\infty dt\,\Pi_t\ ,\ \wh U=\int_0^\infty dt\,\wh\Pi_t$$
and we prove that
\begin{equation} 
\label{duality-pot}
\M(\Phi U\Psi)=\M(\Psi \wh U\Phi).
\end{equation}
This is equivalent to saying that the measure
$$\M(d\mu d\nu) \,U((\mu,\nu),d\mu'd\nu')$$
is invariant under the transformation $(\mu,\nu,\mu',\nu')\rightarrow 
(\nu',\mu',\nu,\mu)$.

To this end, we first derive an explicit expression for the kernel $U$.
By the definition of the kernels $\Pi_t$, we have
$$U\Phi(\mu,\nu)=
E\Big[\int_0^{T_{<\mu,1>}} dt\,\Phi([k_{-I_t}(\mu,\nu),(\rho_t,\eta_t)])\Big].$$
This is computed in a way similar to the proof of Proposition \ref{potker},
using Proposition \ref{invariant-rho-eta} in place of Proposition \ref{mesinv}.
It follows that
\begin{equation}
\label{explicit-pot}
U\Phi(\mu,\nu)=
\int_0^{<\mu,1>}dr\int \M(d\mu'd\nu')\,\Phi([k_r(\mu,\nu),(\mu',\nu')]).
\end{equation}
We then need to get more information about the joint
distribution of $((\mu,\nu),k_r(\mu,\nu))$ under $\M(d\mu d\nu)
1_{[0,{<\mu,1>}]}(r)dr$. Recall the notation $\n$, $\mu_a$, $\nu_a$
introduced before the statement of Proposition \ref{invariant-rho-eta}.
Write
$$\n=\sum_{i\in I} \delta_{(s_i,\ell_i,x_i)}$$
for definiteness, in such a way that
$$(\mu_a,\nu_a)=(\beta m_a+\sum_{s_i\leq a}x_i\,\delta_{s_i},
\beta m_a+\sum_{s_i\leq a}(\ell_i-x_i)\,\delta_{s_i}),$$
where $m_a$ denotes Lebesgue measure on $[0,a]$. Since 
$\M_a$ is the law of $(\mu_a,\nu_a)$, we get
\begin{eqnarray}
\label{duality-tech1}
&&\int \M_a(d\mu d\nu)\int_0^{<\mu,1>} dr\,F((\mu,\nu),k_r(\mu,\nu))\\
&&\qquad=E\Big[\int_0^{<\mu_a,1>} dr\,F((\mu_a,\nu_a),k_r(\mu_a,\nu_a))\Big]
\nonumber\\
&&\qquad=E\Big[\beta\int_0^a ds\,F((\mu_a,\nu_a),(\mu_{a|[0,s]},\nu_{a|[0,s]}))\Big]
\nonumber\\
&&\qquad+E\Big[\sum_{s_i\leq a} \int_0^{x_i} dy\,
F((\mu_a,\nu_a),(\mu_{a|[0,s_i)}+y\delta_{s_i},
\nu_{a|[0,s_i)}+(\ell_i-y)\delta_{s_i}))\Big]\nonumber
\end{eqnarray}
using the definition of $k_r$. 

At this point, we recall the following well-known lemma about 
Poisson measures.

\begin{lemma}
\label{Poisson-Palm}
Let $E$ be a measurable space and let $\Delta$ be a $\sigma$-finite
measure on $E$. Let $\m$ be a Poisson point measure
on $[0,a]\times E$ with intensity $ds\,\Delta(de)$. Then, for any nonnegative
measurable function $\Phi$,
$$E\Big[\int \m(dsde)\,\Phi((s,e),\m)\Big]
=E\Big[\int_0^a ds\int_E \Delta(de)\,\Phi((s,e),\m+\delta_{(s,e)})\Big].$$
\end{lemma}

Thanks to this lemma, the second term in the right side 
of (\ref{duality-tech1}) can be written as
\ba
&&E\Big[\int_0^a ds\int\!\pi(d\ell)\!\int_0^\ell dx\!\int_0^x dy\\
&&\hskip 1cm
F((\mu_a+x\delta_s,\nu_a+(\ell-x)\delta_s),(\mu_{a|[0,s)}+y\delta_s,
\nu_{a|[0,s)}+(\ell-y)\delta_s))\Big].\ea
We now integrate (\ref{duality-tech1}) with respect to $e^{-\alpha a}da$.
After some easy transformations, we get
\ba
&&\int \M(d\mu d\nu)\int_0^{<\mu,1>} dr\,F((\mu,\nu),k_r(\mu,\nu))\\
&&\qquad=\beta\int\M(d\mu_1d\nu_1)\M(d\mu_2d\nu_2)\,F([(\mu_1,\nu_1),(\mu_2,\nu_2)],
(\mu_1,\nu_1))\\
&&\qquad+\int\M(d\mu_1d\nu_1)\M(d\mu_2d\nu_2)\int\pi(d\ell)\int_0^\ell dx\int_0^x dy\\
&&\qquad\qquad\qquad F([(\mu_1,\nu_1),(x\delta_0+\mu_2,(\ell-x)\delta_0+\nu_2)],
[(\mu_1,\nu_1),(y\delta_0,(\ell-y)\delta_0)]).
\ea

Recalling formula (\ref{explicit-pot}) for the potential kernel
$U$, we see that the measure
$$\M(d\mu d\nu) \,U((\mu,\nu),d\mu'd\nu')$$
is the sum of two terms. The first one is the distribution under
$$\beta\M(d\mu_1d\nu_1)\M(d\mu_2d\nu_2)\M(d\mu_3d\nu_3)$$
of the pair 
$$(\mu,\nu)=[(\mu_1,\nu_1),(\mu_2,\nu_2)]\ ,\ 
(\mu',\nu')=[(\mu_1,\nu_1),(\mu_3,\nu_3)].$$
The second one is the distribution under
$$\M(d\mu_1d\nu_1)\M(d\mu_2,d\nu_2)\M(d\mu_3d\nu_3)\pi(d\ell)1_{\{0<y<x<\ell\}}
dx\,dy$$
of the pair
$$(\mu,\nu)=[(\mu_1,\nu_1),(x\delta_0+\mu_2,(\ell-x)\delta_0+\nu_2)]\,,\,
(\mu',\nu')=[(\mu_1,\nu_1),(y\delta_0+\mu_3,(\ell-y)\delta_0+\nu_3)].$$
In this form, it is clear that $\M(d\mu d\nu) \,U((\mu,\nu),d\mu'd\nu')$ has the
desired invariance property. This completes the proof of (\ref{duality-pot}).

Consider now the resolvent kernels
$$U_p((\mu,\nu),d\mu'd\nu')=\int_0^\infty dt\,e^{-pt}\Pi_t((\mu,\nu),d\mu',d\nu').$$
By a standard argument (see e.g. \cite{DM}, p.54), (\ref{duality-pot}) also implies that,
for every $p>0$, $\M(\Phi U_p\Psi)=\M(\Psi \wh U_p\Phi)$, or equivalently
\begin{equation} 
\label{duality-resolvent}
\int_0^\infty dt\,e^{-pt}\,\M(\Phi\Pi_t\Psi)=
\int_0^\infty dt\,e^{-pt}\,\M(\Psi\Pi_t\Phi).
\end{equation}
Recall that our goal is to prove the identity $\M(\Phi\Pi_t\Psi)=\M(\Psi\Pi_t\Phi)$
for every $t>0$.
We may assume that the functions $\Phi$ and $\Psi$ are 
continuous and both dominated by $e^{-a<\mu,1>}$ for some $a>0$. The latter
condition guarantees that $\M(\Phi)<\infty$ and $\M(\Psi)<\infty$. From the
definition of $\Pi_t$ and the right-continuity in probability
of the mapping $t\rightarrow (\rho_t,\eta_t)$ (Lemma \ref{basic-eta}), it is easy to verify
that $t\rightarrow  \Pi_t\Psi(\mu,\nu)$ is right-continuous over $(0,\infty)$.
The same holds for the mapping $t\rightarrow  \M(\Phi\Pi_t\Psi)$, and the statement
of the theorem follows from (\ref{duality-resolvent}). \cq

\medskip
For notational reasons, we make the convention that $\rho_s=\eta_s=0$
if $s<0$.  

\begin{corollary}
\label{reversal}
The process $(\eta_s,s\geq 0)$ has a c\` adl\` ag modification under $N$
or under $P$. Furthermore, the processes $(\rho_s,\eta_s;s\geq 0)$ and
$(\eta_{(\sigma-s)-},\rho_{(\sigma-s)-};s\geq 0)$ have the same
distribution under $N$.
\end{corollary}

A consequence of the corollary is the fact that the processes
$(H_t,t\geq 0)$ and $(H_{(\sigma-t)_+},t\geq 0)$ have the same 
distribution (say in the sense of finite-dimensional marginals
when $H$ is not continuous) under $N$. 
In view of the results of Chapter 2, this is not surprising, as the 
same time-reversal property obviously holds for the discrete contour
process. The more precise statement of the corollary will be
useful in the next sections.

\smallskip
\proof The second part of the corollary is essentially a consequence of the
duality property stated in the previous theorem. Since we have still little 
information about regularity properties of the process $\eta$, we will proceed
with some care. We first introduce the Kuznetsov measure $\K$, which is
the $\sigma$-finite measure on $\R\times \D(\R_+,\R)$ defined by
$$\K(drd\omega)=dr\,N(d\omega).$$
We then define $\gamma(r,\omega)=r$, $\delta(r,\omega)=r+\sigma(\omega)$ 
and, for every $t\in \R$,
$$\ov \rho_t(r,\omega)=\rho_{t-r}(\omega)\ ,\
\ov\eta_t(r,\omega)=\eta_{t-r}(\omega)$$
with the convention explained before the statement of the corollary.
Note that $(\ov\rho_t,\ov\eta_t)\not =(0,0)$ iff $\gamma<t<\delta$. 

It readily follows from Proposition \ref{invariant-rho-eta} that, for every $t\in\R$,
the distribution of $(\ov\rho_t,\ov\eta_t)$ under 
$\K(\cdot \cap \{(\ov\rho_t,\ov\eta_t)\not =(0,0)\})$ is $\M$. Let
$t_1,\ldots,t_p\in\R$ with $t_1<t_2<\cdots<t_p$. Using Proposition
\ref{Markov-eta} and induction on $p$, we easily get that the
restriction to $(M_f(\R_+)^2\backslash \{(0,0)\})^p$ of the distribution
of the $p$-tuple 
$((\ov\rho_{t_1},\ov\eta_{t_1}),\ldots,(\ov\rho_{t_p},\ov\eta_{t_p}))$
is
$$\M(d\mu_1d\nu_1)\,\Pi_{t_2-t_1}((\mu_1,\nu_1),d\mu_2d\nu_2)
\ldots \Pi_{t_p-t_{p-1}}((\mu_{p-1},\nu_{p-1}),d\mu_pd\nu_p).$$
By Theorem \ref{duality-rho}, this measure is equal to
$$\M(d\mu_pd\nu_p)\,\wh\Pi_{t_p-t_{p-1}}((\mu_p,\nu_p),d\mu_{p-1}d\nu_{p-1})
\ldots \wh\Pi_{t_2-t_{1}}((\mu_{2},\nu_{2}),d\mu_1d\nu_1).$$
Hence the two $p$-tuples
$((\ov\rho_{t_1},\ov\eta_{t_1}),\ldots,(\ov\rho_{t_p},\ov\eta_{t_p}))$
and $((\ov\eta_{-t_1},\ov\rho_{-t_1}),\ldots,(\ov\eta_{-t_p},\ov\rho_{-t_p}))$
have the same distribution, in restriction to $(M_f(\R_+)^2\backslash \{(0,0)\})^p$,
under $\K$. Since $(\ov\rho_t,\ov\eta_t)\not =(0,0)$ iff $\gamma<t<\delta$, a simple
argument shows that we can remove the restriction and conclude that these
two $p$-tuples have the same distribution under $\K$. (This distribution is
$\sigma$-finite except for an infinite mass at the point $(0,0)^p$.)
 
In particular, $(\ov\rho_{t_1},\ldots,\ov\rho_{t_p})$ and 
$(\ov\eta_{-t_1},\ldots,\ov\eta_{-t_p})$ have the same distribution
under $\K$. Let $F$ be a bounded continuous function on $M_f(\R_+)^p$,
such that $F(0,\ldots,0)=0$.
Suppose that $0<t_1<t_2<\ldots<t_p$ and let $u<v$.
Then we have
\ba
&&\K(1_{[u,v]}(\gamma)F(\ov\rho_{\gamma+t_1},\ldots,\ov\rho_{\gamma+t_p}))\\
&&\quad=\lim_{\varepsilon\rightarrow  0}\sum_{k\in\Z,\,k\varepsilon\in[u,v]}\K\Big(1_{\{\ov\rho_{k\varepsilon}=0,
\ov\rho_{(k+1)\varepsilon}\not=0\}}F(\ov\rho_{k\varepsilon+t_1},\ldots,\ov\rho_{k\varepsilon+t_p})\Big)\ea
and the similar formula 
\ba
&&\K(1_{[u,v]}(\delta)F(\ov\eta_{\delta-t_1},\ldots,\ov\eta_{\delta-t_p}))\\
&&\quad=\lim_{\varepsilon\rightarrow  0}\sum_{k\in\Z,\,k\varepsilon\in[u,v]}\K\Big(1_{\{\ov\eta_{-k\varepsilon}=0,
\ov\eta_{-(k+1)\varepsilon}\not=0\}}F(\ov\eta_{-k\varepsilon-t_1},\ldots,\ov\eta_{-k\varepsilon-t_p})\Big),
\ea
where we use the right-continuity in $N$-measure of $\eta_t$. Hence the vectors
$(\gamma,\ov\rho_{\gamma+t_1},\ldots,\ov\rho_{\gamma+t_p})$
and $(-\delta,\ov\eta_{\delta-t_1},\ldots,\ov\eta_{\delta-t_p})$ have the same
distribution under $\K$. It follows that the processes $(\rho_t,t\geq 0)$ and $(\eta_{\sigma-t},
t\geq 0)$ have the
same finite-dimensional marginals under $N$. 
Since we already know that $(\rho_t,t\geq 0)$ is c\` adl\` ag, we obtain 
that $(\eta_t,t\geq 0)$ has a c\` adl\` ag modification under $N$. The time-reversal
property of the corollary follows immediately from the previous identification
of finite-dimensional marginals. This property implies in particular that
$\eta_{0+}=\eta_{\sigma-}=0$ $N$ a.e.

\smallskip 
It remains to verify that $(\eta_t,t\geq 0)$ has a c\` adl\` ag modification
under $P$. On each excursion interval of $X-I$ away from $0$,
we can apply the result derived above under the excursion measure $N$. It remains to deal with instants 
$t$ such that $X_t=I_t$, for which $\eta_t=0$. To this end, we note that, for every
$\varepsilon>0$,
$$N\Big(\sup_{s\in[0,\sigma]}\langle\eta_s,1\rangle >\varepsilon\Big)
=N\Big(\sup_{s\in[0,\sigma]}\langle\rho_s,1\rangle >\varepsilon\Big)<\infty.$$
Hence, for any fixed $x>0$, we will have $\langle\eta_s,1\rangle \leq\varepsilon$
for all $s\in[0,T_x]$ except possibly for $s$ belonging to finitely
many excursion intervals of $X-I$. Together with the
continuity of $\eta$ at times $0$ and $\sigma$ under $N$, this
implies that $P$ a.s. for every $t$ such that $X_t=I_t$,
the right and left limits of $\eta_s$ both exist at time $t$
and vanish. \cq

\section{The tree associated with Poissonnian marks}

\subsection{Trees embedded in an excursion}

We first give the definition of the tree associated with a continuous
function $e:[a,b]\rightarrow  \R_+$ and $p$ instants $t_1,\ldots,t_p$
with $a\leq t_1\leq t_2\leq \cdots\leq t_p\leq b$.

Recall from Section 0.1 the definition of a (finite) rooted ordered tree, and the notation
${\bf T}$ for the collection of these trees. If $v$ is an individual (a vertex) in the
tree $\t\in{\bf T}$, the notation $k_v(\t)$ stands for the number of children of $v$. Individuals $v$
without children, i.e. such that $k_v(\t)=0$, are called leaves.
For every $p\geq 1$, we denote by $\TT_p$ the set of all (rooted ordered)
trees with $p$ leaves. 

If $\t^1,\t^2,\ldots,\t^k$ are $k$ trees, the concatenation of $\t^1,\ldots,\t^k$, which is
denoted by $[\t^1,\t^2,\ldots,\t^k]$, is defined in the obvious way: For $n\geq 1$, $(i_1,\ldots, i_n)$
belongs to $[\t^1,\t^2,\ldots,\t^k]$ if and only if $1\leq i_1\leq k$ and $(i_2,\ldots ,i_n)$ belongs
to $\t^{i_1}$.

\smallskip
A {\it marked tree} is a pair
$\theta=(\t,\{h_v,v\in \t\})$, where $h_v\geq 0$ for every $v\in \t$.
The number $h_v$ is interpreted as the lifetime of individual $v$, and $\t$ is 
called the skeleton of $\theta$.
We denote by $\T_p$ the set of all marked trees with $p$ leaves.

Let $\theta^1=(\t^1,\{h^1_v,v\in \t\})\in\T_{p_1},\ldots,\theta^k=(\t^k,\{h^k_v,v\in \t^k\})\in \T_{p_k}$,
and $h\geq 0$.
The concatenation
$[\theta^1,\theta^2,\ldots,\theta^k]_h$
is the element of $\T_{p_1+\ldots+ p_k}$ whose skeleton
is $[\t^1,\t^2,\ldots,\t^k]$ and such that the lifetimes of vertices in $\t^i$, $1\leq i\leq k$
become the lifetimes of the
corresponding vertices in $[\t^1,\t^2,\ldots,\t^k]$, and finally the lifetime
of $\emptyset$ in $[\theta^1,\theta^2,\ldots,\theta^k]_h$ is $h$.

\medskip
Let $e:[a,b]\longrightarrow \R_+$ be a continuous function
defined on a subinterval $[a,b]$ of $\R_+$. For every $a\leq u\leq v\leq b$, we set
$$m(u,v)=\inf_{u\leq t\leq v} e(t).$$
Let $t_1,\ldots,
t_p\in\R_+$ be such that $a\leq t_1\leq t_2\leq \cdots\leq t_p\leq b$. We
will now construct a marked tree
$$\theta(e,t_1,\ldots,t_p)=(\t(e,t_1,\ldots,t_p),\{h_v(e,t_1,\ldots,t_p),
v\in \t\})\in \T_p$$
associated with the function $e$ and the times $t_1,\ldots,t_p$.
We proceed by induction on $p$. If $p=1$, $\t(e,t_1)=\{\emptyset\}$ and 
$h_\emptyset(e,t_1)=e(t_1)$. 

\smallskip
Let $p\geq 2$ and suppose that the tree has been constructed up to 
order $p-1$. Then there exists an integer $k\in\{1,\ldots,p-1\}$
and $k$ integers $1\leq i_1<i_2<\cdots<i_k\leq p-1$ such that 
$m(t_{i},t_{i+1})=m(t_1,t_p)$ iff $i\in\{i_1,\ldots,i_k\}$. For every
$\ell\in\{0,1,\ldots,k\}$, define $e^\ell$ by the formulas
$$\begin{array}{lll}
e^0(t)=e(t)-m(t_1,t_p),&\qquad t\in[a,t_{i_1}],&\\
e^\ell(t)=e(t)-m(t_1,t_p),&\qquad t\in[t_{i_\ell+1},t_{i_{\ell+1}}],&\qquad 1\leq \ell\leq k-1.\\
e^k(t)=e(t)-m(t_1,t_p),&\qquad t\in[t_{i_k+1},b].&
\end{array}
$$
We then set:
$$\theta(e,t_1,\ldots,t_p)=
[\theta(e^0,t_1,\ldots,t_{i_1}),\theta(e^1,t_{i_1+1},\ldots,t_{i_2}),\ldots,
\theta(e^k,t_{i_k+1},\ldots,t_p)]_{m(t_1,t_p)}.$$
This completes the construction of the tree by induction. Note that $k+1$ is the number of children
of $\emptyset$ in the tree $\theta(e,t_1,\ldots,t_p)$, and $m(t_1,t_p)$ is the lifetime
of $\emptyset$.

\subsection{Poissonnian marks}

We consider a standard Poisson process with parameter $\lambda$ defined under the 
probability measure $Q_\lambda$. We denote by $\tau_1\leq \tau_2\leq \cdots$
the jump times of this Poisson process. Throughout this section, we argue 
under the measure $Q_\lambda\otimes N$, which means that we consider 
the excursion measure of $X-I$ together with independent Poissonnian marks with
intensity $\lambda$ on $\R_+$. To simplify notation however, we will 
systematically write $N$ instead of $Q_\lambda\otimes N$.

Set $M=\sup\{i\geq 1:\tau_i\leq \sigma\}$, which represents the number 
of marks that fall in the excursion interval (by convention, $\sup\emptyset=0$). Then,
$$N(M\geq 1)=N(1-e^{-\lambda\sigma})=\psi^{-1}(\lambda),$$
where the second equality follows from the fact that the 
Laplace exponent of the subordinator $T_x$ is $\psi^{-1}(\lambda)$
(see \cite{Be}, Theorem VII.1).

From now on, we assume that the condition $\int_1^\infty {du\over \psi(u)}<\infty$
holds, so that $H$ has continuous sample paths (Theorem \ref{continuityheight}). We can then
use subsection 3.2.1 to define the embedded tree $\theta(H,\tau_1,\ldots,\tau_M)$
under $N(\cdot\mid M\geq 1)$. Our main goal is to determine the law of this
tree.

\begin{theorem}
\label{tree-Poisson}
Under the probability measure $N(\cdot\mid M\geq 1)$, the tree
$\theta(H,\tau_1,\ldots,\tau_M)$ is distributed as the family tree of
a continuous-time Galton-Watson process starting with one
individual at time $0$ and such that:

$\bullet$ Lifetimes of individuals have exponential distributions with
parameter $\psi'(\psi^{-1}(\lambda))$;

$\bullet$ The offspring distribution is the law of the variable $\xi$
with generating function
$$E[r^\xi]=r+{\psi((1-r)\psi^{-1}(\lambda))\over \psi^{-1}(\lambda)
\psi'(\psi^{-1}(\lambda))}.$$
\end{theorem}

\rem As the proof will show, the theorem remains valid without the assumption that 
$H$ has continuous paths. We will leave this extension to the reader. 
Apart from some technical details, it simply
requires the straightforward extension of the construction of subsection 3.2.1 to the case
when the function $e$ is only lower semicontinuous.

\smallskip
The proof of Theorem \ref{tree-Poisson} requires a few intermediate results.
To simplify notation, we will write $\tau=\tau_1$.
We start with an important application of Corollary \ref{reversal}.

\begin{lemma}
\label{rho-expo}
For any nonnegative measurable function $f$ on $M_f(\R_+)$,
$$N(f(\rho_\tau)1_{\{M\geq 1\}})=\lambda\int \M(d\mu d\nu)\,f(\mu)\,e^{-\psi^{-1}(\lambda)\,\langle
\nu,1\rangle}.$$
\end{lemma}

\proof We have
$$N(f(\rho_\tau)1_{\{M\geq 1\}})
=\lambda N\Big(\int_0^\sigma dt\,e^{-\lambda t}\,f(\rho_t)\Big)
=\lambda N\Big(\int_0^\sigma dt\,e^{-\lambda (\sigma-t)}\,f(\eta_t)\Big),$$
using the time-reversal property of Corollary \ref{reversal}. At this point,
we use the Markov property of $X$ under $N$:
$$N\Big(\int_0^\sigma dt\,e^{-\lambda (\sigma-t)}\,f(\eta_t)\Big)
=N\Big(\int_0^\sigma dt\,f(\eta_t)\,E_{X_t}[e^{-\lambda T_0}]\Big).$$
We have already noticed that for $x\geq 0$,
$$E_x[e^{-\lambda T_0}]=E_0[e^{-\lambda T_x}]=e^{-x \psi^{-1}(\lambda)}.$$
Since $X_t=\langle \rho_t,1\rangle$ under $N$, it follows that
$$N(f(\rho_\tau)1_{\{M\geq 1\}})
=\lambda N\Big(\int_0^\sigma\! dt\,f(\eta_t)\,e^{-\langle \rho_t,1\rangle \psi^{-1}(\lambda)}\Big)
=\lambda \!\int\!\! \M(d\mu d\nu)\,f(\nu)\,e^{-\langle \mu,1\rangle \psi^{-1}(\lambda)},$$
using Proposition \ref{invariant-rho-eta}. Since $\M$ is invariant
under the mapping $(\mu,\nu)\rightarrow (\nu,\mu)$, this completes the proof.\cq

\smallskip
We now set
$$K=\left\{\begin{array}{ll}
\inf\{H_s:\tau_1\leq s\leq \tau_M\}\qquad&\hbox{if } M\geq 2\\
\infty&\hbox{if } M\leq 1
\end{array}
\right.
$$
Then $K$ represents the lifetime of the ancestor in the tree
$\theta(H,\tau_1,\ldots,\tau_M)$ (assuming that the event $\{M\geq 2\}$ holds). To give a formula
for the number of children $\xi$ of the ancestor, set
$$\tau_{(K)}=\inf\{t\geq \tau:H_t\leq K\}\ ,\ \tau'_{(K)}=\inf\{t\geq \tau:H_t< K\}.$$
Then, again on the event $\{M\geq 2\}$, $\xi-1$ is the number of excursions
of $H$ above level $K$, on the time interval $[\tau_{(K)},\tau'_{(K)}]$, which
contain at least one of the Poissonnian marks. This 
identification follows readily from the construction of subsection 3.2.1. 

The next proposition gives the joint distribution
of the pair $(K,\xi)$ under $N(\cdot\cap\{M\geq 2\})$. 

\begin{proposition}
\label{ancestor-lifetime}
Let $r\in[0,1]$ and let $\varphi$ be a nonnegative measurable function on
$[0,\infty]$, with $\varphi(\infty)=0$. Then,
\ba
&&N(r^\xi\varphi(K)\mid M\geq 1)\\
&&\qquad=\Big(r\psi'(\psi^{-1}(\lambda))+{\psi((1-r)\psi^{-1}(\lambda))-\lambda\over \psi^{-1}(\lambda)}\Big)
\int_0^\infty db\,\varphi(b)\,e^{-b\psi'(\psi^{-1}(\lambda))}.
\ea
\end{proposition}

The basic idea of the proof
is to apply the Markov property to the process $\rho$ at time $\tau$. To this
end, we need some notation. We write $\PP_\mu$ for the probability measure 
under which $\rho$
starts at an arbitrary measure $\mu\in\M_f(\R^d)$ and is stopped when it hits $0$.
As usual, $H_t=H(\rho_t)$.
Under $\PP_\mu$, the process $X_t=\langle \rho_t,1\rangle$
is the underlying L\'evy process started at $\langle \mu,1\rangle$, and as usual
$T_0=\inf\{t\geq 0:X_t=0\}$. We keep the notation
$I_t$ for the minimum process of $X$. We let $(a_j,b_j)$, $j\in J$ be the
collection of excursion intervals of $X-I$ away from $0$ and before time $T_0$. For
every $j\in J$ we define the corresponding excursion by
$$\omega_j(t)=X_{(a_j+t)\wedge b_j}-I_{a_j}\,,\quad t\geq 0.$$
From excursion theory, we know that the
point measure
$$\sum_{j\in J} \delta_{(I_{a_j},\omega_j)}$$
is Poisson under $\PP_\mu$, with intensity $1_{[0,<\mu,1>]}(u)du\,N(d\omega)$
(cf the proof of Proposition \ref{potker}). On the other hand, by properties 
of the exploration process derived in Chapter 1, we know that $\PP_\mu$
a.s. for every $s\in[0,T_0]$ such that $X_s-I_s=0$ (and in particular for
$s=a_j$, $j\in J$) we have $\rho_s=k_{<\mu,1>-I_{s}}\mu$ and thus
$H_s=H(k_{<\mu,1>-I_{s}}\mu)$. Observe also that the image of the measure
$1_{[0,<\mu,1>]}(u)du$ under the mapping $u\rightarrow  H(k_{<\mu,1>-u}\mu)$
is exactly $\mu(dh)$. By combining these observations, we get:

\smallskip
\noindent (P)
The point measure
$\sum_{j\in J} \delta_{(H_{a_j},\omega_j)}$
is Poisson under $\PP_\mu$, with intensity $\mu(dh)\,N(d\omega)$.

\smallskip
Finally, assume that we are also given
a collection ${\cal P}_\lambda$ of Poisson marks with intensity $\lambda$,
independently of $\rho$ under $\PP_\mu$, and set
\begin{eqnarray*}
&&L=\inf\{H_{a_j}: j\in J,(a_j,b_j)\cap {\cal P}_\lambda\not =\emptyset\},\\
&&\zeta={\rm Card}\{j\in J:H_{a_j}=L\hbox
{ and }(a_j,b_j)\cap {\cal P}_\lambda\not =\emptyset\}.
\end{eqnarray*}

Then the Markov property of the exploration process at time $\tau$
shows that, for any nonnegative measurable function $\varphi$ on $[0,\infty]$,
\begin{equation} 
\label{ancestor-key}
N(1_{\{M\geq 1\}}r^{\xi-1}\varphi(K))=N(1_{\{M\geq
1\}}\EE_{\rho_\tau}[r^\zeta\varphi(L)]).
\end{equation}
To verify this equality, simply observe that those excursions
of $H$ above level $K$ on the time interval $[\tau_{(K)},\tau'_{(K)}]$ that
contain one Poissonnian mark, exactly correspond to those excursions of
the shifted process $\langle\rho_{\tau+t},1\rangle$ above its minimum that
start from the height $K$ and contain one mark.

The next lemma is the key step towards the proof of Proposition
\ref{ancestor-lifetime}.

\begin{lemma}
\label{ances-explor}
Let $a\geq 0$ and let $\mu\in M_f(\R_+)$ be such that $\supp \mu=[0,a]$ and
$\mu(dt)=\beta 1_{[0,a]}(t)dt + \mu_s(dt)$, where $\mu_s$ is
a countable sum of multiples of Dirac point masses at elements of $[0,a]$. Then,  if
$r\in[0,1]$ and $\varphi$ is a nonnegative measurable function on
$[0,\infty]$ such that $\varphi(\infty)=0$,
\begin{eqnarray}
\label{ances-explo}
&&\EE_\mu[r^{\zeta}\varphi(L)]=\beta r \psi^{-1}(\lambda) \int_0^a db\,
e^{-\psi^{-1}(\lambda)\mu([0,b])}\,\varphi(b)\nonumber\\
&&+\sum_{\mu(\{s\})>0} \Big(e^{-(1-r)\mu(\{s\})\psi^{-1}(\lambda)}
-e^{-\mu(\{s\})\psi^{-1}(\lambda)}\Big)e^{-\mu([0,s))\psi^{-1}(\lambda)}\,\varphi(s).
\end{eqnarray}
\end{lemma}

\proof First note that it is easy to derive the law of
$L$ under $\PP_\mu$. Let $b\in [0,a]$. We have by property (P)
\begin{eqnarray*}
\PP_\mu[L>b]&=&\PP_\mu[(a_j,b_j)\cap {\cal P}_\lambda=\emptyset \hbox{ for every }
j\in J\hbox{ s.t. } H_{a_j}\leq b]\\
&=&\EE_\mu[\exp(-\mu([0,b])N(1-e^{-\lambda \sigma}))]\\
&=&\exp(-\mu([0,b])\psi^{-1}(\lambda)).
\end{eqnarray*}
In particular, atoms of the distribution of $L$ exactly correspond to atoms of $\mu$, and the
continuous part of the distribution of
$L$ is the measure
$$\beta\psi^{-1}(\lambda)\exp(-\mu([0,b])\psi^{-1}(\lambda))1_{[0,a]}(b)db.$$

We then need to distinguish two cases:

\smallskip
\noindent(1) Let $s\in[0,a]$ be an atom of $\mu$. By the preceding formula,
$$\PP_\mu[L=s]=(1-e^{-\mu(\{s\})\psi^{-1}(\lambda))})e^{-\mu([0,s))\psi^{-1}(\lambda)}.$$
Note that the excursions $\omega_j$ that start at height $s$ are the atoms of
a Poisson measure with intensity $\mu(\{s\})N$.
Using also the independence properties of Poisson measures, we get that, conditionally on
$\{L=s\}$, $\xi$ is distributed  as a Poisson random variable with intensity
$\mu(\{s\})\psi^{-1}(\lambda)$, conditioned to be greater than or equal to $1$:
$$\EE_\mu[r^\xi\mid L=s]={e^{-(1-r)\mu(\{s\})\psi^{-1}(\lambda))}-e^{-\mu(\{s\})\psi^{-1}(\lambda))}
\over 1-e^{-\mu(\{s\})\psi^{-1}(\lambda))}}.$$

\smallskip\noindent
(2) If $L$ is not an atom of $\mu$, then automatically $\xi=1$. This is so because 
the values $H_{a_j}$ corresponding to indices $j$ such that $\mu(\{H_{a_j}\})=0$
must be distinct, by (P) and standard properties of Poisson measures.

\smallskip
The lemma follows by combining these two cases with the distribution of $L$. \cq

\medskip
\noindent{\bf Proof of Proposition \ref{ancestor-lifetime}.} By
combining (\ref{ancestor-key}), Lemma \ref{rho-expo} and (\ref{ances-explo}), we obtain 
that
$$N(1_{\{M\geq 1\}}r^{\xi-1}\varphi(K))=A_1+A_2$$
where
$$
A_1=\beta r\lambda\psi^{-1}(\lambda)\int_0^\infty da\,e^{-\alpha a}
\int \M_a(d\mu d\nu)e^{-\langle \nu,1\rangle\psi^{-1}(\lambda)}\int_0^a db\,
e^{-\psi^{-1}(\lambda)\mu([0,b])}\,\varphi(b),$$
and 
\ba 
&&A_2=\lambda\int_0^\infty da\,e^{-\alpha a}
\int \M_a(d\mu d\nu)e^{-\langle \nu,1\rangle\psi^{-1}(\lambda)}\\
&&\qquad\qquad\times
\sum_{\mu(\{s\})>0} \Big(e^{-(1-r)\mu(\{s\})\psi^{-1}(\lambda)}
-e^{-\mu(\{s\})\psi^{-1}(\lambda)}\Big)e^{-\mu([0,s))\psi^{-1}(\lambda)}\,\varphi(s).
\ea
To compute $A_1$, we observe that for $u>0$ and $0\leq b\leq a$,
\ba
\M_a(e^{-u(\mu([0,b])+\nu([0,a])})&=&\M_a(e^{-u(\mu+\nu)([0,b])})\M_a(e^{-u\nu((b,a])})\\
&=&e^{-\beta u(a+b)}\exp\Big(-b\int \pi(d\ell)\ell(1-e^{-u\ell})\Big)\\
&&\times\exp\Big(-(a-b)\int \pi(d\ell)\int_0^\ell dx(1-e^{-ux})\Big)\\
&=&e^{\alpha a}\exp(-b\psi'(u)-(a-b){\psi(u)\over u}),
\ea
using the easy formulas
\ba
&&\int \pi(d\ell)\ell(1-e^{-u\ell})=\psi'(u)-\alpha-2\beta u\,,\\
&&\int \pi(d\ell)\int_0^\ell dx(1-e^{-ux})={1\over u}(\psi(u)-\alpha u-\beta u^2)\,.
\ea
It follows that 
\ba A_1&=&\beta r\lambda\psi^{-1}(\lambda)\int_0^\infty da\int_0^a db\, \varphi(b)
e^{-b\psi'(\psi^{-1}(\lambda))- (a-b)(\lambda/\psi^{-1}(\lambda))}\\
&=&\beta r
\psi^{-1}(\lambda)^2\int_0^\infty db\,\varphi(b)e^{-b\psi'(\psi^{-1}(\lambda))}.
\ea

To evaluate $A_2$, first observe that, with the notation
preceding Proposition \ref{invariant-rho-eta}, we have
\ba
A_2&=&\lambda\int_0^\infty da\,e^{-\alpha a}\,E\Big[\int_{\{s\leq a\}}{\cal N}(dsd\ell dx)
\,\varphi(s)\,(e^{-(1-r)x\psi^{-1}(\lambda)}
-e^{-x\psi^{-1}(\lambda)})\\
&&\times \exp\Big(-\psi^{-1}(\lambda)\Big(\int_{\{s'\leq a\}}{\cal N}(ds'd\ell' dx')(\ell'-x')
+\int_{\{s'<s\}}{\cal N}(ds'd\ell' dx')x'\Big)\Big)\Big].
\ea
From Lemma \ref{Poisson-Palm}, it follows that
\ba 
A_2&=&\lambda\int_0^\infty da\,e^{-\alpha a}\int_0^a db\,\varphi(b)\,
\M_a(e^{-(\mu([0,b])+\nu([0,a]))\psi^{-1}(\lambda)})\\
&&\hskip 15mm \times\int\pi(d\ell)\int_0^\ell dx
(e^{-(1-r)x\psi^{-1}(\lambda)}
-e^{-x\psi^{-1}(\lambda)})e^{-(\ell-x)\psi^{-1}(\lambda)}\\
&=&\psi^{-1}(\lambda)\int_0^\infty db\,\varphi(b)\,e^{-b\psi'(\psi^{-1}(\lambda))}\\
&&\hskip 15mm \times\int\pi(d\ell)\int_0^\ell dx
(e^{-(1-r)x\psi^{-1}(\lambda)}
-e^{-x\psi^{-1}(\lambda)})e^{-(\ell-x)\psi^{-1}(\lambda)}
\ea
where the last equality is obtained from the same calculations as those made in
evaluating $A_1$. Furthermore, straightforward calculations give
\ba
&&\int\pi(d\ell)\int_0^\ell dx
(e^{-(1-r)x\psi^{-1}(\lambda)}
-e^{-x\psi^{-1}(\lambda)})e^{-(\ell-x)\psi^{-1}(\lambda)}\\
&&\qquad\qquad=\psi'(\psi^{-1}(\lambda))-\beta r\psi^{-1}(\lambda)+{1\over r\psi^{-1}(\lambda)}
(\psi((1-r)\psi^{-1}(\lambda))-\lambda).
\ea
By substituting this in the previous display and combining with the formula for $A_1$, we
arrive at the result of the proposition. \cq

\medskip
\noindent{\bf Proof of Theorem \ref{tree-Poisson}.} It is convenient to introduce the
random variable $\Lambda$ defined by
$$\Lambda=\left\{\begin{array}{ll}
K\qquad&\hbox{if }M\geq 2\\
H_\tau\qquad&\hbox{if }M=1
\end{array}
\right.
$$
On the event $\{M=1\}$ we also set $\xi=0$.
We can easily compute the law of the pair $(\Lambda,\xi)$. Indeed, by applying the
Markov property at $\tau$ as previously, we easily get
\ba
N(\varphi(\Lambda)1_{\{M=1\}})&=&N(\varphi(H_\tau)1_{\{M=1\}})\\
&=&N(\varphi(H_\tau)\,e^{-<\rho_\tau,1>\psi^{-1}(\lambda)})\\
&=&\lambda\int_0^\infty da\,e^{-\alpha a}\,\varphi(a)\int
\M(d\mu d\nu)\,e^{-(<\mu,1>+<\nu,1>)\psi^{-1}(\lambda)})\\ 
&=&\lambda \int_0^\infty da\,\varphi(a)\,e^{-\psi'(\psi^{-1}(\lambda))a}
\ea
By combining with Proposition \ref{ancestor-lifetime}, we get
\ba
&&N(r^\xi\varphi(\Lambda)1_{\{M\geq 1\}})\\
&&\quad =
\Big(r\psi^{-1}(\lambda)\psi'(\psi^{-1}(\lambda))+{\psi((1-r)\psi^{-1}(\lambda))}\Big)
\int_0^\infty db\,\varphi(b)\,e^{-b\psi'(\psi^{-1}(\lambda))}.
\ea
This formula entails that we have the following properties under $N(\cdot \mid M\geq 1)$:
The variables $\Lambda$ and $\xi$ are independent, $\Lambda$ is exponentially
distributed with parameter $\psi'(\psi^{-1}(\lambda))$, and the generating function
of $\xi$ is as stated in Theorem \ref{tree-Poisson}. To complete the proof, it remains
to verify the ``recursivity property'' of the tree $\theta(H,\tau_1,\ldots,\tau_M)$,
that is to verify that under $N(\cdot\mid M\geq 2)$, the shifted trees corresponding
to each individual in the first generation are independent and distributed as the 
whole tree under $N(\cdot\mid M\geq 1)$. This is a consequence of the
following claim. 

\smallskip\noindent
{\bf Claim.} {\it Let $(\alpha_j,\beta_j)$, $j=1,\ldots,\xi$ be the excursion
intervals of $H$ above level $\Lambda$ that contain at least one mark,
ranked in chronological order, and for 
every $j=1,\ldots,\xi$ let $h_j(s)=H_{(\alpha_j+s)\wedge \beta_j}-\Lambda$
be the corresponding excursion. Then,
conditionally on the pair $(\Lambda,\xi)$, the excursions 
$h_1,\ldots,h_\xi$ are independent and distributed according
to the law of $(H_{s\wedge \sigma},s\geq 0)$ under $N(\cdot\mid M\geq 1)$.}

\smallskip

To verify this property, we first argue under $\PP_\mu$ as previously. Precisely, we 
consider the excursions $\omega_j$ for all $j\in J$ such that $H_{a_j}=L$
and $(a_j,b_j)\cap {\cal P}_\lambda\not =\emptyset$. We denote by $\wt \omega_1,
\ldots,\wt\omega_\zeta$ these excursions, ranked in chronological order.
Then property (P) and familiar properties of Poisson measures 
give the following fact. For every $k\geq 1$, under the measure
$\PP_\mu(\cdot \mid \zeta=k)$, the excursions $\wt \omega_1,
\ldots,\wt\omega_k$ are independent, distributed according to $N(\cdot\mid M\geq 1)$,
and these excursions are also independent of the measure
$$\sum_{H_{a_j}>L} \delta_{(I_{a_j},\omega_j)}\,.$$
Let $\sigma_L:=\inf\{s\geq 0:H_s=L\}$. Excursion theory for $X-I$ allows us to
reconstruct the process $(X_{s\wedge \sigma_L},s\geq 0)$ as a measurable
function of the point measure in the last display. Hence we can also assert
that, under $\PP_\mu(\cdot \mid \zeta=k)$, $\wt \omega_1,
\ldots,\wt\omega_k$ are independent of $(X_{s\wedge \sigma_L},s\geq 0)$. 
In particular, they are
independent of
$L=H(k_{<\mu,1>-I_{\sigma_L}}\mu)$.

We now apply these properties to the shifted process $X_{\tau+\cdot}$
under $N(\cdot\mid M\geq 1)$. We slightly abuse notation
and keep denoting by $\wt \omega_1,
\ldots,\wt\omega_\zeta$ the excursions of $X_{\tau+\cdot}-I_{\tau+\cdot}$
that contain a mark (so that $\zeta=\xi-1$ on the event $M\geq 2$).
By construction, for every $j\in \{2,\ldots,\xi\}$, the function $h_j$
is the height process of $\wt \omega_{j-1}$. Hence it follows from the 
previous properties under $\PP_\mu$ that under $N(\cdot\mid \xi=k)$
(for a fixed $k\geq 2$),
the processes $h_2,\ldots,h_k$ are independent, have the
distribution required in the claim, and are also
independent of the pair $(h_1,\Lambda)$. Hence, for any test functions
$F_1,\ldots,F_k,G$, we get
\ba
&&N(F_1(h_1)\ldots F_k(h_k)G(\Lambda)\mid \xi=k)\\
&&\qquad=N(F_2(H)\mid M\geq 1)\cdots N(F_k(H)\mid M\geq 1)\,N(F_1(h_1)G(\Lambda)\mid \xi=k).
\ea
Now from Corollary \ref{reversal}, we know that the 
time-reversed process 
$(H_{(\sigma-s)_+},s\geq 0)$ has the same distribution under
$N$ as the process $(H_s,s\geq 0)$. Furthermore, this time-reversal operation 
will leave $\Lambda$ and $\xi$ invariant and transform the excursion $h_1$
into the time-reversal of $h_\xi$, denoted by $\wh h_\xi$ (provided we
do simultaneously the similar transformation on the underlying
Poissonnian marks). It follows that
\ba 
N(F_1(h_1)G(\Lambda)\mid \xi=k)
&=&N(F_1(\wh h_k)G(\Lambda)\mid\xi=k)\\
&=&N(F_1(\wh h_k)\mid \xi=k)N(G(\Lambda)\mid\xi=k)\\
&=&N(F_1(H)\mid M\geq 1)N(G(\Lambda)\mid M\geq 1).
\ea
By substituting this equality in the previous displayed formula, we obtain the
claim. This completes the proof of Theorem \ref{tree-Poisson}.

\section{Marginals of stable trees}

We first reformulate Theorem \ref{tree-Poisson} in a way more suitable for
our applications. Recall that $\TT_p$ is the set of all (rooted ordered)
trees with $p$ leaves. If $\t\in \TT_p$ we denote by ${\cal L}_\t$ the set
of all leaves of $\t$, and set $\n_\t=\t\backslash{\cal L}_\t$. Recall the notation
$k_v=k_v(\t)$ for the number of children of an element $v$ of $\t$. 
We write $\TT^*_p$ for the subset of $\TT_p$ composed  of all
trees $\t$ such that $k_v(\t)\geq 2$ for every $v\in\n_\t$. By construction, 
the skeleton of the marked trees $\theta(e,t_1,\ldots,t_p)$ always
belongs to $\TT^*_p$.

\begin{theorem}
\label{Poisson-tree}
Let $p\geq 1$. Then, for any nonnegative measurable function $\Phi$
on $\T_p$, and every $\lambda>0$,
\ba
&&N\Big(e^{-\lambda \sigma}\int_{\{t_1<\cdots<t_p<\sigma\}}dt_1\ldots
dt_p\,\Phi(\theta(H,t_1,\ldots,t_p))\Big)\\
&&\quad=\sum_{\t\in \TT^*_p} \Big(\prod_{v\in\n_\t} {|\psi^{(k_v)}(\psi^{-1}(\lambda))|\over
k_v!}\Big)\\
&&\qquad\qquad
\int \prod_{v\in \t} dh_v\;\exp\Big(-\psi'(\psi^{-1}(\lambda))\sum_{v\in \t}h_v\Big)
\Phi(\t,(h_v)_{v\in \t}).
\ea
\end{theorem}

\proof By elementary properties of the standard Poisson process, the left side
is equal to
$$\lambda^{-p}N(\Phi(\theta(H,\tau_1,\ldots,\tau_M))1_{\{M=p\}})$$
with the notation of the previous section. This quantity can be evaluated
thanks to Theorem \ref{tree-Poisson}: From the generating function
of the offspring distribution, we get
\ba
&&P[\xi=0]={\lambda\over \psi^{-1}(\lambda)\psi'(\psi^{-1}(\lambda))}\\
&&P[\xi=1]=0\\
&&P[\xi=k]={1\over k!}\,{\psi^{-1}(\lambda)^{k-1}|\psi^{(k)}(\psi^{-1}(\lambda))|
\over \psi'(\psi^{-1}(\lambda))},\qquad\hbox{for every }k\geq 2.
\ea
Hence the probability under $N(\cdot\mid M\geq 1)$ that the skeleton of the tree
$\theta(H,\tau_1,\ldots,\tau_M)$ is equal to a given tree $\t\in\TT^*_p$
is 
\ba
&&\Big(\prod_{v\in\n_\t} {1\over k_v!}\,{\psi^{-1}(\lambda)^{k_v-1}|\psi^{(k_v)}(\psi^{-1}(\lambda))|
\over \psi'(\psi^{-1}(\lambda))}\Big)\,\Big(
{\lambda\over \psi^{-1}(\lambda)\psi'(\psi^{-1}(\lambda))}\Big)^p\\
&&\quad={1\over \psi^{-1}(\lambda)}\,{1\over\psi'(\psi^{-1}(\lambda))^{|\t|}}
\,\lambda^p\prod_{v\in\n_\t} \Big({1\over k_v!}\,|\psi^{(k_v)}(\psi^{-1}(\lambda))|\Big).
\ea
Recalling that $N(M\geq 1)=\psi^{-1}(\lambda)$, and 
using the fact that the lifetimes $h_v$, $v\in \t$
are independently distributed according to the exponential distribution with 
parameter $\psi'(\psi^{-1}(\lambda))$, we easily arrive at the formula 
of the theorem. \cq

\medskip
By letting $\lambda\rightarrow  0$ in the preceding theorem, we get the following corollary, which is closely
related to Proposition 3.2 of \cite{LGLJ2}.

\begin{corollary}
\label{uniform-tree}
Suppose that $\int \pi(dr)\,r^p<\infty$. Then, for any nonnegative measurable function $\Phi$
on $\t_p$,
\ba
&&N\Big(\int_{\{t_1<\cdots<t_p<\sigma\}}dt_1\ldots
t_p\,\Phi(\theta(H,t_1,\ldots,t_p))\Big)\\
&&\quad=\sum_{\t\in \TT^*_p} \Big(\prod_{v\in\n_\t} {\beta_{k_v}}\Big)
\int \prod_{v\in \t} dh_v\;\exp\Big(-\alpha\sum_{v\in \t}h_v\Big)
\Phi(\t,(h_v)_{v\in \t})\,,
\ea
where, for every $k=1,\ldots,p$,
$$\beta_k={|\psi^{(k)}(0)|\over k!}=\beta 1_{\{k=2\}}+{1\over k!}\int r^k\,\pi(dr).$$
\end{corollary}

\rem The formula of the corollary still holds without the assumption $\int \pi(dr)\,r^p<\infty$
but it has to be interpreted properly since some of the numbers $\beta_k$
may be infinite.

\medskip
From now on, we concentrate on the stable case $\psi(u)=u^\gamma$
for $1<\gamma<2$. Then the L\'evy process $X$ satisfies the scaling
property 
$$(X_{\lambda t},t\geq 0)\build{=}_{}^{\rm(d)}(\lambda^{1/\gamma}X_t,t\geq 0)$$
under $P$. Thanks to this property, it is possible to choose a regular
version of the conditional probabilities $N_{(u)}:=N(\cdot\mid \sigma=u)$
in such a way that for every $u>0$ and $\lambda>0$, the law of
$(\lambda^{-1/\gamma}X_{\lambda t},t\geq 0)$ under $N_{(\lambda u)}$ is $N_{(u)}$.
Standard arguments then show that the height process $(H_s,s\geq 0)$ is
well defined as a continuous process under the probability measures $N_{(u)}$.
Furthermore, it follows from the approximations of $H_t$ (see Lemma \ref{localstrip})
that the law of $(H_{\lambda s},s\geq 0)$ under $N_{(\lambda u)}$ is equal
to the law of $(\lambda^{1-{1\over \gamma}}H_s,s\geq 0)$ under $N_{(u)}$.

The probability measure $N_{(1)}$ is called the law of the normalized excursion. 

\begin{theorem}
\label{stable-marginal}
Suppose that $\psi(u)=u^\gamma$ for some $\gamma\in(1,2)$.
Then the law of the tree $\theta(H,t_1,\ldots,t_p)$ under the probability measure
$$p!\,1_{\{0<t_1<t_2<\ldots<t_p<1\}}dt_1\ldots dt_p\,N_{(1)}(d\omega)$$
is characterized by the following properties:

\noindent{\rm (i)} The probability of a given skeleton $\t\in\TT^*_p$ is
$${p!\over\displaystyle  \prod_{v\in\n_\t} k_v!}\ 
{\displaystyle\prod_{v\in\n_\t} |(\gamma-1)(\gamma-2)\ldots(\gamma-k_v+1)|\over
(\gamma-1)(2\gamma-1)\ldots ((p-1)\gamma-1)}.$$

\noindent{\rm(ii)} If $p\geq 2$, then conditionally on the skeleton $\t$, the lifetimes 
$(h_v)_{v\in \t}$ have a density with respect to Lebesgue measure on $\R_+^{\t}$
given by 
$${\Gamma(p-{1\over \gamma})\over \Gamma(\delta_\t)}
\,\gamma^{|\t|}\int_0^1 du\,u^{\delta_\t-1}\,q(\gamma\sum_{v\in \t}h_v,1-u)$$
where $\delta_\t=p-(1-{1\over \gamma})|\t|-{1\over \gamma}>0$, and
$q(s,u)$ is the continuous density at time $s$ of the stable
subordinator with exponent $1-{1\over\gamma}$, which is 
characterized by
$$\int_0^\infty du\,e^{-\lambda u}\,q(s,u)=\exp(-s\,\lambda^{1-{1\over \gamma}}).$$
If $p=1$, then $\t=\{\emptyset\}$ and the law of $h_\emptyset$ has density
$$\gamma\,\Gamma(1-{1\over\gamma})\,q(\gamma h,1)$$
with respect to Lebesgue measure on $\R_+$.
\end{theorem}

\proof For every $u>0$, let $\Theta_{(u)}$ be the law of the tree
$\theta(H,s_1,\ldots,s_p)$ under the probability measure
$$p!u^{-p}\,1_{\{s_1<s_2<\cdots<s_p<u\}}\,ds_1\ldots ds_p\,N_{(u)}(d\omega).$$
By the scaling properties of the height process
(see the remarks before Theorem \ref{stable-marginal}), we have, for every $u>0$ and every
$\t\in\TT^*_p$, 
$$\Theta_{(u)}(\{\t\}\times \R_+^\t)=\Theta_{(1)}(\{\t\}\times \R_+^\t).$$
Hence, by conditioning with respect to $\sigma$ in Theorem \ref{Poisson-tree}, we get
\begin{equation} 
\label{tech-stable}
{1\over p!}\,N(\sigma^p e^{-\lambda \sigma})\,\Theta_{(1)}(\{\t\}\times \R_+^\t)=
\psi'(\psi^{-1}(\lambda))^{-|\t|}\,\prod_{v\in\n_\t} {|\psi^{(k_v)}(\psi^{-1}(\lambda))|\over
k_v!}.
\end{equation}
From this, we can compute $\Theta_{(1)}(\{\t\}\times \R_+^\t)$ by observing that, for every $k\geq 1$,
$$\psi^{(k)}(\psi^{-1}(\lambda))=\gamma(\gamma-1)\cdots(\gamma-k+1)\,\lambda^{1-{k\over \gamma}},$$
and
$$N(\sigma^p e^{-\lambda \sigma})=|{d^p\over d\lambda^p}N(1-e^{-\lambda \sigma})|
=|{d^p\over d\lambda^p}\psi^{-1}(\lambda)|=|{1\over \gamma}({1\over \gamma}-1)\cdots
({1\over \gamma}-p+1)|\lambda^{{1\over \gamma}-p}.$$
If we substitute these expressions in (\ref{tech-stable}), the terms in $\lambda$ cancel
and we get part (i) of the theorem.

To prove (ii), fix $\t\in \TT^*_p$, and let $D$ be a bounded Borel subset of $\R_+^\t$. Write
$p_\sigma(du)$ for the law of $\sigma$ under $N$. Then by applying Theorem \ref{Poisson-tree} 
with $\Phi=1_{\{\t\}\times D}$ and $\Phi=1_{\{\t\}\times \R_+^\t}$, we get
\begin{eqnarray}
\label{tech-stable2}
&&\int p_\sigma(du)\,e^{-\lambda u}\,{u^p}\,\Theta_{(u)}(\{h_v\}_{v\in \t}\in D \mid \t)\\
&&=\Big(\int p_\sigma(du)\,e^{-\lambda u}\,{u^p}\Big)
\int_D \prod_{v\in \t} dh_v\,
\psi'(\psi^{-1}(\lambda))^{|\t|}\,\exp\Big(-\psi'(\psi^{-1}(\lambda))\sum_{v\in \t}h_v\Big).\nonumber
\end{eqnarray}
By scaling (or inverting $N(1-e^{-\lambda\sigma})=\lambda^{1/\gamma}$), we have 
$p_\sigma(du)=c\,u^{-1-{1\over \gamma}}\,du$. It follows that
\ba
&&\int_0^\infty du\,e^{-\lambda u}\,u^{p-1-{1\over \gamma}}\,\Theta_{(u)}(\{h_v\}_{v\in \t}\in D
\mid \t)\\
&&\quad={\gamma^{|\t|}\Gamma(p-{1\over \gamma})\over \lambda^{\delta_\t}}
\int_D \prod_{v\in \t} dh_v\,\exp(-\gamma\lambda^{1-{1\over \gamma}}\sum_{v\in \t}h_v),
\ea
where $\delta_\t=p-(1-{1\over \gamma})|\t|-{1\over \gamma}$ as in the theorem.
Suppose first that $p\geq 2$.
To invert the Laplace transform, observe that the right side can be written as
\ba
&&{\gamma^{|\t|}\Gamma(p-{1\over \gamma})\over \Gamma(\delta_\t)}
\int_0^\infty du\,e^{-\lambda u}\,u^{\delta_\t-1}\times
\int_0^\infty du'\,e^{-\lambda u'}\Big(\int_D \prod_{v\in \t} dh_v
\,q(\gamma\sum_{v\in \t}h_v,u')\Big)\\
&&={\gamma^{|\t|}\Gamma(p-{1\over \gamma})\over \Gamma(\delta_\t)}
\int_0^\infty du\,e^{-\lambda u}\Big(\int _D \prod_{v\in \t} dh_v
\int_0^u dr\,r^{\delta_\t-1}\,q(\gamma\sum_{v\in \t}h_v,u-r)\Big).
\ea
The first formula of (ii) now follows. In the case $p=1$, we get 
$$\int_0^\infty du\,e^{-\lambda u}\,u^{-{1\over \gamma}}\,\Theta_{(u)}(h_\emptyset\in
D)=\gamma\Gamma(1-{1\over \gamma})
\int_D dh\,\exp(-\gamma\lambda^{1-{1\over \gamma}}h),$$
and the stated result follows by inverting the Laplace transform.
\cq

\medskip
\rems (a) The previous proof also readily gives the analogue of Theorem 
\ref{stable-marginal} in the case 
$\psi(u)=u^2$, which corresponds to the finite-dimensional marginals of 
Aldous' continuum random tree (see Aldous \cite{Al2}, or
Chapter 3 of \cite{LG99}). In that case, the discrete
skeleton of $\theta(H,t_1,\ldots,t_p)$ is with probability
one a binary tree, meaning that $k_v=2$ for every $v\in \t$. 
The law of $\t(H,t_1,\ldots,t_p)$ is the uniform probability
measure on the set of all binary trees in $\TT^*_p$, 
so that the probability of each possible skeleton is
$${p!\over 2^{p-1}\ (1\times 3\times \cdots\times (2p-3))}.$$
This formula can be deduced informally by letting $\gamma$ tend to $2$
in Theorem \ref{stable-marginal} (i). 

To obtain the analogue of (ii), note that there is an explicit formula
for $q(s,u)$ when $\psi(u)=u^2$:
$$q(s,u)={s\over 2\sqrt{\pi}u^{3/2}}\,e^{-s^2/(4u)}.$$
Observe that when the skeleton is binary, we have always $|\t|=2p-1$. 
It follows that the powers of $\lambda$ cancel in the right side of (\ref{tech-stable2}), and
after straightforward calculations, we obtain that
the density of $(h_v)_{v\in\t}$ on $\R_+^{2p-1}$ is 
$$2^{2p-1}\Gamma(p-{1\over 2})\,q(2\sum h_v,1)
=2^{p}\,(1\times 3\times \cdots\times (2p-3))\ (\sum h_v)\,\exp(-(\sum h_v)^2).$$
Compare with Aldous \cite{Al2} or
Chapter 3 of \cite{LG99}, but note that constants are different because 
$\psi(u)=u^2$ corresponds to a Brownian motion with variance $2t$ (also
the CRT is coded by twice the normalized Brownian excursion in \cite{Al2}).

\smallskip
(b) We could get rid of the factor
$${p!\over\displaystyle  \prod_{v\in\n_\t} k_v!}$$
in Theorem 
\ref{stable-marginal} by considering rooted (unordered) trees with $p$ labelled leaves rather than
rooted ordered trees : See the discussion at the end of Chapter 3 of \cite{LG99}.

 \chapter{The L\'evy snake}

\section{The construction of the L\'evy snake}

Our goal is now to combine the branching structure studied 
in the previous chapters with a spatial displacement
prescribed by a Markov process $\xi$. Throughout this chapter, we assume that 
$H$ has continuous paths (the condition $\int^\infty du/\psi(u)<\infty$
holds) although many of the results can presumably be extended to the general case.

\subsection{Construction of the snake with a 
fixed lifetime process}

We consider a Markov process $\xi$ with c\` adl\` ag
paths and values in a Polish space $E$,
whose topology is defined by a metric $\delta$. For
simplicity, we will assume that $\xi$ is defined on the
canonical space $\D(\R_+,E)$ of c\` adl\` ag functions from
$\R_+$ into $E$. For every $x\in E$, we denote
by $\Pi_x$ the distribution of $\xi$ started at $x$.
It is implicitly assumed in our definition 
of a Markov process that the
mapping $x\rightarrow  \Pi_x$ is measurable.
We also
assume that $\xi$ is continuous in probability
under $\Pi_x$ (equivalently, $\xi$ has no fixed discontinuities,
$\Pi_{x}[\xi_s\not =\xi_{s-}]=0$ for every $s> 0$). On the 
other hand, we do not assume that $\xi$ is strong Markov.

For $x\in E$, we denote by $\W_x$ the space of all 
$E$-valued killed
paths started at $x$. An element of $\W_x$ is a c\` adl\` ag
mapping $\w:[0,\zeta)\la E$
such that $\w(0)=x$. Here $\zeta\in(0,\infty)$ is called the
lifetime of the path. When there is a risk of confusion
we write $\zeta=\zeta_{\w}$. Note that we do not require the
existence of the left limit $\w(\zeta-)$. By
convention, the point
$x$ is also considered as a killed path with lifetime $0$. We
set
$\W=\cup_{x\in E}\W_x$ and equip $\W$ with the distance
$$d(\w,\w')=\delta(\w(0),\w'(0))+|\zeta-\zeta'|
+\int_0^{\zeta\wedge \zeta'}dt\,( d_t(\w_{\leq t},\w'_{\leq
t})\wedge 1),$$
where $d_t$ is the
Skorokhod metric on the space $\D([0,t],E)$, and 
$\w_{\leq t}$ denotes the restriction of $\w$
to the interval $[0,t]$. It is then elementary to check
that the space
$(\W,d)$ is a Polish space. The space $(E,\delta)$ is embedded isometrically
in $\W$ thanks to the previous convention.

Let $x\in E$ and $\w\in\W_x$. If
$a\in [0,\zeta_\w)$ and $b\in[a,\infty)$, we can
define a probability measure $R_{a,b}(\w,d\w')$
on $\W_x$ by requiring that:
\begin{description}
\item{(i)} $R_{a,b}(\w,d\w')$ a.s., $\w'(t)=\w(t)$, $\forall
t\in[0,a)$;
\item{(ii)} $R_{a,b}(\w,d\w')$ a.s., $\zeta_{\w'}=b$;
\item{(iii)} the law of $(\w'(a+t),0\leq t<b-a)$ 
under $R_{a,b}(\w,d\w')$ is the law
of $(\xi_t,0\leq t<b-a)$ under $\Pi_{\w(a-)}$.
\end{description}
In (iii), $\w(0-)=x$ by convention.
In particular, $R_{0,b}(\w,d\w')$ is the law of $(\xi_t,
0\leq t<b)$ under $\Pi_x$, and $R_{0,0}(\w,d\w')=\delta_x(d\w')$.

When $\w(\zeta_\w-)$ exists, we may and will extend the previous
definition to the case 
$a=\zeta_{\w}$.  

\smallskip
We denote by $(W_s,s\geq 0)$ the canonical process on
the product space $(\W)^{\R_+}$. We will abuse notation and
also write $(W_s,s\geq 0)$ for the canonical process on
the set $C(\R_+,\W)$ of all continuous mappings from
$\R_+$ into $\W$. 
Let us fix $x\in E$ and $\w_0\in\W_x$,
and let $h\in C(\R_+,\R_+)$ be such that $h(0)=\zeta_{\w_0}$.
For $0\leq s\leq s'$, we set
$$m_h(s,s')=\inf_{s\leq r\leq s'}h(r).$$
We assume that either $\w_0(\zeta_{\w_0}-)$ exists or
$m_h(0,r)<h(0)$ for every $r>0$. Then, the Kolmogorov
extension theorem can be used to construct the (unique)
probability measure $Q^h_{\w_0}$ on $(\W_x)^{\R_+}$ such
that, for
$0=s_0<s_1<\cdots<s_n$,
\begin{eqnarray*}
&&Q^h_{\w_0}[W_{s_0}\in A_0,\ldots,W_{s_n}\in A_n]\\
&&=1_{A_0}(\w_0)\int_{A_1\times \cdots\times A_n}
R_{m_h(s_0,s_1),h(s_1)}(\w_0,d\w_1)\ldots
R_{m_h(s_{n-1},s_n),h(s_n)}(\w_{n-1},d\w_n).
\end{eqnarray*}
Notice that our assumption on the pair $(\w_0,h)$ is needed
already for $n=1$ to make sense of the measure
$R_{m_h(s_0,s_1),h(s_1)}(\w_0,d\w_1)$.

\smallskip

From the previous definition, it is clear that, for
every $s<s'$, $Q^h_{\w_0}$ a.s.,
$$W_{s'}(t)=W_s(t),\qquad \forall t< m_h(s,s'),$$
and furthermore
$\zeta_{W_s}=h(s)$, $\zeta_{W_{s'}}=h(s')$. Hence,
$$d(W_s,W_{s'})\leq |h(s)-h(s')|+|(h(s)\wedge h(s'))
-m_h(s,s')|=(h(s)\vee h(s'))-m_h(s,s').$$
From this bound, it follows that the mapping $s\la W_s$
is $Q^h_{\w_0}$ a.s. uniformly continuous on the
bounded subsets of $[0,\infty)\cap \Q$. Hence this mapping 
has $Q^h_{\w_0}$ a.s. a continuous extension to the positive
real line. We abuse notation and still denote by $Q^h_{\w_0}$
the induced probability measure on $C(\R_+,\W_x)$. By
an obvious continuity argument, we have
$\zeta_{W_s}=h(s)$, for every $s\geq 0$, $Q^h_{\w_0}$ a.s.,
and
$$W_{s'}(t)=W_s(t),\qquad \forall t< m_h(s,s'),\quad
\forall s<s',\quad Q^h_{\w_0}\ {\rm a.s.}$$
We will refer to this last property as the {\it snake
property}.
The process $(W_s,s\geq 0)$ is under $Q^h_{\w_0}$
a time-inhomogeneous continuous Markov process.

\subsection{The definition of the L\'evy snake} 

Following the remarks of the end of Chapter 1, we now consider the
exploration process $\rho$ as a Markov process with values in 
the set
$$M^0_f=\{\mu\in M_f(\R_+):H(\mu)<\infty \hbox{ and } \supp\mu=[0,H(\mu)]\}\cup\{0\}.$$
We denote by $\PP_\mu$ the law of $(\rho_s,s\geq 0)$
started at $\mu$. We will write indifferently $H(\rho_s)$ or $H_s$. 

We then define $\Theta$ as the set of all
pairs $(\mu,\w)\in M^0_f\times \W$ such that $\zeta_\w=H(\mu)$,
and at least one of the following two properties hold:
\begin{description}
\item{(i)} $\mu(\{H(\mu)\})=0$;
\item{(ii)} $\w(\zeta_{\w}-)$ exists.
\end{description}
We equip
$\Theta$ with the product distance on $M^0_f\times \W$. For every $y\in E$,
we also set $\Theta_y=\{(\mu,\w)\in\Theta:\w(0)=y\}$.

From now on until the end of this section, we fix a point $x\in E$.

\smallskip
Notice that when $H(\mu)>0$ and $\mu(\{H(\mu)\})=0$, 
we have $\inf_{[0,s]}H(\rho_r)<H(\mu)$
for every $s>0$, $\PP_\mu$ a.s. This property
easily follows from (\ref{rhoinitial}) and the fact that $0$ is regular for 
$(-\infty,0)$ for the underlying L\'evy process.

Using the last observation and the previous subsection, we can 
for every $(\mu,\w)\in\Theta_x$ define a probability
measure $\bP_{\mu,\w}$ on $\D(\R_+,M_f(\R_+)\times \W)$ by the formula
$$\bP_{\mu,\w}(d\rho\,dW)=\PP_\mu(d\rho)\,Q^{H(\rho)}_\w(dW),$$
where in the right side $H(\rho)$ obviously stands for
the function $(H(\rho_s),s\geq 0)$, which is continuous
$\PP_\mu$ a.s.

We will write $\bP_x$ instead of $\bP_{0,x}$ when $\mu=0$.

\begin{proposition}
\label{snakedef}
The process $(\rho_s,W_s)$ is under $\bP_{\mu,\w}$
a c\` adl\` ag Markov process in $\Theta_x$.
\end{proposition}

\proof We first verify that $\bP_{\mu,\w}$ a.s. the process
$(\rho_s,W_s)$ does not visit $\Theta_x^c$. We must check that
$W_s(H_s-)$ exists whenever $\rho_s(\{H_s\})>0$. Suppose
thus that $\rho_s(\{H_s\})>0$. Then, 
we have also
$\rho_{s'}(\{H_s\})>0$, and so $H_{s'}\geq H_s$, for 
all $s'>s$ sufficiently close to $s$. In particular, we
can find a rational $s_1>s$ such that $H_{s_1}\geq H_s$
and $\inf_{[s,s_1]}H_r=H_s$, which by the snake 
property implies that $W_{s_1}(t)=W_s(t)$ for
every $t\in [0,H_s)$. However, from the construction
of the measures $Q^h_\w$, it is clear that a.s. for every 
rational $r>0$, the killed path $W_{r}$ must have a
left limit at every 
$t\in(0,H_{r}]$. We conclude that $W_{s}(H_s-)=W_{s_1}(H_s-)$
exists.

The c\` adl\` ag property of paths is obvious by construction.
To obtain the Markov property, we consider 
nonnegative functions $f_1,\ldots,f_n$
on $M^0_f$ and $g_1,\ldots,g_n$ on $\W_x$. Then,
if $0<s_1<\cdots<s_n$,
\begin{eqnarray*}
&&\E_{\mu,\w}[f_1(\rho_{s_1})g_1(W_{s_1})\ldots
f_n(\rho_{s_n})g_n(W_{s_n})]\\
&&=\EE_\mu\Big[f_1(\rho_{s_1})\ldots f_n(\rho_{s_n})
Q^{H(\rho)}_\w[g_1(W_{s_1})\ldots g_n(W_{s_n})]\Big]\\
&&=\EE_\mu\Big[f_1(\rho_{s_1})\ldots f_n(\rho_{s_n})\int R_{m_{H(\rho)}(0,s_1),H(\rho_{s_1})}(\w,d\w_1)\\
&&\quad\ldots
R_{m_{H(\rho)}(s_{n-1},s_n), H(\rho_{s_n})}(\w_{n-1},d\w_n)\,g_1(\w_1)\ldots
g_n(\w_n)\Big]\\
&&=\E_{\mu,\w}\Big[f_1(\rho_{s_1})g_1(W_{s_1})\ldots
f_{n-1}(\rho_{s_{n-1}})g_{n-1}(W_{s_{n-1}})\\
&&\quad \EE_{\rho_{s_{n-1}}}\Big[f_n(\rho_{s_n-s_{n-1}})
\int R_{m_{H(\rho)}(0,s_n-s_{n-1}),H(\rho_{s_n-s_{n-1}})}(W_{s_{n-1}},d\w)
g_n(\w)\Big]\Big],
\end{eqnarray*}
where in the last equality we used the Markov property for
$(\rho_s,s\geq 0)$ at time $s_{n-1}$. We get the desired
result with a transition kernel given by
\begin{eqnarray*}
{\QQ}_rG(\mu,\w)&=&\int \PP_\mu(d\rho)
\int R_{m_{H(\rho)}(0,r),H(\rho_r)}(\w,d\w')\,G(\rho_r,\w')\\
&=&\int \PP_\mu(d\rho)\int Q^{H(\rho)}_\w(dW)\,G(\rho_r,W_r).
\end{eqnarray*}
\par\cq

In what follows we will often use the convenient notation
$\ov W_s=(\rho_s,
W_s)$. 
By our construction, the
conditional distribution under $\bP_{\mu,\w}$ of $(W_s,s\geq
0)$  knowing $(\rho_s,s\geq 0)$ is $Q^{H(\rho)}_\w$. In particular,
if we write $\zeta_s=\zeta_{W_s}$ for the lifetime of $W_s$,
we have
$$\zeta_s=H(\rho_s)=H_s\hbox{ for every }s\geq 0,\ \bP_{\mu,\w}\hbox{ a.s.}$$

\subsection{The strong Markov property}

We denote by $(\f_s)_{s\geq 0}$ the canonical filtration
on $\D(\R_+,M_f(\R_+)\times \W)$. 

\begin{theorem}
\label{strongMarkovsnake}
The process $(\ov W_s,s\geq 0;\bP_{\mu,\w},
(\mu,\w)\in\Theta_x)$ is strong Markov
with respect to the filtration $(\f_{s+})$.
\end{theorem}

\proof Let $(\mu,\w)\in\Theta_x$. It is enough to prove that,
if $T$ is a bounded stopping time of the filtration
$(\f_{s+})$, then, for any bounded $\f_{T+}$-measurable
functional $F$, for any bounded Lipschitz continuous function
$f$ on $\Theta_x$, and for every $t> 0$,
$$\E_{\mu,\w}[F\,f(\ov W_{T+t})]
=\E_{\mu,\w}[F\,\E_{\ov W_T}[f(\ov W_t)]].$$
First observe that
\begin{eqnarray*}
\E_{\mu,\w}[F\,f(\ov W_{T+t})]&=&
\lim_{n\rightarrow \infty}\sum_{k=1}^\infty
\E_{\mu,\w}[F\,1_{\{{k-1\over n}
\leq T<{k\over n}\}}\,f(\ov W_{{k\over n}+t})]\\
&=&\lim_{n\rightarrow \infty}\sum_{k=1}^\infty
\E_{\mu,\w}[F\,1_{\{{k-1\over n}
\leq T<{k\over n}\}}\,\QQ_tf(\ov W_{{k\over n}})].
\end{eqnarray*}
In the first equality, we used the right continuity of paths,
and in the second one the ordinary Markov property.
We see that the desired result follows from the next
lemma.

\begin{lemma}
\label{SMtech}
Let $t>0$, let $T$ be a bounded stopping time of the
filtration $(\f_{s+})$
and let $f$ be a bounded Lipschitz continuous function
on $\Theta_x$. Then the mapping $s\la {\QQ}_tf(\ov W_s)$
is $\bP_{\mu,\w}$ a.s. right-continuous at $s=T$.
\end{lemma}

\noindent{\bf Proof of Lemma \ref{SMtech}.} We use the notation
$Y_s=\big<\rho_s,1\big>$. Recall that $Y$ is distributed under
$\bP_{\mu,\w}$ as the reflected L\'evy process $X-I$
started at $\big<\mu,1\big>$. 
Let $\varepsilon>0$. By the right-continuity 
of the paths of $Y$,
if $s>T$ is sufficiently close to $T$, we have
$$\varepsilon_1(s)=Y_T-\inf_{u\in[T,s]}Y_u<\varepsilon,\qquad
\varepsilon_2(s)=Y_s-\inf_{u\in[T,s]}Y_u<\varepsilon.$$
On the other hand, we know from (\ref{strongMarkov}) 
that $\rho_{T+s}=[k_{\varepsilon_1}\rho_T,\rho^{(T)}_s]$,
and it follows that
$k_{\varepsilon_1}\rho_T= k_{\varepsilon_2}\rho_s$. 
Furthermore, $\inf_{[T,s]}H(\rho_u)=H(
k_{\varepsilon_1}\rho_T)=H(k_{\varepsilon_2}\rho_s)$, and
by the snake property,
$$W_s(u)=W_T(u),\qquad \forall
u\in[0,H(k_{\varepsilon_1}\rho_T)).$$

Let us fix $\ov\w=(\mu,\w)\in\Theta_x$, and set
\begin{eqnarray*}
\v_\varepsilon(\ov \w)&=&
\big\{\ov\w'=(\mu',\w')\in\Theta_x;\,\exists \varepsilon_1,
\varepsilon_2\in[0,\varepsilon),\ k_{\varepsilon_1}\mu
=k_{\varepsilon_2}\mu',\\
&&\qquad{\rm and}\ \w'(u)=\w(u),\ \forall u\in[0,
H(k_{\varepsilon_1}\mu))\big\}.
\end{eqnarray*}
In view of the preceding observations, the proof of Lemma
\ref{SMtech} reduces to checking that
\begin{equation}
\label{SMtech1}
\lim_{\varepsilon\rightarrow  0}
\Big(\sup_{\ov \w'\in\v_\varepsilon(\ov \w)}
\big|\QQ_tf(\ov\w')-\QQ_tf(\ov\w)\big|\Big)=0.
\end{equation}

We will use a coupling argument to obtain (\ref{SMtech1}).
More precisely,
if $\ov\w'\in\v_\varepsilon(\ov\w)$, we will introduce two
(random) variables
$\ov\w_{(1)}$ and $\ov\w_{(2)}$ such that $\ov\w_{(1)}$,
resp. $\ov\w_{(2)}$, is distributed according to
$\QQ_t(\ov\w,\cdot)$, resp. $\QQ_t(\ov\w',\cdot)$, and 
$\ov\w_{(1)}$ and $\ov\w_{(2)}$ are close to each other.
Let us fix $\ov\w'\in\v_\varepsilon(\ov\w)$ and let
$\varepsilon_1,\varepsilon_2\in[0,\varepsilon)$ 
be associated with $\ov\w'$ as in the
definition of $\v_\varepsilon(\ov\w)$. For definiteness
we assume that $\varepsilon_1\leq \varepsilon_2$
(the other case is treated in a symmetric way).
Let $X^{(1)}$ be a copy of the L\'evy process $X$
started at $0$
and let $I^{(1)}$ and $\rho^{(1)}$ be the analogues
of $I$ and $\rho$ for $X^{(1)}$. We can
then define
$\ov\w_{(1)}=(\mu_{(1)},
\w_{(1)})$ by 
\begin{eqnarray*}
&&\mu_{(1)}=[k_{-I^{(1)}_t}\mu,\rho^{(1)}_t]\\
&&\w_{(1)}(r)=\left\{
\begin{array}{ll}
\w(r)&{\rm if}\ r< H(k_{-I^{(1)}_t}\mu),\\
\xi^{(1)}({r-H(k_{-I^{(1)}_t}\mu)})\quad&{\rm if}\
H(k_{-I^{(1)}_t}\mu)\leq r<H(\mu_{(1)}),
\end{array}\right.
\end{eqnarray*}  
where, conditionally on $X^{(1)}$, $\xi^{(1)}=(\xi^{(1)}(t),
t\geq 0)$ is a copy of the spatial motion $\xi$ started
at $\w(H(k_{-I^{(1)}_t}\mu)-)$. Clearly, 
$\ov\w_{(1)}$ is distributed
according to $\QQ_t(\ov\w,\cdot)$.

The definition of $\ov\w_{(2)}$ is analogous but we use 
another copy of the underlying L\'evy process. Precisely,
we let $Z$ be a copy of $X$ independent of the pair
$(X^{(1)},\xi^{(1)})$, and if
$T_{*}(Z):=\inf\{r\geq
0:Z_r=\varepsilon_1-\varepsilon_2\}$, we set
$$X^{(2)}_s=\left\{
\begin{array}{ll}
Z_s&{\rm if}\ 0\leq s\leq
T_{*}(Z),\\
\varepsilon_1-\varepsilon_2+X^{(1)}_{s-
T_{*}(Z)}
&{\rm if}\ s>T_{*}(Z).
\end{array}\right. 
$$
We then take, with an obvious notation,
$$\mu_{(2)}=[k_{-I^{(2)}_t}\mu',\rho^{(2)}_t].$$
The definition of $\w_{(2)}$ is somewhat more intricate.
Let $\tau^{(1)}$ be the (a.s. unique) time of the minimum
of $X^{(1)}$ over $[0,t]$. Consider the event
$$A(\varepsilon,\varepsilon_1,\varepsilon_2)
=\{T_{*}(Z)+\tau^{(1)}<t,
I^{(1)}_t<-\varepsilon\}.$$
Notice that $T_{*}(Z)$
is small in probability when $\varepsilon$ is small,
and $I^{(1)}_t<0$ a.s. It follows that
$P[A(\varepsilon,\varepsilon_1,\varepsilon_2)]\geq
1-\alpha(\varepsilon)$, where the
function $\alpha(\varepsilon)$ satisfies
$\alpha(\varepsilon)\la 0$ as $\varepsilon\rightarrow  0$.
However, on the event
$A(\varepsilon,\varepsilon_1,\varepsilon_2)$,
we have $I^{(2)}_t=\varepsilon_1-\varepsilon_2+I^{(1)}_t$,
and so
$$k_{-I^{(2)}_t}\mu'=
k_{-I^{(1)}_t-\varepsilon_1}k_{\varepsilon_2}\mu'=
k_{-I^{(1)}_t-\varepsilon_1}k_{\varepsilon_1}\mu=
k_{-I^{(1)}_t}\mu.$$
Also recall that from the definition of ${\cal
V}_\varepsilon(\ov\w)$, we have $\w'(r)=\w(r)$
for every $r<H(k_{\varepsilon_1}\mu)$, hence 
for every $r<H(k_{-I^{(1)}_t}\mu)$ when
$A(\varepsilon,\varepsilon_1,\varepsilon_2)$ holds.

We construct $\w_{(2)}$ by imposing that,
on the set $A(\varepsilon,\varepsilon_1,\varepsilon_2)$,
$$\w_{(2)}(r)=\left\{\begin{array}{ll}
\w'(r)=\w(r)&{\rm if}\ r<
H(k_{-I^{(2)}_t}\mu')=H(k_{-I^{(1)}_t}\mu),\\
\xi^{(1)}({r-H(k_{-I^{(1)}_t}\mu)})\quad&{\rm if}\
H(k_{-I^{(1)}_t}\mu)\leq r<H(\mu_{(2)}),
\end{array}
\right.$$
whereas on $A(\varepsilon,\varepsilon_1,\varepsilon_2)^c$,
we take
$$\w_{(2)}(r)=\left\{\begin{array}{ll}
\w'(r)&{\rm if}\ r<
H(k_{-I^{(2)}_t}\mu'),\\
\xi^{(2)}({r-H(k_{-I^{(2)}_t}\mu')})\quad&{\rm if}\
H(k_{-I^{(2)}_t}\mu')\leq r<H(\mu_{(2)}),
\end{array}
\right.$$
where, conditionally on $X^{(2)}$, $\xi^{(2)}$
is 
independent of $\xi^{(1)}$ and distributed according
to the law of $\xi$ started at
$\w'(H(k_{-I_t^{(2)}}\mu')-)$.  Note that, in the first
case, we use the {\it same} process
$\xi^{(1)}$ as in the definition of $\w_{(1)}$.
It is again easy to verify that $\ov\w_{(2)}$ is distributed
according to $\QQ_t(\ov\w',\cdot)$.

To complete the proof, note that the distance in variation
$d_{\rm var}(\mu_{(1)},\mu_{(2)})$ is equal to
$d_{\rm var}(\rho^{(1)}_t,\rho^{(2)}_t)$ on the set
$A(\varepsilon,\varepsilon_1,\varepsilon_2)$. 
Furthermore, from the construction of
$X^{(2)}$, on the set
$A(\varepsilon,\varepsilon_1,\varepsilon_2)$,
we have also
$$\rho^{(2)}_t=
\rho^{(1)}_{t-T_{*}(Z)}.$$
and thus $d_{\rm var}(\mu_{(1)},\mu_{(2)})=
d_{\rm var}(\rho^{(1)}_t,\rho^{(1)}_{t-T_{*}(Z)})$
is small in probability when $\varepsilon$ is small,
because $t$ is a.s. not a discontinuity time of $\rho^{(1)}$.
In addition, again on the set
$A(\varepsilon,\varepsilon_1,\varepsilon_2)$, the paths
$\w_{(1)}$ and
$\w_{(2)}$ coincide on the interval
$[0,H(\mu_{(1)})\wedge H(\mu_{(2)}))$, and so
$$d(\w_{(1)},\w_{(2)})\leq |H(\mu_{(2)})-H(\mu_{(1)})|
=|H(\rho^{(1)}_{t-T_{*}(Z)})-H(\rho^{(1)}_t)|$$
is small in probability when $\varepsilon$ goes to $0$.
The limiting result (\ref{SMtech1}) now follows from these
observations and the fact that 
$P[A(\varepsilon,\varepsilon_1,\varepsilon_2)^c]$
tends to $0$ as $\varepsilon$ goes to $0$. 
\cq

\subsection{Excursion measures}

We know that $\mu=0$ is a regular recurrent point for
the Markov process $\rho_s$, and the associated 
local time is the process $L^0_s$ of Section 1.3. It immediately
follows that
$(0,x)$ is also a regular recurrent point for the
L\'evy snake $(\rho_s,W_s)$, with associated local time
$L^0_s$. We will denote by $\N_x$ the corresponding 
excursion measure. It is straightforward to verify that
\begin{description}
\item{(i)} the law of $(\rho_s,s\geq 0)$ under $\N_x$
is the excursion measure $N(d\rho)$;
\item{(ii)} the conditional distribution of $(W_s,s\geq 0)$
under $\N_x$ knowing $(\rho_s,s\geq 0)$ is $Q^{H(\rho)}_x$.
\end{description}

From these properties and Proposition \ref{mesinv}, we easily
get for any nonnegative measurable function $F$ on $M_f(\R_+)\times \W$,
\begin{equation} 
\label{invar-snake}
\N_x\Big(\int _0^\sigma ds\,F(\rho_s,W_s)\Big)
=\int_0^\infty da\,e^{-\alpha a}\,E^0\otimes \Pi_x[F(J_a,(\xi_r,0\leq r\leq a))]
\end{equation}
Here, as in Chapter 1,
$J_a(dr)$ stands for the measure $1_{[0,a]}(r)
dU_r$, where $U$ is 
under the probability measure $P^0$ a subordinator
with Laplace exponent
$\wt\psi(\lambda)-\alpha$, where $\wt\psi(\lambda)=\psi(\lambda)/\lambda$. 
Note that the right side of (\ref{invar-snake}) gives an invariant measure
for the L\'evy snake $(\rho_s,W_s)$.

\medskip
The strong Markov property of the L\'evy snake can be
extended to the excursion measures in the 
following form. Let $T$ be a stopping time of the 
filtration $(\f_{s+})$ such that $T>0$, $\N_x$ a.e., let
$F$ be a nonnegative $\f_{T+}$-measurable functional
on $\D(\R_+,M_f(\R_+)\times \W)$, and let $G$ be any nonnegative
measurable functional on $\D(\R_+,M_f(\R_+)\times \W)$. Then,
$$\N_x[F\,G(\ov W_{T+s},s\geq 0)]
=\N_x[F\,\E^*_{\ov W_T}[G]],$$
where $\bP^*_{\mu,\w}$ denotes the law under $\bP_{\mu,\w}$
of the process $(\ov W_s,s\geq 0)$ stopped at $\inf\{s\geq 0,
\rho_s=0\}$. This statement follows from Theorem 
\ref{strongMarkovsnake} by standard arguments.

\section{The connection with superprocesses}

\subsection{Statement of the result}

In this section, we state and prove the basic
theorem relating the L\'evy snake with the superprocess with 
spatial motion $\xi$ and branching mechanism
$\psi$. This connection was already obtained
in a less precise form in \cite{LGLJ2}. 

We start with a few simple observations.
Let $\kappa(ds)$ be a random measure on $\R_+$,
measurable with respect to the $\sigma$-field
generated by
$(\rho_s,s\geq 0)$. Then, from the 
form of the conditional distribution
of $(W_s,s\geq 0)$ knowing $(\rho_s,s\geq 0)$,
it is easy to see that, for any 
nonnegative measurable functional $F$
on $\W_x$,
$$\E_{x}\Big[\int \kappa(ds)\,F(W_s)\Big]
=\E_x\Big[\int \kappa(ds)\,\Pi_x[F(\xi_r,0\leq r
<H_s)]\Big],$$
and a similar formula holds under $\N_x$.
This identity implies in particular that
the left limit $W_s(H_s-)$
exists $\kappa(ds)$ a.e., $\bP_x$ a.s.
(or $\N_x$ a.e.). 
We will apply this simple observation
to the random measure $d_sL^a_s$ associated
with the local time of $H$ at level $a$
(cf Chapter 1).
To simplify 
notation, we will write $\wh W_s=W_s(H_s-)$
when the limit exists, and when the limit 
does not exist, we take $\wh W_s=\Delta$,
where $\Delta$ is a cemetery point added to $E$.

In order to state the main theorem of this section,
we denote by $\z_a=\z_a(\rho,W)$
the random
measure on $E$ defined by 
$$\big<\z_a,f\big>=\int_0^{\sigma}
d_sL^a_s\,f(\wh W_s).$$
This definition makes sense under the excursion measures 
$\N_x$.

\begin{theorem}
\label{super}
Let $\mu\in M_f(E)$ and let 
$$\sum_{i\in I} \delta_{(x_i,\rho^i,W^i)}$$
be a Poisson point measure with intensity $\mu(dx)\N_x(d\rho dW)$.
Set $Z_0=\mu$ and for every $a>0$
$$Z_a=\sum_{i\in I} \z_a(\rho^i,W^i).$$
The process $(Z_a,a\geq 0)$ is a
superprocess with spatial motion $\xi$
and branching mechanism $\psi$, started
at $\mu$. 
\end{theorem}

This means that $(Z_a,a\geq 0)$ is a Markov process with values in
$M_f(E)$, whose semigroup is characterized by the following 
Laplace functional. For every $0\leq a\leq b$ and every 
function $f\in \b_{b+}(E)$,
$$E[\exp-\langle Z_b,f\rangle\mid Z_a]=\exp-\langle Z_a,u_{b-a}\rangle$$
where the function $(u_t(y),t\geq 0,y\in E)$ is the unique nonnegative solution
of the integral equation
\begin{equation} 
\label{super2}
u_t(y)+\Pi_y\Big(\int_0^t \psi(u_{t-r}(\xi_{r}))\,dr\Big)=\Pi_y(f(\xi_t)).
\end{equation}
The proof of Theorem \ref{super} is easily reduced to that of the
following proposition.

\begin{proposition}
\label{super-bis}
Let $0<a<b$ and let $f\in \b_{b+}(E)$. Then,
\begin{equation} 
\label{super1}
\N_x(\exp-\langle \z_b,f\rangle\mid (\z_r,0\leq r\leq a))=\exp-\langle \z_a,u_{b-a}\rangle
\end{equation}
where for every $t> 0$ and $y\in E$,
$$u_t(y)=\N_y(1-\exp-\langle \z_t,f\rangle).$$
Furthermore, if we set $u_0(y)=f(y)$, the function $(u_t(y),t\geq 0,y\in E)$ is the unique nonnegative solution
of the integral equation
{\rm (\ref{super2})}.
\end{proposition}

\rem Although $\N_x$ is an infinite measure, the conditioning 
in (\ref{super1}) makes sense because we can restrict our attention
to the set $\{\z_a\not =0\}=\{L^a_\sigma>0\}$ which has finite
$\N_x$-measure (cf Corollary \ref{extinction}). A similar remark applies in several places below, e.g. 
in the statement of Proposition \ref{excursions}.

\medskip

Given Proposition \ref{super-bis}, it is a straightforward exercise to
verify that the process $(Z_a,a\geq 0)$ of Theorem \ref{super} has the finite-dimensional
marginals of the superprocess with spatial motion $\xi$
and branching mechanism $\psi$, started
at $\mu$. In fact the statement of Propostion \ref{super-bis}
means that the laws of $(\z_a,a>0)$ under $\N_y$, $y\in E$
are the canonical measures of the superprocess with spatial motion $\xi$
and branching mechanism $\psi$, and given this fact, Theorem \ref{super}
is just the canonical representation of superprocesses.

The remaining part of this section is devoted
to the proof of Proposition \ref{super-bis}. 
We will proceed in two steps. In the
first one, we introduce a $\sigma$-field 
$\e_a$ that
contains $\sigma(\z_u,0\leq u\leq a)$, and we 
compute
$\N_x(\exp-\big<\z_b,f\big>\mid\e_a)$ in the form
given by (\ref{super1}). In the
second step, we establish the integral equation
(\ref{super2}).

\subsection{First step}

Recall the notation of Section 1.3
$$\tilde \tau^a_s=\inf\{r,\int_0^r du\,1_{\{
H_u\leq a\}}>s\}.$$
Note that $\tilde \tau^a_s<\infty$ for every $s\geq 0$, $\N_x$ a.e.
For $a>0$, we let $\e_a$ be the $\sigma$-field
generated by the right-continuous process
$(\rho_{\tilde \tau^a_s},W_{\tilde
\tau^a_s};s\geq 0)$ and
augmented with the class of all 
sets that are $\N_x$-negligible for every $x\in E$.
From the second approximation of Proposition \ref{LTapprox},
it is easy to verify that $L^a_\sigma$ is measurable with respect
to the $\sigma$-field generated by $(\rho_{\tilde \tau^a_s},s\geq 0)$,
and in particular with respect to $\e_a$
(cf the beginning of the proof of Theorem \ref{RK}). 

We then claim that $\z_a$ is 
$\e_a$-measurable. 
It is
enough to check that, if
$g$ is bounded and continuous on $\W_x$, 
$$\int_0^{\sigma} dL^a_s\,g(W_s)$$
is $\e_a$-measurable. However, by Proposition
\ref{LTapprox}, this integral is the limit
in $\N_x$-measure as $\varepsilon\rightarrow  0$ of
$${1\over \varepsilon}\int_0^{\sigma}
ds\,1_{\{a-\varepsilon<H_s\leq a\}}\,g(W_s).$$
For $\varepsilon<a$, this quantity coincides
with
$${1\over \varepsilon}\int_0^\infty
ds
\,1{\{a-\varepsilon<H_{\tilde\tau^{a}_s}
\leq a\}}\,
g(W_{\tilde\tau^{a}_s}),$$
and the 
claim follows from the definition of $\e_a$.

We then decompose the measure $\z_b$
according to the contributions of the 
different excursions of the process $H$
above level $a$.
Precisely, we let $(\alpha_i,\beta_i)$, $i\in I$ be the
excursion intervals of $H$ above $a$
over the time interval $[0,\sigma]$. We will
use the following simple facts that hold $N$ a.e.:
For every $i\in I$ and every $t>0$, we have 
$$\int_0^{\beta_i+t}1_{\{H_s\leq a\}}ds>\int_0^{\beta_i}1_{\{H_s\leq a\}}ds$$
and $$L^a_{\beta_i+t}>L^a_{\beta_i}.$$ 
The first assertion is
an easy consequence of the strong Markov property of $\rho$,
recalling that $\rho_s(\{a\})=0$ for every $s\geq 0$, $N$ a.e. 
To get the second one, we can use Proposition \ref{reflec}
and the definition of the local time $L^a$ to see that it is enough to
prove that 
$$\int_0^{\beta_i+t}1_{\{H_s> a\}}ds>\int_0^{\beta_i}1_{\{H_s> a\}}ds$$
for every $t>0$ and $i\in I$. Via a time-reversal argument (Corollary \ref{reversal}), it suffices to
verify that, if $\sigma^q_a=\inf\{s>q:H_s>a\}$, we have
$$\int_0^{\sigma^q_a+t}1_{\{H_s\leq a\}}ds>\int_0^{\sigma^q_a}1_{\{H_s\leq a\}}ds$$ for every 
$t>0$ and every
rational $q>0$, $N$ a.e. on the set $\{q<\sigma^q_a<\infty\}$. 
The latter fact is again a consequence of the
strong Markov property of the process $\rho$. 

\smallskip
As was observed in the proof
of Proposition \ref{reflec}, for every $i\in I$, for every
$s\in(\alpha_i,\beta_i)$, the restriction of $\rho_s$
to $[0,a]$ coincides with $\rho_{\alpha_i}=\rho_{\beta_i}$.
Furthermore, the snake
property implies that, for every $i\in I$, the paths $W_s$,
$\alpha_i<s<\beta_i$ take the same value $x_i$ at
time $a$, and this value must be the same
as the limit $\wh W_{\alpha_i}=\wh W_{\beta_i}$
(recall our assumption that $\xi$ has no fixed discontinuities). We can then define the pair $(\rho^i,W^i)
\in\D(\R_+,M_f(\R_+)\times \W)$
by setting
$$\begin{array}{ll}
\big<\rho^i_s,\varphi\big>=\int_{(a,\infty)}
\rho_{\alpha_i+s}(dr)\,\varphi(r-a)
\qquad&{\rm if}\ 0< s<\beta_i-\alpha_i\\
\rho^i_s=0&{\rm if}\ s=0\ {\rm or}\ s\geq\beta_i-\alpha_i,
\end{array}$$
and
$$\begin{array}{ll}
W^i_s(r)=W_{\alpha_i+s}(a+r),\ \zeta_{W^i_s}=
H_{\alpha_i+s}-a
\qquad&{\rm if}\ 0< s<\beta_i-\alpha_i\\
W^i_s=x_i&{\rm if}\ s=0\ {\rm or}\ s\geq\beta_i-\alpha_i.
\end{array}$$

\begin{proposition}
\label{excursions}
Under $\N_x$, conditionally on $\e_a$, the point measure
$$\sum_{i\in I} \delta_{(\rho^i,W^i)}$$
is a Poisson point measure with 
intensity 
$$\int \z_a(dy)\,\N_y(\cdot).$$
\end{proposition}

\proof Let the process $\rho^a_t$ be defined as in Proposition
\ref{reflec}. Note that under $N$ the definition of $\rho^a_t$
only makes sense for $t<\int_0^\sigma ds\,1_{\{H_s>a\}}$. 
For convenience, we take $\rho^a_t=0$ if $t\geq \int_0^\sigma ds\,1_{\{H_s>a\}}$.
We also set 
$$\wt\rho_t=\rho_{\tilde \tau^a_t}\ , \quad  \wt W_t=W_{\tilde \tau^a_t}.$$

With these definitions, the processes $\rho^i$, $i\in I$ are exactly the
excursions of the process $\rho^a$ away from $0$. For every $i\in I$,
introduce the local time at the beginning (or the end) of excursion $\rho^i$:
$$\ell^i=L^a_{\alpha_i}.$$
By Proposition \ref{reflec} and standard excursion theory, we know that
conditionally on the process $\wt \rho_t$, the point measure
$$\sum_{i\in I} \delta_{(\ell^i,\rho^i)}$$
is Poisson with intensity $1_{[0,L^a_\sigma]}(\ell)d\ell\,N(d\rho)$
(recall that
$L^a_\sigma$ is measurable with respect to the
$\sigma$-field generated by $\wt \rho$). Note that Proposition \ref{reflec}
is formulated under $\bP_x$: However, by considering the first excursion of $\rho$
away from $0$ that hits the set $\{\sup H_s>a\}$, we can easily derive the previous
assertion from Proposition \ref{reflec}.

Define $\wt L^a_s=L^a_{\tilde \tau^a_s}$ (note that this is a
continuous process), and let $\gamma^a(r)$ be the right-continuous 
inverse of $\wt L^a$:
$$\gamma^a(r)=\inf\{s\geq 0:\wt L^a_s>r\}.$$
Then, if $f$ is any nonnegative measurable function on $E$,
we have $\N_x$ a.e.
\begin{equation} 
\label{new-exit}
\langle \z_a,f\rangle=\int_0^\infty dL^a_s\,f(\wh W_s)=
\int_0^\infty d\wt L^a_s\,f(\wh{\wt W}_s)=\int_0^{L^a_\sigma} d\ell\,f(\wh{\wt W}_{\gamma^a(\ell)}).
\end{equation}
Notice that both processes $\wt L^a$ and $\gamma^a$ are measurable
with respect to the $\sigma$-field generated by $\wt \rho$ (for $\wt L^a$, this
follows again from Proposition \ref{LTapprox}).

Consider now the processes $\wt W$ and $W^i$, $i\in I$. The following two
properties are straightforward consequences of our construction:

\begin{description}

\item{(i)} The law of $\wt W$ under $Q_x^{H(\rho)}$ is $Q_x^{H(\tilde\rho)}$.

\item{(ii)} Under $Q_x^{H(\rho)}$, conditionally on $\wt W$, the ``excursions''
$W^i$, $i\in I$ are independent and the conditional distribution of $W^i$
is $Q^{H(\rho^i)}_{x_i}$, where $x_i=\wh W_{\beta_i}=\wh{\wt W}_{\gamma^a(\ell^i)}$.

\end{description}
To verify the second expression for $x_i$, note that if $\wt A^a_s=\int_0^s dr\,1_{\{H_r\leq a\}}$,
we have $W_{\beta_i}=\wt W_{\tilde A^a_{\beta_i}}$ (because $\wt A^a_{\beta_i+t}
>\wt A^a_{\beta_i}$ for every $t>0$) and $\wt
A^a_{\beta_i}=\gamma^a(\ell^i)$ (because $L^a_{\beta_i+t}>L^a_{\beta_i}=\ell_i$ for every $t>0$).

As a consequence of (i), the conditional distribution (under $\N_x$) of 
$\wt W$ knowing $\rho$ depends only on $\wt \rho$. Hence, $\wt W$
and the point measure $\sum_{i\in I} \delta_{(\ell^i,\rho^i)}$ are
conditionally independent given $\wt \rho$ under $\N_x$.

We use the previous observations in the following calculation:
\ba
&&\N_x\Big(G(\wt \rho,\wt W)\,\exp(-\sum_{i\in I} F(\rho^i,W^i))\Big)\\
&&\quad=\int N(d\rho)\,Q^{H(\rho)}_x\Big(G(\wt \rho,\wt W)\,\exp(-\sum_{i\in I} F(\rho^i,W^i))\Big)\\
&&\quad=\int N(d\rho)\,Q^{H(\rho)}_x\Big(G(\wt \rho,\wt W)\prod_{i\in I}
Q^{H(\rho^i)}_{\hat{\tilde W}_{\gamma^a(\ell^i)}}(e^{-F(\rho^i,\cdot)})\Big)\\
&&\quad=\N_x\Big(G(\wt \rho,\wt W)\prod_{i\in I}
Q^{H(\rho^i)}_{\hat{\tilde W}_{\gamma^a(\ell^i)}}(e^{-F(\rho^i,\cdot)})\Big)\\
&&\quad=\N_x\Big(G(\wt \rho,\wt W)
\,\exp\Big(-\int_0^{L^a_\sigma} d\ell\int N(d\rho)\,Q^{H(\rho)}_
{\hat{\tilde W}_{\gamma^a(\ell)}}(1-e^{-F(\rho,\cdot)})
\Big)\Big).
\ea
The second equality follows from (ii) above.
In the last one, we used the conditional independence of $\wt W$ and of the point measure
$\sum_{i\in I} \delta_{(\ell^i,\rho^i)}$, given $\wt \rho$, and the fact that the conditional
distribution of this point measure is Poisson with intensity $1_{[0,L^a_\sigma]}(\ell)d\ell\,N(d\rho)$.
Using (\ref{new-exit}), we finally get
\ba 
&&\N_x\Big(G(\wt \rho,\wt W)\,\exp(-\sum_{i\in I} F(\rho^i,W^i))\Big)\\
&&\quad=\N_x\Big(G(\wt \rho,\wt W)
\,\exp\Big(\!\!-\!\int_0^{L^a_\sigma}\! d\ell\,\N_{\hat{\tilde W}_{\gamma^a(\ell)}}(1-e^{-F(\rho,W)})
\Big)\Big)\\
&&\quad=\N_x\Big(G(\wt \rho,\wt W)
\,\exp(-\int \z_a(dy)\,\N_y(1-e^{-F}))\Big).
\ea
This completes the proof. \cq

\smallskip
Let $f$ be a nonnegative measurable function
on $E$, and let $0\leq a<b$. 
With the preceding
notation, it is easy to verify that $\N_x$ a.s.
\begin{eqnarray*}
\big<\z_b,f\big>&=&\int_0^{\sigma}
dL^b_s\,f(\wh W_s)\\
&=&\sum_{i\in I} \int_{\alpha_i}^{\beta_i}
dL^b_s\,f(\wh W_s)\\
&=&\sum_{i\in I} \int_0^\infty
dL^{b-a}_s(\rho^i)\,f(\wh W^i_s)\\
&=&\sum_{i\in I}
\big<\z_{b-a}(\rho^i,W^i),f\big>.
\end{eqnarray*}
As a consequence of Proposition \ref{excursions},
we have then
$$\N_x[\exp-\big<\z_b,f\big>\mid \e_a]
=\exp-\big<\z_a,u_{b-a}\big>,$$
where
$$u_r(y)=\N_y[1-\exp-\big<\z_r,f\big>].$$

\subsection{Second step}

It remains to prove that the function
$(u_r(y),r\geq 0,y\in E)$ introduced at the
end of the first step solves the
integral equation (\ref{super2}). By definition,
we have for $a>0$,
\begin{eqnarray}
\label{supertech1}
u_a(y)&=&\N_y\Big[1-\exp-\int_0^\infty
dL^a_s\,f(\wh W_s)\Big]\nonumber\\
&=&\N_y\Big[\int_0^\infty dL^a_s\,f(\wh W_s)\,
\exp\Big(-\int_s^\infty dL^a_r\,f(\wh W_r)\Big)
\Big]\nonumber\\
&=&\N_y\Big[\int_0^\infty dL^a_s\,f(\wh W_s)\,
\E^*_{\rho_s,W_s}\Big[\exp-\int_0^\infty
dL^a_r\,f(\wh W_r)\Big]\Big]
\end{eqnarray}
where we recall that $\bP^*_{\mu,\w}$ stands for
the law of the L\'evy snake started at $(\mu,\w)$
and stopped when $\rho_s$ first hits $0$. In the
last equality, we replaced 
$\exp-\int_s^\infty dL^a_r\,f(\wh W_r)$
by its optional projection, using the strong
Markov property of the L\'evy snake to identify
this projection. 

We now need to compute
for a fixed $(\mu,\w)\in \Theta_x$, 
$$\E^*_{\mu,\w}\Big[\exp-\int_0^\infty
dL^a_r\,f(\wh W_r)\Big].$$
We will derive this calculation from a more 
general fact, that is also
useful for forthcoming applications. First recall
that
$Y_t=\big<\rho_t,1\big>$ is distributed under
$\bP^*_{\mu,\w}$ as the underlying L\'evy process
started at $\big<\mu,1\big>$ and stopped when
it first hits $0$. We write
$J_t=\inf_{r\leq t} Y_r$, and we denote by
$(\alpha_i,\beta_i)$, ${i\in I}$ the excursion
intervals of $Y-J$ away from $0$.
For every 
$i\in I$, we set $h_i=H_{\alpha_i}=
H_{\beta_i}$.
From the snake property, it is
easy to verify that $W_s(h_i)=\w(h_i-)$
for every $s\in(\alpha_i,\beta_i)$, $i\in I$,
$\bP^*_{\mu,\w}$ a.s.
We then define the pair
$(\rho^i,W^i)$ by the formulas
$$\begin{array}{ll}
\big<\rho^i_s,\varphi\big>=\int_{(h_i,\infty)}
\rho_{\alpha_i+s}(dr)\,\varphi(r-h_i)
\qquad&{\rm if}\ 0\leq s\leq \beta_i-\alpha_i\\
\rho^i_s=0&{\rm if}\ s>\beta_i-\alpha_i,
\end{array}$$
and
$$\begin{array}{ll}
W^i_s(t)=W_{\alpha_i+s}(h_i+t),\quad
\zeta^i_s=H_{\alpha_i+s}-h_i
\qquad&{\rm if}\ 0< s< \beta_i-\alpha_i\\
W^i_s=\w(h_i-)&{\rm if}\ s=0\
{\rm or}\ s\geq\beta_i-\alpha_i.
\end{array}$$

\begin{lemma}
\label{subexcursions}
Let $(\mu,\w)\in\Theta_x$.
The point measure
$$\sum_{i\in I} \delta_{(h_i,\rho^i,W^i)}$$
is under $\bP^*_{\mu,\w}$
a Poisson point measure with intensity
$$\mu(dh)\,\N_{\w(h-)}(d\rho\,dW).$$
\end{lemma}

\proof Consider first the point measure
$$\sum_{i\in I} \delta_{(h_i,\rho^i)}.$$
If $I_s=J_s-\big<\mu,1\big>$, we have
$h_i=H(\rho_{\alpha_i})
=H(k_{-I_{\alpha_i}}\mu)$. Excursion theory
for $Y-J$ ensures that
$$\sum_{i\in I} \delta_{(-I_{\alpha_i},\rho^i)}$$
is under $\bP^*_{\mu,\w}$
a Poisson point measure with intensity
$1_{[0,<\mu,1>]}(u)\,du\,N(d\rho)$.
Since the image measure of
$1_{[0,<\mu,1>]}(u)\,du$ under the 
mapping $u\la H(k_u\mu)$ is
precisely the measure $\mu$, it follows
that
$$\sum_{i\in I} \delta_{(h_i,\rho^i)}$$
is a Poisson point measure with intensity
$\mu(dh)\,N(d\rho)$. To complete the proof,
it remains to obtain the conditional
distribution of $(W^i,i\in I)$ knowing
$(\rho_s,s\geq 0)$. However, the form of the
conditional law $Q^H_\w$ easily implies that
under $Q^H_\w$, the processes $W^i$, $i\in I$
are independent, and furthermore the conditional 
distribution of $W^i$ is $Q^{H^i}_{\w(h_i-)}$,
where $H^i_s=H(\rho^i_s)$. It follows that
$$\sum_{i\in I} \delta_{(h_i,\rho^i,W^i)}$$
is a Poisson measure with intensity
$$\mu(dh)\,N(d\rho)\,Q^{H(\rho)}_{\w(h-)}(dW)
=\mu(dh)\,\N_{\w(h-)}(d\rho dW).$$
This completes the proof.\cq

We apply Lemma \ref{subexcursions} to
a pair $(\mu,\w)$ such that $H(\mu)\leq a$ and
$\mu(\{H(\mu)\})=0$. Then, it is easy to verify
that $\bP^*_{\mu,\w}$ a.s.
$$\int_0^\infty dL^a_r\,f(\wh W_r)
=\sum_{i\in I} \int_{\alpha_i}^{\beta_i}
dL^a_r\,f(\wh W_r)
=\sum_{i\in I}\int_0^\infty dL^{a-h_i}_r(\rho^i)
\,f(\wh W^i_r),$$
and thus, by Lemma \ref{subexcursions},
\begin{equation}
\label{supertech2}
\E^*_{\mu,\w}\Big[\exp-\int_0^\infty
dL^a_r\,f(\wh W_r)\Big]
=\exp\Big(-\int \mu(dh)\,\N_{\w(h-)}
[1-\exp-\big<{\cal Z}_{a-h},f\big>]\Big).
\end{equation}

We now come back to formula (\ref{supertech1}).
As a consequence of Proposition \ref{mesinv}, we
know that 
$\rho_s(\{a\})=0$
for 
every $s\geq 0$, $\N_y$ a.e. We can
thus use  (\ref{supertech2}) to get
\begin{eqnarray}
\label{supertech3}
u_a(y)&=&\N_y\Big[\int_0^\infty dL^a_s\,f(\wh
W_s)\,
\exp\Big(-\int \rho_s(dh)\,
\N_{W_s(h-)}
[1-\exp-\big<{\cal Z}_{a-h},f\big>]\Big)\Big]
\nonumber\\
&=&\N_y\Big[\int_0^\infty dL^a_s\,f(\wh
W_s)\,
\exp\Big(-\int \rho_s(dh)\,
u_{a-h}({W_s(h-)})\Big)\Big].
\end{eqnarray}

\begin{lemma}
\label{moment1super}
For any nonnegative measurable function
$F$ on $\Theta_y$,
$$\N_y\Big[\int_0^\infty dL^a_s\,
F(\rho_s,W_s)\Big]= e^{-\alpha a}\,E^0\otimes
\Pi_y[F(J_a,(\xi_r,0\leq r<a))].$$
\end{lemma}

\proof If $F(\rho_s,W_s)$ depends only on
$\rho_s$, the result follows from
Corollary \ref{localinvar}. In
the general case, we may take $F$ such that
$F(\rho_s,W_s)=F_1(\rho_s)F_2(W_s)$, and we 
use the simple observation of the
beginning of this section.\cq

From (\ref{supertech3}) and Lemma
\ref{moment1super}, we get
\begin{eqnarray*}
u_a(y)&=&e^{-\alpha a}\,E^0\otimes \Pi_y\Big[f(\xi_a)
\exp\Big(-\int
J_a(dh)\,u_{a-h}(\xi_{h-})\Big)\Big]\\
&=&\Pi_y\Big[f(\xi_a)
\exp\Big(-\int_0^a \wt\psi(u_{a-r}(\xi_r))\,dr
\Big)\Big].
\end{eqnarray*}
The proof of (\ref{super2}) is now completed
by routine calculations. We have
\begin{eqnarray*}
u_a(y)&=&\Pi_y[f(\xi_a)]
-\Pi_y\Big[f(\xi_a)\int_0^a
db\,\wt\psi(u_{a-b}(\xi_b))\,\exp\Big(
-\int_b^a \wt\psi(u_{a-r}(\xi_r))\,dr\Big)\Big]\\
&=&\Pi_y[f(\xi_a)]
-\Pi_y\Big[\int_0^a db\,\wt\psi(u_{a-b}(\xi_b))\\
&&\hspace{25mm}\Pi_{\xi_{b}}\Big[f(\xi_{a-b}) \exp\Big(
-\int_0^{a-b} \wt\psi(u_{a-b-r}(\xi_r))\,dr\Big)
\Big]\Big]\\
&=&\Pi_y[f(\xi_a)]
-\Pi_y\Big[\int_0^a db\,\wt\psi(u_{a-b}(\xi_b))\,
u_{a-b}(\xi_b)\Big],
\end{eqnarray*}
which gives (\ref{super2}) and completes the
proof of Proposition \ref{super-bis}.

\section{Exit measures}
\medskip
Throughout this section, we consider an open set
$D\subset E$, and we denote by
$\tau$ the first exit time of $\xi$ from $D$:
$$\tau=\inf\{t\geq 0: \xi_t\notin D\},$$
where $\inf\emptyset=\infty$ as usual.
By abuse of notation, we will also denote by 
$\tau(\w)$ the exit time from $D$ of a killed
path $\w\in \W$,
$$\tau(\w)=\inf\{t\in[0,\zeta_\w) : \w(t)\notin D\}\,.$$
Let $x\in D$.
The next result is much analogous to Proposition \ref{reflec}.

\begin{proposition}
\label{reflec-exit}
Assume that $\Pi_x(\tau<\infty)>0$. Then,
$$\int_0^\infty ds\,1_{\{\tau(W_s)<H_s\}}=\infty\,,\qquad \bP_x \hbox{ a.s.}$$
Furthermore, let
$$\sigma^D_s=\inf\{t\geq 0:\int_0^t dr\,1_{\{\tau(W_r)<H_r\}}>s\},$$
and let $\rho^D_s\in M_f(\R_+)$ be defined by
$$\langle\rho^D_s,f\rangle 
=\int
\rho_{\sigma^D_s}(dr)\,f(r-\tau(W_{\sigma^D_s}))\,1_{\{r>\tau(W_{\sigma^D_s})\}}.$$
Then the process $(\rho^D_s,s\geq 0)$ has the same distribution under
$\bP_x$ as $(\rho_s,s\geq 0)$.
\end{proposition}

\rem We could have considered the more general situation of a
space-time open set $D$ (as a matter of fact, this is not
really more general as we could replace $\xi_t$ by $(t,\xi_t)$). Taking
$D=[0,a)\times E$, we would recover part of the statement of Proposition \ref{reflec}. 
This proposition
contains an independence statement that could also  be extended to the present
setting.

\medskip
\proof To simplify notation, we set
$$A^D_s=\int_0^s dr\,1_{\{\tau(W_r)<H_r\}}.$$
By using (\ref{invar-snake}), excursion theory and our assumption 
$\Pi_x(\tau<\infty)>0$, it is a simple exercise to verify that 
$A^D_\infty=\infty$, $\bP_x$ a.s., and thus the definition of
$\sigma^D_s$ makes sense for every $s\geq 0$, a.s. The arguments then are
much similar to the proof of Proposition \ref{reflec}. For every
$\varepsilon>0$, we introduce the stopping times
$S^k_\varepsilon$,
$T^k_\varepsilon$, $k\geq 1$, defined inductively
by:
\begin{eqnarray*}
&&S^1_\varepsilon=\inf\{s\geq 0:\tau(W_s)<\infty
\hbox{ and }\rho_s((\tau(W_s),\infty))\geq \varepsilon\},\\
&&T^k_\varepsilon=\inf\{s\geq
S^k_\varepsilon:\tau(W_s)=\infty\},\\
&&S^{k+1}_\varepsilon=\inf\{s\geq
T^k_\varepsilon:\tau(W_s)<\infty
\hbox{ and }\rho_s((\tau(W_s),\infty))\geq \varepsilon\}.
\end{eqnarray*}
It is easy to see that these stopping times
are a.s. finite, and $S^k_\varepsilon\ua \infty$,
$T^k_\varepsilon\ua\infty$ as $k\ua \infty$.

From the key formula (\ref{strongMarkov}), we see that for
$$S^k_\varepsilon\leq s<\inf\{r\geq S^k_\varepsilon:
\langle \rho_r,1\rangle\leq \rho_{S^k_\varepsilon}([0,\tau(W_{S^k_\varepsilon})])\}$$
we have $H_s>\tau(W_{S^k_\varepsilon})$, and 
the paths $W_s$ and $W_{S^k_\varepsilon}$ coincide over $[0,\tau(W_{S^k_\varepsilon})]$
(by the snake property),
so that in particular $\tau(W_s)=\tau(W_{S^k_\varepsilon})<\infty$. On the other
hand, for
$$s=\inf\{r\geq S^k_\varepsilon:
\langle \rho_r,1\rangle\leq \rho_{S^k_\varepsilon}([0,\tau(W_{S^k_\varepsilon})])\}$$
the path $W_s$ is the restriction of $W_{S^k_\varepsilon}$
to $[0,\tau(W_{S^k_\varepsilon}))$ and thus $\tau(W_s)=\infty$.
From these observations, we see that
$$T^k_\varepsilon=\inf\{r\geq S^k_\varepsilon:
\langle \rho_r,1\rangle\leq \rho_{S^k_\varepsilon}([0,\tau(W_{S^k_\varepsilon})])\}$$
and that conditionally on
the past up to time $S^k_\varepsilon$, the process
$$Y^{k,\varepsilon}_s=\rho_{(S^k_\varepsilon+s)\wedge T^k_\varepsilon}
((\tau(W_{S^k_\varepsilon}),\infty))$$ is distributed as
the  underlying L\'evy process started at
$\rho_{S^k_\varepsilon}((\tau(W_{S^k_\varepsilon}),\infty))$
and stopped at its first hitting time of $0$.

The same argument as in the proof of (\ref{reflec1}) shows that, for every $t\geq 0$,
\begin{equation}\label{reflec1-snake}
\lim_{\varepsilon\rightarrow  0}
\sup_{\{k\geq 1,S^k_\varepsilon\leq t\}}
\rho_{S^k_\varepsilon}((\tau(W_{S^k_\varepsilon}),\infty))=0,
\qquad{\rm a.s.}
\end{equation}
\par The remaining part of the proof is very similar to the end
of the proof of Proposition \ref{reflec}. Using (\ref{reflec1-snake})
and the observations preceding (\ref{reflec1-snake}), we get by
a passage to the limit $\varepsilon\rightarrow  0$ that the total mass
process $\langle\rho^D_s,1\rangle=\rho_{\sigma^D_s}((\tau(W_{\sigma^D_s}),\infty))$
has the same distribution as the process $\langle \rho_s,1\rangle$. Then the
statement of the proposition follows by an argument similar to the second step
of the proof of Proposition \ref{reflec}. \cq

\medskip
Let $\ell^D=(\ell^D(s),s\geq 0)$ be the local time at $0$ of the
process $\langle \rho^D,1\rangle$. We define the {\it exit local time}
from $D$ by the formula
$$L^D_s=\ell^D(A^D_s)=\ell^D(\int_0^s dr\,1_{\{\tau(W_r)<H_r\}}).$$
Recall from (\ref{invar-snake}) the notation $J_a,\,P^0$.

\begin{proposition}
\label{moment1exit}
For any nonnegative measurable function $\Phi$ on $M_f(\R_+)\times\W$,
$$\N_x\Big(\int_0^\sigma dL^D_s \,\Phi(\rho_s,W_s)\Big)
=E^0\otimes\Pi_x\Big[1_{\{\tau<\infty\}}e^{-\alpha\tau}\,\Phi(J_\tau,(\xi_r,0\leq r<\tau))\Big].$$
\end{proposition}

\medskip
\proof By applying Lemma \ref{heightinfimum} to the reflected L\'evy
process $\langle \rho^D,1\rangle$, we get for every $s\geq 0$,
$$\ell^D(s)=\lim_{\varepsilon\rightarrow  0}{1\over \varepsilon}
\int_0^s dr\,1_{\{0<H(\rho^D_r)\leq \varepsilon\}}$$
in $L^1(\bP_x)$. From a simple monotonicity argument, we have then for every $t\geq
0$
$$\lim_{\varepsilon\rightarrow  0}
\E_x\left[\sup_{s\leq t}\left|\ell^D(s)-{1\over \varepsilon}
\int_0^s dr\,1_{\{0<H(\rho^D_r)\leq \varepsilon\}}\right|\right]=0.$$
Using the formulas $L^D_s=\ell^D(A^D_s)$ and 
$H(\rho_{\sigma^D_r})=\tau(W_{\sigma^D_r})+H(\rho^D_r)$ (the latter
holding on the set $\{H(\rho_{\sigma^D_r})>\tau(W_{\sigma^D_r})\}$, by the definition of
$\rho^D$), we obtain
$$\lim_{\varepsilon\rightarrow  0}
\E_x\left[\sup_{s\leq t}\left|L^D_s-{1\over \varepsilon}
\int_0^s dr\,1_{\{\tau(W_r)<H_r\leq \tau(W_r)+\varepsilon\}}\right|\right]=0.$$
Arguing as in the derivation of (\ref{LTapproxexc}), we get,
for any measurable subset $V$ of $\D(\R_+,M_f(\R_+)\times \W)$ such that $\N_x(V)<\infty$,
\begin{equation} 
\label{tech-momentexit}
\lim_{\varepsilon\rightarrow  0}
\N_x\left(1_V\sup_{s\leq t}\left|L^D_s-{1\over \varepsilon}
\int_0^s dr\,1_{\{\tau(W_r)<H_r\leq \tau(W_r)+\varepsilon\}}\right|\right)=0.
\end{equation}

We then observe that for any bounded measurable function
$F$ on $\R_+\times M_f(\R_+)\times \W$, we have
\begin{equation} 
\label{tech-momentexit2}
\N_x\Big(\int_0^\sigma ds\,F(s,\rho_s,W_s)\Big)
=E\otimes \Pi_x\Big[\int_0^{L_\infty} da\,F(L^{-1}(a),\Sigma_a,(\xi_r,0\leq r< a))\Big]
\end{equation}
where the random measure $\Sigma_a$ is defined under $P$ by
$$\langle \Sigma_a,\varphi\rangle=\int_0^{L^{-1}(a)} dS_s\,\varphi(a-L_s).$$
Indeed, we observe that the special case where $F(s,\mu,\w)$ does not depend
on $\w$,
$$N\Big(\int_0^\sigma ds\,F(s,\rho_s)\Big)
=E\Big[\int_0^{L_\infty} da\,F(L^{-1}(a),\Sigma_a)\Big]$$
is a consequence of Proposition \ref{keyinv} (see the proof of
Proposition \ref{mesinv}), and it then suffices to use 
the conditional distribution of $W$ knowing $(\rho_s,s\geq 0)$.

After these preliminaries, we turn to the proof of the proposition.
We let $F$ be a bounded continuous function
on $\R_+\times M_f(\R_+)\times \W$, and assume in addition that
there exist $\delta>0$ and $A>0$ such that $F(s,\mu,\w)=0$
if $s\leq \delta$ or $s\geq A$. As a consequence of 
(\ref{tech-momentexit}) and (\ref{tech-momentexit2}), we have then
\ba
&&\N_x\Big(\int_0^\sigma dL^D_s \,F(s,\rho_s,W_s)\Big)\\
&&\quad=\lim_{\varepsilon\rightarrow  0}
\N_x\Big({1\over \varepsilon}
\int_0^\sigma dr\,F(r,\rho_r,W_r)\,1_{\{\tau(W_r)<H_r\leq \tau(W_r)+\varepsilon\}}\Big)\\
&&\quad=\lim_{\varepsilon\rightarrow  0}
{1\over \varepsilon} E\otimes \Pi_x\Big[
\int_0^{L_\infty} da\,F(L^{-1}(a),\Sigma_a,(\xi_r,0\leq r<a))\,
1_{\{\tau<a\leq \tau+\varepsilon\}}\Big]\\
&&\quad=E\otimes \Pi_x[1_{\{\tau<L_\infty\}}\,F(L^{-1}(\tau),\Sigma_\tau,(\xi_r,0\leq r<
\tau))].
\ea
From this identity, we easily get
$$\N_x\Big(\int_0^\sigma dL^D_s \,\Phi(\rho_s,W_s)\Big)
=E\otimes\Pi_x[1_{\{\tau<L_\infty\}}\,\Phi(\Sigma_\tau,(\xi_r,0\leq r<\tau))].$$
Recall that $P[L_\infty>a]=e^{-\alpha a}$ and, that conditionally on $\{L_\infty>a\}$,
$\Sigma_a$ has the same distribution as $J_a$. The last formula is thus
equivalent to the statement of the proposition. \par\cq

We now introduce an additional assumption. Namely we assume that
for every $x\in D$, the process $\xi$ is $\Pi_x$ a.s. continuous at
$t=\tau$, on the event $\{\tau<\infty\}$. Obviously this assumption holds 
if $\xi$ has continuous sample paths, but there are other cases
of interest.

Under this assumption, Proposition \ref{moment1exit} ensures that
$\N_x$ a.e. the left limit $\wh W_s$ exists $dL^D_s$ a.e. over $[0,\sigma]$
and belongs to $\partial D$. We define under $\N_x$ the {\it exit measure} $\z^D$ from $D$
by the formula
$$\langle \z^D,\varphi\rangle=\int_0^\sigma dL^D_s\,\varphi(\wh W_s).$$
The previous considerations show that $\z^D$ is a (finite) measure
supported on $\partial D$. As a consequence of
Proposition \ref{moment1exit}, we have for every nonnegative measurable function
$g$ on $\d D$,
$$\N_x(\langle \z^D,g\rangle)=\Pi_x(1_{\{\tau<\infty\}}e^{-\alpha\tau}g(\xi_\tau)).$$

\begin{theorem}
\label{integrexit}
Let $g$ be a bounded nonnegative measurable function on $\d D$. For every $x\in D$ set
$$u(x)=\N_x(1-\exp-\langle \z^D,g\rangle).$$
Then $u$ solves the integral equation
$$u(x)+\Pi_x\Big(\int_0^\tau dt\,\psi(u(\xi_t))\Big)=
\Pi_x(1_{\{\tau<\infty\}}g(\xi_\tau)).$$
\end{theorem}

\proof Several arguments are analogous to the second step of 
the proof of Proposition \ref{super-bis} in Section 4, and so we 
will skip some details. By the definition of $\z^D$, we have
\ba
u(x)&=&\N_x\Big(1-\exp-\int_0^\sigma dL^D_s\,g(\wh W_s)\Big)\\
&=&\N_x\Big(\int_0^\sigma dL^D_s\,g(\wh W_s)\,\exp\Big(-\int_s^\sigma
dL^D_r\,g(\wh W_r)\Big)\Big)\\
&=&\N_x\Big(\int_0^\sigma dL^D_s\,g(\wh W_s)\,
\E^*_{\rho_s,W_s}\Big(\exp -\int_0^\infty
dL^D_r\,g(\wh W_r)\Big)\Big).
\ea
Note that the definition of the random measure $dL^D_r$ makes sense under 
$\bP^*_{\mu,\w}$, provided that $\tau(\w)=\infty$, thanks to 
Lemma \ref{subexcursions} and the 
approximations used in the proof of Proposition \ref{moment1exit}.
Using Lemma \ref{subexcursions} as in subsection 4.2.3, we get
if $(\mu,\w)\in\Theta_x$ is such that $\tau(\w)=\infty$,
\ba
&&\E^*_{\mu,\w}\Big(\exp -\int_0^\infty
dL^D_r\,g(\wh W_r)\Big)\\
&&\quad=
\exp\Big(-\int \mu(dh)\,\N_{\w(h-)}\Big(1-\exp-\int_0^\sigma dL^D_r\,g(\wh W_r)\Big)\Big)\\
&&\quad=
\exp\Big(-\int \mu(dh)\,\N_{\w(h-)}\Big(1-e^{-\langle \z^D,g\rangle}\Big)\Big)\\
&&\quad=
\exp(-\int \mu(dh)\,u(\w(h-))).
\ea
Hence, using also Proposition \ref{moment1exit},
\ba
u(x)&=&\N_x\Big(\int_0^\sigma dL^D_s\,g(\wh W_s)\,
\exp(-\int \rho_s(dh)\,u(W_s(h-))\Big)\\
&=& E^0\otimes \Pi_x\Big(1_{\{\tau<\infty\}}e^{-\alpha \tau} g(\xi_\tau)\,\exp(-\int J_\tau(dh)\,u(\xi_{h-}))\Big)
\\
&=&\Pi_x\Big(1_{\{\tau<\infty\}}\,g(\xi_\tau)\,\exp\Big(-\int_0^\tau dh\,\wt \psi
(u(\xi_h))\Big)\Big).
\ea
The integral equation of the theorem now follows by the same 
routine arguments used in the end of the proof of Proposition \ref{super-bis}. \cq

\section{Continuity properties of the L\'evy snake}

From now on until the end of this chapter we assume that
the underlying spatial motion $\xi$ has continuous
sample paths. The construction of Section 4.1 applies
with the following minor simplification. Rather than
considering c\` adl\` ag paths, we can define $\W_x$
as the set of all $E$-valued killed continuous paths 
started at $x$. An element of $\W$
is thus a continuous mapping $\w:[0,\zeta)\longrightarrow E$,
and the distance between $\w$ and $\w'$ is defined by
\begin{equation} 
\label{newdistance}
d(\w,\w')=\delta(\w(0),\w'(0))+|\zeta-\zeta'|
+\int_0^{\zeta\wedge\zeta'}dt\;(\sup_{r\leq t}\delta(\w(r),\w'(r))\wedge 1).
\end{equation}
Without risk of confusion, we will keep the same notation as 
in Section 4.1. The construction developed there goes through
without change with these new definitions.

Our goal is to provide conditions on $\psi$ and $\xi$
that will ensure that the process $W_s$ is continuous
with respect to a distance finer than $d$, which we now
introduce. We need to consider stopped paths rather than
killed paths. A stopped (continuous) path
is a continuous mapping $\w:[0,\zeta]\longrightarrow E$,
where $\zeta\geq 0$. When $\zeta=0$, we identify $\w$ with $\w(0)\in E$. We denote by $\W^*$
the set of all stopped paths in $E$. The set
$\W^*$ is equipped with the distance
$$d^*(\w,\w')=|\zeta-\zeta'|+\sup_{t\geq 0} \delta(\w(t\wedge \zeta),\w(t\wedge
\zeta')).$$ 
Note that $(\W^*,d^*)$ is a Polish space.

If $\w\in \W$ is a killed path such that $\zeta>0$ and the left limit
$\wh \w=\w(\zeta-)$ exists, we write $\w^*$ for the corresponding
stopped path $\w^*(t)=\w(t)$ if $t<\zeta$, and $\w^*(\zeta)=\wh \w$.
When $\zeta=0$ we make the convention that $x^*=x$. Note that
$W_s^*$ is well defined $\bP_x$ a.s., for every fixed $s\geq 0$.

As in Chapter 1, we set 
$$\gamma=\sup\{a\geq 0:\lim_{\lambda\rightarrow \infty}\lambda^{-a}\psi(\lambda)=\infty\}\geq 1.$$

\begin{proposition}
\label{continuity-snake}
Suppose that there exist three constants $p>0$, $q>0$
and $C<\infty$ such that for every $t>0$ and $x\in E$,
\begin{equation} 
\label{assumspat}
\Pi_x\Big[\sup_{r\leq t} \delta(x,\xi_r)^p\Big]\leq C\,t^q.
\end{equation}
Suppose in addition that
$$q(1-{1\over \gamma})>1.$$
Then the left limit $\wh W_s=W_s(H_s-)$ exists for every $s\geq 0$,
$\bP_x$ a.s. or $\N_x$ a.e. Furthermore the process $(W^*_s,s\geq 0)$
has continuous sample paths with respect to the distance $d^*$, $\bP_x$ a.s.
or $\N_x$ a.e.
\end{proposition}

\rem Only the small values of $t$ are relevant in our assumption (\ref{assumspat})
since we can always replace the distance $\delta$ by $\delta\wedge 1$. 
Uniformity in $x$ could also be relaxed, but we do not strive for
the best conditions.

\smallskip
\proof It is enough to argue under $\bP_x$. Let us fix
$t\geq 0$ and $s\in(0,1)$. Then,
$$ \E_x[d^* (W^*_t , W^*_{t+s } )^p] \leq 2^p\Big( \E_x
[|H_{t+s } -H_t |^p]  + \E_x\Big[
\sup_{r\geq 0 } \delta ( W^*_t(r\wedge H_t) , W^*_{t+s}(r\wedge H_{t+s}) )^p 
\Big]\Big)\; .$$
To simplify notation, set $m =m_H(t,t+s)= \inf_{[t, t+s ] } H_r $.
From the conditional distribution of the process $W$ knowing $H$, we
easily get
\ba
&&\E_x\Big[
\sup_{u\geq 0 } \delta ( W^*_t(u\wedge H_s) , W^*_{t+s}(u\wedge H_{t+s}) )^p 
\Big|\,H_r,r\geq 0\Big]\\
&&\qquad\leq
2^p \Big(
\Pi_x \Big[ \Pi_{\xi_m } 
\Big[\sup_{u\leq H_{t } -m } 
\delta (\xi_0 , \xi_u )^p \Big]\Big]  
+  \Pi_x \Big[ \Pi_{\xi_m } \Big[
\sup_{u\leq H_{t+s} -m } \delta ( \xi_0 , \xi_u ) ^p \Big]
\Big]\Big)  \\
&&\qquad \leq C\,2^p
\left( |H_t -m |^{q} + | H_{t+s } -m |^{q} \right) ,
\ea
using our assumption (\ref{assumspat}) in the last bound.
By combining the previous estimates with Lemma \ref{moments-height},
we arrive at
$$ \E_x[ d^* (W^*_t, W^*_{t+s } )^p ] \leq  
2^{2p+1} \left( C_p \varphi (1/s)^{-p} + 
C\, C_{q} \varphi (1/s)^{- q } \right) , $$
where $\varphi(\lambda)=\lambda/\psi^{-1}(\lambda)$.
Now choose $\alpha\in(0,1-{1\over \gamma})$
such that $q\alpha>1$. Notice that we may also
assume $p\alpha>1$ since by replacing the distance
$\delta$ by $\delta\wedge 1$, we can take $p$ as large as we wish. 
The condition $\alpha<1-{1\over \gamma}$ and the definition
of $\gamma$ imply that $\varphi(\lambda)\geq c\lambda^\alpha$
for every $\lambda\geq 1$, for some constant $c>0$. Hence,
there exists a constant $C'$ independent of $t$ and $s$
such that 
$$\E_x[ d^* (W^*_t, W^*_{t+s } )^p ] 
\leq  
C' (s^{ p\alpha}+s^{q\alpha})
\,. $$
The Kolmogorov lemma then gives the existence of a continuous
modification of the process $(W^*_s,s\geq 0)$ with
respect to the distance $d^*$. The various assertions
of the proposition follow easily, recalling that we already know that the
process $(W_s,s\geq 0)$ has continuous paths for the distance $d$. \cq

\section{The Brownian motion case}

In this section, we concentrate on the case when the
underlying spatial motion $\xi$ is Brownian motion
in $\R^d$. We will give a necessary and sufficient 
condition for the process $W^*$ to have a modification that is continuous
with respect to the distance $d^*$.

To this end, we
introduce the following condition on $\psi$:
\begin{equation} 
\label{keycond}
\int_1^\infty \Big(\int_0^t \psi(u)\,du\Big)^{-1/2} dt<\infty.
\end{equation}
Note that this condition is stronger than the condition 
$\int_1^\infty du/\psi(u)<\infty$ for the
path continuity of $H$. In fact, since $\psi$ is convex, there
exists a positive constant $c$ such $\psi(t)\geq ct$ for every
$t\geq1$. Then, for $t\geq 1$,
$$\int_0^t \psi(u)du\leq t\psi(t)\leq c^{-1}\psi(t)^2$$
and thus
$$\int_1^\infty {du\over\psi(u)}\leq c^{-1/2}
\int_1^\infty \Big(\int_0^t \psi(u)\,du\Big)^{-1/2} dt.$$
Also note that (\ref{keycond}) holds if $\gamma>1$.
On the other hand, it is easy to produce examples where
(\ref{keycond}) does not hold although $H$ has continuous sample
paths.

\smallskip
Condition (\ref{keycond}) was introduced in connection with
solutions of $\Delta u=\psi(u)$ in domains of $\R^d$. We briefly
review the results that will be relevant to our needs
(see \cite{Keller},\cite{Oss} and also Lemma 2.3
in \cite{Sh1}). We denote by $B_r$ the open ball of radius
$r$ centered at the origin in $\R^d$.

\begin{description}
\item{A.} If (\ref{keycond}) holds, then, for every $r>0$, 
there exists a positive solution of the problem
\begin{equation} 
\label{explosionpb}
\left\{\begin{array}{ll}
{1\over 2}\Delta u=\psi(u)\qquad&\hbox{in }B_r\\
u_{|\d B_r}=\infty\,.&
\end{array}
\right.
\end{equation}
Here the condition $u_{|\d B_r}=\infty$ means that $u(x)$ tends to $+\infty$
as $x\rightarrow  y,x\in B_r$, for every $y\in \d B_r$.
\item{B.} If (\ref{keycond}) does not hold, then for every $c>0$,
there exists a positive solution of the problem
\begin{equation} 
\label{pbespace}
\left\{\begin{array}{ll}
{1\over 2}\Delta u=\psi(u)\qquad&\hbox{in }\R^d\\
u(0)=c\,.&
\end{array}
\right.
\end{equation}
\end{description}

Connections between the L\'evy snake and the partial differential equation
$\Delta u=\psi(u)$ follow from Theorem \ref{integrexit}. Note that this is
just a reformulation of the well-known connections involving
superprocesses. We use the notation of Section 4.3. A domain $D$
in $\R^d$ is regular if every point $y$ of $\d D$ is regular for
$D^c$, that is: $\inf\{t>0:\xi_t\notin D\}=0$, $\Pi_y$ a.s.

\begin{proposition}
\label{DirichletPDE}
Assume that $\xi$ is Brownian motion in $\R^d$.
Let $D$ be a bounded regular domain in $\R^d$, and let $g$
be a nonnegative continuous function on $\d D$. Then the function
$$u(x)=\N_x(1-\exp-\langle \z^D,g\rangle)$$
is twice continuously differentiable in $D$ and is the unique nonnegative solution of the problem
\begin{equation} 
\label{Dirichletpb}
\left\{\begin{array}{ll}
{1\over 2}\Delta u=\psi(u)\qquad&\hbox{in }D\\
u_{|\d D}=g\,.&
\end{array}
\right.
\end{equation}
\end{proposition}

\proof This follows from Theorem \ref{integrexit} by standard arguments.
In the context of superprocesses, the result is due to Dynkin \cite{Dy0}.
See e.g. Chapter 5 in \cite{LG99} for a proof in the case $\psi(u)=u^2$, which is
readily extended. \cq

\medskip
We can now state our main result.

\begin{theorem}
\label{compactness-Brownian}
Assume that $\xi$ is Brownian motion in $\R^d$.
The following conditions are equivalent.
\begin{description}
\item{\rm (i)} $\N_0(\z^{B_r}\not =0)<\infty$ for some $r>0$.
\item{\rm (ii)} $\N_0(\z^{B_r}\not =0)<\infty$ for every $r>0$.
\item{\rm (iii)} The left limit $\wh W_s=W_s(\zeta_s-)$ exists for every
$s\geq 0$, $\bP_0$ a.s., and the mapping $s\rightarrow  \wh W_s$
is continuous, $\bP_0$ a.s.
\item{\rm (iv)} The left limit $\wh W_s=W_s(\zeta_s-)$ exists for every
$s\geq 0$, $\bP_0$ a.s., and the mapping $s\rightarrow  W_s^*$
is continuous for the metric $d^*$, $\bP_0$ a.s.
\item{\rm (v)} Condition {\rm (\ref{keycond})} holds.
\end{description}
\end{theorem}

\rem The conditions of Theorem \ref{compactness-Brownian} are also equivalent
to the a.s. compactness of the range of the superprocess with spatial
motion $\xi$ and branching mechanism $\psi$, started at a nonzero initial 
measure $\mu$ with compact support. This fact, that follows from Theorem 5.1
in \cite{Sh1}, can be deduced from the representation of Theorem \ref{super}.

\medskip
\proof The equivalence between (i),(ii) and (v) is easy given facts A. and B. 
recalled above. We essentially reproduce arguments of \cite{Sh1}. 
By fact B.,
if (v) does not hold, then we can for every $c>0$ find a 
nonnegative function $v_c$
such that $v_c(0)=c$ and ${1\over 2}\Delta v_c=\psi(v_c)$ in $\R^d$.
Let $r>0$ and $\lambda>0$. By Proposition \ref{DirichletPDE}, the 
nonnegative function
$$u_{\lambda,r}(x)=\N_x(1-\exp-\lambda\langle \z^{B_r},1\rangle),
\qquad x\in B_r$$
solves ${1\over 2}\Delta u_{\lambda,r}=\psi(u_{\lambda,r})$
in $B_r$ with boundary condition $\lambda$ on $\d B_r$. By choosing
$\lambda$ sufficiently large so that $\sup\{v_c(y),y\in \d B_r\}<\lambda$,
and using the comparison principle for nonnegative solutions
of ${1\over 2}\Delta u=\psi(u)$ (see Lemma V.7 in \cite{LG99}), we see that $v_c\leq u_{\lambda,r}$ in $B_r$.
In particular,
$$c=v_c(0)\leq u_{\lambda,r}(0)\leq \N_0(\z^{B_r}\not =0).$$
Since $c$ was arbitrary, we get $\N_0(\z^{B_r}\not =0)=\infty$
and we have proved that (i) $\rightarrow $ (v). Trivially (ii) $\rightarrow $ (i).

Suppose now that (v) holds. Let $r>0$. By fact A., we can find 
a function $u_{(r)}$ such that ${1\over 2}\Delta u_{(r)}=\psi(u_{(r)})$
in $B_r$ with boundary condition $+\infty$ on $\d B_r$. The 
maximum principle then implies that, for every $\lambda>0$,
$u_{\lambda,r}\leq u_{(r)}$. Hence
$$\N_0(\z^{B_r}\not =0)=\lim_{\lambda\uparrow \infty}\uparrow 
u_{\lambda,r}(0)\leq u_{(r)}(0)<\infty$$
and (ii) holds. We have thus proved the equivalence of (i),(ii) and (v).

Let us prove that (iii) $\rightarrow $ (ii). 
We assume that (iii) holds. Let $r>0$. Then on
the event $\{\z^{B_r}\not =0\}$, there exists $s\in(0,\sigma)$
such that $\wh W_s\in \d B_r$. It follows that
$$\N_0\Big(\sup_{s\in(0,\sigma)}|\wh W_s|\geq r\Big)
\geq \N_0(\z^{B_r}\not =0).$$
Let $T_1=\inf\{s\geq 0:L^0_s>1\}$ (in agreement with the notation
of Chapter 1). The path continuity
of $\wh W_s$ ensures that $\bP_0$ a.s. there are only finitely many
excursions intervals $(\alpha_i,\beta_i)$ of $H_s$
away from $0$, before time $T_1$, such that
$$\sup_{s\in(\alpha_i,\beta_i)}|\wh W_s|\geq r.$$
On the other hand, excursion theory implies that the 
number of such intervals is Poisson with parameter
$$\N_0\Big(\sup_{s\in(0,\sigma)}|\wh W_s|\geq r\Big).$$
We conclude that the latter quantity is finite, and so $\N_0(\z^{B_r}\not =0)<\infty$.

Note that (iv) $\rightarrow $ (iii).
Thus, to complete the proof of Theorem \ref{compactness-Brownian},
it remains to verify that (ii) $\rightarrow $ (iv).
From now on  until the end of the proof, we assume that (ii) holds.

We use the following simple lemma.

\begin{lemma}
\label{compactech}
Let $D$ be a domain in $\R^d$ containing $0$, and let
$$S=\inf\{s\geq 0: W_s(t)\notin D \hbox{ for some } t\in[0,H_s)\}.$$
Then $\N_0(\z^D\not =0)\geq \N_0(S<\infty)$.
\end{lemma}

\proof By excursion theory, we have
$$\bP_0[S\geq T_1]=\exp(-\N_0(S<\infty)).$$
Then, let $\tau$ be as previously the exit time from $D$. If there exists $s<T_1$ such
that
$\tau(W_s)<H_s$, then the same property holds for every $s'>s$ such that $s'-s$ is
sufficiently small, by the continuity of $H$ and the snake property. Hence,
$$\{S< T_1\}\subset\{\int_0^{T_1}ds\,1_{\{\tau(W_s)<H_s\}}>0\},\qquad \bP_0 \hbox{ a.e.}$$
It follows that
$$\bP_0[S\geq T_1]\geq \bP_0\Big[\int_0^{T_1}ds\,1_{\{\tau(W_s)<H_s\}}=0\Big]
=\bP_0[L^D_{T_1}=0],$$
where the second equality is a consequence of the formula 
$$L^D_{T_1}=\ell^D(\int_0^{T_1}ds\,1_{\{\tau(W_s)<H_s\}}),$$ 
together with the
fact that $\ell^D(s)>0$ for every $s>0$, a.s.

Using again excursion theory and the construction of the exit measure 
under $\N_0$, we get
$$\bP_0[S\geq T_1]\geq \bP_0[L^D_{T_1}=0]=\exp(-\N_0(\z^D\not =0)).$$
By comparing with the first formula of the proof, we get the desired inequality. \cq

\smallskip
Let $\varepsilon>0$. We specialize the previous lemma to the case $D=B_\varepsilon$
and write $S=S^\varepsilon_1$. Then, for every $r>0$,
$$\bP_0[S^\varepsilon_1\geq T_r]=\exp(-r\N_0(S^\varepsilon_1<\infty))
\geq \exp(-r\N_0(\z^{B_\varepsilon}\not =0)).$$
From (ii), it follows that $S^\varepsilon_1>0$, $\bP_0$ a.e. Also note that
$S^\varepsilon_1$ is a stopping time of the filtration $({\cal F}_{s+})$ and that
$\{W_{S^\varepsilon_1}(t):0\leq t<H_{S^\varepsilon_1}\}\subset \bar B_\varepsilon$
(if this inclusion were not true, the snake property would contradict the
definition of $S^\varepsilon_1$).

Recall the notation $m_H(s,s')=\inf_{[s,s']}H_r$ for $s\leq s'$. We define 
inductively a sequence $(S^\varepsilon_n)_{n\geq 1}$ of stopping times
(for the filtration $({\cal F}_{s+})$) by setting
$$S^\varepsilon_{n+1}=\inf\{s>S^\varepsilon_n:
|W_s(t)-W_s(t\wedge m_H(S^\varepsilon_n,s))|>\varepsilon \hbox{ for some }
t\in[0,H_s)\}.$$
At this point we need another lemma.

\begin{lemma}
\label{compactMarkov}
Let $T$ be a stopping time of the filtration $({\cal F}_{s+})$, such that
$T<\infty$, $\bP_0$ a.s. For every $s\geq 0$,
define a killed path $\wt W_s$ with lifetime $\wt H_s$ by setting
$$\wt W_s=W_s(m_H(T,T+s)+t)-W_s(m_H(T,T+s)),\quad
0\leq t<\wt H_s:=H_{T+s}-m_H(T,T+s)$$
with the convention that $\wt W_s=0$ if $\wt H_s=0$. Then the process
$(\wt W_s,s\geq 0)$ is independent of ${\cal F}_{T+}$ and
has the same distribution as $(W_s,s\geq 0)$ under $\bP_0$.
\end{lemma}

This lemma follows from the strong Markov property of the L\'evy snake,
together with Lemma \ref{subexcursions}. The translation invariance of
the spatial motion is of course crucial here.

As a consequence of the preceding lemma, we obtain that the
random variables $S^\varepsilon_1,S^\varepsilon_2-S^\varepsilon_1,\ldots,
S^\varepsilon_{n+1}-S^\varepsilon_n,\ldots$ are independent and identically
distributed. Recall that these variables are positive a.s. Also observe that
$$\{W_{S^\varepsilon_{n+1}}(t)-W_{S^\varepsilon_{n+1}}(t\wedge
m_H(S^\varepsilon_n,S^\varepsilon_{n+1})):t\in[0,H_{S^\varepsilon_{n+1}})\}
\subset \bar B_\varepsilon,$$
by the same argument as used previously for $\{W_{S^\varepsilon_1}(t):
t\in[0,H_{S^\varepsilon_1})\}$.

Let $a>0$. We claim that $\N_0$ a.s. we can choose $\delta_1>0$ small enough so
that, for every $s,s'\in [0,T_a]$ such that $s\leq s'\leq s+\delta_1$,
\begin{equation} 
\label{claimcompact}
|W_{s'}(t)-W_{s'}(m_H(s,s')\wedge t)|\leq 3\varepsilon,\qquad
\hbox{for every }t\in[0,H_{s'}).
\end{equation}

Let us verify that the claim holds if we take
$$\delta_1=\inf\{S^\varepsilon_{n+1}-S^\varepsilon_n\,;\,n\geq 1,\,S^\varepsilon_n\leq T_a\}
>0.$$
Consider $s,s'\in [0,T_a]$ with $s\leq s'\leq s+\delta_1$. Then two
cases may occur. 

Either $s,s'$ belong to the same
interval $[S^\varepsilon_n,S^\varepsilon_{n+1}]$. Then, from the
definition of $S^\varepsilon_{n+1}$ we know that
\begin{equation} 
\label{claimtec1}
|W_{s'}(t)-W_{s'}(t\wedge m_H(S^\varepsilon_n,s'))|\leq\varepsilon \hbox{ for every }
t\in[0,H_{s'}).
\end{equation}
Since $m_H(s,s')\geq m_H(S^\varepsilon_n,s')$ we can replace $t$ by $t\wedge m(s,s')$
to get
$$|W_{s'}(t\wedge m_H(s,s'))-W_{s'}(t\wedge m_H(S^\varepsilon_n,s'))|\leq\varepsilon
\hbox{ for every } t\in[0,H_{s'}),$$
and our claim (\ref{claimcompact}) follows by combining this bound with the 
previous one.

Then we need to consider the case where $s\in
[S^\varepsilon_{n-1},S^\varepsilon_{n}]$ and $s'\in
(S^\varepsilon_n,S^\varepsilon_{n+1}]$ for some $n\geq 1$
(by convention $S^\varepsilon_0=0$). 
If $m_H(s,s')=m_H(S^\varepsilon_n,s')$, then the same argument
as in the first case goes through. Therefore we 
can assume that $m_H(s,s')<m_H(S^\varepsilon_n,s')$, which implies
$m_H(S^\varepsilon_{n-1},S^\varepsilon_n)<m_H(S^\varepsilon_n,s')$.
Note that the bound (\ref{claimtec1}) still holds. 
We also know that
\begin{equation} 
\label{claimtec2}
|W_{S^\varepsilon_n}(t)-W_{S^\varepsilon_n}(t\wedge
m_H(S^\varepsilon_{n-1},S^\varepsilon_n))|\leq\varepsilon \hbox{ for every }
t\in[0,H_{S^\varepsilon_n}).
\end{equation}
We replace $t$ by $t\wedge m_H(S^\varepsilon_n,s')$ in this bound, and note that
$W_{S^\varepsilon_n}(t\wedge m_H(S^\varepsilon_n,s'))=W_{s'}(t\wedge
m_H(S^\varepsilon_n,s'))$ for every $t\in[0,H_{S^\varepsilon_n}\wedge H_{s'})$,
by the snake property. It follows that
\begin{equation} 
\label{claimtec3}
|W_{s'}(t\wedge m_H(S^\varepsilon_n,s'))-W_{S^\varepsilon_n}(t\wedge
m_H(S^\varepsilon_{n-1},S^\varepsilon_n))|\leq\varepsilon \hbox{ for every }
t\in[0,H_{s'}).
\end{equation}
Similarly, we can replace $t$ by $t\wedge m_H(s,s')$ in (\ref{claimtec2}),
using again the snake property to write $W_{S^\varepsilon_n}(t\wedge
m_H(s,s'))=W_{s'}(t\wedge m_H(s,s'))$ (note that
$m_H(S^\varepsilon_{n-1},S^\varepsilon_n)\leq m_H(s,s')<m_H(S^\varepsilon_n,s')$). It
follows that
\begin{equation} 
\label{claimtec4}
|W_{s'}(t\wedge m_H(s,s'))-W_{S^\varepsilon_n}(t\wedge
m_H(S^\varepsilon_{n-1},S^\varepsilon_n))|\leq\varepsilon \hbox{ for every }
t\in[0,H_{s'}).
\end{equation}
Our claim (\ref{claimcompact}) is now a consequence of (\ref{claimtec1}),
(\ref{claimtec3}) and (\ref{claimtec4}).

We can already derive from (\ref{claimcompact}) the fact that 
the left limit $\wh W_s$ exists for every $s\in [0,T_a]$, $\bP_0$ a.s.
We know that this left limit exists for every rational $s\in[0,T_a]$,
$\bP_0$ a.s. Let $s\in(0,T_a]$, and let $s_n$ be a sequence of
rationals increasing to $s$. Then the sequence $m_H(s_n,s)$ 
also increases to $H_s$. If $m_H(s_n,s)=H_s$ for some $n$, then
the snake property shows that $W_s(t)=W_{s_n}(t)$ for every 
$t\in[0,H_s)$ and the existence of $\wh W_s$ is an immediate consequence.
Otherwise, (\ref{claimcompact}) shows that for $n$ large enough,
$$\sup_{t\in[0,H_s)}|W_s(t)-W_s(t\wedge m_H(s_n,s))|\leq 3\varepsilon$$
and by applying this to a sequence of values of $\varepsilon$
tending to $0$ we also get the existence of $\wh W_s$. 

We finally use a time-reversal argument. From Corollary \ref{reversal}, we know that
the processes $(H_{t\wedge T_a},t\geq 0)$ and $(H_{(T_a-t)^+},t\geq 0)$
have the same distribution. By considering the
conditional law of $W$ knowing $H$, we immediately obtain that
the processes $(W_{t\wedge T_a},t\geq 0)$ and $(W_{(T_a-t)^+},t\geq 0)$
also have the same distribution. Thanks to this observation and the
preceding claim, we get that $\bP_0$ a.s. there exists $\delta_2>0$ such that
for every $s,s'\in [0,T_a]$ with $s\leq s'\leq s+\delta_2$,
\begin{equation} 
\label{claimcompact2}
|W_{s}(t)-W_{s}(m_H(s,s')\wedge t)|\leq 3\varepsilon,\qquad
\hbox{for every }t\in[0,H_{s}).
\end{equation}
To complete the proof, note that the snake property
implies that 
$$W^*_s(m_H(s,s'))=W^*_{s'}(m_H(s,s')),$$ 
using a continuity
argument in the case $m_H(s,s')=H_s\wedge H_{s'}$. Thus,
if $s\leq s'\leq T_a$ and $s'-s\leq \delta_1\wedge \delta_2$,
\ba
&&\sup_{t\geq 0} |W^*_s(t\wedge H_s)-W^*_{s'}(t\wedge H_{s'})|\\
&&\ \leq\sup_{t\in[0,H_s]}|W^*_s(t)-W^*_s(t\wedge m_H(s,s'))|
+\sup_{t\in[0,H_{s'}]}|W^*_{s'}(t)-W^*_{s'}(t\wedge m_H(s,s'))|\\
&&\ =\sup_{t\in[0,H_s)}|W_s(t)-W_s(t\wedge m_H(s,s'))|
+\sup_{t\in[0,H_{s'})}|W_{s'}(t)-W_{s'}(t\wedge m_H(s,s'))|\\
&&\ \leq 6\varepsilon\,.
\ea
This gives the continuity of the mapping $s\rightarrow  W^*_s$ with respect to the distance
$d^*$, and completes the proof of (iv). \cq

\section{The law of the L\'evy snake at a first exit time}

Our goal in this section is to give explicit formulas 
for the law of the L\'evy snake at its first exit time from
a domain. We keep assuming that $H$ has continuous sample paths
and in addition we suppose that the process $W^*$ 
has continuous sample paths with respect to the metric $d^*$. 
Note that the previous two sections give sufficient
conditions for this property to hold.

\smallskip
Let $D$ be an open set in $E$ and $x\in E$. We slightly
abuse notation by writing
$\tau(\w)=\inf\{t\in[0,\zeta_{\w}]:\w(t)\notin D\}$
for any stopped path $\w\in\W^*$. We also set
$$T_D=\inf\{s> 0:\tau(W_s^*)<\infty\}.$$
The continuity of $W^*$ with respect to the metric $d^*$
immediately implies that $T_D>0$, $\N_x$ a.e. or $\bP_x$ a.e.
Furthermore, on the event $\{T_D<\infty\}$ the path $W^*_{T_D}$ hits the boundary 
of $D$ exactly at its lifetime. 
The main result of this section determines the law of the
pair $(\rho_{T_D},W_{T_D})$ under $\N_x(\cdot\cap\{T_D<\infty\})$.

\smallskip
Before stating this result, we need some notation
and a preliminary lemma. For every
$y\in D$, we set 
$$u(y)=\N_y(T_D<\infty)<\infty.$$
Recall that, for every $a,b\geq 0$, we have defined
$$\gamma_{\psi }(a,b)= \left\{
\begin{array}{ll}
\left(\psi (a) -\psi ( b) \right) / (a-b)
\quad& {\rm if } \quad a\neq b, \\
\psi'(a) \quad& {\rm if } \quad a=b\; .\hfill\\ 
\end{array}
\right. $$
Note that $\gamma_\psi(a,0)=\wt\psi(a)$ (by convention $\wt \psi(0)=\psi'(0)=\alpha$). The following formulas will
be useful: For every
$a,b\geq 0$,
\begin{eqnarray}
\label{tech-hitting1}
&&\int \pi(dr)\int_0^r d\ell \,(1-e^{-a\ell-b(r-\ell)})=\gamma_\psi(a,b)-\alpha-\beta(a+b)\\
\label{tech-hitting2}
&&\int \pi(dr)\int_0^r d\ell \,e^{-a\ell}(1-e^{-b(r-\ell)})=\gamma_\psi(a,b)-\wt\psi(a)-\beta b\,.
\end{eqnarray}
The first formula is easily obtained by
observing that, if $a\neq b$,
$$\int_0^r d\ell(1-e^{-\ell a-(r-\ell)b})={1\over a-b}(r(a-b)+(e^{-ra}-e^{-rb})).$$
The second one is a consequence of the first one and the identity
$$\wt \psi(a)=\alpha+\beta a+\int\pi(dr)\int_0^r d\ell(1-e^{-a\ell}).$$
Recall from Section 3.1 the definition of the probability measures $\M_a$ on $M_f(\R_+)^2$.

\begin{lemma}
\label{lawex-tech}
{\rm (i)} Let $a>0$ and let $F$ be a nonnegative
measurable function on $M_f(\R_+)\times M_f(\R_+)\times\W$. Then,
$$\N_x\Big(\int_0^\sigma dL^a_s\,F(\rho_s,\eta_s,W_s)\Big)
=e^{-\alpha a}\int \M_a(d\mu\,d\nu)\,\Pi_x[F(\mu,\nu,(\xi_r,0\leq r<a))].$$
\noindent {\rm (ii)} Let $f,g$ be two nonnegative measurable functions on $\R_+$. Then,
$$N\Big(\int_0^\sigma dL^a_s\,\exp(-\langle \rho_s,f\rangle-\langle\eta_s,g\rangle)\Big)
=\exp\Big(-\int_0^a \gamma_\psi(f(t),g(t))\,dt\Big).$$
\end{lemma}

\proof (i) As in the proof of Lemma \ref{moment1super}, we may restrict our attention to
a function $F(\rho_s,\eta_s,W_s)=F(\rho_s,\eta_s)$. Then the desired result follows
from Corollary \ref{localinvar} in the same way as Proposition \ref{invariant-rho-eta}
was deduced from Proposition \ref{keyinv}.

\noindent (ii) By part (i) we have
$$N\Big(\int_0^\sigma dL^a_s\,\exp(-\langle \rho_s,f\rangle-\langle\eta_s,g\rangle)\Big)
=e^{-\alpha a} \int \M_a(d\mu\,d\nu)\,\exp(-\langle \mu,f\rangle-\langle\nu,g\rangle).$$
From the definition of $\M_a$ this is equal to
$$\exp\Big(-\alpha a-\beta\int_0^a (f(t)+g(t))dt-\int_0^a dt\int \pi(dr)
\int_0^r d\ell(1-e^{-\ell f(t)-(r-\ell)g(t)})\Big).$$
The stated result now follows from (\ref{tech-hitting1}).
\cq

\begin{theorem}
\label{law-exit}
Assume that $u(x)>0$.
Let $a>0$, let $F$ be a nonnegative measurable function 
on $\W^*_x$ and let $g$ be a nonnegative measurable function 
on $\R_+$ with support contained in $[0,a]$. Then
\begin{eqnarray}
\label{law-hitting1}
&&\N_x\Big(1_{\{T_D<\infty\}}1_{\{a<H_{T_D}\}}F(W_{T_D}(t),0\leq t\leq a)\exp(-\langle \rho_{T_D},g
\rangle)\Big)\nonumber\\
&&\quad=\Pi_x\Big[1_{\{a< \tau\}}u(\xi_a)F(\xi_r,0\leq r\leq a)\exp\Big(-\int_0^a
\gamma_\psi(u(\xi_r),g(r))dr\Big)\Big].
\end{eqnarray}
Alternatively, the law of $W_{T_D}$ under $\N_x(\cdot\cap\{T_D<\infty\})$ is characterized by:
\begin{eqnarray}
\label{law-hitting2}
&&\N_x\Big(1_{\{T_D<\infty\}}1_{\{a<H_{T_D}\}}F(W_{T_D}(t),0\leq t\leq a)\Big)\nonumber\\
&&\quad=\Pi_x\Big[1_{\{a< \tau\}}u(\xi_a)F(\xi_r,0\leq r\leq a)\exp\Big(-\int_0^a
\wt\psi(u(\xi_r))dr\Big)\Big],
\end{eqnarray}
and the conditional law of $\rho_{T_D}$ knowing $W_{T_D}$
is the law of
$$\beta 1_{[0,H_{T_D}]}(r)\,dr+\sum_{i\in I} (v_i-\ell_i)\,\delta_{r_i}$$
where $\sum \delta_{(r_i,v_i,\ell_i)}$ is a Poisson point measure on $\R_+^3$
with intensity 
\begin{equation} 
\label{law-hitting3}
1_{[0,H_{T_D}]}(r)1_{[0,v]}(\ell)e^{-\ell u(W_{T_D}(r))}\,dr\,\pi(dv)d\ell.
\end{equation}
\end{theorem}

\proof We will rely on results obtained in Section 4.2 above. As in subsection 4.2.2,
we denote by $(\rho^i,W^i)$, $i\in I$ the ``excursions'' of the L\'evy snake
above height $a$. We let $(\alpha_i,\beta_i)$ be the time interval
corresponding to the excursion $(\rho^i,W^i)$ and 
$\ell^i=L^a_{\alpha_i}$~.
We also use the obvious notation
$$T_D(W^i)=\inf\{s\geq 0:\tau(W^{i*}_s)<\infty\}.$$ 
For every $s\geq 0$, set
$$G_s=1_{\{s<T_D\}}F(W^*_s)\exp(-\langle \rho_s,g\rangle).$$
Then it is easy to verify that
\begin{equation} 
\label{lawex1}
\sum_{i\in I} G_{\alpha_i}\,1_{\{T_D(W^i)<\infty\}}
=1_{\{T_D<\infty\}}1_{\{a<H_{T_D}\}}F(W_{T_D}(t),0\leq t\leq a)\exp(-\langle \rho_{T_D},g
\rangle).
\end{equation}
In fact, the sum in the left side contains at most one nonzero term, and exactly
one iff $T_D<\infty$ and $a<H_{T_D}$. On this event, $T_D$ belongs to one
excursion interval above height $a$, say $(\alpha_j,\beta_j)$, and then the
restriction of $\rho_{T_D}$ to $[0,a]$ coincides with $\rho_{\alpha_j}$ (see
the second step of the proof of Proposition \ref{reflec}), whereas the snake
property ensures that the paths $W_{T_D}$ and $W_{\alpha_i}$ are the same
over $[0,a)$. Our claim (\ref{lawex1}) follows.

Recall the notation $\wt W,\wt \rho,\gamma^a(\ell)$ introduced in the proof
of Proposition \ref{excursions}. The proof of this proposition shows that
conditionally on the $\sigma$-field $\e_a$, the point measure
$$\sum_{i\in I} \delta_{(\ell^i,\rho^i,W^i)}$$
is Poisson with intensity
$$1_{[0,L^a_\sigma]}(\ell)\,d\ell\,\N_{\hat{\tilde W}_{\gamma^a(\ell)}}(d\rho dW).$$
Note that the statement of Proposition \ref{excursions} is slightly
weaker than this, but the preceding assertion follows readily from the proof.

We now claim that we can find a deterministic function $\Delta$
and an $\e_a$-measurable random variable $Z$ such that, for every $j\in I$, we have
\begin{equation} 
\label{lawex2}
G_{\alpha_j}=\Delta(Z,\ell^j,(\ell^i,W^i)_{i\in I}).
\end{equation}
Precisely, this relation holds if we take for every $\ell \geq 0$,
\ba\Delta(Z,\ell,(\ell^i,W^i)_{i\in I})
&=&\Big(\prod_{i:\ell^i<\ell}1_{\{T_D(W^i)=\infty\}}\Big)\ 1_{\{\tilde W_r(t)\in D,\forall
\,r\in[0,\gamma^a(\ell)],t\in[0,a]\}}\\
&&\qquad\times F(\wt W^*_{\gamma^a(\ell)})
\exp(-\langle \wt\rho_{\gamma^a(\ell)},g\rangle).
\ea
Note that the right side of the last formula depends on $\ell$, on $(\ell^i,W^i)_{i\in I}$
and on the triple $(\wt W,\wt \rho,\gamma^a)$ which is $\e_a$-measurable, and thus
can be written in the form of the left side. Then, to justify (\ref{lawex2}), note that
$$1_{\{\alpha_j<T_D\}}=\Big(\prod_{i:\ell^i<\ell^j}1_{\{T_D(W^i)=\infty\}}\Big)
1_{\{\tilde W_r(t)\in D,\forall
\,r\in[0,\gamma^a(\ell^j)],t\in[0,a]\}},$$
since $\gamma^a(\ell^j)=\int_0^{\alpha_j}dr\,1_{\{H_r\leq a\}}$ as observed in the proof
of Proposition \ref{excursions}. The latter proof also yields the identities
$$W_{\alpha_j}=W_{\beta_j}=\wt W_{\gamma^a(\ell^j)}\ ,\qquad
\rho_{\alpha_j}=\rho_{\beta_j}=\wt \rho_{\gamma^a(\ell^j)}$$
from which (\ref{lawex2}) follows.

Then, by an application of Lemma \ref{Poisson-Palm} to the point measure
$\sum_{i\in I} \delta_{(\ell^i,\rho^i,W^i)}$, which is Poisson
conditional on $\e_a$, we have
\ba
&&\N_x\Big(\sum_{j\in I} G_{\alpha_j}\,1_{\{T_D(W^j)<\infty\}}\,\Big|\,\e_a\Big)\\
&&\qquad=\N_x\Big(\sum_{j\in I} \Delta(Z,\ell^j,(\ell^i,W^i)_{i\in I})\,1_{\{T_D(W^j)<\infty\}}
\,\Big|\,\e_a\Big)\\
&&\qquad=\N_x\Big(\int_0^{L^a_\sigma} d\ell\int
\N_{\hat{\tilde W}_{\gamma^a(\ell)}}(d\rho' dW')\,1_{\{T_D(W')<\infty\}}
\Delta(Z,\ell,(\ell^i,W^i)_{i\in I})\,\Big|\,\e_a\Big)\\
&&\qquad=\N_x\Big(\int_0^{L^a_\sigma} d\ell\, \Delta(Z,\ell,(\ell^i,W^i)_{i\in I})
\,\N_{\hat{\tilde W}_{\gamma^a(\ell)}}(T_D<\infty)
\,\Big|\,\e_a\Big).
\ea
Now use the definition of $\Delta(Z,\ell,(\ell^i,W^i)_{i\in I})$ to get
\begin{eqnarray}
\label{lawex3}
&&\N_x\Big(\sum_{j\in I} G_{\alpha_j}\,1_{\{T_D(W^j)<\infty\}}\Big)\nonumber\\
&&\qquad
=\N_x\Big(\int_0^{L^a_\sigma} d\ell\, u(\wh{\wt W}_{\gamma^a(\ell)})\,
\Big(\prod_{i:\ell^i<\ell}1_{\{T_D(W^i)=\infty\}}\Big)\ 1_{\{T_D(\tilde W)>\gamma^a(\ell)\}}\nonumber\\
&&\qquad\qquad\qquad\times F(\wt
W^*_{\gamma^a(\ell)})
\exp(-\langle \wt\rho_{\gamma^a(\ell)},g\rangle)\Big)\nonumber\\
&&\qquad
=\N_x\Big(\int_0^\sigma dL^a_s\,u(\wh W_s)\,F(W^*_s)\,\exp(-\langle \rho_s,g\rangle)
\,1_{\{s<T_D\}}\Big).
\end{eqnarray}
The last equality is justified by 
the change of variables $\ell=L^a_s$ and the fact that $dL^a_s$ a.e.,
$$\wt W_{\gamma^a(L^a_s)}=\wt W_{\tilde A^a_s}=W_s\ ,\qquad \wt \rho_{\gamma^a(L^a_s)}=\wt \rho_{\tilde
A^a_s}=\rho_s,$$
(where $\wt A^a_s=\int_0^s dr\,1_{\{H_r\leq a\}}$ as previously) and similarly, $dL^a_s$ a.e.,
\ba
\Big(\prod_{i:\ell^i<L^a_s}1_{\{T_D(W^i)=\infty\}}\Big)\ 1_{\{T_D(\wt W)>\gamma^a(L^a_s)\}}
&=&1_{\{W^*_r(t)\in D,\,\forall r\leq s,\;t\in[a,H_r]\}}\,1_{\{T_D(\wt W)>\tilde
A^a_s\}}\\ &=&1_{\{s<T_D\}}.
\ea

To evaluate the right side of (\ref{lawex3}), we use a duality argument. 
It follows from Corollary \ref{reversal} and the construction of the L\'evy snake
that the triples
$$(\rho_s,L^a_s,W_s;0\leq s\leq \sigma)$$
and
$$(\eta_{(\sigma-s)-},L^a_\sigma-L^a_{\sigma-s},W_{\sigma-s};0\leq s\leq \sigma)$$
have the same distribution under $\N_x$. From this we get
\begin{eqnarray}
\label{lawex4}
&&\N_x\Big(\int_0^\sigma dL^a_s\,u(\wh W_s)\,F(W^*_s)\,\exp(-\langle \rho_s,g\rangle)
\,1_{\{s<T_D\}}\Big)\nonumber\\
&&=\N_x\Big(\int_0^\sigma dL^a_s\,u(\wh W_s)\,F(W^*_s)\,\exp(-\langle \eta_s,g\rangle)
\,1_{\{\tau(W^*_r)=\infty,\,\forall r\geq s\}}\Big).
\end{eqnarray}
Now we can use the strong Markov property of the L\'evy snake
(as in the second step of the proof of Proposition \ref{super-bis}),
and then Lemma \ref{subexcursions}, to get
\begin{eqnarray}
\label{lawex5}
&&\N_x\Big(\int_0^\sigma dL^a_s\,u(\wh W_s)\,F(W^*_s)\,\exp(-\langle \eta_s,g\rangle)
\,1_{\{\tau(W^*_r)=\infty,\,\forall r\geq s\}}\Big)\nonumber\\
&&=\N_x\Big(\int_0^\sigma dL^a_s\,u(\wh W_s)\,F(W^*_s)\,\exp(-\langle
\eta_s,g\rangle)\,1_{\{\tau(W^*_s) =\infty\}}
\,\bP^*_{\rho_s,W_s}[T_D=\infty]\Big)\nonumber\\
&&=\N_x\Big(\int_0^\sigma dL^a_s\,u(\wh W_s)\,F(W^*_s)\nonumber\\
&&\hspace{15mm}\times\exp(-\langle
\eta_s,g\rangle)\,1_{\{\tau(W^*_s) =\infty\}}
\,\exp\Big(-\int \rho_s(dt)\,u(W_s(t))\Big)\Big).
\end{eqnarray}
Finally, we use Lemma \ref{lawex-tech}
to write
\ba
&&\hspace{-8mm}\N_x\Big(\int_0^\sigma dL^a_s\,u(\wh W_s)\,F(W^*_s)\,1_{\{\tau(W^*_s)
=\infty\}}\,\exp(-\langle \eta_s,g\rangle)
\,\exp\Big(-\int \rho_s(dt)\,u(W_s(t))\Big)\Big)\\
&&\hspace{-8mm}=e^{-\alpha a}\int \M_a(d\mu d\nu)\,\Pi_x\Big[1_{\{a<\tau\}}u(\xi_a)F(\xi_r,0\leq r<a)
\,e^{-<\nu,g>}\exp(-\int \mu(dr)\,u(\xi_r))\Big]\\
&&\hspace{-8mm}=e^{-\alpha a}\Pi_x\Big[1_{\{a<\tau\}}u(\xi_a)F(\xi_r,0\leq r<a)\,
\int \M_a(d\mu d\nu)\,e^{-<\nu,g>}\exp(-\int \mu(dr)\,u(\xi_r))\Big]\\
&&\hspace{-8mm}=\Pi_x\Big[1_{\{a<\tau\}}u(\xi_a)F(\xi_r,0\leq r\leq a)\exp\Big(-\int_0^a
\gamma_\psi(u(\xi_r),g(r))dr\Big)\Big].
\ea
Formula (\ref{law-hitting1}) follows by combining this equality with
(\ref{lawex1}), (\ref{lawex3}), (\ref{lawex4}) and (\ref{lawex5}). 

Formula (\ref{law-hitting2}) is the special case $g=0$ in (\ref{law-hitting1}).
To prove the last assertion, let $\zeta(W_{T_D},d\mu)$ be the law of the 
random measure 
$$\beta 1_{[0,H_{T_D}]}(r)\,dr+\sum_{i\in I} (v_i-\ell_i)\,\delta_{r_i}$$
where $\sum \delta_{(r_i,v_i,\ell_i)}$ is a Poisson point measure on $\R_+^3$
with intensity 
given by formula (\ref{law-hitting3}). Then, for every $a>0$, we can 
use (\ref{law-hitting2}) to compute
\ba
&&\N_x\Big(F(W_{T_D}(r),0\leq r\leq a)\,1_{\{T_D<\infty\}}1_{\{a<H_{T_D}\}}\int \zeta(W_{T_D},d\mu)\,e^{-\langle
\mu,g\rangle}\Big)\\
&&\ =\Pi_x\Big[1_{\{a<\tau\}}\,F(\xi_r,0\leq r\leq a)u(\xi_a)\,\exp\Big(-\int_0^a \wt\psi(u(\xi_r))\,dr\Big)\\
&&\ \times \exp\Big(-\beta\int_0^a dr\,g(r)-\int_0^a dr\int\pi(dv)\int_0^v d\ell
e^{-\ell u(\xi_r)}(1-e^{-(v-\ell)g(r)})\Big)\Big]\\
&&\ =\Pi_x\Big[1_{\{a<\tau\}}\,F(\xi_r,0\leq r\leq a)u(\xi_a)\,\exp\Big(-\int_0^a
\gamma_\psi(u(\xi_r),g(r))\,dr\Big)\Big],
\ea
using (\ref{tech-hitting2}) in the last equality. 

Set $\N^D_x=\N_x(\cdot\mid T_D<\infty)$ to simplify notation.
By comparing with (\ref{law-hitting1}), we see that for any nonnegative measurable function
$g$ with support in $[0,a]$, we have
$$\N^D_x[e^{-\langle \rho_{T_D},g\rangle}\mid W_{T_D}]=\int \zeta(W_{T_D},d\mu)\,e^{-\langle \mu,g\rangle},$$
a.s. on the set $\{H_{T_D}>a\}$. This is enough to conclude that $\zeta(W_{T_D},d\mu)$
is the conditional distribution of $\rho_{T_D}$ knowing $W_{T_D}$, provided that we already
know that $\rho_{T_D}(\{H_{T_D}\})=0$ a.s. The latter fact however is a simple
consequence of (\ref{invar-snake}). This completes the proof of the theorem.
\cq

\bigskip
{\bf The case of Brownian motion.} Suppose that the spatial motion $\xi$
is $d$-dimensional Brownian motion and that $D$ is a domain in $\R^d$. Then,
it is easy to see that the function $u(x)=\N_x(T_D<\infty)$, $x\in D$ is
of class $C^2$ and solves ${1\over 2}\Delta u=\psi(u)$. In the context 
of superprocesses, this was observed by Dynkin \cite{Dy0}.
We may argue as follows. First note that the set of nonnegative
solutions of ${1\over 2}\Delta u=\psi(u)$ in a domain is closed under pointwise
convergence (for a probabilistic proof, reproduce the arguments of the proof
of Proposition 9 (iii) in \cite{LG99}). Then let $(D_n)$ be a sequence of bounded
regular subdomains of
$D$, such that $\bar D_n\subset D_{n+1}$ and $D=\lim\uparrow D_n$. For every $n\geq 0$, set
$$v_n(x)=\N_x(\z^{D_n}\not =0)\ ,\quad u_n(x)=\N_x(T_{D_n}<\infty)\ ,\quad x\in D_n.$$
From the properties of the exit measure, it is immediate to see that $v_n\leq u_n$.
On the other hand, by writing
$$v_n(x)=\lim_{\lambda\uparrow \infty}\uparrow\,\N_x(1-\exp-\lambda\langle \z^{D_n},1\rangle),$$
we deduce from Proposition \ref{DirichletPDE} and the stability of the set of
nonnegative solutions under pointwise convergence that $v_n$ is of class $C^2$ 
and solves ${1\over 2}\Delta v_n=\psi(v_n)$ in $D$. Since the function
$x\la \N_x(1-\exp-\lambda\langle \z^{D_n},1\rangle)$ has boundary value $\lambda$ on 
$\partial D_n$ (Proposition \ref{DirichletPDE}), we also see that $v_n$ has boundary value
$+\infty$ on $\partial D_n$.

Then, it follows from Lemma \ref{compactech}
and our assumption $\bar D_n\subset D_{n+1}$ that $v_n(x)\geq u_{n+1}(x)$ for $x\in D_n$.
Since it is easy to see that $u_n(x)$ decreases to $u(x)$ as $n\rightarrow  \infty$, for every $x\in D$,
we conclude from the inequalities $u_{n+1}(x)\leq v_n(x)\leq u_n(x)$ that $v_n(x)$
also converges to $u(x)$ pointwise as $n\rightarrow  \infty$. Hence $u$ is a nonnegative solution
of ${1\over 2}\Delta u=\psi(u)$ in $D$. The preceding argument gives more. Let $v$ be
any nonnegative solution of $\Delta v={1\over 2}\psi(v)$ in $D$. Since 
$v_{n|\partial D_n}=+\infty$, the comparison principle (Lemma V.7 in \cite{LG99}) implies that
$v\leq v_n$ in $D_n$. By passing to the limit $n\rightarrow \infty$, we conclude that $v\leq u$. 
Hence $u$ is the maximal nonnegative solution
of ${1\over 2}\Delta u=\psi(u)$ in $D$.

Suppose that $u(x)>0$ for some $x\in D$. It is easy to see that this implies $u(y)>0$ for every $y\in D$
(use a suitable Harnack principle or a probabilistic argument relying on the fact that
$u(\xi_t)\exp(-\int_0^t \wt\psi(u(\xi_r))dr)$ is a martingale).
By applying It\^o's formula to $\log u(\xi_t)$, we see that $\Pi_x$ a.s. on $\{t<\tau\}$,
\ba
\log u(\xi_t)&=&
\log u(x)+\int_0^t {\nabla u\over u}(\xi_r)\cdot d\xi_r+{1\over 2}\int_0^t \Delta(\log u)(\xi_r)dr\\
&=&\log u(x)+\int_0^t {\nabla u\over u}(\xi_r)\cdot d\xi_r+\int_0^t\Big(\wt \psi(u(\xi_r))-{1\over 2}\Big|{\nabla
u\over u}\Big|^2(\xi_r)\Big)dr.
\ea
We can then rewrite (\ref{law-hitting2}) in the form
\ba
&&\N^D_x\Big(1_{\{T_D<\infty\}}1_{\{t<H_{T_D}\}}F(W_{T_D}(r),0\leq r\leq t)\Big)\nonumber\\
&&\qquad=\Pi_x\Big[1_{\{t< \tau\}}\exp\Big(\int_0^t {\nabla u\over u}(\xi_r)\cdot d\xi_r
-{1\over 2}\int_0^t \Big|{\nabla u\over
u}\Big|^2(\xi_r)\,dr\Big)
F(\xi_r,0\leq r\leq
t)\Big].
\ea
An application of Girsanov's theorem then shows that $W_{T_D}$ is distributed
as the solution of the stochastic differential equation
\ba
&&dx_t=dB_t+{\nabla u\over u}(x_t)dt\\
&&x_0=x
\ea
(where $B$ is a standard $d$-dimensional Brownian motion)
which can be defined up to its first hitting time of $\partial D$. See \cite{LG94}
for a discussion and another interpretation of this distribution on
paths in the case $\psi(u)=u^2$.

\section{The reduced tree in an open set}

We keep the notation and assumptions of the previous section. In particular,
we assume that $W^*$ has continuous sample paths with respect to the distance
$d^*$, $D$ is
an open set in $E$, $x\in D$, $T_D=\inf\{s> 0:\tau(W^*_s)<\infty\}$ and
$u(x)=\N_x(T_D<\infty)<\infty$. To avoid trivialities, we assume that $u(x)>0$, and we recall
the notation $\N^D_x=\N_x(\cdot\mid T_D<\infty)$. 
We will assume in addition that 
\begin{equation} 
\label{assum-reduc}
\sup_{y\in K}u(y)<\infty
\end{equation}
for every compact subset $K$ of $D$. This assumption holds in particular when 
$\xi$ is Brownian motion in $\R^d$, under the condition (\ref{keycond})
(use translation invariance and the fact that $u(0)<\infty$ when $D$ is
an open ball centered at the origin).

We also set:
$$L_D=\sup\{s\geq 0:\tau(W^*_s)<\infty\},$$
and 
$$m_D=\inf_{T_D\leq s\leq L_D} H_s.$$
As a consequence of the first lemma below, we will see that $m_D<H_{T_D}$, $\N^D_x$
a.s.

Our goal is to describe the genealogical structure of the paths $W_s$ that exit $D$,
up to their first exit time from $D$, under the probability
measure $\N^D_x$. To be more precise, all paths
$W_s$   such that $\tau(W^*_s)<\infty$ must coincide up to level
$m_D$.
At level $m_D$ there is a branching point with finitely many branches, each 
corresponding to an excursion of $W$ above level $m_D$ that hits $D^c$. In each such
excursion, the  paths $W_s$ that hit $D^c$ will be the same up to a level 
(strictly greater than $m_D$) at which there is another branching point,
and so on. 

We will describe this genealogical structure in a recursive way. We will
first derive the law of the common part to the paths $W_s$ that do exit $D$.
This common part is represented by a stopped path $W^D_0$ in $\W_x$
with lifetime $\zeta_{W^D_0}=m_D$. Then we will obtain the distribution
of the  ``number of branches'' at level $m_D$, that is the number of excursions of $W$
above height $m_D$ that hit $D^c$. Finally, we will see that conditionally
on $W^D_0$, these excursions are independent and distributed according to 
$\N_{\hat W^D_0}(\cdot\mid T_D<\infty)$. This completes our recursive description since we can
apply to each of these excursions the results obtained under $\N^D_x$.
 
Before coming to the main result of this section, we state an important lemma.   

\begin{lemma} 
\label{isolate}
The point $T_D$ is not isolate in $\{s\geq
0:\tau(W^*_s)<\infty\}$, $\N_x$ a.e. on $\{T_D<\infty\}$.
\end{lemma}

\proof We start with some preliminary observations. Let $(\mu,\w)\in\Theta_x$
be such that $\mu(\{H(\mu)\})=0$ and $\w(t)\in D$ for every $t\in[0,H(\mu))$. As an application
of Lemma \ref{subexcursions}, we have
\ba\bP^*_{\mu,\w}[T_D<\infty]&=&1-\exp-\int_{[0,H(\mu))}\N_{\w(t)}(T_D<\infty)\,\mu(dt)\\
&=&1-
\exp-\int_{[0,H(\mu))}u(\w(t))\,\mu(dt).
\ea
By the previous formula, the equality $\bP^*_{\mu,\w}[T_D=0]=1$ can only hold if
\begin{equation} 
\label{regul-cond}
\int_{[0,H(\mu))}u(\w(t))\,\mu(dt)=\infty.
\end{equation}
Conversely, condition (\ref{regul-cond}) also implies that $\bP^*_{\mu,\w}[T_D=0]=1$. 
To see this, first note that our assumption (\ref{assum-reduc}) guarantees that for
every $\varepsilon>0$,
$$\int_{[0,H(\mu)-\varepsilon]}u(\w(t))\,\mu(dt)<\infty,$$
and thus we have also under (\ref{regul-cond})
$$\int_{(H(\mu)-\varepsilon,H(\mu))}u(\w(t))\,\mu(dt)=\infty.$$
Then write $\mu_\varepsilon$ for the restriction of $\mu$ to $[0,H(\mu)-\varepsilon]$, and
set
$$S_\varepsilon=\inf\{s\geq 0:\langle\rho_s,1\rangle=\langle \mu_\varepsilon,1\rangle\}.$$
Lemma \ref{subexcursions} again implies that
$$\bP^*_{\mu,\w}[T_D\leq S_\varepsilon]=1-
\exp-\int_{(H(\mu)-\varepsilon,H(\mu))}u(\w(t))\,\mu(dt)=1.$$
Since $S_\varepsilon\downarrow 0$ as $\varepsilon\da 0$, $\bP^*_{\mu,\w}$ a.s., we get
that $\bP^*_{\mu,\w}[T_D=0]=1$, which was the desired result.

Let us prove the statement of the lemma.
Thanks to the strong Markov property, it is enough to prove
that  
$\bP^*_{\rho_{T_D},W_{T_D}}[T_D=0]=1$, $\N_x$ a.e. on $\{T_D<\infty\}$. Note that we have
$\rho_{T_D}(\{H_{T_D}\})=0$ and $W_{T_D}(t)\in D$ for every $t<H_{T_D}$, 
$\N_x$ a.e. on $\{T_D<\infty\}$. By the preceding observations, it is
enough to prove that
\begin{equation} 
\label{isolate-tech}
\int_{[0,H_{T_D})}u(W_{T_D}(t))\,\rho_{T_D}(dt)=\infty\ ,\quad\hbox{ a.e. on
}\{T_D<\infty\}.
\end{equation} 
To this end, set for every $s>0$,
$$M_s=\N_x(T_D<\infty\mid {\cal F}_s).$$
The Markov property at time $s$ shows that we have for every $s>0$, $\N_x$ a.e.,
\ba
M_s&=&1_{\{T_D\leq s\}}+1_{\{s<T_D\}}\bP^*_{\rho_s,W_s}[T_D<\infty]\\
&=&1_{\{T_D\leq s\}}+1_{\{s<T_D\}}\Big(1-\exp-\int u(W_s(t))\,\rho_s(dt)\Big)
\ea
Since the process $(\rho_s)$ is right-continuous for the variation distance on 
measures, it is easy to verify that the process $1_{\{s<T_D\}}(1-\exp-\int
u(W_s(t))\,\rho_s(dt))$ is right-continuous. Because $(M_s,s>0)$ is a 
martingale with respect to the filtration $({\cal F}_s)$, a standard result 
implies that this process also has left limits at every $s>0$, $\N_x$ a.e.
In particular the left limit at $T_D$
$$\lim_{s\uparrow T_D,s<T_D}M_s=\lim_{s\uparrow T_D,s<T_D}\Big(1-\exp-\int
u(W_s(t))\,\rho_s(dt)\Big)$$
exists $\N_x$ a.e. on $\{T_D<\infty\}$. It is not hard to verify that this limit
is equal to $1$: If $D_n=\{y\in D:{\rm dist}(y,D^c)>n^{-1}\}$ and $T_n=T_{D_n}$, we have
$T_{n}<T_D$ and $T_n\ua T_D$ on $\{T_D<\infty\}$, and
$M_{T_{n}}=\N_x(T_D<\infty\mid{\cal F}_{T_n})$ converges to $1$ as $n\rightarrow  \infty$
on the set $\{T_D<\infty\}$
because $T_D$ is measurable with respect to the $\sigma$-field $\bigvee {\cal F}_{T_n}$.

Summarizing, we have proved that
\begin{equation} 
\label{isolate-tech2}
\lim_{s\uparrow T_D,s<T_D}\int
u(W_s(t))\,\rho_s(dt)=+\infty
\end{equation}
$\N_x$ a.e. on $\{T_D<\infty\}$. Then, for every rational $a>0$, consider 
on the event $\{T_D<\infty\}\cap\{H_{T_D}>a\}$, the number $\alpha_{(a)}$
defined as the left end of the excursion interval of $H$ 
above $a$ that straddles $T_D$. As a consequence of the considerations in subsection 4.2.2, the following two
facts hold on $\{T_D<\infty\}\cap\{H_{T_D}>a\}$:
\ba
&&\rho_{\alpha_{(a)}} \hbox{ is the restriction of } \rho_{T_D} \hbox{ to }[0,a)\\
&&W_{\alpha_{(a)}}(t)=W_{T_D}(t)\ ,\quad\hbox{ for every } t\in[0,a).
\ea
Thus, we have also on the same event
$$\int
u(W_{\alpha_{(a)}}(t))\,\rho_{\alpha_{(a)}}(dt)=\int_{[0,a)} u(W_{T_D}(t))\,\rho_{T_D}(dt).$$
Now on the event $\{T_D<\infty\}$ we can pick a sequence $(a_n)$ of rationals
strictly increasing to $H_{T_D}$. We observe that $\alpha_{(a_n)}$
also converges to $T_D$ (if $S$ is the increasing limit of $\alpha_{(a_n)}$, the snake
property implies that $\wh W_S=\wh W_{T_D}\in D^c$ and so we have $S\geq T_D$,
whereas the other inequality is trivial). Therefore, using (\ref{isolate-tech2}),
$$\infty=\lim_{n\rightarrow \infty}\int
u(W_{\alpha_{(a_n)}}(t))\,\rho_{\alpha_{(a_n)}}(dt)=\lim_{n\rightarrow  \infty}\int_{[0,a_n)}
u(W_{T_D}(t))\,\rho_{T_D}(dt),$$
which yields (\ref{isolate-tech}). \cq

\medskip
Lemma \ref{isolate} implies that $T_D<L_D$, $\N^D_x$ a.s.
Since we know that $\rho_{T_D}(\{H_{T_D}\})=0$, $\N^D_x$ a.s., 
an application of the strong Markov property at time $T_D$
shows that $m_D<H_{T_D}$, $\N^D_x$ a.s. We define $W^D_0$ as the stopped path
which is the restriction of $W_{T_D}$ to $[0,m_{D}]$. Then we define the excursions
of $W$ above level $m_D$ in a way analogous to subsection 4.2.2. If
$$R_D=\sup\{s\leq T_D:H_s=m_D\}\ ,\quad S_D=\inf\{s\geq L_D:H_s=m_D\},$$
we let $(a_j,b_j)$, $j\in J$ be the connected components of the open set
$(R_D,S_D)\cap\{s\geq 0:H_s>m_D\}$. For each $j\in J$, we can then define the process
$W^{(j)}
\in C(\R_+,\W)$
by setting
$$\begin{array}{ll}
W^{(j)}_s(r)=W_{a_j+s}(m_D+r),\ \zeta_{W^{(j)}_s}=
H_{a_j+s}-m_D
\qquad&{\rm if}\ 0< s<b_j-a_j\\
W^{(j)}_s=\wh W^D_0&{\rm if}\ s=0\ {\rm or}\ s\geq b_j-a_j.
\end{array}$$
By a simple continuity argument, the set
$\{j\in J:T_D(W^{(j)})<\infty\}$
is finite a.s., and we set
$$N_D=\card\{j\in J:T_D(W^{(j)})<\infty\}.$$
We write $W^{D,1},W^{D,2},\ldots,W^{D,N_D}$
for the excursions $W^{(j)}$ such that $T_D(W^{(j)})<\infty$, listed in
chronological order.

We are now ready to state our main result.

\begin{theorem}
\label{reduced-domain}
For every $r\geq 0$, set $\theta(r)=\psi'(r)-\wt \psi(r)$. Then the law of $W^D_0$ 
is characterized by the
following formula, valid for any nonnegative measurable function $F$ on $\W^*$:
\begin{eqnarray}
\label{reduc-dom1}
&&\hspace{-12mm}\N_x(1_{\{T_D<\infty\}}F(W^D_0))\nonumber\\
&&\hspace{-12mm}=\int_0^\infty db\,\Pi_x\Big[1_{\{b<\tau\}}u(\xi_b)\,\theta(u(\xi_b))\,
\exp\Big(-\int_0^b \psi'(u(\xi_r))dr\Big) F(\xi_r,0\leq r\leq b)\Big].
\end{eqnarray}
The conditional distribution of $N_D$ knowing $W^D_0$ is given by:
\begin{equation} 
\label{reduc-dom2}
\N^D_x[r^{N_D}\mid W^D_0]=r\,{\psi'(U)-\gamma_\psi(U,(1-r)U)\over
\psi'(U)-\gamma_\psi(U,0)}\ ,\qquad 0\leq r\leq 1,
\end{equation}
where $U=u(\wh W^D_0)$. Finally, conditionally on the pair $(W^D_0,N_D)$, the
processes $W^{D,1},W^{D,2},$ $\ldots,W^{D,N_D}$ are
independent and distributed according to $\N^D_{\hat W^D_0}$.
\end{theorem}

\proof Our first objective is to compute the conditional distribution of $m_D$
knowing $W_{T_D}$.
To this end, we will apply the strong Markov property of the L\'evy snake at time $T_D$.
We have for every $b>0$
$$\N^D_x[m_D>b\mid \rho_{T_D},W_{T_D}]=\bP^*_{\rho_{T_D},W_{T_D}}[\inf_{0\leq s\leq L_D}
H_s>b].$$
By Lemma \ref{subexcursions}, the latter expression is equal to the probability
that in a Poisson point measure with intensity
$$\rho_{T_D}(dh)\,\N_{W_{T_D}(h)}(d\rho,dW)$$
there is no atom $(h_i,\rho^i,W^i)$ such that $h_i\leq b$ and $T_D(W^i)<\infty$. 
We conclude that
\begin{eqnarray}
\label{law-mini}
\N^D_x[m_D>b\mid \rho_{T_D},W_{T_D}]&=&
\exp-\int_{[0,b]}\rho_{T_D}(dh)\,\N_{W_{T_D}(h)}(T_D<\infty)\nonumber\\
&=&
\exp-\int_{[0,b]}\rho_{T_D}(dh)\,u(W_{T_D}(h)).
\end{eqnarray}
Recall that the conditional law of $\rho_{T_D}$ knowing $W_{T_D}$ is given in
Theorem \ref{law-exit}. Using this conditional distribution we see that 
$$\N^D_x[m_D>b\mid W_{T_D}]=
\exp\Big(-\beta\int_0^b da\,u(W_{T_D}(a))\Big)\,E\Big[\exp-\sum_i
(v_i-\ell_i)u(W_{T_D}(r_i))\Big],$$
where $\sum \delta_{(r_i,v_i,\ell_i)}$ is a Poisson point measure with intensity
given by (\ref{law-hitting3}). By exponential formulas for Poisson measures, we have
\ba 
&&E\Big[\exp-\sum_i
(v_i-\ell_i)u(W_{T_D}(r_i))\Big]\\
&&\quad=\exp-\int_0^b dr\int \pi(dv)\int d\ell\, e^{-\ell u(W_{T_D}(r))}(1-
e^{-(v-\ell)u(W_{T_D}(r))}).
\ea
By substituting this in the previous displayed formula, and using (\ref{tech-hitting2}),
we get
\begin{equation} 
\label{reducedom0}
\N^D_x[m_D>b\mid W_{T_D}]=
\exp\Big(-\int_0^b dr\, (\psi'(u(W_{T_D}(r)))-\wt\psi(u(W_{T_D}(r))))\Big).
\end{equation}
Hence, if $\theta(r)=\psi'(r)-\wt \psi(r)$ as in the statement of the theorem, the
conditional law of $m_D$ knowing $W_{T_D}$ has density
$$1_{[0,H_{T_D})}(b)\,\theta(u(W_{T_D}(r)))\,\exp\Big(-\int_0^b
\theta(u(W_{T_D}(r)))\,dr\Big).$$
It follows that 
\ba
&&\N_x\Big(1_{\{T_D<\infty\}}\,F(W^D_0)\Big)\\
&&\ =\N_x\Big(1_{\{T_D<\infty\}}\,F(W_{T_D}(t),0\leq t\leq m_D)\Big)\\
&&\ =\N_x\Big(1_{\{T_D<\infty\}}\,\int_0^{H_{T_D}}
db\,\theta(u(W_{T_D}(b)))\,\exp\Big(-\int_0^b
\theta(u(W_{T_D}(r)))\,dr\Big)\,\\
&&\hskip 6cm \times\ F(W_{T_D}(t),0\leq t\leq b)\Big)\\
&&\ =\int_0^\infty db\,\N_x\Big(1_{\{T_D<\infty\}}1_{\{b<H_{T_D}\}}
\theta(u(W_{T_D}(b)))\,\exp\Big(-\int_0^b
\theta(u(W_{T_D}(r)))\,dr\Big)\,\\
&&\hskip 6cm \times \ F(W_{T_D}(t),0\leq t\leq b)\Big)\\
&&\ =\int_0^\infty db\,\Pi_x\Big[1_{\{b<\tau\}}u(\xi_b)\theta(u(\xi_b))
\exp(-\int_0^b\psi'(u(\xi_r))\,dr)\,F(\xi_r,0\leq r\leq b)\Big],
\ea
using (\ref{law-hitting2}) in the last equality. This gives the first assertion of the
theorem.

We now turn to the distribution of $N_D$. We use again the strong Markov property at time
$T_D$ and Lemma \ref{subexcursions} to analyse the conditional distribution of the pair
$(m_D,N_D)$ knowing $(\rho_{T_D},W_{T_D})$. Conditional on $(\rho_{T_D},W_{T_D})$,
let $\sum \delta_{(h_i,\rho^i,W^i)}$ be a Poisson point measure with intensity
$$\rho_{T_D}(dh)\,\N_{W_{T_D}(h)}(d\rho dW).$$
Set
\ba
&&m=\inf\{h_i:T_D(W^i)<\infty\},\\
&&M=\card\{i:h_i=m\hbox{ and } T_D(W^i)<\infty\}.
\ea
Then Lemma \ref{subexcursions} and the strong Markov property show that the pairs $(m,1+M)$ and
$(m_D,N_D)$ have the same distribution conditional on $(\rho_{T_D},W_{T_D})$. Recall that
the conditional distribution of $m_D$ (or of $m$) is given by (\ref{law-mini}).

Now note that:

\begin{description}
\item{$\bullet$} If $\rho_{T_D}(\{m\})=0$, then $M=1$ because the Poisson measure
$\sum \delta_{(h_i,\rho^i,W^i)}$ cannot have two atoms at a level $h$ such that 
$\rho_{T_D}(\{h\})=0$.
\item{$\bullet$} Let $b\geq 0$ be such that $\rho_{T_D}(\{b\})>0$. The event $\{m=b\}$
occurs with probability
$$\exp\Big(-\int_{[0,b)}
\rho_{T_D}(dh)\,u(W_{T_D}(h))\Big)\,\Big(1-e^{-\rho_{T_D}(\{b\})u(W_{T_D}(b))}\Big).$$
Conditionally on this event, $M$ is distributed as a Poisson variable with parameter
$c=\rho_{T_D}(\{b\})u(W_{T_D}(b))$ and conditioned to be (strictly) positive, whose generating
function is
$${e^{-c(1-r)}-e^{-c}\over 1-e^{-c}}.$$
\end{description}

Since the continuous part of the law of $m$ has density
$$\beta\,u(W_{T_D}(b))\,\exp\Big(-\int_{[0,b)}\rho_{T_D}(dh)\,u(W_{T_D}(h))\Big)$$
we get by combining the previous two cases that
\begin{eqnarray}
\label{reducedom1}
&&\N^D_x[f(m_D)r^{N_D}\mid \rho_{T_D},W_{T_D}]
\nonumber\\
&&\quad =\beta
r^2\int_0^{H_{T_D}}db\,f(b)\,u(W_{T_D}(b))\,
\exp\Big(-\int_{[0,b)}\rho_{T_D}(dh)\,u(W_{T_D}(h))\Big)\nonumber\\
&&\quad\ + r\sum_{\rho_{T_D}(\{b\})>0}f(b)
\exp\Big(-\int_{[0,b)}\rho_{T_D}(dh)\,u(W_{T_D}(h))\Big)\,\nonumber\\
&&\qquad
\times\ \Big(e^{-\rho_{T_D}(\{b\})u(W_{T_D}(b))(1-r)}-e^{-\rho_{T_D}(\{b\})u(W_{T_D}(b))}\Big).
\end{eqnarray}
We now need to integrate the right side of (\ref{reducedom1}) with respect to the conditional
law of $\rho_{T_D}$ knowing $W_{T_D}$. We get
$$\N^D_x[f(m_D)r^{N_D}\mid W_{T_D}]=A_1+A_2$$
where
\begin{eqnarray}
\label{reducedom2}
A_1&=&\beta
r^2\N^D_x\Big[\int_0^{H_{T_D}}db\,f(b)\,u(W_{T_D}(b))\,
\exp\Big(-\int_{[0,b)}\rho_{T_D}(dh)\,u(W_{T_D}(h))\Big)\,\Big|\,W_{T_D}\Big]\nonumber\\
&=&\beta
r^2\int_0^{H_{T_D}}db\,f(b)\,u(W_{T_D}(b))\,
\exp\Big(-\int_{[0,b)}\theta(u(W_{T_D}(h)))\,dh\Big),
\end{eqnarray}
by the calculation used in the proof of (\ref{reducedom0}).
We then compute $A_2$. To this end, let ${\cal N}(dbdvd\ell)$
be (conditionally on $W_{T_D}$) a Poisson point measure in $\R_+^3$ with intensity
$$1_{[0,H_{T_D}]}(b)1_{[0,v]}(\ell)e^{-\ell u(W_{T_D}(b))}\,db\pi(dv)d\ell.$$
From Theorem \ref{law-exit}, we get
\ba
A_2&=&r\,\N^D_x\Big[\int {\cal N}(dbdvd\ell)\,f(b)\,\exp\Big(-\int_{\{a<b\}}{\cal N}(dadv'd\ell')
(v'-\ell')u(W_{T_D}(a))\Big)\\
&&\hskip 2cm\times\
(e^{-(v-\ell)u(W_{T_D}(b))(1-r)}-e^{-(v-\ell)u(W_{T_D}(b))}\Big)\,\Big|\,W_{T_D}\Big]
\ea
From Lemma \ref{Poisson-Palm} and (once again) the calculation used in proving 
(\ref{reducedom0}), we arrive at
\ba
A_2&=&r\int_0^{H_{T_D}} db\,f(b)\,\exp\Big(-\int_{[0,b)}\theta(u(W_{T_D}(a)))\,da\Big)\\
&&\ \times \int\pi(dv)\int_0^v d\ell\,e^{-\ell u(W_{T_D}(b))}
\Big(e^{-(v-\ell)(1-r)u(W_{T_D}(b))}-e^{-(v-\ell)u(W_{T_D}(b))}\Big).
\ea
From (\ref{tech-hitting2}), we have
\ba
&&\int\pi(dv)\int_0^v d\ell\,e^{-\ell u(W_{T_D}(b))}
\Big(e^{-(v-\ell)(1-r)u(W_{T_D}(b))}-e^{-(v-\ell)u(W_{T_D}(b))}\Big)\\
&&\quad=
\psi'(u(W_{T_D}(b)))-\gamma_\psi(u(W_{T_D}(b)),(1-r)u(W_{T_D}(b)))-\beta r u(W_{T_D}(b))).
\ea
By substituting this identity in the previous formula for $A_2$, and then
adding the formula for $A_1$, we arrive at:
\ba
&&\hspace{-8mm}\N^D_x[f(m_D)r^{N_D}\mid W_{T_D}]\\
&&\hspace{-8mm}=r\int_0^{H_{T_D}}db\,f(b)\,\exp\Big(-\int_0^b da\,\theta(u(W_{T_D}(a)))\Big)\\
&&\qquad\qquad\times
\Big(\psi'(u(W_{T_D}(b)))-\gamma_\psi(u(W_{T_D}(b)),(1-r)u(W_{T_D}(b)))\Big)\\
&&\hspace{-8mm}=\N^D_x[f(m_D)\,r\,{\psi'(u(W_{T_D}(m_D)))-\gamma_\psi(u(W_{T_D}(m_D)),(1-r)u(W_{T_D}(m_D)))
\over \psi'(u(W_{T_D}(m_D)))-\gamma_\psi(u(W_{T_D}(m_D)),0)}\,\Big|\,W_{T_D}\Big].
\ea
In the last equality we used the conditional distribution of $m_D$ knowing $W_{T_D}$,
and the fact that $\theta(u)=\psi'(u)-\wt \psi(u)=\psi'(u)-\gamma_\psi(u,0)$.

Finally, if $U=u(W_{T_D}(m_D))=u(\wh W^D_0)$, we have obtained
$$\N^D_x[r^{N_D}\mid W^D_0]=r\,{\psi'(U)-\gamma_\psi(U,(1-r)U)\over
\psi'(U)-\gamma_\psi(U,0)},$$
which is formula (\ref{reduc-dom2}) of the theorem.

It remains to obtain the last assertion of the theorem. Here again, we will rely
on Lemma \ref{subexcursions} and the strong Markov property at time $T_D$. We need
to restate the result of Lemma \ref{subexcursions} in a slightly
different form. Let $(\mu,\w)\in\Theta_x$ with $\mu(\{H(\mu)\}=0$
and $\w(t)\in D$ for every $t<H(\mu)$. Under $\bP^*_{\mu,\w}$, we
write
$Y_t=\langle \rho_t,1\rangle$, $J_t=\inf_{r\leq t}Y_r$ and $I_t=J_t-\langle \mu,1\rangle$. 
If $(\alpha_i,\beta_i)$, $i\in I$ are the excursion intervals of $Y-J$ away from
$0$, we introduce the ``excursions'' $(\rho_i,W^i)$, $i\in I$ as defined before
the statement of Lemma \ref{subexcursions}. The starting height of excursion $(\rho_i,W^i)$
is $h_i=H_{\alpha_i}=H(k_{-I_{\alpha_i}}\mu)$. The proof of Lemma \ref{subexcursions}
shows that the point measure
$$\sum_{i\in I} \delta_{(-I_{\alpha_i},\rho^i,W^i)}$$
is Poisson with intensity $1_{[0,<\mu,1>]}(u)du\,\N_{\w(H(k_u\mu))}(d\rho\,dW)$
(this is slightly more precise than the statement of Lemma \ref{subexcursions}). 

We then write $i_1,i_2,\ldots$ for the indices $i\in I$ such that $T_D(W^i)<\infty$,
ranked in such a way that $I_{\alpha_{i_1}}<I_{\alpha_{i_2}}<\cdots$. Our assumption
(\ref{assum-reduc}) guarantees that this ordering is possible, and
we have clearly $h_{i_1}\leq h_{i_2}\leq \cdots$. By well-known properties of Poisson measures,
the processes $W^{i_1},W^{i_2},\ldots$ are independent conditionally on the sequence
$h_{i_1},h_{i_2},\ldots$, and the conditional distribution of 
$W^{i_\ell}$ is $\N^D_{\w(h_{i_\ell})}$.

If we apply the previous considerations to the shifted process
$(\rho_{T_D+s},W_{T_D+s};s\geq 0)$,  taking $\mu=\rho_{T_D}$ and
$\w=W_{T_D}$ and relying on the strong Markov property at $T_D$, we can easily
identify
\ba &&m_D=h_{i_1}\\
&&N_D=1+\sup\{k\geq 1:h_{i_k}=h_{i_1}\}\\
&&W^{D,N_D}=W^{i_1},\,W^{D,{N_D-1}}=W^{i_2},\ldots,\,W^{D,2}=W^{i_{N_D-1}}.
\ea
By a preceding observation, we know that conditionally on $(m_D,N_D)$, the processes
$W^{i_1},\ldots,$ $W^{i_{N_D-1}}$ are independent and distributed according to $\N^D_{\w(m_D)}$.

Combining this with the strong Markov property at time
$T_D$, we see that, conditionally on $(N_D,W^D_0)$, the
processes
$W^{D,2},\ldots,W^{D,N_D}$ are independent and distributed according
to $\N^D_{\hat W^D_0}$ (recall that $\wh W^D_0=W_{T_D}(m_D)$).
An argument similar to the end of the proof of Theorem \ref{tree-Poisson}
(relying on independence properties of Poisson measures) also
shows that, conditionally on $(N_D,W^D_0)$, the vector $(W^{D,2},\ldots,W^{D,N_D})$
is independent of $W^{D,1}$. Furthermore, denote by $\check W^{D,\ell}$
the time-reversed processes
$$\check W^{D,\ell}_s=W^{D,\ell}_{(\sigma(W^{D,\ell})-s)^+}.$$
The time-reversal property already used in the proof of
Theorem \ref{law-exit} implies that the vectors $(W^{D,1},\ldots,W^{D,N_D})$
and $(\check W^{D,N_D},\ldots,\check W^{D,1})$ have the same conditional
distribution given $(N_D,W^D_0)$. Hence, the conditional distribution
of $\check W^{D,1}$, or equivalently that of $W^{D,1}$, is also 
equal to $\N^D_{\hat W^D_0}$. This completes the proof of Theorem \ref{reduced-domain}.
\cq

\medskip
\rems (i) By considering the special case where the spatial motion $\xi_t$ is deterministic,
$\xi_t=t$, and $E=\R_+$, $x=0$ and $D=[0,T)$ for some fixed $T>0$, we 
obtain an
alternative proof of formulas derived in Theorem \ref{lawreduced}. In particular,
formula (\ref{offspringreduced}) is a special case of (\ref{reduc-dom2}).
Similarly, (\ref{lifetimereduced}) can be seen as a special case of
(\ref{reducedom0}).

(ii) In the stable case $\psi(u)=cu^\alpha$, the variable $N_D$ is independent of
$W^D_0$, and its law is given by
$$\N^D_x[r^{N_D}]={(1-r)^\alpha-1+\alpha r\over \alpha-1}.$$
Of course when $\alpha=2$, we have $N_D=2$.

\bibliographystyle{plain}

\end{document}